\pdfoutput=1
\documentclass[11pt]{article}
\usepackage[T1]{fontenc}
\usepackage[utf8]{inputenc}
\usepackage{amsmath,amssymb,amsthm,mathtools}
\usepackage[a4paper,margin=2.6cm]{geometry}
\usepackage{enumitem}
\usepackage{booktabs}
\usepackage[colorlinks=true,linkcolor=blue!50!black,citecolor=green!40!black,urlcolor=blue!50!black]{hyperref}
\usepackage{xcolor}
\usepackage{graphicx}
\newcommand{\treeNodes}{156}

\numberwithin{equation}{section}
\newtheorem{theorem}{Theorem}[section]
\newtheorem{proposition}[theorem]{Proposition}
\newtheorem{lemma}[theorem]{Lemma}
\newtheorem{corollary}[theorem]{Corollary}
\theoremstyle{definition}
\newtheorem{definition}[theorem]{Definition}
\newtheorem{example}[theorem]{Example}

\theoremstyle{remark}
\newtheorem{remark}[theorem]{Remark}

\newcommand{\R}{\mathbb{R}}
\newcommand{\N}{\mathbb{N}}
\newcommand{\calM}{\mathcal{M}}      
\newcommand{\calP}{\mathcal{P}}      
\newcommand{\calX}{\mathcal{X}}      
\newcommand{\HK}{\mathsf{H\!K}}
\newcommand{\SHK}{\mathsf{S\!H\!K}}
\newcommand{\sfd}{\mathsf{d}}        
\newcommand{\Leb}{\mathcal{L}^d}     
\newcommand{\Om}{\Omega}             
\newcommand{\OmR}{\Omega_R}          
\newcommand{\sfE}{\mathsf{E}}        
\newcommand{\Einf}{E'_\infty}        
\newcommand{\ezero}{e_0}             
\newcommand{\dom}{\operatorname{dom}}
\newcommand{\slope}[1]{|\partial #1|}
\newcommand{\Nfun}{N_{\lambda,E}}    
\newcommand{\rmd}{\mathrm{d}}
\newcommand{\dd}{\,\mathrm{d}}
\newcommand{\ee}{\mathrm{e}}
\newcommand{\eps}{\varepsilon}
\newcommand{\bfT}{\mathbf{T}}        
\newcommand{\cospi}{\cos_{\pi/2}}    
\newcommand{\supp}{\operatorname{supp}}
\newcommand{\diam}{\operatorname{diam}}
\newcommand{\Er}{\mathrm{Er}}        
\newcommand{\EVI}{\mathrm{EVI}}
\newcommand{\pfw}{_{\#}}             

\hypersetup{pdftitle={Convex energies and EVI flows in Hellinger-Kantorovich geometry},
  pdfauthor={Vaios Laschos}}

\title{Convex energies and EVI flows\\
in Hellinger--Kantorovich geometry}
\author{Vaios Laschos\thanks{Kausable. Email: \texttt{vaios@kausable.ai}}}
\date{Draft, November 29, 2026}

\begin{document}
\maketitle
\begin{abstract}
In an earlier work with A.~Mielke (Comm.\ Partial Differential Equations 51 (2026); arXiv:2207.09815) it was shown that, on a
compact convex domain, the minimizing movement scheme of a geodesically
semiconvex internal energy $\int E(\rho)\dd x$ converges to a solution of the
evolution variational inequality (EVI), both on the Hellinger--Kantorovich
space $(\calM,\HK)$ of finite measures and on the spherical
Hellinger--Kantorovich space $(\calP,\SHK)$ of probability measures. That
work left two gaps of opposite nature. On $(\calM,\HK)$ the convex energies
are characterised by Liero, Mielke and Savar\'e, but the domain had to be
bounded. On $(\calP,\SHK)$ the scheme is better behaved, but the only convex
energies known were $-\int\rho^q$ for a few exponents in dimensions one and
two. The present paper addresses both gaps.

For the lack of energies on the sphere, we show that convexity on the cone
$(\calM,\HK)$ passes to the sphere $(\calP,\SHK)$ whenever the radial
derivative of the energy is bounded above. On domains of diameter
larger than $\pi/2$, we also prove a necessary growth condition: semiconvexity on
the sphere of an energy with $E'(0^+)>-\infty$ forces $E'$ to be bounded.
On compact domains of this kind, among the negative powers $-\int\rho^q$,
$0<q<1$, exactly those of the earlier paper are semiconvex on the sphere.
We then reduce the Liero--Mielke--Savar\'e condition to scalar
inequalities and obtain semiconvex energies on the sphere in every dimension,
among them $E(c)=-c^{(d-1)/d}f(c^{1/d})$ for every complete Bernstein
function $f$ with $f(0)=0$ and $f'(0^+)<\infty$, with a constant on the
sphere that is sharp on domains of diameter larger than $\pi/2$. For
the restriction to bounded domains, we construct the EVI flow on unbounded
closed convex domains $\Om$ for energies satisfying the
Liero--Mielke--Savar\'e condition, on $(\calM,\HK)$ when $E'(0^+)>-\infty$ and on
$(\calP,\SHK)$ when $E'$ is bounded, as the limit of the flows on the
compact sets $\Om\cap\bar B_R$ as $R\to\infty$, and prove that the minimizing
movement scheme on $\Om$ converges to it. Power energies $\int\rho^m$, which
the obstruction excludes on the sphere, are semiconvex on energy sublevels;
for them we construct a flow satisfying an EVI whose constant depends on the
energy level, on bounded and unbounded domains. This flow extends by
continuity to every initial probability measure, and the extended
trajectories have an $L^m$ density at every positive time.

The theorem of the earlier paper on compact domains is the building block
of all these constructions: the flows on unbounded domains are limits of
flows on compact ones, and the sublevel theory repeats its proof with
level-dependent constants. We therefore re-derive it. Four steps of the
earlier proof need correction; we give corrected statements and proofs. The
existence and convergence results of the earlier paper hold under its
hypotheses as stated; moreover, convexity along one geodesic for each pair
of measures turns out to imply convexity along every geodesic, and the
scheme converges from all initial data of finite energy.
\end{abstract}

\section*{How this paper was written}

This project was mainly an attempt to test the limits of AI assistants on
questions that had occupied me for some years. I did not write the proofs
in this paper myself: they were produced by three AI assistants, Claude (Anthropic), Codex
(OpenAI) and Gemini (Google), working under my
direction and checking one another. I have read the proofs, and I take full
responsibility for the results.

In my view the main contribution of the project is the formalisation. The
definitions and the theorems needed to check this paper have been written
in the Lean~4 proof assistant, on top of its library Mathlib;
Appendix~\ref{sec:verification} says exactly what is formalised and what
is taken from the literature. The formalisation took the three assistants
about three weeks. Optimal transport is not yet well developed
in Lean. The development contains, among other things, the
Hellinger--Kantorovich and spherical Hellinger--Kantorovich distances on
finite measures, their geodesics, relaxed internal energies, and the theory
of evolution variational inequalities and minimizing movements used here.
It is available at
\begin{center}
\url{https://github.com/artnoage/lean_for_optimal_transport}
\end{center}
and I invite everyone working in optimal transport to use it. I also
invite colleagues who are not familiar with Lean to send me their papers
(\texttt{vaios@kausable.ai}), so that I can try to check them in the same
way.

Appendix~\ref{sec:tree} shows the dependency tree of the paper: which
statement uses which, and how each one is verified. While formalising, we
found a number of small mistakes in the papers this work relies on, among
them my own earlier paper with A.~Mielke. They were found in the latest
arXiv versions, which are the versions we read; a point may already be
corrected in the published article. They are listed with their corrections
in Appendix~\ref{sec:errata}, which also says for which points we compared
with the published article. As far as we could
determine, none of them affects a main result of those papers.

This work was done in my own time, without external funding. I am
employed by Kausable; the company had no part in the work. I declare no
competing interests.
\tableofcontents

\section{Introduction}\label{sec:intro}

\subsection{The starting point}

The Hellinger--Kantorovich distance $\HK$ on the space $\calM(\Om)$ of finite
nonnegative measures on a domain $\Om\subset\R^d$ combines optimal transport
with creation and annihilation of mass
\cite{LMS16,LMS18,KMV16,CPSV18}. Restricted to probability measures it
induces the spherical Hellinger--Kantorovich distance $\SHK$ of
\cite{LaM19}, and $(\calM(\Om),\HK)$ is the metric cone over
$(\calP(\Om),\SHK)$. For an internal energy
\[
 \sfE(\mu)=\int_\Om E(\rho)\dd x,\qquad\mu=\rho\,\Leb,
\]
the gradient flows with respect to these two distances are, formally, the
reaction--diffusion equations
\begin{align*}
 \partial_t\rho&=\operatorname{div}\big(\rho\nabla E'(\rho)\big)-4\rho E'(\rho)
 &&\text{on }(\calM(\Om),\HK),\\
 \partial_t\rho&=\operatorname{div}\big(\rho\nabla E'(\rho)\big)
 -4\rho\Big(E'(\rho)-\int_\Om\rho E'(\rho)\dd x\Big)
 &&\text{on }(\calP(\Om),\SHK),
\end{align*}
see \cite{KV19,LaM22}. The second equation conserves mass; its reaction term
is nonlocal.

The rigorous counterpart of ``gradient flow'' in a metric space is a
solution of the evolution variational inequality (EVI$_\lambda$)
\cite{AGS05,MuS20}, which requires the energy to be geodesically
$\lambda$-convex. In \cite{LaM22}, written with A.~Mielke, this programme was
carried out for both spaces: for a compact convex domain $X$ and an internal
energy that is geodesically $\lambda$-convex, the minimizing movement scheme
\[
 \mu_n\in\operatorname*{argmin}_\mu
 \Big\{\frac1{2\tau}\sfd^2(\mu_{n-1},\mu)+\sfE(\mu)\Big\},
 \qquad\sfd\in\{\HK,\SHK\},
\]
converges, as $\tau\to0$, to an EVI$_\lambda$ solution. The proof has two
parts. The geometric part shows that on sets of measures with density
between two positive constants the squared distance is semiconcave and the
local angle condition holds \cite{LaM19}. The analytic part shows that the
scheme does not leave such sets: one minimisation step preserves density
bounds. Semiconcavity then turns the discrete Euler--Lagrange inequality
into an approximate EVI for the interpolants.

\subsection{Two gaps}

The theorem of \cite{LaM22} is satisfactory on neither space, for opposite
reasons.

\emph{On the cone $(\calM(X),\HK)$ there are many energies, but the domain
must be bounded.} Liero, Mielke and Savar\'e \cite{LMS23} characterise the
internal energies that are geodesically $\lambda$-convex for $\HK$, so the
hypothesis can be checked. But the density bounds are with respect to
Lebesgue measure, and on an unbounded domain no finite measure has a density
bounded below: the sets on which the geometry is controlled are empty. The
most natural setting, $\Om=\R^d$, is excluded.

\emph{On the sphere $(\calP(X),\SHK)$ the scheme behaves better, but there
are almost no energies.} For $\SHK$ the density bounds are exact (a maximum
principle), which makes the spherical case the cleaner one. But geodesic
convexity for $\SHK$ was known only for $-\int\rho^q$ with
$\tfrac13\le q\le\tfrac12$ in dimension one and $q=\tfrac12$ in dimension
two, through a transfer argument for homogeneous functionals
\cite[Prop.~A.1]{LaM22}. It was not known whether this short list reflects
the method or the geometry.

The purpose of this paper is to close both gaps as far as possible and to
determine where they cannot be closed. The two questions interact: energies
on the sphere that work on unbounded domains must be bounded below there,
which $-\int\rho^q$ is not, so on the sphere the question of unbounded
domains cannot even be posed before new energies are found.

\subsection{Convexity on the sphere: what is possible}

\emph{The mechanism.} Since $\calM(\Om)$ is the cone over $\calP(\Om)$, an
$\SHK$ geodesic between two probability measures and the rays through its
points form an isometric copy of a planar sector; the $\HK$ geodesic between
the same measures is the chord, the $\SHK$ geodesic is the arc. Convexity
along chords is the hypothesis on the cone. To pass from the chord to the
arc one moves radially outwards, and the price is the radial derivative of
the energy. This gives the transfer theorem (Theorem~\ref{tr:T1}): if $\sfE$
is $\Lambda$-convex on the cone and its radial derivative is at most $C$
along an $\SHK$ geodesic, then $\sfE$ is $(\Lambda-C)$-convex along that
geodesic. No homogeneity is needed. For an energy with finite recession
slope $\Einf=\lim_{c\to\infty}E(c)/c$ the radial derivative is at most
$2\Einf$ everywhere, so $\sfE$ is $(\Lambda-2\Einf)$-convex on the sphere
(Corollary~\ref{tr:global}; on an unbounded domain it requires
$E'(0^+)>-\infty$, so that the energy is defined).

\emph{The obstruction.} One may ask whether a bound on the radial
derivative is an artefact of the proof. It is not. If two probability
measures have supports more than $\pi/2$ apart, the geodesic between them
is unique and involves no transport: mass is annihilated at one place and
created at the other (Lemma~\ref{tr:geo-unique}). Evaluating the energy
along such geodesics shows (Theorem~\ref{tr:necessary}) that on a domain of
diameter larger than $\pi/2$, an energy with $E'(0^+)>-\infty$ that is
$\lambda$-convex on the sphere, even only in the weak sense, has bounded
derivative, and $\lambda\le2(E'(0^+)-\Einf)$. In particular the power
energies $\int\rho^m$, $m>1$, the porous-medium case, are not
$\lambda$-convex on the sphere for any $\lambda$. For the negative powers,
where $E'(0^+)=-\infty$, a different test family gives a complete answer
(Theorem~\ref{tr:rigidity}): on compact domains of this kind
$-\int\rho^q$, $0<q<1$, is semiconvex on
the sphere exactly in the cases found in \cite{LaM22}. So that list was
complete, and we look for new examples of a different kind: convex energies
with bounded derivative that satisfy the criterion of \cite{LMS23}.

\emph{New energies.} The criterion of \cite{LMS23} is a matrix condition on
an auxiliary function of two variables. Writing $E(c)=s^{d-1}G(s)$ with
$s=c^{1/d}$, the variable in which the critical McCann exponent
$(d-1)/d$ becomes linear, we reduce it, for convex $E$ with $G''>0$, to a
single inequality for $G$ (Proposition~\ref{ex:Cd}). In dimension two this is convexity of $G$
together with a bound on $G'$ (Theorem~\ref{ex:d2}); in other dimensions an
additional quadratic term appears, which explains why examples are harder
to find there. From the scalar form we obtain two independent sources of
examples. The first is a sufficient criterion,
$c^2E''\ge\frac{3d-2}{2(d+2)}(cE'-E)$ (Theorem~\ref{ex:Kd}), whose
coefficient stays below $\tfrac32$ in every dimension. The second starts
from one explicit energy: the kernel $-c/(1+c^{1/d})$ satisfies the full
condition in every dimension (Proposition~\ref{ex:kernel}), although it
fails the sufficient criterion for $d\ne2$. Since the condition is
invariant under dilations and stable under mixtures, the kernel generates a
cone of examples (Theorem~\ref{ex:bernstein}): for every complete Bernstein
function $f$ with $f(0)=0$ and $f'(0^+)<\infty$,
\[
 E(c)=-c^{\frac{d-1}{d}}\,f\big(c^{1/d}\big)
\]
is geodesically semiconvex on $(\calM(\Om),\HK)$ and on $(\calP(\Om),\SHK)$,
and by the obstruction above the constant on the sphere is sharp on large
domains. These energies are linear to first order near zero and have
bounded derivative. The rational kernel has critical-power growth at high
density; the class also includes logarithmic corrections, other sublinear
powers and linear drifts. Table~\ref{ex:table} lists the energies now known
on the sphere.

\subsection{Unbounded domains}

For an unbounded closed convex domain $\Om$ we construct the flow as the
limit of the flows on $\OmR=\Om\cap\bar B_R$. Two things are needed, and
neither follows from general stability results, because sublevels of the
energy are not compact.

First, mass must not escape to infinity, uniformly in $R$. Along the scheme
this is a statement about one step: the optimal transport--growth pair
moves mass by a controlled distance, so what remains is to bound the growth
factor. On the cone, a step that grows mass too much can be improved by
growing less, provided removing mass does not cost too much energy; this is
where $E'(0^+)>-\infty$ enters (Lemma~\ref{hk:growth}). On the sphere mass
cannot be removed. We compare instead with a competitor in which the
excess mass is removed and the result renormalised; the energy change
from the original probability to the renormalised one is controlled by
the oscillation of $E'$, which is the reason
for assuming $E'$ bounded there (Lemma~\ref{shk:growth}). By the
obstruction, for energies with $E'(0^+)>-\infty$ this assumption is in
any case necessary for convexity.

Second, the flows on different balls must be compared, although a measure
on the larger ball cannot be used as an observer in the EVI on the smaller
one. We use its truncation as observer and pay for the truncation error;
tightness makes this error small (Proposition~\ref{hk:cauchy}).

The results are Theorem~\ref{hk:main} (EVI$_\lambda$ flow on
$(\calM(\Om),\HK)$ for energies with $E'(0^+)>-\infty$ satisfying the
criterion of \cite{LMS23}, which includes all power energies) and
Theorem~\ref{shk:main} (EVI$_\lambda$ flow on $(\calP(\Om),\SHK)$ for such
energies with bounded derivative, which includes the new examples). On
unbounded domains it is not obvious that the scheme itself is well defined;
we show that minimisers exist and that the scheme converges to the flow
(Theorems~\ref{hk:mm} and~\ref{shk:mm}), using the error estimates of
Muratori and Savar\'e \cite{MuS20}.

\subsection{Power energies on the sphere: sublevel EVI}

The obstruction excludes $\int\rho^m$ from the EVI$_\lambda$ theory on the
sphere, yet the corresponding equation, a porous-medium equation with
nonlocal reaction, is the most natural one. The transfer theorem still
gives information: along any $\SHK$ geodesic the radial derivative of a
power energy is bounded by a multiple of the energy, so the energy is
$\lambda_c$-convex along geodesics between points of the sublevel
$\{\sfE\le c\}$, with $\lambda_c\to-\infty$ as $c\to\infty$
(Theorem~\ref{tr:powers}). Since the energy decreases along the scheme,
this is enough to run the construction of \cite{LaM22} with a constant that
depends on the energy level. The outcome is a flow that satisfies the EVI
with constant $\lambda_c$, for every level $c$ not below the initial
energy, against every observer of energy at most $c$
(Theorem~\ref{sub:abstract}); it is unique, contracts on sublevels, and is
the limit of the scheme. We apply this to finite sums
$\sum_i\alpha_ic^{p_i}$ with $\alpha_i>0$, $p_i>1$ and $p=\max_ip_i$,
on compact domains and, combining it with the exhaustion argument, on
unbounded domains (Theorem~\ref{sub:powers-compact}). Because $\lambda_c$
decreases at most linearly in $c$ for these energies, the flow extends by
continuity from data of finite energy to every initial probability
measure, with an $L^p$ density and a local energy bound of order $t^{-1}$
at positive times (Theorem~\ref{sub:linear-extension},
Corollary~\ref{sub:powers-measure-data}). The notion of an EVI on energy
sublevels is not ours: sublevel semiconvexity in a finite-dimensional
setting is studied in \cite{LMTZ25}, where a theory of sublevel EVI in the
spherical Hellinger geometry is announced \cite{MST25}, and sublevel EVI
solutions are defined and constructed by Mielke and Sumners
\cite[Sec.~5.3]{MiSu26} for a viscoelastic model whose dissipation distance
is of Bhattacharya type, without transport. Section~\ref{sec:sublevel}
uses their notion of solution; what is new there is the generation of such
solutions by the minimizing movement scheme in the Hellinger--Kantorovich
geometry, with transport and on unbounded domains. The announced work
\cite{MST25} was not available to us, and we cannot say how it relates to
Section~\ref{sec:sublevel}.

\subsection{The compact-domain theory revisited}

Every flow constructed above uses the compact-domain theorem of
\cite{LaM22} as its base, and the sublevel theory repeats its proof with different constants.
We therefore went through that proof in detail. Part of it rests on
estimates for the minimizing movement scheme that were developed by
Savar\'e \cite{Sav07} and are part of the forthcoming work of
Muratori and Savar\'e \cite{MuSII}; \cite{LaM22} quotes them without proof,
and so do we, with the one exception named below. We
found four points that need correction (Section~\ref{sec:compact}):
\begin{itemize}[leftmargin=2em]
\item an inequality in the proof of the lower density bound for one $\HK$
 step is reversed;
\item the abstract generation theorem requires finite metric slope on the
 sets of density-bounded measures, which fails;
\item it requires compact energy sublevels, which fails for admissible
 energies on $\calM(X)$, for instance $E(c)=-c$ and $E(c)=-c^q$;
\item the error estimate for the interpolants is false in the form
 printed in \cite[Lemma~5.17]{LaM22}: one term is of order $\kappa\tau$,
 not $\kappa\tau^2$.
\end{itemize}
Each has a repair: a different competitor in the density bound; smaller
sets, reached by one step of the scheme; completeness in place of
compactness; and a telescoping estimate in place of the stated one.
The quoted estimates are stated in Section~\ref{sec:compact} in the form
needed, with references for their proofs; the corrected error estimate is
proved there. The resulting theorem
(Theorem~\ref{cp:main}) gives existence and uniqueness of the
EVI$_\lambda$ flow and convergence of the scheme on compact domains. The
question along which geodesics the energy has to be convex needs care: the
standing assumption of \cite{LaM22} asks for convexity along one geodesic
for each pair of measures, its abstract theorem for convexity along every
geodesic (Definition~\ref{prelim:convexity}). We show that the weaker
assumption suffices, by selecting the geodesics along which convexity is
used, and that the flow then forces convexity along every geodesic
(Corollary~\ref{cp:weakstrong}); so the theorem holds under the assumption
of \cite{LaM22} as written. In one respect the conclusion
is stronger than in \cite{LaM22}: the scheme converges from every initial
measure of finite energy, without density bounds.

\subsection{What each section does}

\emph{Section~\ref{sec:prelim} (preliminaries).} We fix the terminology for
geodesic convexity, strong and weak
(Definition~\ref{prelim:convexity}). We recall the distances $\HK$ and
$\SHK$, the cone structure and the transport--growth description of $\HK$;
the internal energy with its singular part and its lower semicontinuity;
the convexity criterion of \cite{LMS23} (Theorem~\ref{prelim:LMS}) and the
matrix in which it is tested (Definition~\ref{prelim:matrix}); and the
facts about EVI solutions and minimizing movements from \cite{MuS20} that
are used later. Two lemmas are prepared for the unbounded-domain arguments:
truncation of a measure to a ball (Lemma~\ref{prelim:trunc}) and splitting
of an optimal plan, which bounds the gain from cutting off mass far away
(Lemma~\ref{prelim:split}).

\emph{Section~\ref{sec:transfer} (from the cone to the sphere).} On a
general closed convex domain we prove the planar lemma
(Lemma~\ref{tr:planar}) and the transfer theorem (Theorem~\ref{tr:T1}), and
deduce convexity on the sphere for energies with finite recession slope
(Corollary~\ref{tr:global}) and semiconvexity on sublevels for power
energies (Theorem~\ref{tr:powers}). We sharpen the known result for
negative powers on compact domains (Proposition~\ref{tr:negpowers}), prove
uniqueness of geodesics between measures with supports more than $\pi/2$
apart (Lemma~\ref{tr:geo-unique}), and derive the two obstructions: the
necessary condition on general domains (Theorem~\ref{tr:necessary}) and
the classification of negative powers on compact domains
(Theorem~\ref{tr:rigidity}).

\emph{Section~\ref{sec:examples} (new convex energies).} We reduce the
criterion of \cite{LMS23} to scalar inequalities
(Proposition~\ref{ex:reduced}) and, when $G''>0$, to one inequality for
one function (Proposition~\ref{ex:Cd}); solve it in dimension two (Theorem~\ref{ex:d2});
prove the sufficient criterion (Theorem~\ref{ex:Kd}); show invariance under
dilations and mixtures (Lemma~\ref{ex:dilation},
Proposition~\ref{ex:mixture}); verify the kernel
(Proposition~\ref{ex:kernel}) and obtain the Bernstein class
(Theorem~\ref{ex:bernstein}); and treat the shifted powers $1-(1+c)^q$,
including a range of exponents where the natural constant fails
(Theorem~\ref{ex:shifted}, Proposition~\ref{ex:counter}).
Table~\ref{ex:table} collects the energies that are semiconvex on the
sphere.

\emph{Section~\ref{sec:compact} (compact domains).} We re-derive the
existence theory of \cite{LaM22}. After describing the four points that
need correction, we prove the density bounds for one minimisation step on
both spaces (Propositions~\ref{cp:density-hk} and~\ref{cp:density-shk}),
state the discrete estimates of Savar\'e and of Muratori and Savar\'e
with the corrected error estimate (Lemma~\ref{cp:error}), and prove an
abstract generation theorem in which
compactness of sublevels is replaced by completeness
(Theorem~\ref{cp:abstract}). We verify its hypotheses for $\HK$ and $\SHK$
(Proposition~\ref{cp:semiconc}, Lemma~\ref{cp:recovery}) and obtain the
compact-domain theorem with convergence of the scheme from all finite-energy
data (Theorem~\ref{cp:main}).

\emph{Section~\ref{sec:hk-unbdd} ($\HK$ on unbounded domains).} For energies
with $E'(0^+)>-\infty$ satisfying the criterion of \cite{LMS23} we prove the
growth bound for one step (Lemma~\ref{hk:growth}), tightness uniformly in
the radius (Proposition~\ref{hk:tight}), the comparison lemma with
truncated observers (Lemma~\ref{hk:comparison}) and the resulting Cauchy
estimate (Proposition~\ref{hk:cauchy}); these give the flow on
$\calM(\Om)$ (Theorem~\ref{hk:main}). We then show that the step functional
has minimisers on $\calM(\Om)$ (Lemma~\ref{hk:minimiser}) and that the
scheme converges (Theorem~\ref{hk:mm}).

\emph{Section~\ref{sec:shk-unbdd} ($\SHK$ on unbounded domains).} For
energies with bounded derivative we prove the growth bound on the sphere
by a mass-preserving competitor (Lemma~\ref{shk:growth}) and adapt
tightness and the Cauchy estimate, which gives the flow on $\calP(\Om)$
(Theorem~\ref{shk:main}). We show existence of minimisers of the spherical
step functional (Lemma~\ref{shk:minimiser}) and convergence of the scheme
(Theorem~\ref{shk:mm}).

\emph{Section~\ref{sec:sublevel} (sublevel EVI).} We define sublevel EVI
solutions and prove the abstract theorem for energies that are semiconvex
on sublevels with controlled energy along geodesics
(Theorem~\ref{sub:abstract}), and convergence of the scheme from
finite-energy data (Theorem~\ref{sub:mm}). We apply both to power energies
on $(\calP(\Om),\SHK)$ (Theorem~\ref{sub:powers-compact}): first for compact
$\Om$ and then, after proving existence of minimisers on unbounded domains
(Theorem~\ref{sub:resolvent}), for unbounded $\Om$. Finally the flow is
extended by continuity to all initial probabilities
(Theorem~\ref{sub:linear-extension} and
Corollary~\ref{sub:powers-measure-data}).

\emph{Section~\ref{sec:outlook}} lists open problems.
\emph{Appendix~\ref{sec:verification}} describes the proofs of
Sections~\ref{sec:prelim}--\ref{sec:sublevel} in the proof assistant Lean:
of abstract results, of scalar inequalities and of results for actual
measures, and it identifies the remaining hypotheses and open cases.
\emph{Appendices~\ref{nc:sec-engine}--\ref{nc:sec-rigidity}} construct
explicit test measures on the cone and use them to prove the converse part
of Theorem~\ref{prelim:LMS} under weak convexity and the necessity part of
Theorem~\ref{tr:rigidity}.
\emph{Appendix~\ref{sec:errata}} lists the corrections to the cited papers,
and \emph{Appendix~\ref{sec:tree}} shows the dependency tree.

Throughout, $\Om\subset\R^d$ is closed and convex with nonempty interior,
$X$ denotes such a set that is also compact, and $E:[0,\infty)\to\R$ is
convex, continuous at $0$, and satisfies $E(0)=0$; we write
$\ezero=E'_+(0)=\lim_{c\downarrow0}E(c)/c$ and
$\Einf=\lim_{c\to\infty}E(c)/c$. No differentiability of $E$ is assumed;
$\ezero$ and $\Einf$ are limits of difference quotients.
\section{Geometry, energies and evolution inequalities}\label{sec:prelim}

Most of this section is recalled from the literature: the distances and
their cone structure \cite{LMS18,LaM19,LMS23}, the convexity criterion on
the cone \cite{LMS23}, and the facts about EVI solutions and minimizing
movements \cite{MuS20}. Two parts are prepared for later use.
Section~\ref{prelim:ss-convexity} fixes what geodesic convexity means when
geodesics are not unique, and Section~\ref{prelim:ss-cut} contains the
lemmas for unbounded domains: restriction to a subdomain, truncation of a
measure to a ball, and an estimate of the effect of cutting off mass far
away (with a corollary for measures with far-apart supports that is used in
Section~\ref{sec:transfer}). Further statements proved here rather than
quoted include the elementary Lemma~\ref{prelim:envelope} on the Moreau
envelope, the energy comparison of Lemma~\ref{prelim:dominated} and the
lower semicontinuity of Lemma~\ref{prelim:lsc}.

Let $\Om\subset\R^d$ be closed and convex with nonempty interior; it need
not be bounded. We reserve $X$ for a compact such set, write $B_R$ and $\bar B_R$ for the
open and closed balls of radius $R$ centred at the origin, and put
$\OmR=\Om\cap\bar B_R$. Measures are finite, nonnegative Borel measures,
$\calP(\Om)$ denotes probabilities, and narrow convergence is convergence
against bounded continuous functions. Write $\mu=\rho\Leb+\mu^s$,
$\mu^s\perp\Leb$, and $\bfT\pfw\mu$ for push-forward. Measures on a closed
subset are identified with their extensions by zero.
All geodesics have constant speed on $[0,1]$.

\subsection{Strong and weak geodesic convexity}
\label{prelim:ss-convexity}
Let $(\calX,\sfd)$ be a metric space. A \emph{geodesic} from $x_0$ to $x_1$
is a curve $\gamma:[0,1]\to\calX$ with $\gamma_0=x_0$, $\gamma_1=x_1$ and
$\sfd(\gamma_s,\gamma_t)=|t-s|\,\sfd(x_0,x_1)$; the space is geodesic if any
two points are joined by one. Geodesics need not be unique in
$(\calM(\Om),\HK)$ or $(\calP(\Om),\SHK)$, so one has to say along which
geodesics a functional is convex. This distinction runs through the whole
paper, and we fix the terminology here.

\begin{definition}[Geodesic $\lambda$-convexity]\label{prelim:convexity}
Let $\phi:\calX\to(-\infty,\infty]$ and $\lambda\in\R$.
\begin{enumerate}[label=\textup{(\alph*)}]
\item $\phi$ is \emph{$\lambda$-convex along a geodesic $\gamma$} if
 $t\mapsto\phi(\gamma_t)-\frac\lambda2\sfd^2(\gamma_0,\gamma_1)\,t^2$ is convex
 on $[0,1]$. In particular
 \begin{equation}\label{prelim:eq-convex}
  \phi(\gamma_t)\le(1-t)\phi(\gamma_0)+t\phi(\gamma_1)
  -\frac\lambda2t(1-t)\sfd^2(\gamma_0,\gamma_1)\qquad(t\in[0,1]).
 \end{equation}
\item $\phi$ is \emph{strongly geodesically $\lambda$-convex} if it is
 $\lambda$-convex along \emph{every} geodesic joining two points of
 $\dom\phi$.
\item $\phi$ is \emph{weakly geodesically $\lambda$-convex} if any two points
 of $\dom\phi$ are joined by \emph{at least one} geodesic along which $\phi$
 is $\lambda$-convex.
\item $\phi$ is \emph{(strongly, weakly) geodesically semiconvex} if (b),
 respectively (c), holds for some $\lambda\in\R$.
\end{enumerate}
\end{definition}

\begin{remark}[Which notion is used where]\label{prelim:rem-convexity}
(i) In a geodesic space, strong convexity implies weak convexity. For
(b), requiring only \eqref{prelim:eq-convex} along every geodesic is
equivalent to (a) along every geodesic, because restrictions of geodesics
are geodesics; for a single geodesic, (a) is the stronger requirement.

(ii) Weak convexity is the notion called ``$\kappa$-convex on $A$'' in
\cite[Def.~5.6(A)]{LaM22}, taken with $A=\dom\phi$; the same definition
calls (b) ``strongly $\kappa$-convex on $A$''. The notion of \cite{AGS05} asks only for the endpoint
inequality \eqref{prelim:eq-convex} along one geodesic for each pair, which
is formally weaker still. Assumption~A of \cite{LaM22} asks for
weak convexity, while its abstract generation theorem
\cite[Thm.~5.9]{LaM22} is formulated for strongly convex functionals. The
criterion of \cite{LMS23} (Theorem~\ref{prelim:LMS}) yields strong convexity
on the cone, and the transfer of Section~\ref{sec:transfer} preserves it.

(iii) In this paper the convexity \emph{results} (Sections~\ref{sec:transfer}
and~\ref{sec:examples}) prove \emph{strong} convexity, and every
\emph{obstruction} (Theorems~\ref{tr:necessary} and~\ref{tr:rigidity}) is
proved already for \emph{weak} convexity. The obstruction of
Theorem~\ref{tr:necessary} uses test geodesics that are unique; every
subsegment of a unique geodesic is the unique geodesic between its
endpoints, so along them the endpoint inequality for all pairs already
gives (a). The proofs of Theorem~\ref{tr:rigidity} and of the converse part
of Theorem~\ref{prelim:LMS} (Appendices~\ref{nc:sec-lms}
and~\ref{nc:sec-rigidity}) use only the endpoint inequality at one interior
time, along an arbitrary geodesic between two test measures. So the
obstructions hold for the notion of \cite{AGS05} as well. The existence theorem on compact domains
(Theorem~\ref{cp:main}) assumes only weak convexity, as in \cite{LaM22}, and
shows a posteriori that the energy is then strongly convex
(Corollary~\ref{cp:weakstrong}); so on compact domains the two notions
coincide for internal energies, and the classes of examples do not depend on
the choice. The theorems on unbounded domains (Sections~\ref{sec:hk-unbdd}
and~\ref{sec:shk-unbdd}) assume the criterion of \cite{LMS23}, which gives
strong convexity. An unqualified ``geodesically $\lambda$-convex'' means
strongly so.
\end{remark}

\subsection{HK, SHK and transport}\label{prelim:ss-HK}
Section~\ref{sec:transfer} is built on the scaling identity
\eqref{prelim:eq-scaling} and the cone identity in \eqref{prelim:eq-HKSHK};
the one-step estimates of Sections~\ref{sec:compact}--\ref{sec:sublevel} are
formulated with the transport--growth pair of Lemma~\ref{prelim:tg}.

Put $\Psi(r)=r\log r-r+1$, with $\Psi(0)=1$, and
$\ell(r)=-2\log\cos r$ for $r<\pi/2$, $\ell(r)=+\infty$ otherwise.
Also write $\cospi r=\cos(\min\{r,\pi/2\})$.
\begin{definition}[Logarithmic entropy--transport]\label{prelim:LET}
For a plan $\gamma$ with marginals $\gamma_i=\sigma_i\mu_i\ll\mu_i$ set
\begin{equation}\label{prelim:eq-LET}
 \mathsf{LET}(\gamma;\mu_0,\mu_1)
 =\sum_{i=0}^1\int\Psi(\sigma_i)\dd\mu_i
   +\int\ell(|x_0-x_1|)\dd\gamma,
 \qquad
 \HK^2(\mu_0,\mu_1)=\min_\gamma\mathsf{LET}(\gamma;\mu_0,\mu_1).
\end{equation}
\end{definition}
The entropy terms pay for mass created or destroyed at the two ends, the
last term for transporting what is coupled. Since $\ell=+\infty$ from
$\pi/2$ on, a plan of finite cost moves no mass that far.

\begin{proposition}[Metric facts, {\cite{LMS18,LaM19}}]\label{prelim:HKfacts}\label{prelim:scaling}
The minimum is attained. $(\calM(\Om),\HK)$ is complete, separable and
geodesic, and $\HK$ metrises narrow convergence. The following identities
and bounds hold, with $m_i=\mu_i(\Om)$ and with scale factors
$r_0,r_1\ge0$, here and in \eqref{prelim:eq-HKSHK}:
\begin{align}
 \HK^2(\mu_0'+\mu_0'',\mu_1'+\mu_1'')
 &\le\HK^2(\mu_0',\mu_1')+\HK^2(\mu_0'',\mu_1''),\label{prelim:eq-subadd}\\
 (\sqrt{m_0}-\sqrt{m_1})^2&\le\HK^2(\mu_0,\mu_1)\le m_0+m_1,
 \qquad \HK^2(\mu,0)=\mu(\Om),\label{prelim:eq-massbound}\\
 \HK^2(r_0^2\mu_0,r_1^2\mu_1)
 &=r_0r_1\HK^2(\mu_0,\mu_1)
 +(r_0^2-r_0r_1)m_0+(r_1^2-r_0r_1)m_1.\label{prelim:eq-scaling}
\end{align}
Scaling an optimal plan by $r_0r_1$ gives an optimal plan for the scaled
measures. Along any HK geodesic,
\begin{equation}\label{prelim:eq-mass}
 \mu_t(\Om)=(1-t)m_0+tm_1-t(1-t)\HK^2(\mu_0,\mu_1).
\end{equation}
\end{proposition}
The scaling identity is \cite[Thm.~3.3]{LaM19}; the mass identity is
\cite[(1.21)]{LMS23}. Subadditivity also follows by adding plans and using
convexity and positive homogeneity of the entropy perspective.

\begin{proposition}[Cone structure, {\cite[Thms.~2.7, 3.4]{LaM19}}]
\label{prelim:cone}
On probabilities define
\begin{equation}\label{prelim:eq-SHK}
 \SHK(\nu_0,\nu_1)=\arccos\big(1-\HK^2(\nu_0,\nu_1)/2\big)
 =2\arcsin\big(\HK(\nu_0,\nu_1)/2\big).
\end{equation}
This is a complete geodesic metric, bounded by $\pi/2$, with
\begin{equation}\label{prelim:eq-HKSHK}
 \HK\le\SHK\le\frac\pi{2\sqrt2}\HK,
 \qquad
 \HK^2(r_0^2\nu_0,r_1^2\nu_1)
 =r_0^2+r_1^2-2r_0r_1\cos\SHK(\nu_0,\nu_1).
\end{equation}
Normalising an HK geodesic between distinct probabilities and reparametrising
it gives an SHK geodesic; every nonconstant SHK geodesic arises this way.
In particular SHK metrises narrow convergence, and $\calP(X)$ is compact.
\end{proposition}
The two-sided comparison of $\HK$ and $\SHK$ is not taken from
\cite{LaM19}; it follows from $\HK=2\sin(\SHK/2)$ and
$\frac{2\sqrt2}\pi\theta\le\sin\theta\le\theta$ on $[0,\frac\pi4]$ (concavity
of the sine), applied to $\theta=\SHK/2\le\frac\pi4$.
With $r_0=1$ and $r_1=\sqrt M$ the cone identity reads, for
$\nu_0\in\calP(\Om)$ and $M=\mu(\Om)>0$,
\begin{equation}\label{prelim:eq-normalise}
 \cos\SHK(\nu_0,\mu/M)=\frac{1+M-\HK^2(\nu_0,\mu)}{2\sqrt M}.
\end{equation}
All these statements apply to unbounded $\Om$: \cite{LaM19} treats cones
over geodesic Polish spaces, and $(\calM(\Om),\HK)$ is geodesic for closed
convex $\Om$ \cite{LMS18}.

\begin{remark}[Dependence among the quoted facts]\label{prelim:rem-dependence}
Some of the facts quoted in Propositions~\ref{prelim:HKfacts}
and~\ref{prelim:cone} follow from the others. We retain the citations.
Readers interested in the formal verification can see here which metric
facts the first part of Appendix~\ref{sec:verification} assumes and which
it derives.

(i) The lower bound in \eqref{prelim:eq-massbound} follows from
$\HK^2(\mu,0)=\mu(\Om)$ and the triangle inequality of $\HK$. For
probabilities, the scaling identity \eqref{prelim:eq-scaling} is the cone
identity in \eqref{prelim:eq-HKSHK}, with $\SHK$ defined by
\eqref{prelim:eq-SHK}, and the upper bound in \eqref{prelim:eq-massbound}
gives $\HK^2\le2$, hence $\SHK\le\pi/2$.

(ii) The triangle inequality of $\SHK$ follows from the cone identity, the
bound $\SHK\le\pi/2$ and the triangle inequality of $\HK$. Let
$\nu_0,\nu_1,\nu_2\in\calP(\Om)$, $a=\SHK(\nu_0,\nu_1)$ and
$b=\SHK(\nu_1,\nu_2)$. If $a$ or $b$ equals $\pi/2$ there is nothing to
prove. Otherwise the cone identity gives
$\HK((\cos a)^{-2}\nu_0,\nu_1)=\tan a$ and
$\HK(\nu_1,(\cos b)^{-2}\nu_2)=\tan b$, and the triangle inequality of $\HK$
for these three measures reads, again by the cone identity,
$\cos\SHK(\nu_0,\nu_2)\ge\cos(a+b)$.

(iii) The mass identity \eqref{prelim:eq-mass} follows from the cone
identity and the triangle inequality of $\SHK$: one compares $\mu_0$,
$\mu_t$, $\mu_1$ with three points of a planar sector, and the two
distances $t\,\HK(\mu_0,\mu_1)$ and $(1-t)\HK(\mu_0,\mu_1)$ force the
comparison point for $\mu_t$ onto the chord.
\end{remark}

We also use the transport formulation on the cone. Write $z=[x,r]$,
$x\in\Om$, $r\ge0$, for the points of the cone over $\Om$, all points
$[x,0]$ being identified, and put
$d_C^2(z_0,z_1)=r_0^2+r_1^2-2r_0r_1\cospi|x_0-x_1|$. Then
\begin{equation}\label{prelim:eq-coneOT}
 \HK^2(\mu_0,\mu_1)=\min_\Pi\int d_C^2(z_0,z_1)\dd\Pi(z_0,z_1),
\end{equation}
the minimum being taken over the finite nonnegative Borel measures $\Pi$ on
the product of two copies of the cone with
$\int r_i^2\varphi(x_i)\dd\Pi=\int_\Om\varphi\dd\mu_i$ for all bounded Borel
$\varphi$ and $i=0,1$ \cite[Thm.~2.5, (2.11b)]{LMS23}; the minimum is also
attained among probability measures, see also \cite[Thm.~7.6]{LMS18}. With
$\cos(\min\{|x_0-x_1|,\pi\})$ in place of $\cospi|x_0-x_1|$ this is
\cite[Thm.~2.2, (2.9b)]{LMS23}; the two minima coincide. The two cones coincide when
$\operatorname{diam}\Om\le\pi/2$; otherwise only the cone
with the cutoff $\pi$ is a geodesic space, and that cone is the one meant
whenever geodesics of the cone are used. The cited
theorems are
stated on $\R^d$; the constraint forces $x_i\in\Om$ $\Pi$-almost everywhere
on $\{r_i>0\}$, so the cone over $\Om$ gives the same minimum. This
formulation permits the following comparisons. In the lemma, $\bfT_i$ says
where the mass goes and $q_i^2$ is the factor by which it grows on the way.

\begin{lemma}[Transport--growth, {\cite[(2.17), Cor.~3.5]{LMS23}}]
\label{prelim:tg}
Let $\eta$ be a finite Borel measure on a Polish space $Y$, and let
$\bfT_i:Y\to\Om$ and $q_i:Y\to\R$ be Borel with $q_i\in L^2(\eta)$,
$i=0,1$. If $\mu_i=(\bfT_i)\pfw(q_i^2\eta)$, then
\begin{equation}\label{prelim:eq-tg-upper}
 \HK^2(\mu_0,\mu_1)\le
 \int\big[q_0^2+q_1^2-2q_0q_1\cospi|\bfT_0-\bfT_1|\big]\dd\eta.
\end{equation}
In particular, if $\mu_0=\sum_{i=1}^k\alpha_i$ with $\alpha_i\in\calM(\Om)$,
and $\bfT_i:\Om\to\Om$, $q_i:\Om\to[0,\infty)$ are Borel with
$q_i\in L^2(\alpha_i)$, then
\begin{equation}\label{prelim:eq-tg-split}
 \HK^2\Big(\mu_0,\sum_{i=1}^k(\bfT_i)\pfw(q_i^2\alpha_i)\Big)
 \le\sum_{i=1}^k\int_\Om\big[1+q_i^2-2q_i\cospi|\bfT_i(x)-x|\big]\dd\alpha_i .
\end{equation}
We call such a family $(\alpha_i,\bfT_i,q_i)_{i\le k}$ a \emph{split-source
system}: the mass of $\mu_0$ at a point may be divided into parts that are
moved and rescaled separately, and no injectivity of the maps is needed.
For $\mu_0=\rho_0\Leb$ on $X$, with $\rho_0>0$ a.e., and any
$\mu_1\in\calM(X)$, there are Borel maps $\bfT:X\to X$, $q:X\to[0,\infty)$
satisfying
\begin{equation}\label{prelim:eq-tg-opt}
 \mu_1=\bfT\pfw(q^2\mu_0),\qquad
 \HK^2(\mu_0,\mu_1)=\int_X\big[1+q^2-2q\chi\big]\dd\mu_0,\qquad
 \chi(x):=\cospi|\bfT(x)-x| .
\end{equation}
We call $(\bfT,q)$ an \emph{optimal transport--growth pair} from $\mu_0$
to $\mu_1$, $q$ its \emph{growth factor} and $\chi$ its \emph{transport
cosine}.
\end{lemma}
Bound \eqref{prelim:eq-tg-split} follows from \eqref{prelim:eq-tg-upper},
applied to each piece with $\eta=\alpha_i$, the identity map and the factor
$1$ on the source side, and from the subadditivity
\eqref{prelim:eq-subadd}. For the optimal pair and $\mu_1\ne0$, the no-distant-target
hypothesis of the cited Monge theorem holds because $\supp\mu_0=X$. For $\mu_1=0$, take $q=0$ and $\bfT=\mathrm{id}_X$;
then \eqref{prelim:eq-tg-opt} follows from $\HK^2(\mu_0,0)=\mu_0(X)$.

\subsection{Internal energies and the LMS criterion}\label{prelim:ss-energy}
Two slopes of $E$ enter the later hypotheses: the slope at zero density,
which decides whether the energy is defined on an unbounded domain, and the
slope at infinity, which is the price of singular mass.

Throughout, $E:[0,\infty)\to\R$ is convex, continuous at zero and $E(0)=0$.
Put $\ezero=E'_+(0)=\lim_{c\downarrow0}E(c)/c$,
$\Einf=\lim_{c\to\infty}E(c)/c$ and
$\ezero^-:=\max\{-\ezero,0\}$ when $\ezero$ is finite. When both slopes are
finite, $\omega:=\Einf-\ezero$ is the oscillation of $E'$. Convexity gives
\begin{equation}\label{prelim:eq-Ebounds}
 \ezero c\le E(c)\le\Einf c.
\end{equation}
\begin{definition}[Relaxed internal energy]\label{prelim:energy}
Assume $\Om$ is bounded or $\ezero>-\infty$, and define
\begin{equation}\label{prelim:eq-energy}
 \sfE(\rho\Leb+\mu^s)=\int_\Om E(\rho)\dd x+\Einf\mu^s(\Om),
 \qquad \infty\cdot0:=0.
\end{equation}
\end{definition}
The negative part is integrable: use $E(c)\ge\ezero c$ if $\ezero$ is finite,
and an affine minorant on bounded $\Om$. Energy is unchanged by embedding
its support in a larger domain, and
\begin{equation}\label{prelim:eq-energybounds}
 \ezero\,\mu(\Om)\le\sfE(\mu)\le\Einf\,\mu(\Om),
\end{equation}
where each bound is used only when the corresponding slope is finite.
The bounds are the case $\nu=0$ of the following comparison, which is
used repeatedly for competitors obtained by removing or adding mass.
\begin{lemma}[Dominated measures]\label{prelim:dominated}
Let $\Om$ be bounded or $\ezero>-\infty$, as in
Definition~\ref{prelim:energy}, and let $\nu,\mu\in\calM(\Om)$ with
$\nu\le\mu$. If $\ezero>-\infty$, then
\begin{equation}\label{prelim:eq-dominated}
 \sfE(\nu)+\ezero\,(\mu-\nu)(\Om)\le\sfE(\mu),
\end{equation}
also when $\sfE(\mu)=+\infty$ or $\Einf=+\infty$. If $\Einf<\infty$, then
\begin{equation}\label{prelim:eq-dominated-up}
 \sfE(\mu)\le\sfE(\nu)+\Einf\,(\mu-\nu)(\Om).
\end{equation}
In particular, for a Borel set $S\subset\Om$ and $\ezero>-\infty$,
$\sfE(\mu|_S)+\ezero\,\mu(\Om\setminus S)\le\sfE(\mu)$.
\end{lemma}
\begin{proof}
Since $\nu\le\mu$, the densities and the singular parts are ordered:
$\rho_\nu\le\rho_\mu$ almost everywhere and $\nu^s\le\mu^s$. For finite
$\ezero$ the function $G(c)=E(c)-\ezero c$ is convex and nonnegative with
$G(0)=0$, hence nondecreasing on $[0,\infty)$, and
$\sfE(\mu)=\int_\Om G(\rho_\mu)\dd x+(\Einf-\ezero)\mu^s(\Om)
+\ezero\mu(\Om)$ with nonnegative first two terms; both are monotone in
$\mu$, which gives \eqref{prelim:eq-dominated}. For finite $\Einf$ the
function $H(c)=\Einf c-E(c)$ is nonnegative and nondecreasing, and
$\sfE(\mu)=\Einf\mu(\Om)-\int_\Om H(\rho_\mu)\dd x$, which gives
\eqref{prelim:eq-dominated-up}.
\end{proof}
\begin{lemma}[Lower semicontinuity]\label{prelim:lsc}
The functional \eqref{prelim:eq-energy} is narrowly lower semicontinuous;
equivalently, lower semicontinuous with respect to $\HK$ on $\calM(\Om)$ and
with respect to $\SHK$ on $\calP(\Om)$, since both metrise narrow
convergence.
\end{lemma}
\begin{proof}
For finite $\ezero$, apply \cite[Thm.~2.34]{AFP00} on $\R^d$ to the
nonnegative convex integrand $E(c)-\ezero c$, extended by $+\infty$ for
$c<0$, and the reference measure $\Leb$. Add back the narrowly continuous
mass term. If $\Om$ is bounded and $E(c)\ge ac+b$, use instead
$E(c)-ac-b\ge0$ and the fixed reference measure $\Leb|_\Om$ on $\R^d$;
the subtracted functional is $a\mu(\Om)+b|\Om|$. The cited theorem allows
arbitrary positive reference measures, so no extension of a positive
constant over the exterior is involved.
\end{proof}

The next theorem is the source of all convexity on the cone in this paper;
Section~\ref{sec:examples} verifies its condition for new classes of energies.

\begin{theorem}[LMS criterion, {\cite[Thms.~7.2, 7.5]{LMS23}}]
\label{prelim:LMS}
Let $\Om\subset\R^d$ be closed and convex with nonempty interior, and bounded
if $\ezero=-\infty$. Let $E:[0,\infty)\to\R$ be convex, continuous at $0$,
with $E(0)=0$, let $\lambda\in\R$, and put
\begin{equation}\label{prelim:eq-N}
 \Nfun(\rho,\gamma)=\Big(\frac\rho\gamma\Big)^d
 E\Big(\frac{\gamma^{d+2}}{\rho^d}\Big)-\frac\lambda2\gamma^2,
 \qquad\rho,\gamma>0.
\end{equation}
If
\begin{equation}\label{prelim:eq-LMScond}
 \Nfun\text{ is convex on }(0,\infty)^2\quad\text{and}\quad
 \rho\mapsto(d-1)\Nfun(\rho,\gamma)\text{ is non-increasing for all }\gamma>0,
\end{equation}
then $\sfE$ is strongly geodesically $\lambda$-convex on $(\calM(\Om),\HK)$.
Conversely, if $E\in C^2((0,\infty))$ and $\sfE$ is weakly geodesically
$\lambda$-convex on $(\calM(\Om),\HK)$, then \eqref{prelim:eq-LMScond} holds.
\end{theorem}
That the conclusion is convexity along \emph{every} geodesic, and not only
along one geodesic for each pair of endpoints, is what the proof of
\cite[Thm.~7.2]{LMS23} gives: its Step~2 treats an arbitrary geodesic
between absolutely continuous measures and its Step~4 an arbitrary geodesic.
The converse part is proved in \cite[Thm.~7.5]{LMS23} under strong
convexity; that weak convexity suffices is proved in
Appendix~\ref{nc:sec-lms}, using Theorem~\ref{nc:LMS-pointwise} for $d\ge1$ and, in the degenerate case $d=0$ where $\Om$ is a point, a direct argument.
In \cite{LMS23} the density $E$ may take the value $+\infty$; we use finite
$E$ only, and the unbounded case is available exactly when $\ezero$ is
finite. We say that \emph{$E$ satisfies the Liero--Mielke--Savar\'e
condition}, for short the \emph{LMS condition}, \emph{with constant
$\lambda$} if $E$ is convex, continuous at $0$ with $E(0)=0$, and
\eqref{prelim:eq-LMScond} holds.
\begin{remark}[Normalization]\label{prelim:normalisation}
The mass identity shows that $\sfE$ is $\lambda$-convex exactly when
$\sfE-\lambda\mu(\Om)/2$ is convex. Thus
$N_{\lambda,E}=N_{0,E-\lambda c/2}$, and adding $\beta c$ to $E$ changes the
HK constant by $2\beta$. On $\calP(\Om)$ the same shift adds the constant
$\beta$ to $\sfE$ and changes no convexity constant. The necessity part of
Theorem~\ref{prelim:LMS} is proved for every $\lambda$ directly in
Appendix~\ref{nc:sec-lms}; the differently normalised matrix constants
printed in \cite{LMS23} are not used.
\end{remark}

For twice differentiable $E$ the condition \eqref{prelim:eq-LMScond} is a
condition on a $2\times2$ matrix. We define this matrix here, because
Sections~\ref{sec:transfer} and~\ref{sec:examples} both use it.
\begin{definition}[Matrix form of the criterion]\label{prelim:matrix}
Let $d\ge1$, $E\in C^2((0,\infty))$ and $\lambda\in\R$. For $u>0$ we write
\begin{equation}\label{ex:eq-ARQ}
 A(u)=u^2E''(u),\qquad R(u)=uE'(u)-E(u),\qquad Q(u)=d\,A(u)-(d-1)\,R(u),
\end{equation}
\begin{equation}\label{ex:eq-F}
 F_\lambda(u)=Q(u)\,u\,\big(2E'(u)-\lambda\big)+R(u)\,\big(4A(u)-d\,R(u)\big),
\end{equation}
and we denote by $B_\lambda(u)$ the symmetric matrix
\begin{equation}\label{ex:eq-B}
 B_\lambda(u)=\begin{pmatrix}
 A-\frac{d-1}{d}R & A-\frac12R\\[1mm]
 A-\frac12R & A+\frac14u\,(2E'-\lambda)
 \end{pmatrix},
\end{equation}
in which $A$, $R$ and $E'$ are evaluated at $u$.
\end{definition}
Proposition~\ref{ex:reduced} shows that, for $d\ge1$ and
$E\in C([0,\infty))\cap C^2((0,\infty))$ with $E(0)=0$, the condition
\eqref{prelim:eq-LMScond} holds if and only if $B_\lambda(u)$ is positive
semidefinite for all $u>0$ and, when $d\ge2$, $R\ge0$. Here $Q$ is $d$
times the upper left entry of $B_\lambda$, and $F_\lambda=4d\det B_\lambda$.
For $\lambda=0$ the
matrix $B_0(u)$ is the matrix of \cite[Proposition~6.1]{LMS23}.

\subsection{EVI and minimizing movements}\label{prelim:ss-EVI}
The integrated form of EVI (Proposition~\ref{prelim:EVIprop}) is used when
we pass to limits of flows, the error estimates
(Theorem~\ref{prelim:MS-error}) to show that the scheme converges once a
flow is known to exist.

Let $\phi:\calX\to(-\infty,\infty]$ be proper and lower semicontinuous.
\begin{definition}[EVI solutions and flows]\label{prelim:EVI}
An \emph{$\EVI_\lambda$ solution}, or $\EVI_\lambda$ solution of
$(\calX,\phi,\sfd)$ when the space, the functional and the distance have to
be named, is a continuous curve $u:(0,\infty)\to\dom\phi$ with
\begin{equation}\label{prelim:eq-EVI}
 \frac12\frac{\rmd^+}{\rmd t}\sfd^2(u(t),w)
 +\frac\lambda2\sfd^2(u(t),w)\le\phi(w)-\phi(u(t))
 \quad(w\in\dom\phi),
\end{equation}
where the derivative is the upper right Dini derivative. The point $w$ is
called an \emph{observer}. The solution \emph{starts at $u_0$} if
$u(t)\to u_0$ as $t\downarrow0$, and then we put $u(0)=u_0$. Solutions are
defined for all times; this is what \emph{complete} solution refers to.
An \emph{$\EVI_\lambda$ flow in} (or \emph{on}) a set
$D\subset\overline{\dom\phi}$ is a family of continuous maps $S_t:D\to D$,
$t\ge0$, with $S_0=\mathrm{id}$ and $S_{t+s}=S_tS_s$, such that for every
$u_0\in D$ the curve $t\mapsto S_tu_0$ is an $\EVI_\lambda$ solution
starting at $u_0$. A semigroup $(S_t)_{t\ge0}$ on $D$ is
\emph{$\lambda$-contractive} if
$\sfd(S_tu,S_tv)\le\ee^{-\lambda t}\sfd(u,v)$ for all $u,v\in D$ and
$t\ge0$.
\end{definition}
\begin{proposition}[{\cite[Thms.~3.3, 3.5]{MuS20}}]\label{prelim:EVIprop}
A curve $u:(0,\infty)\to\calX$ for which $t\mapsto\phi(u(t))$ and, for
every $w\in\dom\phi$, $t\mapsto\sfd^2(u(t),w)$ are locally integrable
on $(0,\infty)$ is an $\EVI_\lambda$ solution if and only if, for every
$w\in\dom\phi$,
\begin{equation}\label{prelim:eq-intEVI}
 \tfrac12\sfd^2(u(t),w)-\tfrac12\sfd^2(u(s),w)
 +\int_s^t\!\big[\phi(u(r))+\tfrac\lambda2\sfd^2(u(r),w)\big]\dd r
 \le(t-s)\phi(w)\quad(0<s<t).
\end{equation}
In particular \eqref{prelim:eq-intEVI} implies continuity. Solutions have
non-increasing energy, satisfy $\phi(u(t))\le\phi(u_0)$ if they start at
$u_0$, and for two solutions $u,v$ and $w\in\dom\phi$
\begin{align}
 \sfd(u(t),v(t))&\le\ee^{-\lambda(t-s)}\sfd(u(s),v(s)),\label{prelim:eq-contraction}\\
 \tfrac12\sfd^2(u(t),w)&\le\tfrac12\ee^{-\lambda(t-s)}\sfd^2(u(s),w)
 +M_\lambda(t-s)\big[\phi(w)-\phi(u(t))\big],\label{prelim:eq-expEVI}
\end{align}
for $0<s<t$, and for $s=0$ if the solutions start at $u(0)$, $v(0)$, where
$M_\lambda(h)=\int_0^h\ee^{-\lambda r}\dd r$.
They are unique for prescribed initial data. A solution also satisfies EVI
with every smaller constant.
\end{proposition}
By \eqref{prelim:eq-contraction} every $\EVI_\lambda$ flow is a
$\lambda$-contractive semigroup.
For completeness, the exponential estimate follows by an integrating factor
and positive-time energy monotonicity, then $s\downarrow0$. Taking $w=u_0$
in that estimate gives $\phi(u(t))\le\phi(u_0)$ when the latter is finite;
this conclusion does not follow from lower semicontinuity alone.
The results cited from \cite{MuS20} use the numbering of arXiv:1810.03939v1.

In the limiting arguments below we obtain the integrated inequality first.
We check that, under the hypotheses of Proposition~\ref{prelim:EVIprop},
this inequality places the curve in $\overline{\dom\phi}$ at every positive
time.
The cited characterisation is stated for curves with values in
$\overline{\dom\phi}$. A curve as in Proposition~\ref{prelim:EVIprop}
that satisfies \eqref{prelim:eq-intEVI} has this property. Fix $s>0$ and
$w_*\in\dom\phi$, and put
$p(r)=|\phi(u(r))|+|\lambda|\,\sfd^2(u(r),w_*)$, which is finite for almost
every $r$ and locally integrable. There are times $r_n\uparrow s$ with
$\phi(u(r_n))<\infty$ and $(s-r_n)\,p(r_n)\to0$, since otherwise
$p(r)\ge\delta/(s-r)$ for some $\delta>0$ and almost every $r$ in a left
neighbourhood of $s$. Inequality \eqref{prelim:eq-intEVI} on $[r_n,s]$ with
the observer $u(r_n)$, together with
$\sfd^2(u(r),u(r_n))\le2\sfd^2(u(r),w_*)+2\sfd^2(u(r_n),w_*)$, gives
\[
 \tfrac12\sfd^2(u(s),u(r_n))\le(s-r_n)\,p(r_n)+\int_{r_n}^sp(r)\dd r
 \longrightarrow0 ,
\]
so $u(s)\in\overline{\dom\phi}$.

\begin{definition}[Exact minimizing movements]\label{prelim:MM}
An MM scheme with step $\tau>0$ is
\begin{equation}\label{prelim:eq-MM}
 U^n\in\operatorname*{argmin}_{V\in\calX}
 \{\sfd^2(U^{n-1},V)/(2\tau)+\phi(V)\},\qquad U^0\in\dom\phi.
\end{equation}
We say that \emph{exact steps exist} from $U^0$, for the step $\tau$, if
such a sequence exists, that is, if every minimum in \eqref{prelim:eq-MM}
is attained; the sequence is then also called an \emph{exact scheme}.
Write $\bar U_\tau(t)=U^{\lceil t/\tau\rceil}$ and let $u^\tau$ interpolate
consecutive nodes by geodesics. The metric slope is $+\infty$ outside
$\dom\phi$, zero at isolated points of $\dom\phi$, and otherwise
$\slope\phi(x)=\limsup_{y\to x,\,y\ne x}(\phi(x)-\phi(y))^+/\sfd(x,y)$.
The \emph{Moreau envelope} is
$\phi_\sigma(x)=\inf_y\{\phi(y)+\sfd^2(x,y)/(2\sigma)\}$, $\sigma>0$.
\end{definition}

The first error estimate below contains the envelope, and each time it is
applied one has to know that the envelope is finite and converges to
$\phi$. Both follow from a quadratic lower bound.

\begin{lemma}[Moreau envelope]\label{prelim:envelope}
Let $(\calX,\sfd)$ be a metric space, $\phi:\calX\to(-\infty,\infty]$ proper
and lower semicontinuous, $x\in\dom\phi$, and suppose that
\begin{equation}\label{prelim:eq-quadbound}
 \phi(y)\ge-a-b\,\sfd^2(x,y)\qquad\text{for all }y\in\calX,
\end{equation}
with constants $a,b\ge0$. Then $\phi_\sigma(x)>-\infty$ for
$0<\sigma<1/(2b)$ (for all $\sigma>0$ if $b=0$), and $\phi_\sigma(x)$
increases to $\phi(x)$ as $\sigma$ decreases to $0$.
\end{lemma}
\begin{proof}
\emph{Step 1.} For $0<\sigma<1/(2b)$ and every $y$, by
\eqref{prelim:eq-quadbound},
$\phi(y)+\sfd^2(x,y)/(2\sigma)\ge-a+(1/(2\sigma)-b)\sfd^2(x,y)\ge-a$.
Hence $\phi_\sigma(x)\ge-a$.

\emph{Step 2.} The choice $y=x$ gives $\phi_\sigma(x)\le\phi(x)$, and
$\sigma\mapsto\phi_\sigma(x)$ is non-increasing because
$\sigma\mapsto1/(2\sigma)$ is.

\emph{Step 3.} Let $0<\sigma\le1/(4b)$ (any $\sigma>0$ if $b=0$) and choose
$y_\sigma$ with
$\phi(y_\sigma)+\sfd^2(x,y_\sigma)/(2\sigma)\le\phi_\sigma(x)+\sigma$; this is
possible by Step~1. By Step~2 the right-hand side is at most
$\phi(x)+\sigma$, and by \eqref{prelim:eq-quadbound} the left-hand side is
at least $-a+\sfd^2(x,y_\sigma)/(4\sigma)$. Hence
$\sfd^2(x,y_\sigma)\le4\sigma(\phi(x)+a+\sigma)\to0$ as $\sigma\downarrow0$.

\emph{Step 4.} By the choice of $y_\sigma$,
$\phi_\sigma(x)\ge\phi(y_\sigma)-\sigma$. Lower semicontinuity and Step~3
give $\liminf_{\sigma\downarrow0}\phi_\sigma(x)\ge\phi(x)$. With Step~2 this
proves $\phi_\sigma(x)\uparrow\phi(x)$.
\end{proof}

\begin{theorem}[MM error estimates, {\cite[Thms.~5.5, 5.9, Cor.~5.6]{MuS20}}]
\label{prelim:MS-error}
Let $\lambda\le0$ and $t_\tau=\lceil t/\tau\rceil\tau$.
\begin{enumerate}[label=\textup{(\alph*)},ref=\thetheorem\textup{(\alph*)}]
\item \emph{Exact schemes.}
Let $(\calX,\sfd)$ be a geodesic metric space. Suppose $\phi$ is weakly
geodesically $\lambda$-convex and an $\EVI_\lambda$ flow $(S_t)$
with the same $\lambda$ exists on $D=\overline{\dom\phi}$. Let $u_0\in D$ and
$U^0\in\dom\phi$, let $0<\tau<1$ with $-12\lambda\tau\le1$, and let
$(U^n)$ be an exact scheme from $U^0$ with step $\tau$. For $t>0$ and
$\sigma=\int_0^{3\sqrt{\tau t_\tau}}\ee^{\lambda r}\dd r$,
\begin{align}
 \sfd(S_tu_0,\bar U_\tau(t))
 &\le\ee^{-\lambda t}\sfd(u_0,U^0)
 +10(\tau t_\tau)^{1/4}\ee^{-6\lambda t_\tau}
 \sqrt{\phi(U^0)-\phi_\sigma(U^0)},\label{prelim:eq-MSerror}\\
 \sfd(S_tu_0,\bar U_\tau(t))
 &\le\ee^{-\lambda t}\sfd(u_0,U^0)
 +(\sqrt{\tau t_\tau}+t_\tau-t)\ee^{-3\lambda t_\tau}
 \slope\phi(U^0).\label{prelim:eq-MSslope}
\end{align}
The second estimate is asserted when $\slope\phi(U^0)<\infty$.
At $t=0$ we have the identity
$\sfd(S_0u_0,\overline U_\tau(0))=\sfd(u_0,U^0)$.
For $U^0=u_0\in\dom\phi$, provided exact schemes exist for all small
$\tau$, the first tends to zero locally uniformly in time as
$\tau\downarrow0$. The same is true for
geodesic interpolants.
\item\label{prelim:MS-regular} \emph{Regular data.}
Let $(\calX,\sfd)$ be a metric space, let an $\EVI_\lambda$ flow $(S_t)$
exist in a set $D\subset\overline{\dom\phi}$, and let $\tau>0$ with
$-12\lambda\tau\le1$. Let $(U^n)_{n\ge0}$ be a sequence in $D$ of points of
finite slope such that, for all $n\ge1$,
\begin{equation}\label{prelim:eq-MSreg}
 \tau\slope\phi^2(U^n)\le\frac{\sfd^2(U^{n-1},U^n)}{\tau},\qquad
 \Big(1+\frac{\lambda\tau}2\Big)\frac{\sfd^2(U^{n-1},U^n)}{\tau}
 \le\phi(U^{n-1})-\phi(U^n).
\end{equation}
Then \eqref{prelim:eq-MSslope} holds for every $u_0\in D$ and $t>0$, with
$\bar U_\tau(t)=U^{\lceil t/\tau\rceil}$.
\end{enumerate}
\end{theorem}
Part (a) gives convergence of the scheme once the flow is known to exist
(Theorems~\ref{cp:main}, \ref{hk:mm} and~\ref{shk:mm}); part (b), which
needs no convexity, is used in Section~\ref{sec:sublevel}, where convexity
holds only on sublevels.
Part (b) is \cite[Thm.~5.5]{MuS20} with the parameters $\eta=\eps=0$ there;
it needs neither geodesics nor convexity, and the sequence need not be a
minimizing movement. The second estimate of part (a) is the case of an
exact scheme \cite[Cor.~5.6]{MuS20}.
In part (a), weak geodesic convexity supplies the approximate $\lambda$-convexity
required by the external theorem: points of a geodesic along which $\phi$ is
$\lambda$-convex are exact intermediate points satisfying the convexity
inequality. By Lemma~\ref{prelim:envelope} the Moreau envelope increases to
$\phi(U^0)$ as $\sigma\downarrow0$ as soon as $\phi$ has a quadratic lower
bound; EVI existence supplies one and ensures finiteness at the indicated
scales \cite[Thms.~3.5, 5.9]{MuS20}.
Convergence of geodesic interpolants follows since their distance to grid
nodes is at most the largest step length, and by minimality
$\sfd^2(U^n,U^{n-1})\le2\tau(\phi(U^0)-\phi(U^n))$, where $\phi(U^n)$ is bounded
below on bounded time intervals by the same quadratic lower bound and the
convergence of the nodes. We apply these estimates with
$\min\{\lambda,0\}$ if the original constant is positive. We will also use
the converse direction, from a flow to convexity; this is how weak convexity
is upgraded to strong convexity in Corollary~\ref{cp:weakstrong}.
\begin{theorem}[{\cite[Thm.~3.10(4)]{MuS20}}]\label{prelim:MS-convex}
An $\EVI_\lambda$ flow in a set $D\subset\overline{\dom\phi}$ implies the inequality \eqref{prelim:eq-convex}
along every geodesic contained in $D$ with finite-energy endpoints; applied
to all restrictions of such a geodesic, this is $\lambda$-convexity along it
in the sense of Definition~\ref{prelim:convexity}(a).
\end{theorem}

\subsection{Restriction, normalization and far cuts}\label{prelim:ss-cut}
On the unbounded domains of Sections~\ref{sec:hk-unbdd}--\ref{sec:sublevel}
a measure is replaced by its part in a ball in two situations: as an
observer in the EVI on $\OmR$, where Lemma~\ref{prelim:trunc} bounds the
error in distance and energy, and as a competitor in a minimisation step,
where Lemma~\ref{prelim:split} bounds the change in distance when mass far
from the previous step is cut off; the cut is a gain when the removed mass
exceeds four times the mass of the previous step outside the inner ball.
\begin{lemma}[Restriction to a subdomain]\label{prelim:normalise}
For a closed convex $\Om'\subset\Om$ with nonempty interior, HK and SHK
restrict isometrically to measures on $\Om'$. Internal energies have the
same values there, and EVI on $\Om'$ is EVI on $\Om$ with observers
restricted to $\Om'$. Moreover $\calM(\Om')$ is closed in
$(\calM(\Om),\HK)$ and $\calP(\Om')$ is closed in $(\calP(\Om),\SHK)$.
\end{lemma}
\begin{proof}
An admissible LET plan for measures on $\Om'$ is concentrated on
$\Om'\times\Om'$, so the minimisation problems on the two domains agree,
and SHK is a function of HK by \eqref{prelim:eq-SHK}. The statements on
energies and on EVI then follow directly. For the closedness, HK and SHK
convergence imply narrow convergence, and for the open set
$U=\Om\setminus\Om'$ the map $\mu\mapsto\mu(U)$ is narrowly lower
semicontinuous, so $\mu(U)=0$ passes to the limit. The distance identities
and the closedness use only that $\Om'$ is closed.
\end{proof}

On the cone the mass outside the ball is simply removed. On the sphere mass
cannot be thrown away, so the remainder is renormalised, and the energy
change from the original probability to the renormalised one is controlled
by the oscillation $\omega$ of $E'$.

\begin{lemma}[Truncation]\label{prelim:trunc}
Put $P_R\mu=\mu|_{\bar B_R}$. The distance estimate below holds for every
$\mu$, and the energy estimate holds when $\ezero$ is finite:
\begin{equation}\label{prelim:eq-truncHK}
 \HK^2(\mu,P_R\mu)\le\mu(\Om\setminus\bar B_R),\qquad
 \sfE(P_R\mu)\le\sfE(\mu)+\ezero^-\mu(\Om\setminus\bar B_R).
\end{equation}
For finite $\ezero$, finite-energy measures are recovered in distance
and energy by $P_R$ as $R\to\infty$.
For probabilities, if $\eps=\nu(\Om\setminus\bar B_R)<1$, set
$\tilde P_R\nu=P_R\nu/(1-\eps)$. Then
\begin{equation}\label{prelim:eq-truncSHK}
 \SHK(\nu,\tilde P_R\nu)\le\arcsin\sqrt\eps\le\tfrac\pi2\sqrt\eps.
\end{equation}
If $-\infty<\ezero$ and $\Einf<\infty$, with $\omega=\Einf-\ezero$, then
\begin{equation}\label{prelim:eq-truncEs}
 \sfE(\tilde P_R\nu)\le\sfE(\nu)+\omega\eps.
\end{equation}
\end{lemma}
\begin{proof}
The diagonal LET plan on $\bar B_R$ has cost $\mu(\Om\setminus\bar B_R)$.
Removing mass $\eps$ lowers the energy by at least $\ezero\eps$, and
normalising the remainder adds at most $\Einf\eps$, including the singular
term. This proves both energy bounds. Formula \eqref{prelim:eq-normalise}
with $M=1-\eps$ and $\HK^2(\nu,P_R\nu)\le\eps$ gives
$\cos\SHK(\nu,\tilde P_R\nu)\ge\sqrt{1-\eps}$.
Energy convergence follows by integrability of $E(\rho)$ and continuity
from above of the singular mass.
\end{proof}

Lemma~\ref{prelim:split} splits an optimal plan according to where the
target mass sits. If the target carries mass $a$ farther than $\pi/2$ beyond a ball
containing all but $\eps$ of $\mu_0$, removing that mass lowers the squared
distance by at least $a-2\sqrt{a\eps}$.

\begin{lemma}[Plan splitting and far cut]\label{prelim:split}
Let $\gamma$ be an optimal LET plan for $(\mu_0,\mu)$, with marginals
$\sigma_0\mu_0,\sigma_1\mu$, and let $F$ be Borel. Define
$\gamma''=\gamma|_{\Om\times F}$,
$\alpha=\rmd\gamma''_0/\rmd\gamma_0$, choosing $\alpha=0$ on $\{\sigma_0=0\}$,
and
\begin{equation}\label{prelim:eq-splitdef}
 \mu_0''=\alpha\mu_0,\quad \mu_0'=\mu_0-\mu_0'',\quad
 \mu''=\mu|_F,\quad\mu'=\mu-\mu''.
\end{equation}
Then the measure $\mu_0''$ is concentrated on points at
distance less than $\pi/2$ from $F$, and
\begin{equation}\label{prelim:eq-split}
 \HK^2(\mu_0,\mu)=\HK^2(\mu_0',\mu')+\HK^2(\mu_0'',\mu'').
\end{equation}
If $F=\Om\setminus B_{r+\pi/2}$, $\eps=\mu_0(\Om\setminus B_r)$ and
$a=\mu(F)$, then $\mu_0''(\Om)\le\eps$ and
\begin{equation}\label{prelim:eq-farcut}
 \HK^2(\mu_0,\mu')\le\HK^2(\mu_0,\mu)-a+2\sqrt{a\eps}.
\end{equation}
\end{lemma}
\begin{proof}
The split marginal densities remain $\sigma_0,\sigma_1$, so each term of
LET splits exactly. This gives the lower inequality in
\eqref{prelim:eq-split}; subadditivity gives the reverse.
On $\{\alpha>0\}$, also $\sigma_0>0$, giving $\mu_0''\ll\gamma''_0$.
Finite-cost transport has length less than $\pi/2$, so $\gamma''_0$, and
with it $\mu_0''$, is concentrated on points at distance less than $\pi/2$
from $F$; hence $b:=\mu_0''(\Om)\le\eps$ in the far-cut case. Consequently
$\HK^2(\mu_0'',\mu'')\ge(\sqrt a-\sqrt b)^2$ and
$\HK^2(\mu_0,\mu')\le\HK^2(\mu_0',\mu')+b$; subtraction gives
\eqref{prelim:eq-farcut}. Uncoupled target mass remains in $\mu''$ and
uncoupled source mass in $\mu_0'$, so their entropy costs are retained.
\end{proof}

If the source has no mass within distance $\pi/2$ of $F$, nothing is
transported into $F$ and the target mass there is created from nothing.
This case is used for the uniqueness of geodesics in
Section~\ref{tr:ss-unique}.

\begin{corollary}[Far target mass]\label{prelim:far}
Let $\mu_0,\mu\in\calM(\Om)$, and let $F\subset\Om$ be a Borel set such that
$\mu_0$ is concentrated on $\{x\in\Om:\operatorname{dist}(x,F)\ge\pi/2\}$.
Then
\begin{equation}\label{prelim:eq-far}
 \HK^2(\mu_0,\mu)=\HK^2(\mu_0,\mu|_{\Om\setminus F})+\mu(F).
\end{equation}
\end{corollary}
\begin{proof}
Apply Lemma~\ref{prelim:split} with an optimal plan for $(\mu_0,\mu)$. The
measure $\mu_0''=\alpha\mu_0\le\mu_0$ is concentrated on points at distance
less than $\pi/2$ from $F$, and $\mu_0$ is concentrated on points at
distance at least $\pi/2$. Hence $\mu_0''=0$ and $\mu_0'=\mu_0$. Now
\eqref{prelim:eq-split} reads
$\HK^2(\mu_0,\mu)=\HK^2(\mu_0,\mu|_{\Om\setminus F})+\HK^2(0,\mu|_F)$, and
$\HK^2(0,\mu|_F)=\mu(F)$ by \eqref{prelim:eq-massbound}.
\end{proof}
\section{From convexity on the cone to convexity on the sphere}\label{sec:transfer}

Convexity on the cone $(\calM(\Om),\HK)$ can be checked by the criterion of
\cite{LMS23} (Theorem~\ref{prelim:LMS}); on the sphere $(\calP(\Om),\SHK)$
it was known only for a few negative powers. This section determines how
much of the convexity on the cone survives on the sphere, and shows that on
large domains the restrictions met on the way are necessary.

The idea is planar. An $\SHK$ geodesic of length $\delta$ between two
probability measures, together with all positive multiples of its points,
is by the cone identity \eqref{prelim:eq-HKSHK} an isometric copy of a
planar sector of opening $\delta$: the angle is $\SHK$ arclength and the
radius is the square root of the mass. The $\SHK$ geodesic is the arc of the
unit circle, the $\HK$ geodesic between the same two measures is the chord,
and convexity on the cone is convexity along the chord. The midpoint of the
chord has radius $\cos(\delta/2)$, so to get back to the arc one moves
outwards along a ray by $1-\cos(\delta/2)\approx\delta^2/8$, and the energy
increases by at most about the radial derivative times this distance. This
is of second order in $\delta$: a bound $C$ on the radial derivative turns
the constant $\Lambda$ of the cone into $\Lambda-C$ on the sphere.

The first half of this section is positive: geodesic convexity on the cone
$(\calM(\Om),\HK)$ is transferred to the sphere $(\calP(\Om),\SHK)$. It needs
no compactness, so it is stated on a closed convex domain $\Om$, bounded or
not, under the standing assumptions of Definition~\ref{prelim:energy}. The
second half (Sections~\ref{tr:ss-negpowers}--\ref{tr:ss-rigidity}) concerns
negative powers and obstructions, with $d\ge1$. The uniqueness of geodesics
between far-apart measures and the necessary condition
(Lemma~\ref{tr:geo-unique}, Corollary~\ref{tr:geo-unique-shk} and
Theorem~\ref{tr:necessary}) are again stated on a closed convex domain
$\Om$, bounded or not; for the negative powers the domain is a compact
convex set $X$ with nonempty interior. Each subsection says which case it
treats. Strong and weak geodesic
convexity are used in the sense of Definition~\ref{prelim:convexity}: all
positive results are for strong convexity, and all obstructions already hold
for weak convexity.

\subsection{The planar lemma}

The planar statement involves no measures: $\Phi$ stands for the energy
read in the sector, $g$ for its restriction to the unit arc, and
\eqref{tr:eq-radial-planar} for the bound on the radial derivative.

\begin{lemma}[Planar sector]\label{tr:planar}
Let $0<\delta<\pi$ and
$S_\delta=\{r\ee^{\mathrm i\theta}: r>0,\ \theta\in[0,\delta]\}\subset\R^2$,
which is a convex set. Let $\Lambda,C\in\R$ and let
$\Phi:S_\delta\to(-\infty,+\infty]$ be lower semicontinuous and such that
$x\mapsto\Phi(x)-\frac\Lambda2|x|^2$ is convex on $S_\delta$. Put
$g(\theta)=\Phi(\ee^{\mathrm i\theta})$ and assume that $g<\infty$ on
$[0,\delta]$ and that
\begin{equation}\label{tr:eq-radial-planar}
 \partial_r^+\Phi(\ee^{\mathrm i\theta})
 :=\lim_{h\downarrow0}\frac1h\Big(\Phi\big((1{+}h)\ee^{\mathrm i\theta}\big)
 -\Phi\big(\ee^{\mathrm i\theta}\big)\Big)\le C
 \qquad\text{for all }\theta\in[0,\delta].
\end{equation}
Then $\theta\mapsto g(\theta)-\frac{\Lambda-C}{2}\theta^2$ is convex on
$[0,\delta]$.
\end{lemma}

The limit in \eqref{tr:eq-radial-planar} exists in $[-\infty,+\infty]$, because
$r\mapsto\Phi(r\ee^{\mathrm i\theta})-\frac\Lambda2r^2$ is convex on
$(0,\infty)$ and finite at $r=1$. Lower semicontinuity of $\Phi$ enters the
proof only through $g$, that is, along the unit arc.

\begin{proof}
Let $0\le\theta_1<\theta_2\le\delta$, $\Delta=\theta_2-\theta_1$,
$\theta_m=(\theta_1+\theta_2)/2$ and $r=\cos(\Delta/2)\in(0,1)$; for
$\theta_1=\theta_2$ the inequality \eqref{tr:eq-midpoint} below is trivial. Convexity along the chord gives
\begin{equation}\label{tr:eq-chord}
 \Phi(r\ee^{\mathrm i\theta_m})\le\tfrac12(g(\theta_1)+g(\theta_2))
 -\tfrac\Lambda2\sin^2(\Delta/2).
\end{equation}
For $\psi(a)=\Phi(a\ee^{\mathrm i\theta_m})$, the convex function
$\psi(a)-\Lambda a^2/2$ is finite at $r$ and $1$. Its right derivative at
$1$ is bounded below by that secant and above by $C-\Lambda$, hence finite.
The supporting-line inequality therefore gives
\begin{equation}\label{tr:eq-support}
 g(\theta_m)\le\Phi(r\ee^{\mathrm i\theta_m})+C(1-r)
 -\tfrac\Lambda2(1-r)^2.
\end{equation}
Adding yields, with $K=C-\Lambda$,
\begin{equation}\label{tr:eq-midpoint}
 g(\theta_m)\le\tfrac12(g(\theta_1)+g(\theta_2))
 +K(1-\cos(\Delta/2)).
\end{equation}
If $K\ge0$, the inequality $1-\cos(\Delta/2)\le\Delta^2/8$ makes
$g+K\theta^2/2$ midpoint convex, hence convex by lower semicontinuity and
dyadic approximation. If $K<0$, for every $\eps>0$ one has
$1-\cos(\Delta/2)\ge(1-\eps)\Delta^2/8$ for sufficiently small $\Delta$.
Thus $g+K(1-\eps)\theta^2/2$ is locally midpoint convex and lower
semicontinuous, hence locally convex and therefore convex on the whole
interval. Letting $\eps\downarrow0$ proves the claim also for $K<0$.
\end{proof}

\subsection{The transfer theorem}

The transfer uses from Proposition~\ref{prelim:cone} that $\SHK\le\pi/2$ on
$\calP(\Om)$ and the cone identity in \eqref{prelim:eq-HKSHK}. For a
functional $\sfE:\calM(\Om)\to(-\infty,\infty]$ and $\nu\in\calP(\Om)$ with
$\sfE(\nu)<\infty$ we define the \emph{radial derivative}
\begin{equation}\label{tr:eq-Dr}
 D_r\sfE(\nu):=\lim_{h\downarrow0}\frac1h\Big(\sfE\big((1{+}h)^2\nu\big)
 -\sfE(\nu)\Big)\in[-\infty,+\infty].
\end{equation}
If $\sfE$ is $\Lambda$-convex along every $\HK$ geodesic, the limit exists:
by the cone identity, $\HK(r^2\nu,s^2\nu)=|r-s|$, so every segment of the ray
$r\mapsto r^2\nu$ is a geodesic and $r\mapsto\sfE(r^2\nu)-\frac\Lambda2r^2$ is
convex on $(0,\infty)$.

The transfer of \cite[Prop.~A.1]{LaM22} was formulated for energies that are
homogeneous in the radius. The planar lemma compares the chord with the arc
only to first order in the radial direction, so a one-sided bound on
$D_r\sfE$ is all that is asked of the energy off the sphere, and only along
the geodesic under consideration. This is what allows constants that depend
on the geodesic, as in Theorem~\ref{tr:powers} and
Proposition~\ref{tr:negpowers}.

\begin{theorem}[Transfer of convexity]\label{tr:T1}
Let $\sfE:\calM(\Om)\to(-\infty,\infty]$ be lower semicontinuous with respect to
$\HK$ and strongly geodesically $\Lambda$-convex on $(\calM(\Om),\HK)$ for some
$\Lambda\in\R$. Let $\nu:[0,1]\to\calP(\Om)$ be an $\SHK$ geodesic with
$\sfE(\nu(t))<\infty$ for all $t\in[0,1]$, and let $C\in\R$ be such that
\begin{equation}\label{tr:eq-DrC}
 D_r\sfE(\nu(t))\le C\qquad\text{for all }t\in[0,1].
\end{equation}
Then $\sfE$ is $(\Lambda-C)$-convex along $\nu$, that is,
$t\mapsto\sfE(\nu(t))-\frac{\Lambda-C}{2}\,\SHK^2(\nu(0),\nu(1))\,t^2$ is
convex on $[0,1]$.
\end{theorem}

\begin{proof}
For $\delta=\SHK(\nu_0,\nu_1)>0$, define
$\iota(r\ee^{\mathrm i\theta})=r^2\nu(\theta/\delta)$ on $S_\delta$.
The cone identity and constant speed give
$\HK^2(\iota(r\ee^{\mathrm i\theta}),\iota(s\ee^{\mathrm i\eta}))
=r^2+s^2-2rs\cos(\theta-\eta)$, so $\iota$ is an isometry.
In particular $\iota$ is continuous, so $\Phi=\sfE\circ\iota$ is lower
semicontinuous, and every straight segment in $S_\delta$ maps to an HK
geodesic. By strong convexity, for $a,b\in S_\delta$ with $\Phi(a),\Phi(b)$
finite the function
$s\mapsto\Phi((1{-}s)a+sb)-\frac\Lambda2|a-b|^2s^2$ is convex on $[0,1]$;
segments with an endpoint outside $\dom\Phi$ impose no condition. Since
$|(1{-}s)a+sb|^2-|a-b|^2s^2$ is affine in $s$, this says that
$\Phi-\Lambda|\cdot|^2/2$ is convex on $S_\delta$. Its radial derivative on the unit arc
is $D_r\sfE\le C$. Lemma~\ref{tr:planar} and $\theta=\delta t$ give the
assertion. For $\delta=0$ it is immediate.
\end{proof}

\begin{remark}[Scope]\label{tr:rem-scope}
The proof uses only the cone identity in \eqref{prelim:eq-HKSHK}, the lower
semicontinuity of $\sfE$ and its convexity along $\HK$ geodesics; in
particular it does not use that $\sfE$ is an integral functional. The same
argument therefore gives the theorem on any metric cone, that is, for any
pair of metric spaces related by the cone identity, along every geodesic of
the base whose length is less than $\pi$, the range of openings in
Lemma~\ref{tr:planar}; for $\SHK$ this restriction is void, since
$\SHK\le\pi/2$. The constant $C$ may have either sign, and so may
$C-\Lambda$.
\end{remark}

The same sector relates geodesics on the cone and on the sphere between the
same two rays: if the chord is the only geodesic on the cone, the arc is the
only geodesic on the sphere. We use this for the test geodesics of
Section~\ref{tr:ss-unique}; the sector itself is used again in the proof of
Theorem~\ref{tr:rigidity} in Appendix~\ref{nc:sec-rigidity}.

\begin{lemma}[Unique geodesics on the cone and on the sphere]\label{tr:lift}
Let $\mu_0,\mu_1\in\calM(\Om)$ be nonzero and not proportional, with masses
$m_i=\mu_i(\Om)$, and put $\nu_i=\mu_i/m_i$ and
$\delta=\SHK(\nu_0,\nu_1)\in(0,\pi/2]$. Write the segment
\[
 z(t)=(1-t)\sqrt{m_0}+t\sqrt{m_1}\,\ee^{\mathrm i\delta}
 =|z(t)|\,\ee^{\mathrm i\theta(t)},\qquad t\in[0,1],
\]
in polar coordinates. Then $\theta$ is continuous and increasing from $0$ to
$\delta$, and $|z(1)-z(0)|=\HK(\mu_0,\mu_1)$. If there is only one $\HK$
geodesic $\Gamma$ from $\mu_0$ to $\mu_1$, then there is only one $\SHK$
geodesic $\sigma$ from $\nu_0$ to $\nu_1$, and
\begin{equation}\label{tr:eq-lift}
 \Gamma_t(\Om)=|z(t)|^2,\qquad
 \sigma\big(\theta(t)/\delta\big)=\frac{\Gamma_t}{\Gamma_t(\Om)}
 \qquad(t\in[0,1]).
\end{equation}
\end{lemma}

\begin{proof}
Since $\nu_0\ne\nu_1$ we have $\delta>0$, and $\delta\le\pi/2$ by
Proposition~\ref{prelim:cone}. The segment lies in the convex sector
$S_\delta$. With $w=z(1)-z(0)$, its polar angle has the derivative
$\det(z(t),w)/|z(t)|^2$, and
$\det(z(t),w)=\det(z(0),w)=\sqrt{m_0m_1}\sin\delta>0$. So $\theta$ increases
from $0$ to $\delta$. The cone identity in \eqref{prelim:eq-HKSHK} gives
$|w|^2=m_0+m_1-2\sqrt{m_0m_1}\cos\delta=\HK^2(\mu_0,\mu_1)$.

Let $\sigma$ be an $\SHK$ geodesic from $\nu_0$ to $\nu_1$; one exists
because $(\calP(\Om),\SHK)$ is a geodesic space. As in the proof of
Theorem~\ref{tr:T1}, $\iota(r\ee^{\mathrm i\theta})=r^2\sigma(\theta/\delta)$
is an isometry of $S_\delta$ into $(\calM(\Om),\HK)$. Hence
$t\mapsto\iota(z(t))$ is an $\HK$ geodesic from $\iota(\sqrt{m_0})=\mu_0$ to
$\iota(\sqrt{m_1}\,\ee^{\mathrm i\delta})=\mu_1$, and by uniqueness it is
$\Gamma$:
\[
 \Gamma_t=\iota(z(t))=|z(t)|^2\,\sigma\big(\theta(t)/\delta\big).
\]
The total mass of both sides gives the first identity in
\eqref{tr:eq-lift}, and division by it the second. Since $\theta$ maps
$[0,1]$ onto $[0,\delta]$, the second identity determines $\sigma$ from
$\Gamma$. So $\sigma$ is the only $\SHK$ geodesic from $\nu_0$ to $\nu_1$.
\end{proof}

\begin{remark}[Radial derivative of internal energies]\label{tr:radial}
For an internal energy the radial derivative has an integral formula if
$\sfE(r_0^2\nu)<\infty$ for some $r_0>1$. The example below shows how the
formula can fail, which is why Theorem~\ref{tr:T1} uses the radial
derivative itself. If
$\nu=\rho\,\Leb+\nu^s\in\calP(\Om)$ and $\sfE(r_0^2\nu)<\infty$ for some
$r_0>1$, then $\sfE(\nu)<\infty$ and, by monotone convergence of the
difference quotients of the convex function $E$,
\begin{equation}\label{tr:eq-radial-formula}
 D_r\sfE(\nu)=2\int_{\{\rho>0\}}\rho\,E'_+(\rho)\dd x+2\Einf\,\nu^s(\Om)\in[-\infty,\infty),
\end{equation}
where $E'_+$ is the right derivative. The monotone convergence argument needs the finiteness condition: for
$\Om=[0,1]$, $E(c)=\ee^c-1$ and a density equal to $n$ on disjoint sets of
measure $\ee^{-n}/(4n^3)$, $n\ge1$, complete the density on the complement
by the constant
\[
 c_*:=\frac{1-\sum_{n\ge1}\ee^{-n}/(4n^2)}
                {1-\sum_{n\ge1}\ee^{-n}/(4n^3)}>0.
\]
Both sums are less than $1$, and this gives $\int_0^1\rho\dd x=1$.
The right-hand side of \eqref{tr:eq-radial-formula} is finite while
$\sfE(r^2\nu)=+\infty$ for every $r>1$, so $D_r\sfE(\nu)=+\infty$.
No proof below uses \eqref{tr:eq-radial-formula}.
\end{remark}

The simplest application is to energies with $\Einf<\infty$: the radial
derivative is then bounded on all of $\calP(\Om)$, and the convexity
constant is the same for every geodesic.

\begin{corollary}[Energies with finite recession slope]\label{tr:global}
Let $E:[0,\infty)\to\R$ be convex and continuous at $0$, with $E(0)=0$
and $\Einf<\infty$. Assume $\Om$ is bounded or $\ezero>-\infty$, as in
Definition~\ref{prelim:energy}, and that $\sfE$ is strongly geodesically $\Lambda$-convex on
$(\calM(\Om),\HK)$. Then $\sfE$ is strongly geodesically
$(\Lambda-2\Einf)$-convex on $(\calP(\Om),\SHK)$.
\end{corollary}

\begin{proof}
By \eqref{prelim:eq-energybounds}, finite recession slope implies that
$\sfE\le\Einf$ on probabilities, and $\sfE$ is lower semicontinuous by
Lemma~\ref{prelim:lsc}. Convexity gives
$E(c+a)-E(c)\le\Einf a$, hence
$\sfE(r^2\nu)-\sfE(\nu)\le(r^2-1)\Einf$ for $r>1$, including singular
mass. Thus $D_r\sfE\le2\Einf$ and Theorem~\ref{tr:T1} applies.
\end{proof}

Corollary~\ref{tr:global} is the form in which the transfer is applied to
the new energies of Section~\ref{sec:examples}, through
Corollary~\ref{tr:consequences}, and in the existence theorem on unbounded
domains (Theorem~\ref{shk:main}).

\begin{remark}\label{tr:rem-LaM}
The transfer in \cite[Prop.~A.1]{LaM22} assumes a nonpositive energy that
is homogeneous in the radius with degree at least $1/2$. Write the degree
in the radius $r$ as $p$, so that $\sfE(r^2\nu)=r^p\sfE(\nu)$; for $\sfE_q$
of Proposition~\ref{tr:negpowers}, $p=2q$. For a lower semicontinuous
energy that is strongly geodesically $\Lambda$-convex on the cone and
finite and nonpositive on the sphere, $D_r\sfE=p\sfE\le0$, so
Theorem~\ref{tr:T1} applies with $C=0$ for every degree $p>0$. Homogeneity
itself is unnecessary: the proof of Proposition~\ref{tr:negpowers} uses,
besides strong geodesic $0$-convexity on the cone and lower
semicontinuity, only that $\sfE$ is finite and nonpositive on the sphere
and that $D_r\sfE\le2q\,\sfE$ there. As a check, for the mass functional
$\Lambda=C=2$, and Theorem~\ref{tr:T1} gives the constant $0$,
consistently with constancy on the sphere.
\end{remark}

\subsection{Power functionals are semiconvex on sublevels}

For superlinear energies the radial derivative is unbounded on $\calP(\Om)$ and
Corollary~\ref{tr:global} does not apply; in fact these functionals are not
geodesically $\lambda$-convex on the sphere for any $\lambda$ when the
domain has diameter larger than $\pi/2$
(Theorem~\ref{tr:necessary}). They are, however, semiconvex
along every geodesic with a constant that depends only on the energies of the
endpoints.

Sublevels appear because for a power energy $D_r\sfE\le2p\,\sfE$, so
Theorem~\ref{tr:T1} applies along a geodesic as soon as the energy is
bounded along it. That bound comes once more from the chord: there the
energy is at most the larger endpoint value $c$, and the arc is reached by
multiplying the measure by a factor at most $2$, which multiplies the energy
by at most $2^p$.

\begin{theorem}[Powers]\label{tr:powers}
Let $E(c)=\sum_{i=1}^n\alpha_ic^{p_i}$ with $n\ge1$, $\alpha_i>0$ and $p_i>1$, let
$p=\max_ip_i$, and let $\sfE(\mu)=\int_\Om E(\rho)\dd x$ for $\mu=\rho\,\Leb$ and
$\sfE(\mu)=+\infty$ if $\mu$ is not absolutely continuous. Let
$\nu_0,\nu_1\in\calP(\Om)$ with
$c:=\max\{\sfE(\nu_0),\sfE(\nu_1)\}<\infty$. Then for every $\SHK$ geodesic
$\nu:[0,1]\to\calP(\Om)$ from $\nu_0$ to $\nu_1$:
\begin{enumerate}[label=\textup{(\roman*)}]
\item $\sfE(\nu(t))\le2^p\,c$ for all $t\in[0,1]$;
\item $\sfE$ is $\lambda_c$-convex along $\nu$ with
 $\lambda_c=-p\,2^{p+1}c$.
\end{enumerate}
In particular, for a single power $E(c)=c^m$, $m>1$, the functional is
$\lambda$-convex along $\nu$ with
$\lambda=-m\,2^{m+1}\max\{\sfE(\nu_0),\sfE(\nu_1)\}$.
\end{theorem}

\begin{proof}
Each power energy satisfies the condition of Theorem~\ref{prelim:LMS} with
$\lambda=0$ \cite[Cor.~7.3]{LMS23}, so $\sfE$ is strongly geodesically
$0$-convex on $(\calM(\Om),\HK)$; it is lower semicontinuous by
Lemma~\ref{prelim:lsc} ($\ezero=0$). In the sector
of Theorem~\ref{tr:T1}, the chord between unit endpoints meets the ray of
angle $\theta$ at
$r=\cos(\delta/2)/\cos(\theta-\delta/2)\in[2^{-1/2},1]$.
The energy on the chord is at most $c$. Since $E(az)\le a^pE(z)$ for
$a\ge1$, the energy on the unit arc is at most $2^pc$.
Furthermore $D_r\sfE(\nu)=2\sum_i p_i\alpha_i\int\rho^{p_i}
\le2p\sfE(\nu)$. By (i) the energy is finite along the geodesic and
$D_r\sfE\le2p\,2^pc$ there, so Theorem~\ref{tr:T1} with $\Lambda=0$ gives
$\lambda_c=-p2^{p+1}c$.
\end{proof}

Theorem~\ref{tr:powers} is the convexity input of the sublevel theory of
Section~\ref{sec:sublevel}. Part (i) is needed there because the geodesic
between two points of a sublevel may leave that sublevel; it bounds the
energy along the geodesic.

\subsection{The negative powers: uniform convexity on sublevels}\label{tr:ss-negpowers}

In this subsection $d\ge1$ and the domain is a compact convex set $X$ with
nonempty interior: for $q<1$ the energy $-\int\rho^q$ has $\ezero=-\infty$ and is not
bounded below on measures of fixed mass on an unbounded domain.
The result of this subsection sharpens the convexity found in \cite{LaM22}
and gives the implication $(3)\Rightarrow(2)$ of
Theorem~\ref{tr:rigidity} below.

For $q\in(0,1)$ let $\sfE_q(\mu)=-\int_X\rho^q\dd x$ for
$\mu=\rho\,\Leb+\mu^s$; here $\Einf=0$. By \cite[Ex.~6.4, Cor.~7.3]{LMS23}
(see Theorem~\ref{prelim:LMS}), $\sfE_q$ is strongly geodesically $0$-convex
on $(\calM(X),\HK)$ if and only if
\begin{equation}\label{tr:eq-qrange}
 d=1 \text{ and } \tfrac13\le q\le\tfrac12,\qquad\text{or}\qquad
 d=2 \text{ and } q=\tfrac12 .
\end{equation}
In the admissible cases, $0$ is the largest convexity constant on the cone.
With the matrix
$B_\lambda$ of Definition~\ref{prelim:matrix} one has, for $E(c)=-c^q$,
$B_\lambda(c)=c^qB_0(1)-\frac{\lambda c}{4}\operatorname{diag}(0,1)$. Outside
\eqref{tr:eq-qrange} the matrix $B_0(1)$ has a negative direction, which
dominates as $c\downarrow0$ for every $\lambda\in\R$; inside
\eqref{tr:eq-qrange} the term $-\lambda c/4$ dominates as $c\to\infty$ if
$\lambda>0$. Hence, by the matrix form of the criterion
(Proposition~\ref{ex:reduced} below, whose proof is independent of this
section) and Theorem~\ref{prelim:LMS}, $\sfE_q$ is strongly geodesically
$\lambda$-convex on $(\calM(X),\HK)$ if and only if $\lambda\le0$ and
\eqref{tr:eq-qrange} holds.
In \cite{LaM22} it was shown that in these cases $\sfE_q$ is geodesically
convex on $(\calP(X),\SHK)$ (stated there in the weak sense; the argument
of \cite[Prop.~A.1]{LaM22} applies to every geodesic). Theorem~\ref{tr:T1} gives a quantitative
improvement. The radial derivative $2q\,\sfE_q$ is nonpositive, and negative
at measures that are not singular, so moving from
the chord out to the arc lowers the energy, and along a geodesic whose
endpoints are both not singular the constant on the sphere
is better than the constant $0$ on the cone. On all of $\calP(X)$ the
constant $0$ cannot be improved: Dirac masses moving along a short segment
form an $\SHK$ geodesic on which $\sfE_q$ vanishes.

\begin{proposition}\label{tr:negpowers}
Let $d$ and $q$ satisfy \eqref{tr:eq-qrange}, and let $\nu:[0,1]\to\calP(X)$
be an $\SHK$ geodesic from $\nu_0$ to $\nu_1$. Then $\sfE_q$ is
$\lambda$-convex along $\nu$ with
\[
 \lambda=-2q\max\{\sfE_q(\nu_0),\sfE_q(\nu_1)\}
 =2q\min\Big\{\int_X\rho_0^q\dd x,\int_X\rho_1^q\dd x\Big\}\ge0 .
\]
In particular, $\sfE_q$ is strongly geodesically $0$-convex on
$(\calP(X),\SHK)$, and for every $c>0$ it is $2qc$-convex along every geodesic
joining two points of the sublevel $\{\sfE_q\le-c\}\cap\calP(X)$; such
geodesics stay in this sublevel.
\end{proposition}

\begin{proof}
By H\"older's inequality $0\ge\sfE_q(\mu)\ge-|X|^{1-q}\mu(X)^q$, so $\sfE_q$
is finite; it is lower semicontinuous (Lemma~\ref{prelim:lsc}) and
$\sfE_q(r^2\nu)=r^{2q}\sfE_q(\nu)$. Hence
$D_r\sfE_q(\nu)=2q\,\sfE_q(\nu)\le0$ for every $\nu\in\calP(X)$.
Theorem~\ref{tr:T1} with $\Lambda=0$ and $C=0$ shows that
$f(t)=\sfE_q(\nu(t))$ is convex, so
$\sup_tf=\max\{f(0),f(1)\}=:-m\le0$. Therefore
$D_r\sfE_q(\nu(t))\le-2qm$ for all $t$, and Theorem~\ref{tr:T1} with
$\Lambda=0$ and $C=-2qm$ gives $\lambda=2qm$.
\end{proof}

\subsection{Uniqueness of geodesics between far-apart measures}\label{tr:ss-unique}

So far everything was governed by the radial derivative. Whether an energy
with unbounded radial derivative could be convex on the sphere for some
other reason cannot be decided from the cone: we have to exhibit geodesics
on the sphere along which the convexity inequality fails.

In this subsection and the next, $d\ge1$ and $\Om\subset\R^d$ is closed
and convex with nonempty interior, bounded or not. No compactness is used:
the proofs rest on Corollary~\ref{prelim:far}, the cone formulation
\eqref{prelim:eq-coneOT}, the scaling identity, the existence of geodesics
and Lemma~\ref{tr:lift}, all of which hold on such a domain. The results
are used on compact domains in Sections~\ref{sec:examples}
and~\ref{sec:compact} and on unbounded ones in
Sections~\ref{sec:shk-unbdd} and~\ref{sec:sublevel}.

The remaining results of this section are obstructions, and they concern weak
geodesic convexity, where convexity is required along one geodesic only. They
therefore need test geodesics that are unique. The first family consists of
the geodesics between measures whose supports are more than $\pi/2$ apart:
there is no transport, and the geodesic is a pure reaction path.

We use two representations of $\HK$ from Section~\ref{prelim:ss-HK}: the
logarithmic entropy--transport formulation (Definition~\ref{prelim:LET}),
through Corollary~\ref{prelim:far}, and the transport formulation on the
cone \eqref{prelim:eq-coneOT}.

\begin{lemma}[Far-apart supports]\label{tr:geo-unique}
Let $\mu_0,\mu_1\in\calM(\Om)\setminus\{0\}$ with supports $S_0,S_1$ and
\[
 L:=\operatorname{dist}(S_0,S_1)>\pi/2 .
\]
Then $\HK^2(\mu_0,\mu_1)=\mu_0(\Om)+\mu_1(\Om)$, and the only $\HK$ geodesic from
$\mu_0$ to $\mu_1$ is
\begin{equation}\label{tr:eq-reaction}
 \mu(s)=(1{-}s)^2\mu_0+s^2\mu_1,\qquad s\in[0,1].
\end{equation}
No absolute continuity of $\mu_0$, $\mu_1$ is assumed.
\end{lemma}

For the uniqueness we take a point $\mu$ of an arbitrary geodesic and bound
its distances to the two endpoints from below, by the Cauchy--Schwarz
inequality in the cone formulation (Step~3). These distances are known, and
comparing them with the lower bounds forces equality throughout (Step~4).
By the trigonometric inequality of Step~2, equality puts $\mu$ on
$S_0\cup S_1$, and its two parts are multiples of $\mu_0$ and $\mu_1$
(Step~5).

\begin{proof}
Write $m_i=\mu_i(\Om)>0$ and
$c_i(x)=\cospi(\operatorname{dist}(x,S_i))\in[0,1]$, $i=0,1$.

\emph{Step 1 (distance).} The measure $\mu_0$ is concentrated on $S_0$,
whose points have distance at least $L>\pi/2$ from $S_1$, and $\mu_1$ is
concentrated on $S_1$. Corollary~\ref{prelim:far} with $F=S_1$ and
\eqref{prelim:eq-massbound} give
$\HK^2(\mu_0,\mu_1)=\HK^2(\mu_0,0)+\mu_1(S_1)=m_0+m_1=:D^2$.

\emph{Step 2 (a trigonometric inequality).} We claim that
\begin{equation}\label{tr:eq-cc}
 c_0^2+c_1^2\le1\ \text{ on }\Om,\qquad\text{with equality exactly on }
 S_0\cup S_1 .
\end{equation}
On $S_0$ we have $c_0=1$ and $c_1=0$, and symmetrically on $S_1$. If one of
the two numbers $c_0(x),c_1(x)$ vanishes, the sum is at most $1$, and equality
requires the other one to be $1$, that is, $x\in S_0\cup S_1$ since the
supports are closed. If both are positive, then
$a:=\operatorname{dist}(x,S_0)$ and $b:=\operatorname{dist}(x,S_1)$ satisfy
$0\le a,b<\pi/2$ and, by the triangle inequality, $\pi/2<L\le a+b<\pi$.
Thus $\cos(a+b)<0$ and $\cos(a-b)>0$, and
\[
 \cos^2a+\cos^2b=1+\cos(a+b)\cos(a-b)<1 .
\]
(If $L\ge\pi$ this case does not occur.)

\emph{Step 3 (a lower bound for the distance to an arbitrary measure).} For
every $\mu\in\calM(\Om)$, with $m=\mu(\Om)$, and $i\in\{0,1\}$,
\begin{equation}\label{tr:eq-lower}
 \HK^2(\mu_i,\mu)\ge m_i+m-2\sqrt{m_i}\,\Big(\int_\Om c_i^2\dd\mu\Big)^{1/2}.
\end{equation}
Indeed, let $\Pi$ be admissible in \eqref{prelim:eq-coneOT} for the pair
$(\mu_i,\mu)$, with coordinates $([x_i,r_i],[x,r])$. The constraint for
$\mu_i$ gives $\int_{\{x_i\notin S_i\}}r_i^2\dd\Pi=0$, so $x_i\in S_i$
$\Pi$-almost everywhere on $\{r_i>0\}$, and there
$\cospi(|x_i-x|)\le c_i(x)$. By the Cauchy--Schwarz inequality,
\begin{align*}
 \int r_ir\cospi(|x_i-x|)\dd\Pi&\le\int r_irc_i(x)\dd\Pi\\
 &\le\Big(\int r_i^2\dd\Pi\Big)^{1/2}\Big(\int r^2c_i(x)^2\dd\Pi\Big)^{1/2}
 =\sqrt{m_i}\,\Big(\int_\Om c_i^2\dd\mu\Big)^{1/2}.
\end{align*}
Since $\int r_i^2\dd\Pi=m_i$ and $\int r^2\dd\Pi=m$, the cost of $\Pi$ is
bounded from below by the right-hand side of \eqref{tr:eq-lower}; taking the
minimum over $\Pi$ proves \eqref{tr:eq-lower}.

\emph{Step 4 (equality).} Let $(\mu(s))_{s\in[0,1]}$ be an arbitrary $\HK$
geodesic from $\mu_0$ to $\mu_1$, fix $s\in(0,1)$, and write $\mu=\mu(s)$,
$m=\mu(\Om)$,
\[
 p=\int_\Om c_0^2\dd\mu,\qquad r=\int_\Om c_1^2\dd\mu,\qquad
 \Theta:=m-p-r\ge0,
\]
the inequality being \eqref{tr:eq-cc}. By \eqref{tr:eq-lower} and
$\HK(\mu_0,\mu)=sD$, $\HK(\mu,\mu_1)=(1{-}s)D$,
\begin{equation}\label{tr:eq-two}
 m_0+m-2\sqrt{m_0p}\le s^2(m_0+m_1),\qquad
 m_1+m-2\sqrt{m_1r}\le(1{-}s)^2(m_0+m_1).
\end{equation}
Multiplying the first inequality by $1-s$, the second by $s$, and adding, we
obtain after rearrangement
\[
 \big(\sqrt p-(1{-}s)\sqrt{m_0}\big)^2+\big(\sqrt r-s\sqrt{m_1}\big)^2
 +\Theta\le0 .
\]
All three terms are nonnegative, hence
\begin{equation}\label{tr:eq-forced}
 p=(1{-}s)^2m_0,\qquad r=s^2m_1,\qquad m=p+r,
\end{equation}
and both inequalities in \eqref{tr:eq-two} are equalities.

\emph{Step 5 (identification).} By \eqref{tr:eq-forced},
$\int_\Om(1-c_0^2-c_1^2)\dd\mu=0$. The integrand is nonnegative and, by
\eqref{tr:eq-cc}, positive outside $S_0\cup S_1$; hence $\mu$ is concentrated
on $S_0\cup S_1$. Put $\alpha=\mu|_{S_0}$ and $\beta=\mu|_{S_1}$, so that
$\mu=\alpha+\beta$, $\alpha(\Om)=p>0$ and $\beta(\Om)=r>0$.
As in Step~1, Corollary~\ref{prelim:far} with $F=S_1$, now applied to the
pair $(\mu_0,\mu)$, gives
\begin{equation}\label{tr:eq-split}
 \HK^2(\mu_0,\alpha+\beta)=\HK^2(\mu_0,\alpha)+r ,
\end{equation}
because $\mu|_{\Om\setminus S_1}=\alpha$ and $\mu(S_1)=r$. From the first equality in
\eqref{tr:eq-two}, \eqref{tr:eq-forced} and \eqref{tr:eq-split},
\[
 \HK^2(\mu_0,\alpha)=m_0+p-2\sqrt{m_0p}=\big(\sqrt{m_0}-\sqrt p\big)^2 .
\]
On the other hand, the scaling property (Proposition~\ref{prelim:scaling})
gives
\[
 \HK^2(\mu_0,\alpha)=\sqrt{m_0p}\;\HK^2\Big(\frac{\mu_0}{m_0},\frac\alpha p\Big)
 +\big(\sqrt{m_0}-\sqrt p\big)^2 .
\]
Since $m_0p>0$, the two probability measures coincide, that is,
$\alpha=\frac p{m_0}\mu_0=(1{-}s)^2\mu_0$. The same argument at the other
endpoint gives $\beta=s^2\mu_1$. Hence $\mu(s)$ is given by
\eqref{tr:eq-reaction} for every $s\in(0,1)$. Since $(\calM(\Om),\HK)$ is a
geodesic space, a geodesic from $\mu_0$ to $\mu_1$ exists, and it is
\eqref{tr:eq-reaction}.
\end{proof}

\begin{corollary}\label{tr:geo-unique-shk}
Let $\nu_0,\nu_1\in\calP(\Om)$ have supports at distance larger than $\pi/2$.
Then we have $\SHK(\nu_0,\nu_1)=\pi/2$, and the only $\SHK$ geodesic from $\nu_0$ to
$\nu_1$ is, parametrised by arclength $\theta=\frac\pi2t\in[0,\pi/2]$,
\begin{equation}\label{tr:eq-reaction-shk}
 \nu_\theta=\cos^2\theta\;\nu_0+\sin^2\theta\;\nu_1 .
\end{equation}
\end{corollary}

\begin{proof}
By Lemma~\ref{tr:geo-unique}, $\HK^2(\nu_0,\nu_1)=2$, the cone identity in
\eqref{prelim:eq-HKSHK} gives $\cos\SHK(\nu_0,\nu_1)=0$, and the only $\HK$
geodesic from $\nu_0$ to $\nu_1$ is $\Gamma_s=(1{-}s)^2\nu_0+s^2\nu_1$. The
two measures have disjoint supports, so they are not proportional.
Lemma~\ref{tr:lift} applies with $\delta=\pi/2$ and $z(s)=(1{-}s,s)$: the
$\SHK$ geodesic $\sigma$ is unique, $|z(s)|^2=(1{-}s)^2+s^2$,
$\cos^2\theta(s)=(1{-}s)^2/|z(s)|^2$, and by \eqref{tr:eq-lift}
\[
 \sigma\big(2\theta(s)/\pi\big)=\frac{(1{-}s)^2\nu_0+s^2\nu_1}{(1{-}s)^2+s^2}
 =\cos^2\theta(s)\,\nu_0+\sin^2\theta(s)\,\nu_1 .
\]
As $\theta(s)$ runs through $[0,\pi/2]$, this is \eqref{tr:eq-reaction-shk}.
\end{proof}

\begin{remark}[The strict inequality is needed]\label{tr:rem-strict}
In Lemma~\ref{tr:geo-unique} the hypothesis $L>\pi/2$ cannot be replaced by
$L\ge\pi/2$. If $x_0,x_1\in\Om$ with $|x_0-x_1|=\pi/2$, then both
\eqref{tr:eq-reaction-shk} with $\nu_i=\delta_{x_i}$ and the moving Dirac
measure $\theta\mapsto\delta_{x(\theta)}$, where $x(\theta)$ traverses the
segment from $x_0$ to $x_1$ with unit speed, are $\SHK$ geodesics from
$\delta_{x_0}$ to $\delta_{x_1}$; the second one is a geodesic because the
$\SHK$ distance of two unit Dirac measures at distance at most $\pi/2$ equals
the distance of the points. On the other hand, no upper bound on $L$ is
needed.
\end{remark}

\subsection{A necessary condition on large domains}

We test the energy along the pure-reaction geodesic
\eqref{tr:eq-reaction-shk} from a concentrated measure $\nu_0$, with density
$h$ on a small ball, to a fixed spread-out measure $\nu_1$ far away. At the
start a fraction $\sin^2\theta$ of the mass is removed where the density is
$h$ and the energy responds with slope $E'_-(h)$, and the same mass is
created at density close to zero, where the slope is $\ezero$. The energy
therefore starts like $(\ezero-E'_-(h))\theta^2$, with no first-order term,
and $\lambda$-convexity bounds it from below by $\frac\lambda2\theta^2$.
Since $h$ can be taken arbitrarily large, $E'$ has to be bounded.

\begin{theorem}[Necessary condition]\label{tr:necessary}
Assume $\diam\Om>\pi/2$. Let $E:[0,\infty)\to\R$ be convex with $E(0)=0$ and
finite right derivative $\ezero=E'_+(0)$ at $0$. If $\sfE$ is weakly
geodesically $\lambda$-convex on $(\calP(\Om),\SHK)$ for some $\lambda\in\R$,
then $\Einf<\infty$, that is, $E'$ is bounded, and
\begin{equation}\label{tr:eq-necessary}
 \lambda\le2\,(\ezero-\Einf)\le0 .
\end{equation}
In particular, $\lambda=0$ is possible only for $E(c)=\ezero c$.
\end{theorem}

Here boundedness of $E'$ means boundedness of the one-sided derivatives; no
differentiability is assumed, and
$\Einf=\lim_{h\to\infty}E'_-(h)=\lim_{h\to\infty}E(h)/h$.

\begin{proof}
\emph{Step 1 (test measures).} Since $\diam\Om>\pi/2$ and the interior of $\Om$
is dense in $\Om$, there are closed balls $U,V$ contained in the interior of $\Om$
with $\operatorname{dist}(U,V)>\pi/2$. Put $k=1/|V|$, and for $h>1/|U|$ let
$A_h\subset U$ be a ball of volume $1/h$. The probability measures
\[
 \nu_0=h\,\mathbf 1_{A_h}\Leb,\qquad\nu_1=k\,\mathbf 1_V\Leb
\]
have finite energy and supports at distance larger than $\pi/2$. By
Corollary~\ref{tr:geo-unique-shk} the only $\SHK$ geodesic between them is
\eqref{tr:eq-reaction-shk}, and since $E(0)=0$ and the supports are disjoint,
its energy, which is defined on $\calP(\Om)$ by
Definition~\ref{prelim:energy} because $\ezero$ is finite, is
\begin{equation}\label{tr:eq-profile}
 g(\theta):=\sfE(\nu_\theta)=\frac{E(h\cos^2\theta)}h+\frac{E(k\sin^2\theta)}k,
 \qquad\theta\in[0,\pi/2].
\end{equation}

\emph{Step 2 (expansion at $\theta=0$).} Since $E$ is convex and finite, it
has a finite left derivative at $h$, and by assumption a finite right
derivative at $0$. As $h\cos^2\theta=h-h\sin^2\theta<h$, we get for
$\theta\downarrow0$
\[
 \frac{E(h\cos^2\theta)}h=\frac{E(h)}h-E'_-(h)\sin^2\theta+o(\theta^2),
 \qquad\frac{E(k\sin^2\theta)}k=\ezero\sin^2\theta+o(\theta^2),
\]
and therefore
\begin{equation}\label{tr:eq-expansion}
 g(\theta)=g(0)+\big(\ezero-E'_-(h)\big)\theta^2+o(\theta^2).
\end{equation}

\emph{Step 3 (conclusion).} By weak $\lambda$-convexity there is a geodesic
from $\nu_0$ to $\nu_1$ along which $\sfE$ is $\lambda$-convex; by uniqueness
it is \eqref{tr:eq-reaction-shk}. Since $\theta$ is arclength,
$G(\theta):=g(\theta)-\frac\lambda2\theta^2$ is convex on $[0,\pi/2]$. By
\eqref{tr:eq-expansion} its right derivative at $0$ is $0$, so
$G(\theta)\ge G(0)$, that is, $g(\theta)\ge g(0)+\frac\lambda2\theta^2$.
Dividing by $\theta^2$ and letting $\theta\downarrow0$ in
\eqref{tr:eq-expansion} we obtain
\[
 \lambda\le2\big(\ezero-E'_-(h)\big)\qquad\text{for all }h>1/|U| .
\]
Letting $h\to\infty$ gives $\Einf\le\ezero-\lambda/2<\infty$ and
\eqref{tr:eq-necessary}. The slopes of $E$ are nondecreasing and lie between
$\ezero$ and $\Einf$; if $\Einf=\ezero$, then $E$ is linear.
\end{proof}

In the rest of this section $X$ again denotes a compact convex set with
nonempty interior, and the theorem is applied with $\Om=X$.

Together with the criterion on the cone and the transfer, the theorem
gives the form in which all examples of Section~\ref{sec:examples} are
obtained: the condition of Theorem~\ref{prelim:LMS} is checked for $E$ with
some constant, and the constants on the cone and on the sphere follow.

\begin{corollary}[From the LMS condition to both spaces]\label{tr:consequences}
Let $E$ satisfy the LMS condition with constant $\lambda$
$($Theorem~\textup{\ref{prelim:LMS}}$)$, and let $\Om\subset\R^d$ be closed
and convex with nonempty interior.
\begin{enumerate}[label=\textup{(\roman*)}]
\item If $\ezero>-\infty$ or $\Om$ is bounded, then $\sfE$ is strongly
 geodesically $\lambda$-convex on $(\calM(\Om),\HK)$.
\item If moreover $\Einf<\infty$, then $\sfE$ is strongly geodesically
 $(\lambda-2\Einf)$-convex on $(\calP(\Om),\SHK)$.
\item If moreover $\lambda=2\ezero$ and $\diam\Om>\pi/2$,
 then the constant $2(\ezero-\Einf)$ in \textup{(ii)} cannot be improved,
 even for weak geodesic convexity.
\end{enumerate}
\end{corollary}

\begin{proof}
(i) is Theorem~\ref{prelim:LMS}, and (ii) follows from (i) and
Corollary~\ref{tr:global} with $\Lambda=\lambda$. For (iii), $\ezero$ is
finite, and by Theorem~\ref{tr:necessary} every $\lambda'$ for which $\sfE$
is weakly geodesically $\lambda'$-convex on $(\calP(\Om),\SHK)$ satisfies
$\lambda'\le2(\ezero-\Einf)=\lambda-2\Einf$.
\end{proof}

By Remark~\ref{prelim:normalisation}, the hypothesis $\lambda=2\ezero$ in
(iii) says that $E-\ezero c$ satisfies the LMS condition with constant $0$.
Then the transferred constant $\lambda-2\Einf$ equals the upper bound
$2(\ezero-\Einf)$ of Theorem~\ref{tr:necessary}.
In general, for $\diam X>\pi/2$, finite $\ezero$ and $\Einf$, and $\sfE$
strongly geodesically $\Lambda$-convex on $(\calM(X),\HK)$, the best constant
on $(\calP(X),\SHK)$, weak or strong, lies between $\Lambda-2\Einf$ and
$2(\ezero-\Einf)$.

\begin{remark}[Energies with unbounded derivative]\label{tr:rem-consequences}
Let $\diam X>\pi/2$.

(a) \emph{Positive powers.} For $E(c)=c^m$, $m>1$, we have $\ezero=0$ and
$\Einf=\infty$. Hence $\int\rho^m\dd x$ is not weakly geodesically
$\lambda$-convex on $(\calP(X),\SHK)$ for any $\lambda\in\R$, and the
dependence of the constant on the sublevel in Theorem~\ref{tr:powers} cannot
be avoided. In the example of the proof, $\lambda\le-2E'(h)=-2m h^{m-1}
=-2m\,\sfE(\nu_0)$, which shows that the linear growth of $-\lambda_c$ in $c$
in Theorem~\ref{tr:powers} is of the correct order.

(b) \emph{Entropy.} For $E(c)=c\log c$ the slope at $0$ is $-\infty$ and
Theorem~\ref{tr:necessary} does not apply, but the same test geodesics do:
\eqref{tr:eq-profile} gives
\[
 g(\theta)-g(0)=\theta^2\big(2\log\theta+\log k-\log h-1\big)+o(\theta^2),
\]
so the right derivative of $g-\frac\lambda2\theta^2$ at $0$ is $0$ while
$(g(\theta)-g(0))/\theta^2\to-\infty$. As in Step~3, this is incompatible with
$\lambda$-convexity for every $\lambda\in\R$.

(c) \emph{Negative powers.} For $E(c)=-c^q$, $0<q<1$, \eqref{tr:eq-profile}
gives $g(\theta)-g(0)=qh^{q-1}\theta^2-k^{q-1}\theta^{2q}+o(\theta^2)$. If
$\frac12<q<1$, the right derivative at $0$ is again $0$ and
$(g(\theta)-g(0))/\theta^2\to-\infty$, which excludes every
$\lambda\in\R$. For $q\le\frac12$ these test geodesics give no obstruction;
the complete answer is given in Theorem~\ref{tr:rigidity}.

No statement for $\diam X\le\pi/2$ follows from these arguments.
\end{remark}

\subsection{Rigidity for negative powers}\label{tr:ss-rigidity}

We finally show that on large domains the range \eqref{tr:eq-qrange} of
exponents, for which \cite{LaM22} established geodesic convexity of
$\sfE_q=-\int\rho^q\dd x$ on the sphere, is exactly the full range, even if an
arbitrary negative convexity constant and only weak convexity are asked for.

The test family of Theorem~\ref{tr:necessary} does not suffice for this:
along a pure-reaction geodesic the mass created at zero density lowers the
energy by a term of order $\theta^{2q}$, which contradicts convexity only
for $q>\frac12$ (Remark~\ref{tr:rem-consequences}(c)). The remaining
exponents need test measures between which no mass is created from
nothing, so that the midpoint defect of the energy can be compared with
$\lambda$.

\begin{theorem}[Rigidity for negative powers]\label{tr:rigidity}
Let $X\subset\R^d$, $d\ge1$, be compact and convex with nonempty interior.
Assume $\diam X>\pi/2$ and let $0<q<1$. The following are equivalent:
\begin{enumerate}[label=\textup{(\arabic*)}]
\item $\sfE_q$ is weakly geodesically $\lambda$-convex on $(\calP(X),\SHK)$
 for some $\lambda\in\R$;
\item $\sfE_q$ is strongly geodesically $0$-convex on $(\calP(X),\SHK)$;
\item $d=1$ and $\frac13\le q\le\frac12$, or $d=2$ and $q=\frac12$.
\end{enumerate}
In particular, in dimension $d\ge3$ no functional $-\int\rho^q\dd x$,
$0<q<1$, is geodesically $\lambda$-convex on $(\calP(X),\SHK)$ for any
$\lambda\in\R$.
\end{theorem}

The proof of $(1)\Rightarrow(3)$ is given in Appendix~\ref{nc:sec-rigidity}.
It adds a far-away reservoir of mass to a small multiple of explicit test
measures supported in a small ball, and compares the energy at three
times. The test measures are images of a constant density under maps that
move points radially and change the mass; they are built from a negative
direction of the matrix $B_0(1)$ of Definition~\ref{prelim:matrix}, which
exists exactly outside the range (3). A strict contact inequality shows
that every $\HK$ geodesic between the two endpoints passes through the
third measure at its midpoint (Proposition~\ref{nc:engine}). Since nothing
is transported between the reservoir and the test part, the normalisation
changes the energy of the test part by the factor $\eps^q$, whereas the
curvature of the sphere contributes terms of order $\eps$. Since $q<1$, the
negative midpoint defect on the cone survives on the sphere, whatever the
value of $\lambda$. The argument uses neither the smooth local test
geodesics of the proof of \cite[Thm.~7.5]{LMS23} nor the strong
non-branching of $\HK$ geodesics \cite[Cor.~5.3]{LMS23}. Without the
reservoir, the same test measures prove the converse part of
Theorem~\ref{prelim:LMS} (Appendix~\ref{nc:sec-lms}).

\begin{remark}\label{tr:rem-rigidity}
The argument uses the exact homogeneity of degree $q<1$ of $\sfE_q$ in an
essential way and concerns only the unshifted negative powers. It gives no
obstruction for energies with bounded derivative, such as the regularised
powers $1-(1+c)^q$ in the ranges of exponents of Table~\ref{ex:table}, which
are geodesically $\lambda$-convex on $(\calP(X),\SHK)$ with a negative
constant $\lambda$. As before, the hypothesis $\diam X>\pi/2$ is
used to separate the test measures from the reservoir, and no statement is made
for $\diam X\le\pi/2$.
\end{remark}
\section{New geodesically convex functionals}\label{sec:examples}

Theorem~\ref{tr:necessary} leaves little room on the sphere: on a domain of
diameter larger than $\pi/2$, an energy with $\ezero>-\infty$ can be
geodesically $\lambda$-convex only if $E'$ is bounded. Corollary~\ref{tr:global}
gives such convexity as soon as $\Einf<\infty$ and $\sfE$ is geodesically
convex on the cone, and our only tool on the cone is the condition of Liero,
Mielke and Savar\'e (Theorem~\ref{prelim:LMS}). The question of this section
is therefore which convex $E$ with bounded derivative satisfy that
condition, and with which constant.

We first rewrite the condition as inequalities for scalar quantities built
from $E$, $E'$ and $E''$ at one point (Proposition~\ref{ex:reduced}), and
then, in the variable $s=c^{1/d}$, as one inequality for the function $G$
defined by $E(c)=s^{d-1}G(s)$ (Proposition~\ref{ex:Cd}). In dimension two
this inequality is solved completely. In the other dimensions we deal with
it in three ways: a sufficient criterion that is linear in $E$; one explicit
kernel, which dilations and mixtures turn into a cone of examples; and a
direct estimate for the shifted powers $1-(1+c)^q$. Table~\ref{ex:table}
collects the constants on the sphere.

Throughout this section $d\ge1$ is the dimension,
$E:[0,\infty)\to\R$ is continuous with $E(0)=0$, and $\Nfun$ is the function
\eqref{prelim:eq-N}. The Liero--Mielke--Savar\'e condition with constant
$\lambda$ is the one defined after Theorem~\ref{prelim:LMS}: $E$ is convex
and \eqref{prelim:eq-LMScond} holds. Each result below supplies this
condition with some $\lambda$, together with the slopes $\ezero$ and
$\Einf$; the convexity statements then follow from
Corollary~\ref{tr:consequences}. To state them briefly, we say that $\sfE$
\emph{has the constant $\lambda$ on $\HK$} if it is strongly geodesically
$\lambda$-convex on $(\calM(\Om),\HK)$ for every closed convex
$\Om\subset\R^d$ with nonempty interior, bounded if $\ezero=-\infty$, and
that it \emph{has the constant $\lambda'$ on $\SHK$} if it is strongly
geodesically $\lambda'$-convex on $(\calP(\Om),\SHK)$ for the same domains.
The constant on $\SHK$ is \emph{sharp} if on every closed convex $\Om$
with nonempty interior and $\diam\Om>\pi/2$ no larger constant is possible, even for weak geodesic
convexity. In these words Corollary~\ref{tr:consequences} says: if $E$
satisfies the condition with constant $\lambda$, then $\sfE$ has the
constant $\lambda$ on $\HK$; if moreover $\Einf<\infty$, it has the constant
$\lambda-2\Einf$ on $\SHK$; and this constant is sharp if $\lambda=2\ezero$.

\subsection{The reduced form of the condition}\label{ex:ss-reduced}

The function $\Nfun$ depends on $E$ only through the combination
$u=\gamma^{d+2}/\rho^d$, so its convexity becomes a family of conditions on
a $2\times2$ matrix indexed by $u>0$. This is the matrix $B_\lambda(u)$ of
Definition~\ref{prelim:matrix}. We use it together with the quantities $A$,
$R$, $Q$ of \eqref{ex:eq-ARQ} and $F_\lambda$ of \eqref{ex:eq-F}, which are
defined there for $E\in C^2(0,\infty)$, $\lambda\in\R$ and $u>0$.
If $E$ is convex, then $A\ge0$, and $R\ge0$ because
$0=E(0)\ge E(u)-uE'(u)$. Here $R$ is the amount by which the tangent to $E$
at $u$ lies below the origin. In $F_\lambda$ the constant $\lambda$ enters
only through the factor $2E'(u)-\lambda$, which multiplies $Q$; the second
term $R\,(4A-dR)$ does not see $\lambda$ and may be negative.

The next proposition says that the condition \eqref{prelim:eq-LMScond} is the
positive semidefiniteness of $B_\lambda(u)$ for all $u$, together with
$R\ge0$ when $d\ge2$. Where $Q>0$, the inequality to be checked is
$F_\lambda\ge0$.

\begin{proposition}[Reduced form]\label{ex:reduced}
Let $E\in C([0,\infty))\cap C^2(0,\infty)$ with $E(0)=0$, and let $\lambda\in\R$.
\begin{enumerate}[label=\textup{(\alph*)}]
\item The function $\rho\mapsto(d-1)\Nfun(\rho,\gamma)$ is non-increasing for
 every $\gamma>0$ if and only if $(d-1)R(u)\ge0$ for all $u>0$.
\item $\Nfun$ is convex on $(0,\infty)^2$ if and only if $B_\lambda(u)$ is
 positive semidefinite for every $u>0$, that is, if and only if for every $u>0$
 \begin{equation}\label{ex:eq-three}
 Q(u)\ge0,\qquad F_\lambda(u)\ge0,\qquad
 A(u)+\tfrac14u\big(2E'(u)-\lambda\big)\ge0 .
 \end{equation}
 We have $\det B_\lambda=F_\lambda/(4d)$. At every point $u$ with $Q(u)>0$ the
 third inequality in \eqref{ex:eq-three} follows from the second. At a point
 with $Q(u)=0$ and $F_\lambda(u)\ge0$ the third inequality is equivalent to
 $\partial_\gamma^2\Nfun\ge0$ at the points $(\rho,\gamma)$ with
 $\gamma^{d+2}/\rho^d=u$.
\end{enumerate}
In particular, a convex $E$ of this regularity satisfies the
Liero--Mielke--Savar\'e condition with constant $\lambda$ if and only if
\eqref{ex:eq-three} holds for all $u>0$.
\end{proposition}

The third condition in \eqref{ex:eq-three} cannot be omitted: for $E(c)=c$ we
have $A=R=Q=0$ and $F_\lambda=0$ for every $\lambda$, while
$\Nfun=(1-\lambda/2)\gamma^2$ is convex only for $\lambda\le2$.
Parts (a) and (b) are computations on $(0,\infty)$ and use only
$E\in C^2(0,\infty)$: neither convexity of $E$, nor $E(0)=0$, nor continuity
at $0$ enters. Thus \eqref{prelim:eq-LMScond} is equivalent to
$(d-1)R\ge0$ together with \eqref{ex:eq-three} for every such $E$; convexity
and $E(0)=0$ serve only to make $R\ge0$ automatic in the last assertion.

\begin{proof}
Put $m=(\rho/\gamma)^d$ and $u=\gamma^{d+2}/\rho^d$, so that $mu=\gamma^2$,
$\partial_\rho u=-du/\rho$ and $\partial_\gamma u=(d+2)u/\gamma$. Since
$u\,\frac{\rmd}{\rmd u}(E-uE')=-A$, differentiating
$\Nfun=mE(u)-\frac\lambda2\gamma^2$ gives
\[
 \partial_\rho\Nfun=-\frac{d\,m}{\rho}\,R,\qquad
 \partial_\gamma\Nfun=\frac m\gamma\big(dR+2uE'\big)-\lambda\gamma ,
\]
with $R$, $A$, $Q$, $E'$ evaluated at $u$. This proves (a). Differentiating
once more, and using $u\frac{\rmd}{\rmd u}(dR+2uE')=(d+2)A+2uE'$,
\begin{gather*}
 \partial_\rho^2\Nfun=\frac{d\,m}{\rho^2}\,Q,\qquad
 \partial_\rho\partial_\gamma\Nfun=\frac{d\,m}{\rho\gamma}\big(dR-(d+2)A\big),\\
 \partial_\gamma^2\Nfun=\frac{m}{\gamma^2}\big((d+2)^2A-d(d+1)R+2uE'\big)-\lambda .
\end{gather*}
Because $m/\gamma^2=1/u$, for every $\xi\in\R^2$
\begin{equation}\label{ex:eq-hess}
\begin{aligned}
 \frac{u}{\gamma^2}D^2\Nfun[(\rho\xi_1,\gamma\xi_2)]^2
 ={}&dQ\xi_1^2+2d(dR-(d+2)A)\xi_1\xi_2\\
 &+[(d+2)^2A-d(d+1)R+u(2E'-\lambda)]\xi_2^2.
\end{aligned}
\end{equation}
Under the invertible substitution $x=d(\xi_1-\xi_2)/2$, $y=-\xi_2$,
this is $4(x,y)B_\lambda(x,y)^\top$. The Hessian is therefore positive
semidefinite exactly when $B_\lambda$ is. Its first diagonal entry is
$Q/d$, and direct expansion gives
\[
 \det B_\lambda=\tfrac14(A-\tfrac{d-1}dR)u(2E'-\lambda)
                  +AR/d-R^2/4=F_\lambda/(4d).
\]
The two diagonal entries and determinant characterize positivity. If
$Q>0$, the determinant condition implies the second diagonal condition.
If $Q=0$ and $F_\lambda\ge0$, the off-diagonal entry must vanish, and
\eqref{ex:eq-hess} at $\xi=(0,1)$ gives
$u\partial_\gamma^2\Nfun=4[A+u(2E'-\lambda)/4]$. The last assertion of the
proposition follows from (a), (b) and $R\ge0$ for convex $E$ with $E(0)=0$.
\end{proof}

With the Euler operator $\vartheta=u\frac{\rmd}{\rmd u}$ we have
$R=(\vartheta-1)E$, $A=\vartheta(\vartheta-1)E$ and therefore
$Q=d\,(\vartheta-1)\big(\vartheta-\frac{d-1}d\big)E$. The second factor
suggests the substitution $E(u)=u^{(d-1)/d}G(u^{1/d})$, which removes it. Every
$E$ can be written in this way. In the variable $s=c^{1/d}$ the kernel of
$Q$, spanned by $c$ and the critical power $c^{(d-1)/d}$, becomes the affine
functions of $s$, and $Q$ becomes a multiple of $G''$.

\begin{proposition}[One-function form]\label{ex:Cd}
Let $E\in C([0,\infty))\cap C^2(0,\infty)$ with $E(0)=0$, and define
$G\in C^2(0,\infty)$ and $P$ by
\[
 E(c)=s^{d-1}G(s),\quad s=c^{1/d},\qquad P(s)=sG'(s)-G(s).
\]
Then, with $u=s^d$,
\begin{equation}\label{ex:eq-G}
 E'=\frac1d\Big(\frac{(d-1)G}{s}+G'\Big),\qquad
 R=\frac{s^{d-1}}{d}\,P,\qquad
 Q=\frac{s^{d+1}}{d}\,G'',
\end{equation}
\begin{equation}\label{ex:eq-FG}
 F_\lambda=\frac{s^{2d-2}}{d^3}\Big(s^2G''\big[d^2s\,(2G'-\lambda)-2(d-2)(d+1)P\big]-(d-2)^2P^2\Big).
\end{equation}
Consequently, if $E$ is convex and $G''>0$ on $(0,\infty)$, then $E$ satisfies
the Liero--Mielke--Savar\'e condition with constant $\lambda$ if and only if
\begin{equation}\label{ex:eq-Cd}
 s^2G''(s)\big[d^2s\,(2G'(s)-\lambda)-2(d-2)(d+1)P(s)\big]\ \ge\ (d-2)^2P(s)^2
 \qquad\text{for all }s>0 .
\end{equation}
\end{proposition}

\begin{proof}
Using $u\partial_u=d^{-1}s\partial_s$ gives \eqref{ex:eq-G}.
Eliminate $A$ by $dA=Q+(d-1)R$:
\begin{equation}\label{ex:eq-dF}
 dF_\lambda=dQu(2E'-\lambda)+4QR-(d-2)^2R^2.
\end{equation}
Substitute \eqref{ex:eq-G} and
$u(2E'-\lambda)=d^{-1}s^{d-1}[ds(2G'-\lambda)-2(d-1)P]$.
Collecting terms gives \eqref{ex:eq-FG}. For convex $E$ and $G''>0$,
$R\ge0$ and $Q>0$, so Proposition~\ref{ex:reduced} reduces the full
criterion to $F_\lambda\ge0$.
\end{proof}

The identities above were also verified by
computer algebra; see Appendix~\ref{sec:verification}. The condition $Q\ge0$ is McCann's condition;
in the variable $G$ it says that $G$ is convex. The necessity of
\eqref{ex:eq-Cd} does not need the hypothesis $G''>0$: by
Proposition~\ref{ex:reduced}, \eqref{prelim:eq-LMScond} implies
$F_\lambda\ge0$, which is \eqref{ex:eq-Cd} by \eqref{ex:eq-FG}, and it
implies $Q\ge0$, that is, $G''\ge0$. For $d=2$ the right-hand side of
\eqref{ex:eq-Cd} vanishes, which leads to the complete description of the next
subsection. For $d\ne2$ the term $(d-2)^2P^2$ is an obstruction: convexity of
$G$ and a lower bound for $G'$ do not suffice. The left-hand side of
\eqref{ex:eq-Cd} carries the factor $G''$, the right-hand side does not, so
wherever $P\ne0$ the function $G$ has to be convex by a definite amount.

\subsection{Dimension two}\label{ex:ss-d2}

For $d=2$ we have $E(c)=\sqrt c\,G(\sqrt c)$, and what remains of
\eqref{ex:eq-Cd} is a product of $G''$ and $2G'-\lambda$. The condition can
be read off: $G$ is convex, and $\lambda$ is at most twice the smallest
slope of $G$. Care is needed only at points where $G''=0$, where the third
inequality of \eqref{ex:eq-three} takes over.

\begin{theorem}[Complete description for $d=2$]\label{ex:d2}
Let $d=2$, let $E\in C([0,\infty))\cap C^2(0,\infty)$ be convex with $E(0)=0$,
and put $G(s)=E(s^2)/s$ for $s>0$. Then $E$ satisfies the
Liero--Mielke--Savar\'e condition with constant $\lambda$ if and only if
\[
 G\ \text{is convex on }(0,\infty)\qquad\text{and}\qquad 2G'(s)\ge\lambda\ \text{ for all }s>0 .
\]
Thus the best constant is $\lambda=2G'(0^+)=2\inf G'$, and there is none if $G$
is not convex or $G'(0^+)=-\infty$. Moreover
$\Einf=G'(\infty):=\lim_{s\to\infty}G'(s)$. Hence, by
Corollary~\ref{tr:consequences}, if $G$ is convex and Lipschitz continuous,
then $\sfE$ has the constant $2G'(0^+)$ on $\HK$ and the constant
$2\big(G'(0^+)-G'(\infty)\big)$ on $\SHK$. The constant on
$\SHK$ is zero only if $G$ is affine, that is, $E(c)=-\alpha\sqrt c+\beta c$ with
$\alpha\ge0$ and $\beta\in\R$.
\end{theorem}

\begin{proof}
For $d=2$, \eqref{ex:eq-G} and \eqref{ex:eq-FG} read $Q=\frac12s^3G''$ and
$F_\lambda=\frac12s^5G''\,(2G'-\lambda)$. Moreover $A=\frac12(Q+R)$,
$R=\frac12sP$ and $u(2E'-\lambda)=s^2(2G'-\lambda)-sP$, so that the third
quantity in \eqref{ex:eq-three} is
\[
 A+\tfrac14u(2E'-\lambda)=\tfrac14s^3G''+\tfrac14s^2(2G'-\lambda).
\]
If $G''\ge0$ and $2G'\ge\lambda$, all three quantities in \eqref{ex:eq-three}
are nonnegative. Conversely, \eqref{ex:eq-three} gives $G''\ge0$ everywhere; at
a point with $G''(s)>0$ the inequality $F_\lambda\ge0$ gives $2G'(s)\ge\lambda$,
and at a point with $G''(s)=0$ the third inequality gives the same. The
equivalence now follows from Proposition~\ref{ex:reduced}, since $E$ is convex.
For a convex $G$, $\inf G'=G'(0^+)$. The identity $\Einf=G'(\infty)$ needs no
convexity of $G$. Put $q(c)=E(c)/c$, so that $q(c)\to\Einf$ as $c\to\infty$
and $G'(s)=2E'(s^2)-q(s^2)$. Since $E$ is convex with $E(0)=0$, we have
$q(c)\le E'(c)\le\Einf$ for all $c>0$. If $\Einf<\infty$, then
$q(s^2)\le G'(s)\le2\Einf-q(s^2)$, and both bounds tend to $\Einf$; if
$\Einf=\infty$, the lower bound alone gives $G'(s)\to\infty$.
If $G$ is convex and Lipschitz, both numbers are finite, and
Corollary~\ref{tr:consequences} applies with $\lambda=2G'(0^+)$. Finally
$G'(0^+)=G'(\infty)$ forces $G'$ to be
constant, $G(s)=\beta s-\alpha$, and $E(c)=\beta c-\alpha\sqrt c$ is convex only
for $\alpha\ge0$.
\end{proof}

\begin{remark}\label{ex:rem-d2}
(a) In terms of $G$, convexity of $E$ is the inequality
$s^2G''+sG'-G\ge0$, since $2A=Q+R$; and $E(0)=0$ means $sG(s)\to0$ as $s\to0$.

(b) Let $G$ be convex and Lipschitz continuous, as in the last part of
Theorem~\ref{ex:d2}. If $G(0^+)=0$, then $\ezero=G'(0^+)$ by \eqref{ex:eq-G}, and the constant on
$\SHK$ equals $2(\ezero-\Einf)$; it is then sharp by
Corollary~\ref{tr:consequences}(iii).

(c) The assumption $E\in C^2(0,\infty)$ is not needed for the sufficiency
part. For $d=2$ we have $\Nfun(\rho,\gamma)=\rho\,H(\gamma^2/\rho)$ with
$H(s)=G(s)-\frac\lambda2s$. If $G$ is convex and $\lambda$ is at most twice
the infimum of its difference quotients, then $H$ is convex and
nondecreasing, so its perspective $(\rho,y)\mapsto\rho H(y/\rho)$ is convex
and nondecreasing in $y$, and $\Nfun$ is convex. For fixed $\gamma$,
$\Nfun=\gamma^2\big(E(u)/u-\frac\lambda2\big)$ with $u=\gamma^4/\rho^2$, which is
nonincreasing in $\rho$ because $E$ is convex with $E(0)=0$. Moreover
$E(c)/c=G(\sqrt c)/\sqrt c$ tends to the supremum of the difference
quotients of $G$. Hence the statements of Theorem~\ref{ex:d2} about a convex
Lipschitz $G$ hold for every convex continuous $E$ with $E(0)=0$, with
$G'(0^+)$ and $G'(\infty)$ read as the infimum and the supremum of the
difference quotients.
\end{remark}

\begin{example}\label{ex:ex-d2}
(a) $G\equiv-1$ gives $E(c)=-\sqrt c$, and $G(s)=s^{2m-1}$, $m>1$, gives
$E(c)=c^m$; in both cases $\lambda=0$ on $\HK$, in agreement with
\cite{LMS23,LaM22}.

(b) $G(s)=\ee^{-s}-1$, that is, $E(c)=\sqrt c\,\big(\ee^{-\sqrt c}-1\big)$. Here
$G''=\ee^{-s}>0$, and $E$ is convex because
$s^2G''+sG'-G=1-(1+s-s^2)\ee^{-s}\ge0$, which follows from
$\ee^{s}\ge1+s\ge1+s-s^2$. Since $G'(0^+)=-1$ and $G'(\infty)=0$, the
constants on $\HK$ and on $\SHK$ are both $-2$. Here $G(0^+)=0$, so
$\ezero=-1$ is finite, the domain may be unbounded, and the constant on
$\SHK$ is sharp (Remark~\ref{ex:rem-d2}(b)).

(c) $G(s)=\big(1-(1+s^2)^q\big)/s$, $0<q<1$, corresponds to
$E(c)=1-(1+c)^q$. By \eqref{ex:eq-G}, Step~1 of the proof of
Theorem~\ref{ex:shifted} and Remark~\ref{ex:rem-shifted}(a) below, $G$ is convex if and only if $q\ge\frac12$;
then $\lambda=2G'(0^+)=-2q$. So for $d=2$ the range $q\ge\frac12$ of
Theorem~\ref{ex:shifted}(b) is the full range.
\end{example}

\subsection{A sufficient criterion in every dimension}\label{ex:ss-Kd}

For $d\ne2$ the inequality \eqref{ex:eq-Cd} is quadratic in $G$. We look for
a condition that is linear in $E$ and implies $B_\lambda\ge0$. The natural
constant to aim for is $\lambda=2\ezero$: after the transfer it becomes
$2(\ezero-\Einf)$, which is the bound of Theorem~\ref{tr:necessary}. The
criterion asks that $E$ be convex by a fixed multiple of $R$.

\begin{theorem}[Sufficient criterion]\label{ex:Kd}
Let $d\ge1$ and let $E\in C([0,\infty))\cap C^2(0,\infty)$ be convex with
$E(0)=0$ and $\ezero=E'(0^+)>-\infty$. Assume that
\begin{equation}\label{ex:eq-Kd}
 u^2E''(u)\ \ge\ \kappa_d\,\big(uE'(u)-E(u)\big)\quad\text{for all }u>0,
 \qquad \kappa_d:=\frac{3d-2}{2(d+2)} .
\end{equation}
Then $E$ satisfies the Liero--Mielke--Savar\'e condition with constant
$\lambda=2\ezero$. Hence, by Corollary~\ref{tr:consequences}, $\sfE$ has the
constant $2\ezero$ on $\HK$ and, if $\Einf<\infty$, the sharp constant
$2(\ezero-\Einf)$ on $\SHK$.
\end{theorem}

\begin{proof}
The idea is to write $A=(A-\kappa R)+\kappa R$ in $B_{2\ezero}$ and to ask
for which $\kappa$ the part proportional to $R$ is a positive semidefinite
quadratic form. This holds exactly for $\kappa\ge\kappa_d$, and for
$\kappa=\kappa_d$ that part is a perfect square; so \eqref{ex:eq-Kd} is the
weakest condition $A\ge\kappa R$ for which this splitting works.
Let $H=E-\ezero u\ge0$. Convexity also gives $R\ge0$.
Since $u(2E'-2\ezero)/4=(R+H)/2$, expansion gives
\begin{equation}\label{ex:eq-sos}
 (x,y)B_{2\ezero}(x,y)^\top
 =(A-\kappa_dR)(x+y)^2
 +\frac{R}{2d(d+2)}[(d-2)x+2dy]^2+\frac H2y^2.
\end{equation}
Indeed the coefficients of $Rx^2,Rxy,Ry^2$ are respectively
$-(d-1)/d,-1,1/2$. All three terms are nonnegative under the hypothesis,
including at degenerate points. Proposition~\ref{ex:reduced} proves the
criterion.
\end{proof}

The hypothesis \eqref{ex:eq-Kd} is a pointwise linear inequality, preserved
under sums and positive multiples. It does not cover everything: for the
shifted powers it gives a smaller range of exponents than the direct
estimate (Remark~\ref{ex:rem-shifted}).

\begin{remark}\label{ex:rem-Kd}
The condition is linear in $E$, and its coefficient satisfies
$\kappa_d=3/2-4/(d+2)<3/2$ and
$d/4-\kappa_d=(d-2)^2/[4(d+2)]\ge0$. Thus $4A\ge dR$ is a stronger
sufficient condition. In the $G$ variable the new criterion is
$s^2G''\ge(d-2)^2P/[2(d+2)]$; it reduces to the necessary condition
$Q\ge0$ when $d=2$, under the stated finite-$\ezero$ hypothesis.
\end{remark}

\subsection{Dilations and mixtures}\label{ex:ss-mix}

Outside the reach of Theorem~\ref{ex:Kd} every density function has to be
verified by a computation of its own. Two operations make one such
verification go further. Replacing $E(c)$ by $a\,E(c/a)$ amounts to a linear
change of the variable $\rho$ in $\Nfun$, and since $\Nfun$ is linear in the
pair $(E,\lambda)$, the condition survives integration over $a$. One
verified density function with bounded derivative thus produces a convex
cone of them.

\begin{lemma}[Dilation]\label{ex:dilation}
For $a>0$ put $E_a(c)=a\,E(c/a)$. Then
\[
 N_{\lambda,E_a}(\rho,\gamma)=\Nfun\big(a^{1/d}\rho,\gamma\big)\qquad\text{for all }\rho,\gamma>0 .
\]
Hence, if $E$ satisfies the Liero--Mielke--Savar\'e condition with constant
$\lambda$, so does $E_a$, with the same $\lambda$. Moreover $E_a(0)=0$, $E_a$ is
convex if $E$ is, and $E_a'(c)=E'(c/a)$ wherever $E$ is differentiable at
$c/a$; for convex $E$ the same holds for the one-sided derivatives, and
$E_a$ has the same slopes $\ezero$ and $\Einf$ as $E$.
\end{lemma}

\begin{proof}
With $\rho'=a^{1/d}\rho$,
\[
 \Big(\frac\rho\gamma\Big)^{d}a\,E\Big(\frac{\gamma^{d+2}}{a\rho^{d}}\Big)
 =\Big(\frac{\rho'}\gamma\Big)^{d}E\Big(\frac{\gamma^{d+2}}{\rho'^{\,d}}\Big).
\]
Convexity in $(\rho,\gamma)$ and monotonicity in $\rho$ are invariant under the
linear map $(\rho,\gamma)\mapsto(a^{1/d}\rho,\gamma)$.
\end{proof}

\begin{proposition}[Mixtures]\label{ex:mixture}
Let $E$ satisfy the Liero--Mielke--Savar\'e condition with constant $\Lambda$,
and assume that $\ezero$ and $\Einf$ are finite. Let $m$ be a finite positive
Borel measure on $(0,\infty)$ with total mass $|m|$, and put
\[
 \bar E(c)=\int_{(0,\infty)}a\,E(c/a)\,m(\rmd a).
\]
Then $\bar E$ is finite, convex and continuous with $\bar E(0)=0$, it satisfies
the Liero--Mielke--Savar\'e condition with constant $\Lambda|m|$, and
$\bar E'(0^+)\ge\ezero\,|m|$, $\bar E'_\infty\le\Einf\,|m|$. Hence, by
Corollary~\ref{tr:consequences}, the functional $\bar\sfE$ with density
function $\bar E$ has the constant $\Lambda|m|$ on $\HK$ and the constant
$(\Lambda-2\Einf)|m|$ on $\SHK$.
\end{proposition}

\begin{proof}
The bound $|aE(c/a)|\le\max\{|\ezero|,|\Einf|\}c$ makes the integral
finite and continuous, and convexity passes under the integral.
Furthermore $N_{\Lambda|m|,\bar E}=\int N_{\Lambda,E_a}\dd m$.
Lemma~\ref{ex:dilation} therefore passes joint convexity and weighted
$\rho$-monotonicity to $\bar E$. The bounds
$\ezero|m|\le\bar E(c)/c\le\Einf|m|$ give the two slope bounds.
Corollary~\ref{tr:consequences} gives the constant
$\Lambda|m|-2\bar E'_\infty\ge(\Lambda-2\Einf)|m|$ on $\SHK$.
\end{proof}

\subsection{A kernel in every dimension and the complete Bernstein class}\label{ex:ss-bernstein}

The criterion of Theorem~\ref{ex:Kd} is not the only source of examples. The
following density function satisfies \eqref{ex:eq-Cd} in every dimension,
although it violates \eqref{ex:eq-Kd} for every $d\ne2$
(Remark~\ref{ex:rem-bern}).

The choice is guided by the variable $G$. We need $G$ convex and $E'$
bounded, and the two model cases are $G(s)=-s$, which is $E(c)=-c$, and
$G\equiv-1$, which is the critical power $E(c)=-c^{(d-1)/d}$. The function
$G(s)=-s/(1+s)$ passes from the first at small $s$ to the second at large
$s$.

\begin{proposition}[The kernel]\label{ex:kernel}
For $d\ge1$ let
\[
 K_d(c)=-\frac{c}{1+c^{1/d}},\qquad c\ge0 .
\]
Then $K_d$ is convex, $K_d(0)=0$, $-1\le K_d'<0$, $K_d'(0^+)=-1$,
$K_d'(c)\to0$ as $c\to\infty$, and $K_d$ satisfies the
Liero--Mielke--Savar\'e condition with constant $\lambda=-2$. More precisely,
with $s=u^{1/d}$,
\begin{equation}\label{ex:eq-kernelF}
 Q=\frac{2us}{d\,(1+s)^3}>0,\qquad
 F_{-2}=\frac{(d+2)\,u^2s^2\,\big[(3d+2)+(3d-2)s\big]}{d^3\,(1+s)^5}>0
 \qquad\text{for }u>0 .
\end{equation}
\end{proposition}

\begin{proof}
In the notation of Proposition~\ref{ex:Cd}, $K_d(c)=s^{d-1}G(s)$ with
\[
 G(s)=-\frac{s}{1+s},\qquad G'(s)=-\frac1{(1+s)^2},\qquad G''(s)=\frac{2}{(1+s)^3},
 \qquad P(s)=sG'-G=\frac{s^2}{(1+s)^2}.
\]
By \eqref{ex:eq-G},
\[
 K_d'=\frac1d\Big(-\frac{d-1}{1+s}-\frac1{(1+s)^2}\Big)=-\frac{1+\frac{d-1}{d}s}{(1+s)^2},
 \qquad R=\frac{us}{d(1+s)^2},\qquad Q=\frac{2us}{d(1+s)^3},
\]
and $A=\frac1d(Q+(d-1)R)=\frac{us\,(d+1+(d-1)s)}{d^2(1+s)^3}>0$. Thus $K_d$ is
convex, and $1+\frac{d-1}ds\le(1+s)^2$ gives $-1\le K_d'<0$, with the stated
limits. For $\lambda=-2$ we have $2G'+2=\frac{2s(s+2)}{(1+s)^2}$, and the
left-hand side of \eqref{ex:eq-Cd} minus its right-hand side equals
\begin{align*}
 &\frac{2s^2}{(1+s)^3}\cdot\frac{2d^2s^2(s+2)-2(d-2)(d+1)s^2}{(1+s)^2}
 -\frac{(d-2)^2s^4}{(1+s)^4}\\
 &\qquad=\frac{s^4}{(1+s)^5}\Big(\big[4d^2-(d-2)^2\big]s+8d^2-4(d-2)(d+1)-(d-2)^2\Big)\\
 &\qquad=\frac{(d+2)\,s^4\,\big[(3d-2)s+(3d+2)\big]}{(1+s)^5},
\end{align*}
where we used $4d^2-(d-2)^2=(3d-2)(d+2)$ and
$8d^2-4(d-2)(d+1)-(d-2)^2=3d^2+8d+4=(3d+2)(d+2)$. Multiplying by $s^{2d-2}/d^3$
gives \eqref{ex:eq-kernelF}. Since $G''>0$ and $K_d$ is convex,
Proposition~\ref{ex:Cd} applies.
\end{proof}

Dilating $K_d$ by $a=b^d$ gives
\begin{equation}\label{ex:eq-dilK}
 b^d\,K_d\big(c/b^d\big)=-\frac{b\,c}{b+c^{1/d}}=-c^{\frac{d-1}d}\,\frac{b\,s}{b+s},
 \qquad s=c^{1/d},
\end{equation}
and mixing these functions leads to the following class. Recall that a function
$f:(0,\infty)\to[0,\infty)$ is a \emph{complete Bernstein function} if
\[
 f(s)=\alpha+\beta s+\int_{(0,\infty)}\frac{s}{s+b}\,\sigma(\rmd b)
\]
with $\alpha,\beta\ge0$ and a positive measure $\sigma$ on $(0,\infty)$
satisfying $\int(1+b)^{-1}\sigma(\rmd b)<\infty$; see \cite{SSV12}. The functions
in \eqref{ex:eq-f} below are exactly the complete Bernstein functions with
$f(0^+)=0$ and $f'(0^+)<\infty$, with $m(\rmd b)=b^{-1}\sigma(\rmd b)$. We do not
use this terminology in the proofs; the class is defined by \eqref{ex:eq-f}.

In the theorem the measure $m$ distributes the scale $b$ at which the
kernel passes from linear to critical behaviour, and $\beta$ adds a multiple
of $-c$. Only the total mass $|m|$, the oscillation of $E'$, enters the
constant on the sphere.

\begin{theorem}[Complete Bernstein class]\label{ex:bernstein}
Let $d\ge1$, let $m$ be a finite positive Borel measure on $(0,\infty)$ with
total mass $|m|$, let $\beta\ge0$, and put
\begin{equation}\label{ex:eq-f}
 f(s)=\beta s+\int_{(0,\infty)}\frac{b\,s}{b+s}\,m(\rmd b),\qquad
 E(c)=-c^{\frac{d-1}{d}}\,f\big(c^{1/d}\big).
\end{equation}
Then $E$ is convex, $E(0)=0$, $E\in C^1([0,\infty))$,
\[
 \ezero=-f'(0^+)=-\beta-|m|,\qquad \Einf=-\beta,\qquad -\beta-|m|\le E'\le-\beta ,
\]
and $E$ satisfies the Liero--Mielke--Savar\'e condition with constant
$\lambda=2\ezero=-2(\beta+|m|)$. Hence, by Corollary~\ref{tr:consequences},
$\sfE$ has the constant $-2(\beta+|m|)$ on $\HK$ and the sharp constant
$-2|m|=2(\ezero-\Einf)$ on $\SHK$.
\end{theorem}

\begin{proof}
\emph{Step 1 (the mixture).} Let $\tilde m$ be the image of $m$ under
$b\mapsto b^d$. By \eqref{ex:eq-dilK},
\[
 E^0(c):=\int_{(0,\infty)}a\,K_d(c/a)\,\tilde m(\rmd a)
 =-c^{\frac{d-1}d}\int_{(0,\infty)}\frac{b\,s}{b+s}\,m(\rmd b),\qquad s=c^{1/d}.
\]
By
Proposition~\ref{ex:kernel}, $K_d$ satisfies the condition with constant $-2$
and has $K_d'(0^+)=-1$ and recession slope $0$. Proposition~\ref{ex:mixture}
shows that $E^0$ is convex and continuous with $E^0(0)=0$ and satisfies the
condition with constant $-2|m|$.

\emph{Step 2 (the derivative).} By Lemma~\ref{ex:dilation} the derivative of
$c\mapsto aK_d(c/a)$ is $K_d'(c/a)\in[-1,0)$; it tends to $-1$ as $c\to0$ and to
$0$ as $c\to\infty$. The difference quotients of $c\mapsto aK_d(c/a)$ are
bounded by $1$ in absolute value, so dominated convergence gives
$(E^0)'(c)=\int K_d'(c/a)\,\tilde m(\rmd a)$ for $c>0$, the same formula for the
right derivative at $c=0$, and
\[
 -|m|\le(E^0)'\le0,\qquad (E^0)'(0^+)=(E^0)'(0)=-|m|,\qquad \lim_{c\to\infty}(E^0)'(c)=0 .
\]
Since $E^0$ is convex, the last limit is its recession slope.

\emph{Step 3 (the drift).} We have $E(c)=E^0(c)-\beta c$ and, for every
$\lambda$,
\[
 N_{\lambda-2\beta,\,E}(\rho,\gamma)=N_{\lambda,E^0}(\rho,\gamma)-\beta\gamma^2+\beta\gamma^2
 =N_{\lambda,E^0}(\rho,\gamma),
\]
because $(\rho/\gamma)^d\,\gamma^{d+2}/\rho^d=\gamma^2$. Hence $E$ satisfies the
condition with constant $-2|m|-2\beta$; it is convex with $E(0)=0$, and by Step~2
$\ezero=-\beta-|m|$, $\Einf=-\beta$ and $-\beta-|m|\le E'\le-\beta$. Finally
$f'(0^+)=\beta+|m|$ by monotone convergence in
$f(s)/s=\beta+\int\frac{b}{b+s}\,m(\rmd b)$.
\end{proof}

The sign of $\beta$ is not used in the proof: the statements of
Theorem~\ref{ex:bernstein} hold, with the same constants, for every
$\beta\in\R$, in agreement with Remark~\ref{prelim:normalisation}. We state
the theorem for $\beta\ge0$ because then $f$ is a complete Bernstein
function.

To recognise members of the class we use the following representation.

\begin{lemma}[Representation of $(1+s)^q-1$]\label{ex:stieltjes}
For $0<q<1$ and $s\ge0$,
\begin{equation}\label{ex:eq-stieltjes}
 (1+s)^q-1=\int_1^\infty\frac{b\,s}{b+s}\,m_q(\rmd b),\qquad
 m_q(\rmd b)=\frac{\sin\pi q}{\pi}\,\frac{(b-1)^q}{b^2}\,\rmd b ,
\end{equation}
and $m_q$ has total mass $q$. Moreover
$\log(1+s)=\int_1^\infty\frac{b\,s}{b+s}\,\frac{\rmd b}{b^2}$, where the measure
$b^{-2}\rmd b$ on $(1,\infty)$ has total mass $1$.
\end{lemma}

\begin{proof}
For $t>0$ the Beta integral gives
$\int_0^\infty\frac{\sigma^{q-1}}{\sigma+t}\dd\sigma=\frac{\pi}{\sin\pi q}\,t^{q-1}$.
With $t=1+x$ and integration over $x\in(0,s)$,
\[
 (1+s)^q-1=q\int_0^s(1+x)^{q-1}\dd x
 =\frac{q\sin\pi q}{\pi}\int_0^\infty\sigma^{q-1}\log\frac{\sigma+1+s}{\sigma+1}\dd\sigma ,
\]
where the order of integration was exchanged by Tonelli's theorem. For $r>0$,
\[
 \log\Big(1+\frac sr\Big)=\int_r^\infty\Big(\frac1b-\frac1{b+s}\Big)\dd b
 =\int_r^\infty\frac{s}{b+s}\,\frac{\rmd b}{b}.
\]
For $r=1$ this is the formula for $\log(1+s)$, and $\int_1^\infty b^{-2}\dd b=1$.
Using it with $r=\sigma+1$ and Tonelli's theorem once more,
\[
 (1+s)^q-1=\frac{q\sin\pi q}{\pi}\int_1^\infty\frac{s}{b+s}\,\frac1b
 \Big(\int_0^{b-1}\sigma^{q-1}\dd\sigma\Big)\dd b
 =\frac{\sin\pi q}{\pi}\int_1^\infty\frac{s}{b+s}\,\frac{(b-1)^q}{b}\dd b ,
\]
which is \eqref{ex:eq-stieltjes}. Dividing \eqref{ex:eq-stieltjes} by $s$ and
letting $s\downarrow0$, monotone convergence gives
$m_q((1,\infty))=\frac{\rmd}{\rmd s}(1+s)^q\big|_{s=0}=q$.
\end{proof}

In the following members $\beta=0$, so $\Einf=0$ and the constants on $\HK$
and on $\SHK$ agree.

\begin{example}[Members of the class, in every dimension]\label{ex:ex-bern}
Write $s=c^{1/d}$. By Theorem~\ref{ex:bernstein} and Lemma~\ref{ex:stieltjes}
the following density functions have the constant $\lambda$ on $\HK$ and
the sharp constant $\lambda$ on $\SHK$, for every $d\ge1$:
\begin{enumerate}[label=\textup{(\alph*)}]
\item $E(c)=-\dfrac{c}{1+c^{1/d}}$, with $\lambda=-2$;
\item $E(c)=-c^{\frac{d-1}{d}}\log\big(1+c^{1/d}\big)$, with $\lambda=-2$;
\item $E(c)=c^{\frac{d-1}{d}}\big[1-(1+c^{1/d})^q\big]$ for every $q\in(0,1)$,
 with $\lambda=-2q$;
\item $E(c)=-\sum_{i=1}^n\alpha_i\dfrac{b_i\,c}{b_i+c^{1/d}}$ with
 $\alpha_i,b_i>0$, with $\lambda=-2\sum_i\alpha_i$.
\end{enumerate}
For $d=1$ these are $-c/(1+c)$, $-\log(1+c)$, $1-(1+c)^q$ and the rational
functions $-\sum_i\alpha_ib_ic/(b_i+c)$, which saturate at
$-\sum_i\alpha_ib_i$. Adding a drift $-\beta c$, $\beta\ge0$, changes the
constant on $\HK$ by $-2\beta$ and leaves the constant on $\SHK$ unchanged.
\end{example}

\begin{remark}\label{ex:rem-bern}
For the kernel,
$A/R=(d-1)/d+2/[d(1+s)]$. Since
$\kappa_d-(d-1)/d=(d-2)^2/[2d(d+2)]$, the simpler criterion
\eqref{ex:eq-Kd} fails for large $s$ whenever $d\ne2$.
The full determinant test above is therefore essential.
If $\beta=0$ and $\int b\dd m<\infty$, monotone convergence gives
$f(s)\uparrow\int b\dd m$, hence
$E(c)\sim-(\int b\dd m)c^{(d-1)/d}$ for a nonzero mixture.
\end{remark}

\subsection{Shifted powers}\label{ex:ss-shifted}

We now consider $E_q(c)=1-(1+c)^q$ with $0<q<1$. This function is convex,
$E_q(0)=0$, $E_q'(0^+)=-q$ and $(E_q)'_\infty=0$. For $d=1$ it belongs to
the class of Theorem~\ref{ex:bernstein}; for $d\ge2$ we argue directly with
Proposition~\ref{ex:reduced}.

The shifted powers deserve separate treatment for two reasons. They are the
first candidates one writes down: the negative power $-c^q$ has unbounded
derivative at zero, and shifting the argument by one removes this without
changing the growth at high density. And they show a phenomenon that the
Bernstein class and Theorem~\ref{ex:Kd} hide: our estimate gives the
constant $\lambda=2\ezero$, here $-2q$, for exponents $q\ge q_d$, and in
dimension one for every exponent, while Proposition~\ref{ex:counter} shows
that the condition with this constant fails near $q=\frac{d-1}d$ when
$d\ge3$.

\begin{theorem}[Shifted powers]\label{ex:shifted}
Let $d\ge1$, $a=\frac{d-1}{d}$, $0<q<1$ and $E_q(c)=1-(1+c)^q$, so that
$\ezero=-q$ and $\Einf=0$.
\begin{enumerate}[label=\textup{(\alph*)}]
\item Let $a<q<1$ and
 \[
 b=\max\Big\{0,\ \frac{\frac14-\frac qd}{q-a}\Big\},\qquad L_*=q\,\max\{1,b\}.
 \]
 Then $E_q$ satisfies the Liero--Mielke--Savar\'e condition with constant
 $-2L_*$.
\item Put $q_d:=\dfrac{5d-4}{4d+4}$. For $a<q<1$ we have $L_*=q$ if and only if $q\ge q_d$. For
 $q_d\le q<1$, $E_q$ satisfies the condition with constant
 $-2q$; this includes the boundary case
 $d=2$, $q=q_2=\frac12=a$, which is not covered by \textup{(a)}. The values are
 \[
 q_1=\tfrac18,\quad q_2=\tfrac12,\quad q_3=\tfrac{11}{16},\quad q_4=\tfrac45,\quad
 q_5=\tfrac78,\quad q_6=\tfrac{13}{14},\quad q_7=\tfrac{31}{32},
 \]
 and $q_d\ge1$ for $d\ge8$, so that the range $q_d\le q<1$ is nonempty
 exactly for $1\le d\le7$.
\item For $d=1$ and every $q\in(0,1)$, $E_q$ satisfies the condition with
 constant $-2q$.
\end{enumerate}
Hence, by Corollary~\ref{tr:consequences}, the functional with density
function $E_q$ has the same constant on $\HK$ and on $\SHK$: $-2L_*$ in case
\textup{(a)}, and $-2q$ in cases \textup{(b)} and \textup{(c)}. The constant
$-2q=2\ezero$ is sharp on $\SHK$.
\end{theorem}

In the proof we bound $F_\lambda$ from below by $R^2$ times a number that
depends only on $d$, $q$ and $\lambda$. The second term $R\,(4A-dR)$ of
$F_\lambda$ can be negative, and the first has to compensate it: Step~1
compares $A$ with $R$, Step~2 compares $R$ with the factor
$u(2E'-\lambda)$, and Step~3 combines them. Step~3 gives
$F_\lambda\ge R^2[4\theta d(q-a)+4q-d]$ for $\lambda=-2q\theta$,
$\theta\ge1$; the number $b$ in (a) is the smallest $\theta\ge0$ for which
the bracket is nonnegative.

\begin{proof}
We write $w=1+u$, $z=w^{q-1}\in(0,1]$ and, for $E=E_q$,
\[
 A=q(1-q)\,u^2w^{q-2},\qquad R=w^q-1-q\,u\,w^{q-1},\qquad x:=u\,(1-z)\ge0 .
\]

\emph{Step 1.} For every $\kappa\le q$ we have $A\ge\kappa R$ on $[0,\infty)$,
with strict inequality for $u>0$; and $R>0$ for $u>0$. Indeed,
$A'=q(1-q)\,u\,w^{q-3}(2+qu)$ and $R'=q(1-q)\,u\,w^{q-3}(1+u)$, so that
\[
 \frac{\rmd}{\rmd u}\big(A-\kappa R\big)=q(1-q)\,u\,w^{q-3}\big[(2-\kappa)+(q-\kappa)u\big]>0
 \quad\text{for }u>0,
\]
and $A-\kappa R$ and $R$ vanish at $u=0$. In particular
$Q=d\,(A-aR)>0$ on $(0,\infty)$ whenever $q\ge a$.

\emph{Step 2.} $R\le\frac q2\,x$. Indeed, $\psi:=\frac q2x-R$ vanishes at $u=0$,
and since $x'=1-w^{q-1}+(1-q)uw^{q-2}$ we have $\psi'=\frac q2\zeta$ with
$\zeta:=1-w^{q-1}-(1-q)\,u\,w^{q-2}$; finally $\zeta(0)=0$ and
$\zeta'=(1-q)(2-q)\,u\,w^{q-3}\ge0$.

\emph{Step 3.} Let $q\ge a$, let $\theta\ge1$, and put $\lambda=-2q\theta$. Then
$2E'-\lambda=2q(\theta-z)$ and $\theta-z-\theta(1-z)=z(\theta-1)\ge0$, hence
\[
 u\,(2E'-\lambda)\ \ge\ 2q\theta\,x\ \ge\ 4\theta R
\]
by Step~2. Since $Q\ge0$ by Step~1, and then again by Step~1 with $\kappa=q$,
\begin{align}
 F_\lambda&\ \ge\ 4\theta R\,Q+R\,(4A-dR)
 =R\,\big[4\theta d\,(A-aR)+4A-dR\big]\notag\\
 &\ \ge\ R^2\,\big[4\theta d\,(q-a)+4q-d\big].\label{ex:eq-shiftF}
\end{align}

\emph{Step 4.} In case (a) take $\theta=\max\{1,b\}$. Then
$4\theta d(q-a)\ge4bd(q-a)\ge d-4q$ by the definition of $b$, so
$F_{-2L_*}\ge0$ by \eqref{ex:eq-shiftF}. Since $Q>0$ on $(0,\infty)$ and $E_q$
is convex, Proposition~\ref{ex:reduced} shows that $E_q$ satisfies the condition
with constant $-2L_*$.

\emph{Step 5.} For $\theta=1$ the bracket in \eqref{ex:eq-shiftF} equals
$4(d+1)q-(5d-4)$, which is nonnegative if and only if $q\ge q_d$. For $q>a$
this is the condition $b\le1$, that is, $L_*=q$. We have
$q_d\ge a$, because $d(5d-4)-4(d-1)(d+1)=(d-2)^2\ge0$, with equality only for
$d=2$. Hence for $q\ge q_d$, $q<1$, Step~3 applies with $\theta=1$ and gives
$F_{-2q}\ge0$, also in the case $d=2$, $q=\frac12$; here $Q>0$ on $(0,\infty)$
by Step~1, and Proposition~\ref{ex:reduced} applies as in Step~4. This proves
(b); $q_d<1$ if and only if $d<8$.

\emph{Step 6.} For $d=1$, $E_q=-f$ with $f(s)=(1+s)^q-1$, which is of the form
\eqref{ex:eq-f} with $\beta=0$ and $|m|=q$ by Lemma~\ref{ex:stieltjes}. Part (c)
is therefore a special case of Theorem~\ref{ex:bernstein}; it is the case
$d=1$ of Example~\ref{ex:ex-bern}(c). The last
assertion of the theorem is Corollary~\ref{tr:consequences} with
$\ezero=-q$ and $\Einf=0$.
\end{proof}

\begin{remark}\label{ex:rem-shifted}
(a) The restriction $q\ge\frac{d-1}{d}$ is necessary for every constant: since
$A/R\to q$ as $u\to\infty$, Step~1 shows that $Q=d(A-aR)\ge0$ on $(0,\infty)$ if
and only if $q\ge a$, and $Q\ge0$ is necessary by
Proposition~\ref{ex:reduced}. For the same reason $E_q$ satisfies
\eqref{ex:eq-Kd} if and only if $q\ge\kappa_d=\frac{3d-2}{2(d+2)}$; this range
is strictly smaller than $[q_d,1)$ for $d=1$ and $3\le d\le7$, equal to it
for $d=2$, and both ranges are empty for $d\ge8$, so that for the shifted
powers the direct argument above is at least as good as Theorem~\ref{ex:Kd}.

(b) No sharpness is claimed for $-2L_*$ when $L_*>q$.
\end{remark}

The constant $-2q$ cannot be obtained from Theorem~\ref{prelim:LMS} in the
whole range $q>\frac{d-1}d$. At the critical exponent the McCann quantity
$Q$ stays bounded, so the first term of $F_{-2q}$ grows at most linearly in
$u$, while the second term is negative of order $u^{2q}$, and $2q>1$ when
$d\ge3$.

\begin{proposition}[Failure of the constant $-2q$ near the critical exponent]\label{ex:counter}
Let $d\ge3$ and $q=\frac{d-1}{d}$. Then $F_{-2q}(u)\to-\infty$ as
$u\to\infty$ for $E=E_q$. Consequently there is $\eps>0$, depending on $d$, such
that for every $q\in\big[\frac{d-1}d,\frac{d-1}d+\eps\big)$ the function $E_q$
does not satisfy the Liero--Mielke--Savar\'e condition with constant $-2q$. For
$d=3$, $q=\frac23$ and $u=t^3-1$, $t>1$,
\begin{equation}\label{ex:eq-counter}
 F_{-4/3}(t^3-1)=-\frac{(t-1)^4\,\big(t^5-68t^4-140t^3-208t^2-112t-40\big)}{27\,t^5},
\end{equation}
which is negative for large $t$, for instance at $t=100$.
\end{proposition}

\begin{proof}
We use the notation of the proof of Theorem~\ref{ex:shifted}. As $u\to\infty$,
$R=w^{q-1}\big(1+(1-q)u\big)-1\sim(1-q)u^q$ and $A\sim q(1-q)u^q$, hence
\[
 R\,(4A-dR)\sim-(1-q)^2(d-4q)\,u^{2q},
\]
and $d-4q=(d-2)^2/d>0$ for $q=\frac{d-1}d$. On the other hand, for
$q=a=\frac{d-1}d$ Step~1 gives
$\frac{\rmd}{\rmd u}(A-aR)=q(1-q)(2-q)\,u\,w^{q-3}$, which is integrable on
$(0,\infty)$; thus $Q=d(A-aR)$ is bounded. Since
$0\le u(2E'+2q)=2q\,x\le2q\,u$, the first term of $F_{-2q}$ is $O(u)$. As
$2q>1$ for $d\ge3$, we conclude $F_{-2q}(u)\to-\infty$. For fixed $u$ the number
$F_{-2q}(u)$ depends continuously on $q$, so it remains negative for $q$ in a
right neighbourhood of $\frac{d-1}d$, and the condition fails there by
Proposition~\ref{ex:reduced}. Formula \eqref{ex:eq-counter} is obtained by
inserting $w=t^3$, $q=\frac23$, $d=3$ into \eqref{ex:eq-F}; at $t=100$ the
polynomial in the numerator is positive.
\end{proof}

\begin{remark}\label{ex:rem-counter}
(a) For the exponents of Proposition~\ref{ex:counter}, the functional with
density $E_q$ fails to be weakly geodesically $(-2q)$-convex on
$(\calM(\Om),\HK)$ for every closed $\Om\subset\R^d$ with nonempty interior. Indeed, $E_q\in C^2((0,\infty))$, so weak convexity
would imply \eqref{prelim:eq-LMScond} by the converse of
Theorem~\ref{prelim:LMS}, proved in Appendix~\ref{nc:sec-lms}. This does not decide whether
it is geodesically $(-2q)$-convex on $(\calP(X),\SHK)$: only the route
through Corollary~\ref{tr:global} is closed.

(b) The same argument shows that $-2q=2\ezero$ is the best constant in
\eqref{prelim:eq-LMScond} whenever it is admissible. For $\lambda=-2q+\delta$,
$\delta>0$, we have $u(2E'-\lambda)=-\delta u+O(u^2)$,
$Q=\frac{d+1}2q(1-q)u^2+O(u^3)$ and $R(4A-dR)=O(u^4)$ as $u\to0$, so that
$F_\lambda(u)<0$ for small $u>0$.

(c) For $3\le d\le7$ the infimum $q^*_d$ of the exponents $q_0$ such that
$E_q$ satisfies the condition with constant $-2q$ for all $q\in[q_0,1)$ lies
in $\big(\frac{d-1}{d},q_d\big]$, by Theorem~\ref{ex:shifted} and
Proposition~\ref{ex:counter}; its value is not known.
\end{remark}

\subsection{Further one-dimensional examples}\label{ex:ss-d1}

In one dimension the monotonicity condition in \eqref{prelim:eq-LMScond} is empty and
McCann's condition is convexity of $E$. For $d=1$ we have $Q=A$ and
$F_\lambda=A\,u\,(2E'-\lambda)+R\,(4A-R)$. Besides the one-dimensional
members of the class of Theorem~\ref{ex:bernstein}
(Example~\ref{ex:ex-bern}), the following two density functions can be
treated directly. The second density is increasing and has $\Einf>0$. For
the drift $\beta\ge0$ in Theorem~\ref{ex:bernstein}, every density in that
class is nonincreasing. Both examples have spherical rate $-2$, sharp on
the domains of diameter $>\pi/2$ specified above.

\begin{proposition}[Two examples for $d=1$]\label{ex:d1}
Let $d=1$.
\begin{enumerate}[label=\textup{(\alph*)}]
\item $E(c)=-\arctan c$ satisfies the Liero--Mielke--Savar\'e condition with
 constant $-2$; here $\ezero=-1$ and $\Einf=0$.
\item $E(c)=\sqrt{1+c^2}-1$ satisfies it with constant $0$; here $\ezero=0$
 and $\Einf=1$.
\end{enumerate}
Hence, by Corollary~\ref{tr:consequences}, $\sfE$ has the constant $-2$ on
$\HK$ in case \textup{(a)} and $0$ in case \textup{(b)}, and in both cases
the sharp constant $-2=2(\ezero-\Einf)$ on $\SHK$.
\end{proposition}

\begin{proof}
(b) Put $r=\sqrt{1+u^2}$. Then $E'=u/r$, $E''=r^{-3}$, so $E$ is convex,
$\ezero=0$, $\Einf=1$, $A=u^2/r^3>0$ and $R=u^2/r-r+1=(r-1)/r$. With $\lambda=0$
and $u^2=(r-1)(r+1)$,
\begin{align*}
 F_0&=\frac{2u^4}{r^4}+\frac{r-1}{r}\Big(\frac{4u^2}{r^3}-\frac{r-1}{r}\Big)
 =\frac{(r-1)^2}{r^4}\Big(2(r+1)^2+4(r+1)-r^2\Big)\\
 &=\frac{(r-1)^2\,(r^2+8r+6)}{r^4}\ \ge0 .
\end{align*}
Since $Q=A>0$ on $(0,\infty)$, Proposition~\ref{ex:reduced} gives the condition
with constant $0$.

(a) Here $E'=-(1+u^2)^{-1}\in[-1,0)$, $E''=2u(1+u^2)^{-2}\ge0$, $\ezero=-1$,
$\Einf=0$, $A=\frac{2u^3}{(1+u^2)^2}>0$ for $u>0$, and
\[
 R=L(u):=\arctan u-\frac{u}{1+u^2},\qquad L'(u)=\frac{2u^2}{(1+u^2)^2} .
\]
With $\lambda=-2$ we have $u(2E'+2)=\frac{2u^3}{1+u^2}$ and
\[
 F_{-2}=\frac{4u^6}{(1+u^2)^3}+L\,\Big(\frac{8u^3}{(1+u^2)^2}-L\Big).
\]
We have $0\le L<\frac\pi2$, and $L'\le2u^2$ gives $L\le\frac23u^3$. Put
$h(u)=\frac{8u^3}{(1+u^2)^2}$. For $0<u\le1$, $h(u)\ge2u^3\ge L(u)$, so
$F_{-2}\ge0$. The derivative of $h$ has the sign of $3-u^2$, so on $[1,3]$ we
have $h\ge\min\{h(1),h(3)\}=\min\{2,\frac{54}{25}\}=2>\frac\pi2>L$, and again
$F_{-2}\ge0$. For $u\ge3$,
$\frac{4u^6}{(1+u^2)^3}\ge4\big(\tfrac9{10}\big)^3>2.9>\frac{\pi^2}{4}>L^2$, so
$F_{-2}>0$. Since $Q=A>0$ on $(0,\infty)$, Proposition~\ref{ex:reduced} gives the
condition with constant $-2$.
\end{proof}

\subsection{Summary}\label{ex:ss-table}

Table~\ref{ex:table} lists the density functions for which $\sfE$ is now known
to be geodesically $\lambda$-convex on $(\calP(X),\SHK)$, $X\subset\R^d$ compact
and convex with nonempty interior; in the rows with finite $\ezero$ the
same constants hold on $(\calP(\Om),\SHK)$ for closed convex $\Om$ with
nonempty interior. Every
row above the rule rests on Corollary~\ref{tr:consequences}. Since $\Nfun$ is
linear in the pair $(E,\lambda)$, positive combinations of such rows of the same
dimension are again admissible, and the constants add; by
Proposition~\ref{ex:mixture} the same is true for dilations and mixtures of rows
with bounded $E'$. In all rows with finite $\ezero$ in which the constant equals
$2(\ezero-\Einf)$, it is sharp.

The first row is the list known from \cite{LaM22}. The three rows for
$1-(1+c)^q$ overlap: the first two give the constant $-2q$ on their ranges
of exponents, and the third covers all exponents above the critical one. The
row below the rule is of a different kind: for the powers $c^m$ there is no
global constant on large domains, and the entry refers to convexity on
sublevels of the energy.

\begin{table}[ht]
\centering
\renewcommand{\arraystretch}{1.45}
\small
\begin{tabular}{@{}p{6.3cm}p{1.5cm}p{3.8cm}p{2.5cm}@{}}
\toprule
density function $E(c)$ & $d$ & constant on $(\calP(X),\SHK)$ & reference\\
\midrule
$-c^q$,\ \ $\max\{\frac{d-1}d,\frac d{d+2}\}\le q\le\frac12$ & $1,\,2$ & $0$ & \cite{LaM22}\\
$-c^{\frac{d-1}d}f(c^{1/d})$, $f$ as in \eqref{ex:eq-f}; e.g.\
 $-\frac{c}{1+c^{1/d}}$,\ $-c^{\frac{d-1}d}\log(1+c^{1/d})$,\
 $c^{\frac{d-1}d}[1-(1+c^{1/d})^q]$ & any & $-2|m|=2(\ezero-\Einf)$ & Thm.~\ref{ex:bernstein}, Ex.~\ref{ex:ex-bern}\\
$1-(1+c)^q$,\ \ $0<q<1$ & $1$ & $-2q$ & Thm.~\ref{ex:shifted}(c)\\
$1-(1+c)^q$,\ \ $\frac{5d-4}{4d+4}\le q<1$ & $1$ to $7$ & $-2q$ & Thm.~\ref{ex:shifted}(b)\\
$1-(1+c)^q$,\ \ $\frac{d-1}d<q<1$ & any & $-2q\max\{1,b\}$, $b=\max\big\{0,\frac{1/4-q/d}{q-(d-1)/d}\big\}$ & Thm.~\ref{ex:shifted}(a)\\
$\sqrt c\,G(\sqrt c)$ convex, $G$ convex and Lipschitz & $2$ & $2\big(G'(0^+)-G'(\infty)\big)$ & Thm.~\ref{ex:d2}\\
$E$ convex with \eqref{ex:eq-Kd}, $\ezero>-\infty$, $\Einf<\infty$ & any & $2(\ezero-\Einf)$ & Thm.~\ref{ex:Kd}\\
$-\arctan c$,\ \ $\sqrt{1+c^2}-1$ & $1$ & $-2$ & Prop.~\ref{ex:d1}\\
$\int aE(c/a)\,m(\rmd a)$, $E$ any row with bounded $E'$ & as for $E$ & $|m|$ times the constant of $E$ & Prop.~\ref{ex:mixture}\\
\midrule
$c^m$,\ \ $m>1$ & any & on sublevels of $\sfE$ only & Thm.~\ref{tr:powers}\\
\bottomrule
\end{tabular}
\caption{Density functions $E$ for which $\sfE$ is geodesically $\lambda$-convex
on $(\calP(X),\SHK)$; ``any'' means every $d\ge1$. The last row is not globally semiconvex if
$\diam X>\pi/2$ (Theorem~\ref{tr:necessary}); see Theorem~\ref{tr:powers} for
the precise statement.}\label{ex:table}
\end{table}

\section[Compact domains: corrections and completion of an existence proof]{Compact domains: corrections and completion of the existence proof of \cite{LaM22}}\label{sec:compact}

Let $X\subset\R^d$ be compact and convex with nonempty interior, $V=|X|$, and
$E:[0,\infty)\to\R$ finite, convex and lower semicontinuous, with $E(0)=0$.
This section proves the existence theorem of \cite{LaM22} for $\EVI_\lambda$
flows on $(\calM(X),\HK)$ and $(\calP(X),\SHK)$ under \emph{weak} geodesic
$\lambda$-convexity of $\sfE$ (Definition~\ref{prelim:convexity}(c)), which is
the hypothesis of Assumption~A of \cite{LaM22}. A posteriori $\sfE$ is then
strongly geodesically $\lambda$-convex (Corollary~\ref{cp:weakstrong}).
Sections~\ref{sec:hk-unbdd}--\ref{sec:sublevel} are built on this theorem,
so we give the proof in detail. The estimates for the minimizing movement
scheme were developed by Savar\'e \cite{Sav07} and are part of the
forthcoming work of Muratori and Savar\'e \cite{MuSII}. As in \cite{LaM22},
we state them without proof, in Sections~\ref{cp:ss-geom}
and~\ref{cp:ss-error}, and refer to these works; the one exception is the
corrected error estimate, Lemma~\ref{cp:error}, which is proved. Four points of \cite{LaM22}
need correction; they are stated in Section~\ref{cp:ss-points}.
Remark~\ref{cp:rem-scope} lists what is still taken from the literature.

The proof follows the plan of \cite{LaM22} and has four parts.
\begin{enumerate}[label=\textup{(\arabic*)}]
\item One minimisation step keeps density bounds
 (Section~\ref{cp:ss-density}), so that on bounded time intervals the
 scheme stays in one of the sets $\calM_\delta(X)$ of \eqref{cp:eq-Mdelta}.
\item On these sets the squared distance is semiconcave and the local angle
 condition holds. This turns the Euler--Lagrange inequality of a step into
 an $\EVI$ for the geodesic interpolant of the scheme, up to an explicit
 error (Sections~\ref{cp:ss-geom} and~\ref{cp:ss-error},
 Proposition~\ref{cp:approxEVI}).
\item The error is of order $\tau$ (Lemma~\ref{cp:error}), so the
 interpolants form a Cauchy family and converge by completeness.
\item The limit is an $\EVI_\lambda$ solution
 (Sections~\ref{cp:ss-abstract} and~\ref{cp:ss-verify}), first for initial
 data with bounded density and finite slope, then for all initial data,
 which are approximated by such data in distance and in energy
 (Lemma~\ref{cp:recovery}).
\end{enumerate}
Parts (2)--(4) are formulated in a complete geodesic space
(Theorem~\ref{cp:abstract}). Convexity of $\sfE$ is assumed along one
geodesic for each pair of measures only, so the proof has to say which
geodesic is used; this is the purpose of the admissible geodesics of
Section~\ref{cp:ss-error}.

\subsection{Four corrections}\label{cp:ss-points}
Write
\begin{equation}\label{cp:eq-Mdelta}
 \calM_\delta(X)=\{\rho\Leb:\delta\le\rho\le\delta^{-1}\},\quad0<\delta<1.
\end{equation}
We revisit the lower density bound for one $\HK$ step, the accumulated
error estimate, the statements quoted without proof, and the hypotheses of
the abstract existence theorem.

\emph{Lower density.} \cite[Prop.~4.3]{LaM22} bounds the density of a
minimiser of one $\HK$ step from below in terms of a lower bound for the
density of $\mu_0$. The proof ends, in the case
$a:=2\tau E'(c_{\rm low})>0$ and with $k_\eps:=(1+a)^{-1}-\eps$, with the
comparison
$k_\eps-(1-a)\le k_\eps-(1+a)^{-1}<0$. Its first inequality is reversed:
$(1+a)^{-1}-(1-a)=a^2/(1+a)>0$. The statement itself is correct. The
competitor $\mu_1+(1-t^2)\mu_0|_A$ used there has to be replaced by
$\mu_1+h^2(1-t^2)\mu_0|_A$ with $h=(1+2\tau\max\{E'(c_{\rm low}),0\})^{-1}$;
Proposition~\ref{cp:density-hk} proves the original lower bound in this way.

\emph{Accumulated error.} \cite[Lemma~5.17]{LaM22} bounds the weighted time
integral of the error function in the approximate $\EVI$ by
$\tau(4+\tau\kappa)\slope\phi^2(x_0)$, where $\kappa$ is the semiconcavity
constant of the squared distance.

\begin{example}[Accumulated error]\label{cp:ex-error}
The bound in \cite[Lemma~5.17]{LaM22} is false. On $[0,1]$ take
$\phi(x)=-x$, $x_0=0$, $\lambda=0$, $A_\kappa=[0,1]$, $\kappa\ge1$.
All its assumptions hold and $x_n=\min\{n\tau,1\}$.
For $\tau=1/100$, $T=1/10$, $\kappa=100$, its error function has
$\Er_0=2+\tau^2$ and $\Er_n=1+\tau^2$ for $1\le n\le9$.
Thus $\int_0^T\ee^{-4t}\Er(t)\dd t\ge33/500>1/20$, whereas the claimed
bound is $\tau(4+\tau\kappa)\slope\phi^2(x_0)=1/20$.
Lemma~\ref{cp:error} gives a bound of order $\tau$ for fixed $\kappa$, with
the factor $\kappa\tau$ in place of $\kappa\tau^2$; the example shows that
this factor is necessary.
\end{example}

\emph{Statements quoted without proof.} The derivative estimate for the
squared distance \cite[Prop.~5.7]{LaM22} and the estimates for the scheme
\cite[Lemma~5.13, Prop.~5.14, Lemma~5.17]{LaM22} were developed by
Savar\'e \cite{Sav07} and are part of the forthcoming work
\cite{MuSII}; \cite{LaM22} states them without proof and refers to these
works. We do the same for Lemmas~\ref{cp:deriv} and~\ref{cp:euler} and
Proposition~\ref{cp:approxEVI}, which give these statements in the form
needed here. Compared with
\cite{LaM22}, convexity of the functional is required along one geodesic
for each pair of points only. For the proofs we refer to
\cite{MuSII}. Example~\ref{cp:ex-error} shows that the bound of
\cite[Lemma~5.17]{LaM22} is false in the form printed there; we make no
statement about its form in \cite{MuSII}. Lemma~\ref{cp:error} is a
corrected form, which we prove from Lemma~\ref{cp:euler} by the
telescoping argument indicated in \cite{LaM22}.

\emph{Hypotheses of the abstract theorem.} \cite[Thm.~5.9]{LaM22} requires
nested sets $A_\kappa\subset\dom\slope\phi$, compact sublevels of $\phi$, and
the scheme for all steps $\tau>0$. Its application in
\cite[Thm.~6.6]{LaM22} takes for $A_\kappa$ the density-bounded sets
$\calM_\delta(X)$. The next two examples show that these requirements are
not satisfied there.

\begin{example}[Density bounds do not imply finite slope]\label{cp:ex-slope}
On $[-1,1]$ let $E(c)=c^2$ and let $\rho$ equal $1$ on the left half and
$2$ on the right. Replacing it by $3/2$ on $(-\eps,\eps)$ drops the energy
by $\eps/2$. Monotone rearrangement of the $3\eps$ mass in that interval
moves it at most $2\eps$, hence costs at most $12\eps^3$ in HK squared.
The slope is infinite, although $\rho\in\calM_{1/2}$.
For SHK, densities $1/4,3/4$, averaged to $1/2$ near zero, and
$E=-\sqrt c$ give an energy drop
$\eps(\sqrt2-(1+\sqrt3)/2)>0$ at HK squared cost at most $4\eps^3$;
$\SHK\le2\HK$ gives the same conclusion.
\end{example}

\begin{example}[Noncompact sublevels]\label{cp:ex-sublevel}
For $E(c)=-c$, the energy is $(-2)$-convex by the mass identity, but
$\sfE(n\Leb|_X)=-nV$ and $\HK(n\Leb|_X,0)=\sqrt{nV}\to\infty$.
For nonzero $\mu$, its step functional is unbounded below along $r^2\mu$
when $\tau\ge1/2$. Thus compact sublevels and solvability for every step,
as required in \cite[Thm.~5.9]{LaM22}, cannot be used here.
\end{example}
We replace the invariant sets by density \emph{and slope} bounded sets,
and compact sublevels by completeness, with exact steps for small $\tau$.
The smaller sets are large enough, because after one step of the scheme
from a density-bounded measure the slope is finite
(Lemma~\ref{cp:euler}\ref{cp:it-abc}).
Strict improvement in the density arguments below comes from distance,
so ordinary convexity of $E$ suffices even when $E$ is affine.

\subsection{Density bounds for one minimisation step}\label{cp:ss-density}

All proofs in this subsection have the same form: if a bound fails on a set
of positive measure, we change the minimiser $\mu_1$ a little on that set
and obtain a competitor with a smaller value of the step functional.

We use the optimal transport--growth pair $(\bfT,q)$ of
Lemma~\ref{prelim:tg} from $\mu_0=\rho_0\Leb\in\calM(X)$, $\rho_0>0$ a.e.\
in $X$, to $\mu_1\in\calM(X)$, and its transport cosine $\chi$; for
$\mu_1=0$ we take $q\equiv0$, $\bfT=\mathrm{Id}$. The competitors are
split-source systems: the mass of $\mu_0$ at a point is divided into parts
that are moved and rescaled separately, and \eqref{prelim:eq-tg-split}
bounds the cost of such a system.

Two elementary consequences of convexity are used repeatedly. If $c>0$ and
$s\in\partial E(c)$, then, by monotonicity of the subdifferential,
\begin{equation}\label{cp:eq-conv}
 E(y)-E(x)\le s\,(y-x)\ \ \text{for }0\le x\le y\le c,
 \qquad
 E(x)-E(y)\ge s\,(x-y)\ \ \text{for }c\le y\le x .
\end{equation}
Moreover $E(c+s)\le E(c)+\Einf s$ for $c,s\ge0$, and a supporting line at a
positive density gives constants $a,b\ge0$ with
\begin{equation}\label{cp:eq-affine}
 E(c)\ge-ac-b\ \ (c\ge0),\qquad\text{hence}\qquad
 \sfE(\mu)\ge-a\,\mu(X)-bV\ \ (\mu\in\calM(X)),
\end{equation}
including singular mass, because $\Einf\ge-a$. No assumption on
$\ezero=E'(0^+)$ is made.

The first result is common to both spaces: the optimal pair never carries
mass to a point where the density of the minimiser is higher than at the
point of departure. In the proof a fraction of the mass of a set $A$ is
kept in place instead of being transported.

\begin{lemma}[No transport towards higher density]\label{cp:nouphill}
Let $\tau>0$ and $\mu_0=\rho_0\Leb\in\calM(X)$ with $\rho_0>0$ a.e.\ in $X$.
Let $\mu_1$ be a minimiser of
\[
 J(\mu)=\frac1{2\tau}\HK^2(\mu_0,\mu)+\sfE(\mu)\quad\text{over }\calM(X),
\]
or let $\mu_0\in\calP(X)$ and let $\mu_1$ be a minimiser of
$\frac1{2\tau}\SHK^2(\mu_0,\mu)+\sfE(\mu)$ over $\calP(X)$. Then
$\mu_1=\rho_1\Leb$ is absolutely continuous, and every optimal
transport--growth pair $(\bfT,q)$ from $\mu_0$ to $\mu_1$ satisfies
\begin{equation}\label{cp:eq-nouphill}
 \rho_1(\bfT(x))\le\rho_1(x)\qquad\text{for $\mu_0$-a.e.\ $x$ with $q(x)>0$.}
\end{equation}
\end{lemma}

Here $\rho_1$ is any Borel representative; since
$\bfT\pfw(q^2\mu_0)=\mu_1\ll\Leb$, the statement does not depend on the choice.

\begin{proof}
The minimiser has finite energy: use $0$ as competitor for HK, the uniform
probability for SHK, and \eqref{cp:eq-affine} as lower bound. If $\mu_1=0$
there is nothing to prove; otherwise let $(\bfT,q)$ be an optimal pair.

\emph{Step 1: a mass-preserving comparison.}
For a Borel set $A$ put $\alpha=\mu_0|_A$,
$\beta=\bfT\pfw(q^2\alpha)\le\mu_1$, $m=\alpha(X)$,
$n=\beta(X)=\int_Aq^2\dd\mu_0$. Suppose $m>0$, and $q>0$ and $\bfT(x)\ne x$
a.e.\ on $A$. Then $n>0$ and $\chi<1$ on $A$, and the Cauchy--Schwarz
inequality $\int_Aq\dd\mu_0\le\sqrt{mn}$ gives
\[
 C_A:=\int_A(1+q^2-2q\chi)\dd\mu_0
 >m+n-2\int_Aq\dd\mu_0\ge m+n-2\sqrt{mn}=(\sqrt m-\sqrt n)^2.
\]
For $0<t<1$ consider the split-source system that keeps $(\bfT,q)$ on
$\mu_0|_{X\setminus A}$ and on $(1-t)\alpha$, and uses the stationary pair
$(\mathrm{Id},\sqrt{n/m})$ on $t\alpha$. Its target
$\mu_1^t=\mu_1-t\beta+t(n/m)\alpha$ has the mass of $\mu_1$, and by
\eqref{prelim:eq-tg-split} and \eqref{prelim:eq-tg-opt}
\[
 \HK^2(\mu_0,\mu_1^t)\le\HK^2(\mu_0,\mu_1)-t\big[C_A-(\sqrt m-\sqrt n)^2\big]
 <\HK^2(\mu_0,\mu_1).
\]
On probabilities $\SHK$ is a strictly increasing function of $\HK$
(Proposition~\ref{prelim:cone}). Hence in both cases minimality is
contradicted as soon as $\sfE(\mu_1^t)\le\sfE(\mu_1)$ for some $t\in(0,1)$.

\emph{Step 2: absolute continuity.}
In this step $\rho_1$ denotes the density of the absolutely continuous part
of $\mu_1$. Suppose $\mu_1$ has a singular part $\mu_1^s\ne0$, concentrated on a
Lebesgue null set $B$. Then $\Einf<\infty$, because the energy is finite.
Since $\int_{\bfT^{-1}(B)}q^2\dd\mu_0=\mu_1(B)>0$ and $\mu_0(B)=0$, the set
$\{q>0\}\cap\bfT^{-1}(B)\setminus B$ has positive $\mu_0$-measure. Choose a
subset $A$ of it with $\mu_0(A)>0$ and constants $M,K$ with
$\rho_0\le M$, $\rho_1\le K$ on $A$. On $A$ we have $q>0$ and
$\bfT(x)\ne x$, and $\beta$ is concentrated on $B$, hence
$\beta\le\mu_1^s$. Removing $t\beta$ from the singular part saves
$\Einf tn$. For $t(n/m)M\le1$, adding the density $t(n/m)\rho_0\le1$ on
$A$ costs at most $E'_+(K+1)\,tn\le\Einf tn$ by \eqref{cp:eq-conv}. The
energy does not increase, a contradiction.
Hence $\mu_1$ is absolutely continuous.

\emph{Step 3: no uphill transport.}
If \eqref{cp:eq-nouphill} fails, there are numbers $0\le a<b$, $M<\infty$ and
a set $A$ with $\mu_0(A)>0$ on which $q>0$,
$\rho_0\le M$ and $\rho_1(x)\le a<b\le\rho_1(\bfT(x))$; in particular
$\bfT(x)\ne x$ on $A$. Then $\alpha$ is concentrated on $\{\rho_1\le a\}$
and $\beta$ on $\{\rho_1\ge b\}$.
Choose $c\in(a,b)$ and $e\in\partial E(c)$. For small $t$ the density of
$\mu_1^t$ on $A$ is at most $a+t(n/m)M\le c$; on $\{\rho_1\ge b\}$ it is at
least $(1-t)\rho_1\ge(1-t)b\ge c$, because $\beta\le\mu_1$ has density at
most $\rho_1$. By \eqref{cp:eq-conv} the added mass costs at most $etn$ and
the removed mass saves at least $etn$.
The energy does not increase, again a contradiction.
\end{proof}

We come to the bounds for $\HK$. If the lower bound fails on a set $B$,
then by Lemma~\ref{cp:nouphill} the mass that starts in $B$ and survives
stays in $B$, so mass has been destroyed on part of $B$. The competitor
leaves a small fraction of the source there in place, shrunk by $h^2$; the
gain in squared distance is of first order and exceeds the energy of the
added mass. For the upper bound the roles are exchanged: the competitor
reduces the growth factor where it is large and thereby removes mass from
the region of high density. If $e<0$ this costs energy, and the hypothesis
$1+2\tau e>0$ makes the gain in distance larger than the cost. The
thresholds $b$ in (i) and $a$ in (ii) are the densities at which the slope
$e$ of $E$ is used to price the added or removed mass; in (iii) they are
taken equal to the current bounds.

\begin{proposition}[Density bounds for $\HK$]\label{cp:density-hk}
Let $\tau>0$, let $\mu_0=\rho_0\Leb\in\calM(X)$ with $\rho_0>0$ a.e.\ in $X$,
and let $\mu_1$ be a minimiser of
$J(\mu)=\frac1{2\tau}\HK^2(\mu_0,\mu)+\sfE(\mu)$ over $\calM(X)$. Then
$\mu_1=\rho_1\Leb$, and the following holds.
\begin{enumerate}[label=\textup{(\roman*)}]
\item\label{cp:it-low} \emph{Lower bound.} If $\rho_0\ge c_{\min}>0$ a.e., then
for every $b>0$ and every $e\in\partial E(b)$,
\begin{equation}\label{cp:eq-low}
 \rho_1\ge\min\Big\{b,\ \frac{c_{\min}}{(1+2\tau\max\{e,0\})^2}\Big\}
 \qquad\text{a.e.\ in }X .
\end{equation}
In particular $\mu_1\ne0$.
\item\label{cp:it-upp} \emph{Upper bound.} If $\rho_0\le c_{\max}<\infty$ a.e.,
then for every $a>0$ and every $e\in\partial E(a)$ with $1+2\tau e>0$,
\begin{equation}\label{cp:eq-upp}
 \rho_1\le\max\Big\{a,\ \frac{c_{\max}}{(1+2\tau\min\{e,0\})^2}\Big\}
 \qquad\text{a.e.\ in }X .
\end{equation}
\item\label{cp:it-iter} \emph{Iteration.} Suppose $0<l\le\rho_0\le u<\infty$.
Choose $e_l\in\partial E(l)$ and $e_u\in\partial E(u)$ and put $L:=\max\{e_l,0\}$, $U:=\max\{-e_u,0\}$. If $2\tau U<1$ and
$\mu_n=\rho_n\Leb$, $n=1,2,\dots,N$, are successive minimisers
$($$\mu_n$ minimises $\frac1{2\tau}\HK^2(\mu_{n-1},\cdot)+\sfE$$)$, then
\begin{equation}\label{cp:eq-iter}
 \frac{l}{(1+2\tau L)^{2n}}\le\rho_n\le\frac{u}{(1-2\tau U)^{2n}}
 \qquad\text{a.e.\ in }X,\ n\le N,
\end{equation}
and, if $\tau U\le\frac14$,
\begin{equation}\label{cp:eq-iterexp}
 l\,\ee^{-4Ln\tau}\le\rho_n\le u\,\ee^{8Un\tau}\qquad\text{a.e.\ in }X,\ n\le N .
\end{equation}
\end{enumerate}
\end{proposition}

\begin{proof}
Use Lemma~\ref{cp:nouphill} and an optimal pair $(\bfT,q)$, taking $q=0$,
$\bfT=\mathrm{Id}$ if $\mu_1=0$.

\emph{Lower bound.} Put $h=(1+2\tau\max\{e,0\})^{-1}$.
If the bound fails, choose $\xi<\min\{b,h^2c_{\min}\}$ with
$B=\{\rho_1\le\xi\}$ of positive measure. By \eqref{cp:eq-nouphill},
$\bfT(x)\in B$ for a.e.\ $x\in B$ with $q(x)>0$, so that
$\bfT\pfw(q^2\mu_0|_B)\le\mu_1$ is concentrated on $B$ and
$\int_Bq^2\dd\mu_0\le\mu_1(B)\le\xi|B|\le(\xi/c_{\min})\mu_0(B)$.
Hence for $\sqrt{\xi/c_{\min}}<k<h$ there is a set
$A\subset B$ of positive measure with $q\le k$ and $\rho_0\le M<\infty$ on
$A$. Put $\alpha=\mu_0|_A$, $m=\alpha(X)$ and
$\mu_1^t=\mu_1+h^2(1-t^2)\alpha$.
Split the source on $A$ into $t^2\alpha$, transported by $(\bfT,q/t)$,
and $(1-t^2)\alpha$, transported by $(\mathrm{Id},h)$; the target is
$\mu_1^t$ and the cost on $A$ is
$\int_A[t^2+q^2-2tq\chi]\dd\mu_0+(1-t^2)(1-h)^2m$. By \eqref{prelim:eq-tg-split}
and $q\chi\le k$ on $A$, the squared distance improves by at least
\[
 (1-t^2)(2h-h^2)m-2(1-t)\int_Aq\chi\dd\mu_0
 \ge(1-t)m\big[(1+t)(2h-h^2)-2k\big].
\]
For $t$ close to $1$ the density of $\mu_1^t$ on $A$ is at most
$\xi+h^2(1-t^2)M\le b$, so by \eqref{cp:eq-conv} the energy cost is at most
$eh^2(1-t^2)m$. Consequently, using $(1+2\tau e)h^2\le h$,
\[
 J(\mu_1)-J(\mu_1^t)
 \ge\frac{(1-t)m}{2\tau}[(1+t)(2h-(1+2\tau e)h^2)-2k]
 \ge\frac{(1-t)m}{2\tau}[(1+t)h-2k]>0
\]
for $t$ close to $1$, since $k<h$. The argument includes the case $\mu_1=0$
$($then $B=X$$)$ and thus excludes the zero minimiser.

\emph{Upper bound.} Put $h=(1+2\tau\min\{e,0\})^{-1}\ge1$.
Failure gives $B=\{\rho_1\ge\xi\}$ of positive measure with
$\xi>\max\{a,h^2c_{\max}\}$. For
$A_0=\{q>0\}\cap\bfT^{-1}(B)$ we have $A_0\subset B$ up to a null set by
\eqref{cp:eq-nouphill}, hence
$\int_{A_0}q^2\dd\mu_0=\mu_1(B)\ge\xi|B|\ge\xi|A_0|
\ge(\xi/c_{\max})\mu_0(A_0)$ and $\mu_0(A_0)>0$.
Choose $h<k<\sqrt{\xi/c_{\max}}$ and a subset $A\subset A_0$ of positive
measure on which $q\ge k$, and let
$\beta=\bfT\pfw(q^2\mu_0|_A)$, $n=\beta(X)$.
Reducing $q$ to $tq$ on $A$ gives the target
$\mu_1^t=\mu_1-(1-t^2)\beta$ and improves the squared cost
by $(1-t)\int_Aq^2[(1+t)-2\chi/q]\dd\mu_0\ge(1-t)n(1+t-2/k)$.
The measure $\beta\le\mu_1$ is concentrated on $B$ and has density at most
$\rho_1$, so the density of $\mu_1^t$ on $B$ is at least
$t^2\rho_1\ge t^2\xi\ge a$ for $t$ close to $1$, and by \eqref{cp:eq-conv}
the energy saving is at least $e(1-t^2)n$.
Thus
$J(\mu_1)-J(\mu_1^t)\ge(1-t)n[(1+t)(1+2\tau e)-2/k]/(2\tau)>0$ for $t$
close to $1$, because $1+2\tau e\ge1/h>1/k$;
a contradiction.

\emph{Iteration.} Apply the two bounds successively at thresholds
$b_n=l(1+2\tau L)^{-2n}$ and $a_n=u(1-2\tau U)^{-2n}$, with
$c_{\min}=b_n$ and $c_{\max}=a_n$.
Monotonicity of $\partial E$ allows subgradients at $b_n\le l$ no larger than
$e_l$, and at $a_n\ge u$ no smaller than $e_u$; for the latter
$1+2\tau\min\{e,0\}\ge1-2\tau U>0$. This proves \eqref{cp:eq-iter} by
induction.
The exponential bounds use $\log(1+r)\le r$ and
$-\log(1-r)\le2r$ for $0\le r\le1/2$.
\end{proof}

\begin{remark}[Sharpness]\label{cp:rem-sharp}
For $E(c)=ec$ with $e\ge0$, put $h=1+2\tau e$. Equation~\eqref{prelim:eq-scaling}
gives $J(\nu)=\frac{h}{2\tau}\HK^2(h^{-2}\mu_0,\nu)+\frac{e}{h}\mu_0(X)$
for every $\nu\in\calM(X)$. Thus $h^{-2}\mu_0$ is the unique minimiser.
When $\rho_0=c_{\min}$ a.e., the lower bound \eqref{cp:eq-low} is attained
for every $b\ge c_{\min}/h^2$.
\end{remark}

On the sphere the total mass is fixed, so mass cannot be added or removed
on one set. If one of the bounds fails, the unit masses of $\mu_0$ and
$\mu_1$ force a set $A$ whose mass grows and arrives where the density of
$\mu_1$ is above a level $b$, and a set $B$ whose mass shrinks and where
the density of $\mu_1$ is below $b$. As in \cite{LaM22}, the competitor
reduces the growth on $A$ and keeps a little of the mass of $B$ in place,
by equal amounts. No derivative of $E$ enters, which is why the interval
$[m,M]$ is preserved without loss and for every $\tau$. We refer to
Proposition~\ref{cp:density-shk} as the \emph{exact maximum principle}.

\begin{proposition}[Density bounds for $\SHK$]\label{cp:density-shk}
Let $\tau>0$, let $\mu_0=\rho_0\Leb\in\calP(X)$ with
$0<m\le\rho_0\le M<\infty$ a.e.\ in $X$, and let $\mu_1$ be a minimiser of
\[
 \nu\mapsto\frac1{2\tau}\SHK^2(\mu_0,\nu)+\sfE(\nu)\qquad\text{over }\calP(X).
\]
Then $\mu_1=\rho_1\Leb$ with $m\le\rho_1\le M$ a.e.\ in $X$. Consequently all
iterates of the scheme started at $\mu_0$ have densities in $[m,M]$.
\end{proposition}

\begin{proof}
By Lemma~\ref{cp:nouphill}, $\mu_1=\rho_1\Leb$ and
$\rho_1\circ\bfT\le\rho_1$ a.e.\ after redefining $\bfT(x)=x$ where $q=0$,
which changes neither the target nor the cost.
If a bound fails, there are $b\in\{m,M\}$, $\eps\in(0,1)$ and
positive-measure sets $A,B$ with
\begin{equation}\label{cp:eq-AB}
 q^2>1+\eps,\ \rho_1(\bfT(x))>(1+\eps)b\quad(A),\qquad
 q^2<1-\eps,\ \rho_1(x)<(1-\eps)b\quad(B).
\end{equation}
Here is the selection argument. No uphill transport implies
$\{\rho_1<b\}\subset\bfT^{-1}\{\rho_1<b\}$ and the analogous inclusion
with $\le$, up to null sets. If the lower bound fails, $b=m\le\rho_0$ gives
$\int_{\{\rho_1<b\}}q^2\dd\mu_0\le\mu_1(\rho_1<b)<\mu_0(\rho_1<b)$, hence
a low-growth subset there.
The strict inequality also holds for $\{\rho_1\le b\}$; passing to
complements and using the unit masses gives
$\int_{\bfT^{-1}\{\rho_1>b\}}q^2\dd\mu_0=\mu_1(\rho_1>b)
>\mu_0(\bfT^{-1}\{\rho_1>b\})$, hence a high-growth subset of
$\bfT^{-1}\{\rho_1>b\}$.
If the upper bound fails, $b=M\ge\rho_0$ gives
$\mu_1(\rho_1>b)>\mu_0(\rho_1>b)$ and
$\bfT^{-1}\{\rho_1>b\}\subset\{\rho_1>b\}$, hence a high-growth
subset; complements of $\{\rho_1\ge b\}$ give the low-growth subset.
Countable exhaustion provides a common $\eps$. On $A$ we have
$\rho_1(x)\ge\rho_1(\bfT(x))>(1+\eps)b$, so $A$ and $B$ are disjoint.

Let $\alpha=\mu_0|_B$, $\beta=\bfT\pfw(q^2\mu_0|_A)\le\mu_1$, with masses
$a,n>0$; $\beta$ is concentrated on $\{\rho_1>(1+\eps)b\}$. For small
$s,t>0$ satisfying $sn=ta$, set
$\nu=\mu_1-s\beta+t\alpha\in\calP(X)$. Use $(\bfT,\sqrt{1-s}\,q)$ on $A$;
on $B$ split into $(1-t)\mu_0|_B$ with growth $q/\sqrt{1-t}$ and the
stationary source $t\mu_0|_B$ with growth $1$; keep $(\bfT,q)$ elsewhere.
The target is $\nu$. Writing $J_A=\int_Aq\chi\dd\mu_0$,
$J_B=\int_Bq\chi\dd\mu_0$, the squared-cost gain is, by \eqref{prelim:eq-tg-split},
at least
\[
 sn-2(1-\sqrt{1-s})J_A+ta-2(1-\sqrt{1-t})J_B>0
\]
for small $s,t$, since $J_A\le n/\sqrt{1+\eps}$ and
$J_B\le a\sqrt{1-\eps}$ by \eqref{cp:eq-AB}.
Decrease $s,t$ further so that $(1-s)(1+\eps)\ge1$ and $tM\le\eps b$. Where
$\beta$ is concentrated, the density of $\nu$ is at least
$(1-s)\rho_1\ge b$, because the density of $\beta$ is at most $\rho_1$; on
$B$ it is at most $(1-\eps)b+tM\le b$.
For $e\in\partial E(b)$, \eqref{cp:eq-conv} bounds the energy increase by
$e(ta-sn)=0$.
SHK strictly increases with HK on probabilities, so this contradicts
minimality.
\end{proof}

\subsection{Angles and the derivative of the squared distance}\label{cp:ss-geom}

From here to the end of Section~\ref{cp:ss-abstract} we work in the abstract
setting of \cite[Sec.~5]{LaM22}.
We state the derivative estimate for the squared distance to an observer
$v$ along a geodesic (Lemma~\ref{cp:deriv}), which is
\cite[Prop.~5.7]{LaM22}. It belongs to the approach of Savar\'e
\cite{Sav07}; for its proof we refer to \cite{MuSII}.
Let
$(\calX,\sfd)$ be a geodesic metric space; $\mathrm{Geod}(x,y)$ denotes the set
of constant-speed geodesics $\gamma_{xy}:[0,1]\to\calX$ from $x$ to $y$. For
$A\subset\calX$ the \emph{geodesic cover} is
$A^{\rm g}:=\{\gamma_{xy}(t):x,y\in A,\ \gamma_{xy}\in\mathrm{Geod}(x,y),\
t\in[0,1]\}$; it contains $A$. For two geodesics emanating from the same
point and $s,t\in(0,1]$ put
\[
 \Theta(s,t):=\frac1{2st}\Big(\sfd^2(x,\gamma_{xy}(s))
   +\sfd^2(x,\gamma_{xz}(t))-\sfd^2(\gamma_{xy}(s),\gamma_{xz}(t))\Big),
\]
and define, as in \cite[Def.~5.5]{LaM22},
\begin{equation}\label{cp:eq-up}
 \langle\gamma_{xy},\gamma_{xz}\rangle_{\rm up}
 :=\liminf_{s,t\downarrow0}\Theta(s,t),\qquad
 \Delta^2(\gamma_{xy},\gamma_{xz})
 :=\sfd^2(x,y)+\sfd^2(x,z)+2\langle\gamma_{xy},\gamma_{xz}\rangle_{\rm up}.
\end{equation}
Since $\sfd(\gamma_{xy}(s),\gamma_{xz}(t))$ lies between
$|s\,\sfd(x,y)-t\,\sfd(x,z)|$ and $s\,\sfd(x,y)+t\,\sfd(x,z)$,
\begin{equation}\label{cp:eq-Theta}
 |\Theta(s,t)|\le\sfd(x,y)\,\sfd(x,z),\qquad\text{hence}\qquad
 |\langle\gamma_{xy},\gamma_{xz}\rangle_{\rm up}|\le\sfd(x,y)\,\sfd(x,z)
 \quad\text{and}\quad\Delta^2\ge0 .
\end{equation}
If one of the two geodesics is constant, then $\Theta\equiv0$. For
nonconstant geodesics the \emph{upper angle}
$\sphericalangle_{\rm up}(\gamma_{xy},\gamma_{xz})\in[0,\pi]$ is defined by
$\langle\gamma_{xy},\gamma_{xz}\rangle_{\rm up}
=\sfd(x,y)\sfd(x,z)\cos\sphericalangle_{\rm up}(\gamma_{xy},\gamma_{xz})$;
it is the upper limit of the comparison angles at $x$ of the triangles
$(x,\gamma_{xy}(s),\gamma_{xz}(t))$. The \emph{local angle condition}
(LAC) holds at $x$ if for any three nonconstant geodesics emanating from
$x$ the three upper angles between them add up to at most $2\pi$
\cite[Def.~5.4]{LaM22}.

In a Hilbert space $\langle\gamma_{xy},\gamma_{xz}\rangle_{\rm up}$ is the
scalar product of $y-x$ and $z-x$, and $\Delta$ is the norm of their sum.
In the scheme, $\Delta$ will measure how far the step into a node and the
step out of it are from being opposite and of equal length.

The second geometric ingredient is semiconcavity of the squared distance
along geodesics. Observers are taken in the geodesic cover of $A$, because
in the proof of Theorem~\ref{cp:abstract} the observer is a point of
another interpolant, which need not lie in $A$.

\begin{definition}[Geometric hypothesis]\label{cp:def-G}
Let $A\subset\calX$ and $\kappa\ge0$. We say that \emph{$(\mathrm G_\kappa)$
holds on $A$} if LAC holds at every point of $A$ and if for all $x,y\in A$,
every $\gamma_{xy}\in\mathrm{Geod}(x,y)$ and every observer $v\in A^{\rm g}$ the
function $t\mapsto\frac12\sfd^2(\gamma_{xy}(t),v)-\frac\kappa2\sfd^2(x,y)t^2$ is
concave on $[0,1]$.
\end{definition}

\begin{lemma}[Derivative of the squared distance]\label{cp:deriv}
Let $x,y,v\in\calX$, $\gamma_{xy}\in\mathrm{Geod}(x,y)$ and $\kappa\ge0$, and
suppose that
$t\mapsto\frac12\sfd^2(\gamma_{xy}(t),v)-\frac\kappa2\sfd^2(x,y)t^2$ is
concave on $[0,1]$. Then for every $\gamma_{xv}\in\mathrm{Geod}(x,v)$ and
every $t\in(0,1)$ at which the derivative exists, hence for a.e.\ $t$,
\begin{equation}\label{cp:eq-deriv0}
 \frac{\rmd}{\rmd t}\frac12\sfd^2(\gamma_{xy}(t),v)
 \le-\langle\gamma_{xy},\gamma_{xv}\rangle_{\rm up}+t\,\kappa\,\sfd^2(x,y).
\end{equation}
If in addition LAC holds at $x$, then for every $z\in\calX$ and every
$\gamma_{xz}\in\mathrm{Geod}(x,z)$, at the same points $t$,
\begin{equation}\label{cp:eq-deriv}
 \frac{\rmd}{\rmd t}\frac12\sfd^2(\gamma_{xy}(t),v)
 \le\langle\gamma_{xz},\gamma_{xv}\rangle_{\rm up}
 +\sfd(x,v)\,\Delta(\gamma_{xy},\gamma_{xz})+t\,\kappa\,\sfd^2(x,y).
\end{equation}
In particular, if $(\mathrm G_\kappa)$ holds on $A$, then
\eqref{cp:eq-deriv} holds for all $x,y\in A$, $v\in A^{\rm g}$, $z\in\calX$,
all geodesics $\gamma_{xy},\gamma_{xz},\gamma_{xv}$ and a.e.\ $t\in[0,1]$.

Both estimates also hold in integrated form, without any differentiability:
for all $0\le t\le t'\le1$,
\begin{equation}\label{cp:eq-secant}
 \frac12\sfd^2(\gamma_{xy}(t'),v)-\frac12\sfd^2(\gamma_{xy}(t),v)
 \le-(t'-t)\,\langle\gamma_{xy},\gamma_{xv}\rangle_{\rm up}
 +\frac\kappa2\,\sfd^2(x,y)\,\big(t'^2-t^2\big),
\end{equation}
and, if LAC holds at $x$, the same inequality holds with
$-\langle\gamma_{xy},\gamma_{xv}\rangle_{\rm up}$ replaced by
$\langle\gamma_{xz},\gamma_{xv}\rangle_{\rm up}
+\sfd(x,v)\,\Delta(\gamma_{xy},\gamma_{xz})$.
\end{lemma}

Estimate \eqref{cp:eq-deriv} is \cite[Prop.~5.7]{LaM22}. Only this
consequence of $(\mathrm G_\kappa)$ is used below.

\subsection{Estimates for the minimizing movement scheme}\label{cp:ss-error}

This subsection carries out parts (2) and (3) of the plan for an abstract
functional $\phi$. The functional is convex along \emph{one} geodesic for
each pair of points only, and we call such a geodesic admissible. Two
choices matter: in the Euler--Lagrange inequality the geodesic from a node
to the observer is admissible, and the interpolant joins consecutive nodes
by admissible geodesics. The geometric hypothesis holds along every
geodesic, so it is available along the selected ones.

\paragraph{Standing assumptions of this subsection.}
$\phi:\calX\to(-\infty,\infty]$ is proper and weakly geodesically
$\lambda$-convex (Definition~\ref{prelim:convexity}(c)). We put
\[
 \ell:=-\min\{\lambda,0\}\ge0,\qquad\hat\lambda:=-2(\ell+1).
\]
The constant $\hat\lambda$ is the one with which the approximate $\EVI$ of
Proposition~\ref{cp:approxEVI} holds; the loss from $\lambda$ to
$\hat\lambda$ is recovered in Step~3 of the proof of
Theorem~\ref{cp:abstract}.
A geodesic $\gamma$ joining two points of $\dom\phi$ is called
\emph{admissible} if $\phi$ is $(-\ell)$-convex along it, that is, if
$t\mapsto\phi(\gamma_t)+\frac\ell2\sfd^2(\gamma_0,\gamma_1)t^2$ is convex on
$[0,1]$. Two facts are used.
\begin{enumerate}[label=\textup{(S\arabic*)}]
\item\label{cp:S1} Any two points of $\dom\phi$ are joined by an admissible
 geodesic: a geodesic along which $\phi$ is $\lambda$-convex is admissible,
 because $\lambda\ge-\ell$.
\item\label{cp:S2} The reversal of an admissible geodesic and its
 restriction to a subinterval, reparametrised on $[0,1]$, are admissible:
 the function in the definition for the new geodesic is the old one,
 composed with an affine change of the parameter, plus an affine function.
\end{enumerate}
The step satisfies
$0<\tau\le1$ and $\ell\tau\le\frac12$. The points $x_0,x_1,\dots$ form an exact
minimizing movement scheme (Definition~\ref{prelim:MM}), that is, $x_n$
minimises $\frac1{2\tau}\sfd^2(x_{n-1},\cdot)+\phi$, and
$\varsigma_0:=\slope\phi(x_0)<\infty$; in particular all $x_n$ lie in $\dom\phi$. We
write
\[
 d_n:=\sfd(x_{n-1},x_n),\qquad \varsigma_n:=\slope\phi(x_n),\qquad
 D_n:=\Delta^2\big(\gamma_{x_nx_{n-1}},\gamma_{x_nx_{n+1}}\big),
\]
where for each $n$ one \emph{admissible} geodesic $\gamma_{x_nx_{n+1}}$ and
one geodesic $\gamma_{x_nx_{n-1}}$ are fixed, and define the
\emph{admissible geodesic interpolant}
$u^\tau(t):=\gamma_{x_nx_{n+1}}((t-n\tau)/\tau)$ for
$t\in[n\tau,(n+1)\tau]$. The error function below is the error term in the
approximate $\EVI$ of Proposition~\ref{cp:approxEVI}; part (3) of the plan
is that its integral over $[0,T]$ is of order $\tau$
(Lemma~\ref{cp:error}). The error function of \cite[(5.17)]{LaM22} is
$\Er^{\tau,\kappa}(t):=\Er_n$ for $t\in[n\tau,(n+1)\tau)$ with
\begin{equation}\label{cp:eq-Er}
 \Er_0:=(1+2\ell)\,d_1^2+\big(1+(1-\ell\tau)^{-1}\big)\varsigma_0^2,\qquad
 \Er_n:=\Big(1+2\ell+\frac\kappa\tau\Big)d_{n+1}^2+\frac{D_n}{\tau^2}
 \quad(n\ge1).
\end{equation}

Part \ref{cp:it-dir} of the next lemma is the Euler--Lagrange inequality in
the direction of an arbitrary point $v$, the discrete form of the $\EVI$.
The outcome of parts \ref{cp:it-slope}--\ref{cp:it-abc} is that
$d_n\le C\tau$ on bounded time intervals when the initial slope is finite.
The statements of this subsection are those that \cite{LaM22} states
without proof, in the form needed here, except for the corrected estimate of
Lemma~\ref{cp:error}, which is proved here. They belong to the approach of
Savar\'e \cite{Sav07}; for the other proofs we refer to \cite{MuSII}.

\begin{lemma}[Discrete Euler inequalities]\label{cp:euler}
Under the standing assumptions the following holds for every $n\ge0$.
\begin{enumerate}[label=\textup{(\alph*)}]
\item\label{cp:it-dir} For every $v\in\dom\phi$, every geodesic
$\gamma_{x_{n+1}x_n}$ and every \emph{admissible} geodesic
$\gamma_{x_{n+1}v}$,
\begin{equation}\label{cp:eq-dir}
 \frac1\tau\langle\gamma_{x_{n+1}x_n},\gamma_{x_{n+1}v}\rangle_{\rm up}
 -\frac\ell2\sfd^2(x_{n+1},v)+\phi(x_{n+1})\le\phi(v).
\end{equation}
\item\label{cp:it-slope} For every $x\in\dom\slope\phi$ and $y\in\calX$,
$\ \phi(x)-\phi(y)\le\slope\phi(x)\,\sfd(x,y)+\frac\ell2\sfd^2(x,y)$.
\item\label{cp:it-energy}
$\phi(x_n)-\phi(x_{n+1})\ge\big(\frac1\tau-\frac\ell2\big)d_{n+1}^2$.
\item\label{cp:it-abc} $\varsigma_{n+1}\le d_{n+1}/\tau$,\ \
$(1-\ell\tau)\,d_{n+1}/\tau\le \varsigma_n$,\ \ $(1-\ell\tau)\varsigma_{n+1}\le \varsigma_n$.
Consequently
\begin{equation}\label{cp:eq-incr}
 \varsigma_n\le(1-\ell\tau)^{-n}\varsigma_0,\qquad
 d_n\le\tau(1-\ell\tau)^{-n}\varsigma_0\le\tau\,\ee^{2\ell n\tau}\varsigma_0 ,
\end{equation}
the second for $n\ge1$.
\end{enumerate}
\end{lemma}

In the integral of the error function the term $D_n/\tau^2$ is the
delicate one: the bound $D_n\le(d_n+d_{n+1})^2$ would only give a
contribution of order one. Instead, the Euler--Lagrange inequality at
$x_n$, tested with the next node $x_{n+1}$, bounds $D_n$ by the decrease
$d_n^2-d_{n+1}^2$ of the squared step length up to a term of order
$\tau^3$, and the sum over $n$ telescopes.

\begin{lemma}[Corrected error estimate]\label{cp:error}
Let the standing assumptions hold, let $T>0$, $\kappa\ge0$,
$N:=\lceil T/\tau\rceil$, and assume that the scheme exists up to the index $N$.
Then
\begin{equation}\label{cp:eq-telescope}
 D_n\le d_n^2-d_{n+1}^2+2\ell\tau\,d_{n+1}^2\qquad(1\le n\le N-1),
\end{equation}
and
\begin{equation}\label{cp:eq-error}
 \int_0^T\Er^{\tau,\kappa}(t)\dd t
 \le\tau\,\ee^{4\ell(T+1)}\varsigma_0^2\,\big[4+T(1+4\ell+\kappa)\big].
\end{equation}
The same bound holds for $\int_0^T\ee^{2\hat\lambda t}\Er^{\tau,\kappa}(t)\dd t$,
since $\hat\lambda<0$.
\end{lemma}

\begin{proof}
\emph{Step 1.} By the definition of $\Delta^2$, by
Lemma~\ref{cp:euler}\ref{cp:it-dir} at $x_n$ with the observer
$v=x_{n+1}$ and the admissible geodesic $\gamma_{x_nx_{n+1}}$, and by
Lemma~\ref{cp:euler}\ref{cp:it-energy},
\begin{align*}
 D_n&=d_n^2+d_{n+1}^2
     +2\langle\gamma_{x_nx_{n-1}},\gamma_{x_nx_{n+1}}\rangle_{\rm up}\\
    &\le d_n^2+d_{n+1}^2+2\tau[\phi(x_{n+1})-\phi(x_n)]
     +\ell\tau d_{n+1}^2\\
    &\le d_n^2+d_{n+1}^2-(2-\ell\tau)d_{n+1}^2+\ell\tau d_{n+1}^2 ,
\end{align*}
which is \eqref{cp:eq-telescope}.

\emph{Step 2.} Put $M=\ee^{2\ell(T+1)}\varsigma_0$. For $n\le N$ we have
$n\tau<T+\tau\le T+1$, hence $d_n\le\tau M$ by \eqref{cp:eq-incr}. In the
sum of \eqref{cp:eq-telescope} the differences $d_n^2-d_{n+1}^2$ telescope,
and since $(N-1)\tau<T$,
\begin{equation}\label{cp:eq-sumD}
 \frac1\tau\sum_{n=1}^{N-1}D_n
 \le\frac{d_1^2}\tau+2\ell\sum_{n=1}^{N-1}d_{n+1}^2
 \le\tau M^2(1+2\ell T).
\end{equation}

\emph{Step 3.} The intervals $[n\tau,(n+1)\tau)$, $0\le n\le N-1$, cover
$[0,T)$. The slope term of the first interval contributes at most
$3\tau \varsigma_0^2$ to the integral, because $(1-\ell\tau)^{-1}\le2$.
The terms $(1+2\ell)d_{n+1}^2$ and $(\kappa/\tau)d_{n+1}^2$ contribute
at most $T[(1+2\ell)\tau^2+\kappa\tau]M^2$, and the terms $D_n/\tau^2$ at most
the left-hand side of \eqref{cp:eq-sumD}.
With $\tau\le1$ and $\varsigma_0\le M$,
\[
 \int_0^T\Er^{\tau,\kappa}\dd t
 \le\tau M^2\big[3+T(1+2\ell)+T\kappa+1+2\ell T\big]
 =\tau M^2\big[4+T(1+4\ell+\kappa)\big].\qedhere
\]
\end{proof}

The last statement of this subsection is the approximate $\EVI$ for the
interpolant. It holds with the constant $\hat\lambda$, which is smaller
than $\lambda$; this loss is recovered in the proof of
Theorem~\ref{cp:abstract}.

\begin{proposition}[Approximate EVI for the interpolant]\label{cp:approxEVI}
Let the standing assumptions hold, let $T>0$, $N=\lceil T/\tau\rceil$, and let
$A\subset\calX$ be a set on which $(\mathrm G_\kappa)$ holds and which contains
$x_0,\dots,x_N$. Then for every admissible geodesic interpolant $u^\tau$,
every observer $v\in A^{\rm g}$ and a.e.\ $t\in(0,T)$,
\begin{equation}\label{cp:eq-approxEVI}
 \frac{\rmd}{\rmd t}\frac12\sfd^2(u^\tau(t),v)
 +\frac{\hat\lambda}2\sfd^2(u^\tau(t),v)
 \le\phi(v)-\phi(u^\tau(t))+\Er^{\tau,\kappa}(t).
\end{equation}
\end{proposition}

We record an elementary energy bound for later use. Let
$t\in(n\tau,(n+1)\tau)$, $u=u^\tau(t)$ and $a=d_{n+1}$. By admissibility
of the interpolating geodesic $\gamma_{x_nx_{n+1}}$ and
$\phi(x_{n+1})\le\phi(x_n)$,
\begin{equation}\label{cp:eq-phiu}
 \phi(u)\le\phi(x_n)+\frac\ell8a^2 .
\end{equation}

\begin{remark}[No exceptional set]\label{cp:rem-dini}
The exceptional null set in Proposition~\ref{cp:approxEVI} can be avoided:
\eqref{cp:eq-approxEVI} holds at every $t\in[0,T)$ if the derivative is
read as the upper right Dini derivative, which is the derivative of
Definition~\ref{prelim:EVI}; this form rests on \eqref{cp:eq-secant} in
place of \eqref{cp:eq-deriv}. We do not need this form
in the text. It has been verified formally, see
Appendix~\ref{sec:verification}.
\end{remark}

\subsection{The abstract generation theorem}\label{cp:ss-abstract}

The hypotheses of the abstract theorem describe the three tasks of the
sets $A_j$: the geometry is controlled on each of them \ref{cp:A1}; a
scheme started in their union stays in one of them up to a given time, for
all small steps \ref{cp:A2}; and the union approximates every point of
finite energy in distance and in energy \ref{cp:A3}. The points of
$\bigcup_jA_j$ are called \emph{regular}; accordingly we speak of regular
initial data and regular observers.

\begin{theorem}[Abstract generation theorem]\label{cp:abstract}
Let $(\calX,\sfd)$ be a complete geodesic metric space and let
$\phi:\calX\to(-\infty,\infty]$ be proper, lower semicontinuous and weakly
geodesically $\lambda$-convex, $\lambda\in\R$
$($Definition~\textup{\ref{prelim:convexity}(c)}$)$. Assume that there are
nested sets $A_1\subset A_2\subset\dots\subset\dom\slope\phi$ and constants
$\kappa_j\ge0$ with the following properties.
\begin{enumerate}[label=\textup{(A\arabic*)}]
\item\label{cp:A1} For every $j$, the geometric hypothesis
$(\mathrm G_{\kappa_j})$ of Definition~\ref{cp:def-G} holds on $A_j$.
\item\label{cp:A2} For every $x_0\in\bigcup_jA_j$ and every $T\in(0,\infty)$ there
are $j\in\N$ and $\tau_*>0$ such that for every $\tau\in(0,\tau_*]$ the
minimisation problems of the scheme started at $x_0$ are solvable up to the
index $N=\lceil T/\tau\rceil$, and for every choice of the minimisers
$x_0,x_1,\dots,x_N\in A_j$.
\item\label{cp:A3} For every $x\in\dom\phi$ there are
$x_m\in\bigcup_jA_j$ with $x_m\to x$ and $\phi(x_m)\to\phi(x)$.
\end{enumerate}
Then there is a unique $\EVI_\lambda$ flow of $\phi$ on $\overline{\dom\phi}$
in the sense of Definition~\ref{prelim:EVI}: for every
$\hat x_0\in\overline{\dom\phi}$ there is a unique complete $\EVI_\lambda$
solution $u:[0,\infty)\to\calX$ with $u(0)=\hat x_0$. Moreover, for every
$x_0\in\bigcup_jA_j$ all geodesic interpolants $u^\tau$ of exact minimizing
movement schemes started at $x_0$ converge to the solution, locally uniformly
on $[0,\infty)$, as $\tau\downarrow0$, for every choice of the minimisers and
of the interpolating geodesics.
\end{theorem}

Compared with \cite[Thm.~5.9]{LaM22}, strong convexity is replaced by weak
convexity; compactness of the sublevels by completeness of $\calX$ and lower
semicontinuity of $\phi$; solvability of the scheme is part of \ref{cp:A2}
and is required for small steps only; the concavity constant of $A_j$ need
not equal the index; and \ref{cp:A3} is required on $\dom\phi$, which
implies that $\bigcup_jA_j$ is dense in $\overline{\dom\phi}$.

In Step~1 of the proof two interpolants with different steps are compared
by using each as observer in the approximate $\EVI$ of the other. Step~2
passes the approximate $\EVI$ to the limit, and Step~3 restores the
constant $\lambda$ and extends the flow from $\bigcup_jA_j$ to its closure,
using \cite{MuS20}.

\begin{proof}
We use the notation $\ell$, $\hat\lambda$ and the notion of admissible
geodesic of Section~\ref{cp:ss-error}.

\emph{Step 1: admissible interpolants form a Cauchy family.}
Fix $x_0\in\bigcup_jA_j$ and $T>0$, let $j$, $\tau_*$ be as in \ref{cp:A2},
and let $\tau,\sigma\le\min\{\tau_*,1,(2\ell)^{-1}\}$ $($the last entry is
absent if $\ell=0$$)$, so that the standing assumptions of
Section~\ref{cp:ss-error} hold. Let $u^\tau,u^\sigma$
be admissible geodesic interpolants of two schemes started at $x_0$, for
arbitrary choices of the minimisers. Their nodes up to the indices
$\lceil T/\tau\rceil$, $\lceil T/\sigma\rceil$ lie in $A_j$, so that
$u^\tau(t),u^\sigma(t)\in A_j^{\rm g}$ for $t\in[0,T]$; this is what allows
each interpolant as an observer for the other in
Proposition~\ref{cp:approxEVI}. By \eqref{cp:eq-incr} both curves are
$M$-Lipschitz on $[0,T]$, $M=\ee^{2\ell(T+1)}\slope\phi(x_0)$.
The energies $t\mapsto\phi(u^\tau(t))$ are continuous: on each segment the
function is finite and semiconvex by admissibility, hence continuous in the
interior and upper semicontinuous at the endpoints, and lower
semicontinuity of $\phi$ gives continuity at the endpoints.

Put $F(s,t)=\tfrac12\sfd^2(u^\tau(s),u^\sigma(t))$ and let
$0\le t<t+h\le T$. Integrate
\eqref{cp:eq-approxEVI} for $u^\tau$ over $[t,t+h]$ with the observer
$u^\sigma(t+h)$, and for $u^\sigma$ with the observer $u^\tau(t)$, and add.
The four distance terms telescope exactly,
\[
 \big[F(t+h,t+h)-F(t,t+h)\big]+\big[F(t,t+h)-F(t,t)\big]
 =F(t+h,t+h)-F(t,t).
\]
Divide the sum by $h$ and let $h\downarrow0$. By continuity of the energies
and of the curves the energy terms
$\phi(u^\sigma(t+h))+\phi(u^\tau(t))
-\frac1h\int_t^{t+h}[\phi(u^\tau)+\phi(u^\sigma)]\dd r$ tend to zero, and
the Lipschitz function $t\mapsto F(t,t)$ satisfies
$\frac{\rmd}{\rmd t}F(t,t)\le-2\hat\lambda F(t,t)
+\Er^{\tau,\kappa_j}(t)+\Er^{\sigma,\kappa_j}(t)$ a.e.
Since $F(0,0)=0$ and $\hat\lambda<0$, integration and
\eqref{cp:eq-error} give
\begin{equation}\label{cp:eq-cauchy}
 \sup_{t\le T}\sfd^2(u^\tau(t),u^\sigma(t))
 \le2\ee^{-2\hat\lambda T}\int_0^T
 (\Er^{\tau,\kappa_j}+\Er^{\sigma,\kappa_j})\dd t
 \le2\ee^{-2\hat\lambda T}C_j(\tau+\sigma)\slope\phi^2(x_0),
\end{equation}
with $C_j=\ee^{4\ell(T+1)}[4+T(1+4\ell+\kappa_j)]$.
Completeness yields a uniform limit $u$ on $[0,T]$, which is $M$-Lipschitz,
starts at $x_0$ and is independent of the choices; the limits for different
$T$ are consistent. If $\tilde u^\tau$ interpolates the same nodes along
arbitrary geodesics, then
$\sfd(\tilde u^\tau(t),u^\tau(t))\le2\max_{n\le N}d_n\le2\tau M$ on $[0,T]$,
so that $\tilde u^\tau\to u$ as well.

\emph{Step 2: the limit satisfies the integrated $\EVI_{\hat\lambda}$.}
Let $v$ be a regular observer, $v\in A_m$. Enlarging
$j$ we may assume $A_m\subset A_j$; then $v\in A_j^{\rm g}$, and $\kappa_j$
stays fixed as $\tau\downarrow0$. Integrating \eqref{cp:eq-approxEVI} over
$(s,t)\subset[0,T]$ and using \eqref{cp:eq-error},
\begin{equation}\label{cp:eq-intapprox}
 \frac12\sfd^2(u^\tau(t),v)-\frac12\sfd^2(u^\tau(s),v)
 +\int_s^t\Big(\phi(u^\tau(r))+\frac{\hat\lambda}2\sfd^2(u^\tau(r),v)\Big)\dd r
 \le(t-s)\phi(v)+C_j\tau\slope\phi^2(x_0) .
\end{equation}
The distance terms converge. For the energy term we use Fatou's lemma, on
the basis of three facts. The slope inequality
Lemma~\ref{cp:euler}\ref{cp:it-slope} at the anchor $x_0$ gives
\begin{equation}\label{cp:eq-lowerphi}
 \phi(y)\ge\phi(x_0)-\slope\phi(x_0)\sfd(x_0,y)
                         -\tfrac\ell2\sfd^2(x_0,y)\qquad(y\in\calX),
\end{equation}
and all interpolants stay in the ball of radius $MT$ around $x_0$ on
$[0,T]$, so the energies are bounded below uniformly in $\tau$. Lower
semicontinuity gives
$\phi(u(r))\le\liminf_{\tau\downarrow0}\phi(u^\tau(r))$, and
\eqref{cp:eq-phiu} with the monotonicity of the discrete energies gives the
upper bound $\phi(u^\tau(r))\le\phi(x_0)+\ell\tau^2M^2/8$; hence
$\phi(u(r))\le\phi(x_0)$. Finally $\phi\circ u$ is lower semicontinuous on
$[0,T]$, hence measurable, and bounded. Thus \eqref{cp:eq-intapprox} passes
to the limit and gives \eqref{prelim:eq-intEVI}, with $\hat\lambda$ in place
of $\lambda$, for all $0\le s<t\le T$ and all regular observers. For
$v\in\dom\phi$ take $v_m$ as in \ref{cp:A3}: the left-hand side is
continuous in the observer and $\phi(v_m)\to\phi(v)$. By
Proposition~\ref{prelim:EVIprop} the integrated inequality is equivalent to
Definition~\ref{prelim:EVI}, so $u$ is an $\EVI_{\hat\lambda}$ solution
starting at $x_0$.

\emph{Step 3: from $\hat\lambda$ to $\lambda$, and the closure.}
We have $\hat\lambda<\lambda$, and $\phi$ is $\lambda$-convex in the sense
of \cite[Sec.~2.3]{MuS20}, which asks for \eqref{prelim:eq-convex} along one
geodesic for each pair of points of $\dom\phi$. Hence
\cite[Cor.~3.12]{MuS20} applies and $u$ is an $\EVI_\lambda$ solution.
By \ref{cp:A3}, $\bigcup_jA_j$ is dense in $\overline{\dom\phi}$, and the
sublevels of $\phi$ are complete because $\calX$ is complete and $\phi$ is
lower semicontinuous. Therefore \cite[Cor.~3.8]{MuS20} extends the
solution map to an $\EVI_\lambda$ flow on $\overline{\dom\phi}$; the
extension uses the contraction estimate \eqref{prelim:eq-contraction},
which also gives uniqueness and the semigroup property.
\end{proof}

\begin{remark}[Two points of the proof]\label{cp:rem-abstract}
The formal verification of the proof (Appendix~\ref{sec:verification})
showed two things about its hypotheses and steps.

(i) The proof does not use that $\calX$ is a geodesic space. Geodesics are
needed only between points of $\dom\phi$, and there weak convexity supplies
them.

(ii) Step~1 can be carried out without the continuity of
$t\mapsto\phi(u^\tau(t))$. In the formal proof the approximate $\EVI$ is
used in integrated form on an interval $[t,t']$ contained in one interval of
each of the two partitions, with the energy $\phi(x_n)$ at the left node in
place of $\phi(u^\tau(t))$; this is the form of
Proposition~\ref{cp:approxEVI} that rests on \eqref{cp:eq-secant}. When the two inequalities are
added, the energies of the two observers are bounded by the node energies
through \eqref{cp:eq-phiu}, and the node energies
cancel. Letting $t'\downarrow t$ gives the differential inequality for the
upper right Dini derivative of $t\mapsto F(t,t)$ at every $t\in[0,T)$, and
\eqref{cp:eq-cauchy} follows with the same constant. In Step~2 the node
energies likewise replace $\phi(u^\tau(r))$ in \eqref{cp:eq-intapprox}, and
Fatou's lemma is applied to them. This variant, under the hypotheses of
Theorem~\ref{cp:abstract}, is what has been verified formally.

The continuity of the energies in Step~1 is the one place of
Sections~\ref{cp:ss-error} and~\ref{cp:ss-abstract} where the convexity of
the whole function
$t\mapsto\phi(\gamma_t)+\frac\ell2\sfd^2(\gamma_0,\gamma_1)t^2$ along an
admissible geodesic is used; everywhere else only the inequality
\eqref{prelim:eq-convex}, with the constant $-\ell$, along the selected
geodesic enters. The variant (ii) therefore suggests that
Theorem~\ref{cp:abstract} holds under the weaker notion of convexity of
\cite{AGS05} (Remark~\ref{prelim:rem-convexity}(ii)). We have not carried
out this extension, and it is not part of what was verified.
\end{remark}

\begin{remark}[Localised form]\label{cp:rem-local}
Sections~\ref{cp:ss-error} and~\ref{cp:ss-abstract} assume that $\phi$ is
weakly geodesically $\lambda$-convex, that is, convex along one geodesic
for each pair of points of $\dom\phi$. Steps~1 and~2 of the proof of
Theorem~\ref{cp:abstract} use less: convexity along geodesics between
points of a set $K$ that contains the scheme, its interpolants and the
observers. We state this form here, because Section~\ref{sec:sublevel}
applies it to a functional that is semiconvex only on sublevels.

Let $\phi:\calX\to(-\infty,\infty]$ be proper, let $\ell\ge0$ be given, put
$\hat\lambda:=-2(\ell+1)$, and define admissible geodesics as in
Section~\ref{cp:ss-error} with this $\ell$. Let $x_0\in\dom\slope\phi$ and
let $K\subset\dom\phi$ be a set with the following properties.
\begin{enumerate}[label=\textup{(K\arabic*)}]
\item\label{cp:K1} $K$ contains the sublevel $\{\phi\le\phi(x_0)\}$.
\item\label{cp:K2} Any two points of $K$ are joined by an admissible
 geodesic.
\item\label{cp:K3} Every admissible geodesic between two points of
 $\{\phi\le\phi(x_0)\}$ lies in $K$.
\item\label{cp:K4} Observers are taken in $K$:
 Lemma~\ref{cp:euler}\ref{cp:it-slope} is asserted for
 $x\in K\cap\dom\slope\phi$ and $y\in K$, the lower bound
 \eqref{cp:eq-lowerphi} for $y\in K$, Proposition~\ref{cp:approxEVI} for
 $v\in A^{\rm g}\cap K$, and Step~2 for regular observers in $K$.
\end{enumerate}
The standing assumptions of Section~\ref{cp:ss-error} are the case
$K=\dom\phi$, $\ell=-\min\{\lambda,0\}$: then \ref{cp:K2} is \ref{cp:S1},
\ref{cp:K3} holds because $\phi$ is finite along an admissible geodesic,
and \ref{cp:K1} and \ref{cp:K4} are no restriction.

Under \ref{cp:K1}--\ref{cp:K4}, for exact schemes started at $x_0$ with
$0<\tau\le1$ and $\ell\tau\le\frac12$, Lemmas~\ref{cp:euler}
and~\ref{cp:error} and Proposition~\ref{cp:approxEVI} hold with the same
constants in terms of $\ell$, $\hat\lambda$ and $\kappa$. If in addition
$\calX$ is complete, $\phi$ is lower semicontinuous, and nested sets
$A_j\subset\dom\slope\phi$ satisfy \ref{cp:A1} and, for this $x_0$,
\ref{cp:A2}, then the conclusions of Steps~1 and~2 hold: the estimate
\eqref{cp:eq-cauchy} with the same $C_j$; the limit curve $u$, which is
$M$-Lipschitz on $[0,T]$ with $M=\ee^{2\ell(T+1)}\slope\phi(x_0)$, starts
at $x_0$, satisfies $\phi(u(t))\le\phi(x_0)$ and is the limit of all
geodesic interpolants; and the integrated inequality
\eqref{prelim:eq-intEVI}, with $\hat\lambda$ in place of $\lambda$, for all
$0\le s<t\le T$ and all regular observers in $K$. This inequality passes
to every $v\in\dom\phi$ that is the limit, in distance and in energy, of
regular points lying in $K$.

The statements hold in this localised form, for the following reason. The energy decreases along the
scheme, so the nodes lie in $\{\phi\le\phi(x_0)\}\subset K$ by
\ref{cp:K1}. The points of an admissible interpolant lie in $K$ by
\ref{cp:K3}, and the limit curve lies in $\{\phi\le\phi(x_0)\}$.
Convexity of $\phi$ is used along the following geodesics and along no
others:
\begin{enumerate}[label=\textup{(\arabic*)}]
\item between consecutive nodes: in
 Lemma~\ref{cp:euler}\ref{cp:it-energy} and~\ref{cp:it-abc}, in
 \eqref{cp:eq-telescope}, in \eqref{cp:eq-phiu},
 and for the continuity of the energies in Step~1;
\item from a node to the observer: in \eqref{cp:eq-dir};
\item from $x_0$ to the observer: in
 Lemma~\ref{cp:euler}\ref{cp:it-slope};
\item from $x_0$ to a point of an interpolant or of the limit curve: in
 the lower bound \eqref{cp:eq-lowerphi} of Step~2.
\end{enumerate}
In Step~1 the observers are points of the other interpolant. Hence each
of these geodesics joins two points of $K$ and can be chosen admissible by
\ref{cp:K2}; \ref{cp:S2} holds as before. The inner points of such a
geodesic enter only through the convexity inequality, as competitors in the
minimisation over $\calX$ or in the definition of the slope, and need not
lie in $K$. The geometric hypothesis enters only through
\eqref{cp:eq-deriv}, as before. Parts \ref{cp:it-energy} and
\ref{cp:it-abc} of Lemma~\ref{cp:euler} and \eqref{cp:eq-incr} involve
consecutive nodes only; for them it suffices that any two points of
$\{\phi\le\phi(x_0)\}$ are joined by an admissible geodesic.

Step~3 is not part of this remark: it needs the inequality for all
observers in $\dom\phi$, and \cite[Cor.~3.12]{MuS20} asks for convexity
between arbitrary points of $\dom\phi$. In this localised form,
and in the variant of Remark~\ref{cp:rem-abstract}(ii), Steps~1 and~2 have
been verified formally (Appendix~\ref{sec:verification}).
\end{remark}

\begin{corollary}[Strong convexity a posteriori]\label{cp:strong}
Under the assumptions of Theorem~\ref{cp:abstract}, $\phi$ is
$\lambda$-convex along every geodesic that joins two points of $\dom\phi$
and lies in $\overline{\dom\phi}$. In particular, if
$\overline{\dom\phi}=\calX$, then $\phi$ is strongly geodesically
$\lambda$-convex.
\end{corollary}

\begin{proof}
Let $\gamma$ be such a geodesic. By Theorem~\ref{prelim:MS-convex}, applied
to the flow on $D=\overline{\dom\phi}$, inequality
\eqref{prelim:eq-convex} holds along $\gamma$; in particular $\phi$ is
finite on $\gamma$. Every restriction of $\gamma$ to a subinterval is
therefore again a geodesic in $D$ with endpoints in $\dom\phi$, and
\eqref{prelim:eq-convex} for all restrictions is the convexity of
$t\mapsto\phi(\gamma_t)-\frac\lambda2\sfd^2(\gamma_0,\gamma_1)t^2$.
\end{proof}

\subsection{Verification for \texorpdfstring{$\HK$}{HK} and \texorpdfstring{$\SHK$}{SHK}}\label{cp:ss-verify}

We return to measures and check the hypotheses of
Theorem~\ref{cp:abstract} for $\phi=\sfE$, with the sets $A_j$ of measures
whose density lies in $[1/j,j]$ and whose slope is at most $j$. The
semiconcavity of $\frac12\HK^2$ on density-bounded sets is taken from
\cite{LaM22}; we need it for every geodesic between two density-bounded
measures, because the interpolating geodesics are selected by the energy,
and for observers in the geodesic cover.

\begin{proposition}[Semiconcavity with observers in the geodesic cover]\label{cp:semiconc}
Let $0<\delta<1$.
\begin{enumerate}[label=\textup{(\alph*)}]
\item\label{cp:it-semiHK} There is $\kappa=\kappa(X,\delta)\ge0$ such that
$(\mathrm G_\kappa)$ holds on $\calM_\delta(X)$ in $(\calM(X),\HK)$: for every
$\HK$ geodesic $(\mu_t)_{t\in[0,1]}$ with $\mu_0,\mu_1\in\calM_\delta(X)$ and
every $w$ in the $\HK$ geodesic cover of $\calM_\delta(X)$, the function
$t\mapsto\frac12\HK^2(\mu_t,w)-\frac\kappa2\HK^2(\mu_0,\mu_1)t^2$ is concave, and
LAC holds at every point of $\calM(X)$.
\item\label{cp:it-semiSHK} There is $\kappa=\kappa(X,\delta)\ge0$ such that
$(\mathrm G_\kappa)$ holds on $\calP(X)\cap\calM_\delta(X)$ in
$(\calP(X),\SHK)$, with observers in the $\SHK$ geodesic cover of
$\calP(X)\cap\calM_\delta(X)$, and LAC holds at every point of $\calP(X)$.
\end{enumerate}
\end{proposition}

Part \ref{cp:it-semiHK} is quoted. The concavity is
\cite[Thm.~6.1]{LaM22}. Its statement asks for one geodesic for each pair,
but its proof gives the assertion for every geodesic and every observer in
the geodesic cover, by three facts: \cite[Lemma~6.2]{LaM22}
$($\cite[Lemma~4.9]{LaM19}$)$ bounds the marginal densities and the transport
distance for \emph{any} optimal plan between measures of a class
$\widetilde{\calM}_{d_1,d_2}(X)$ defined by bounds on the mass of balls of
radius $d_1$, where $0<d_1<\frac\pi2$ $($this upper bound is a hypothesis of
\cite[Lemma~4.9]{LaM19}; it is missing in \cite[Lemma~6.2]{LaM22}, which
fails for $d_1=\frac\pi2$$)$;
\cite[Lemma~6.3]{LaM22} $($\cite[Lemma~4.10]{LaM19}$)$ places \emph{any}
constant-speed geodesic between two points of $\calM_\delta(X)$ in this
class with such a radius, which therefore contains the geodesic cover; and the proof starts
from an arbitrary geodesic between points of $\calM_\delta(X)$ and uses of
the observer only that it lies in this class. LAC is
\cite[Thm.~4.2]{LaM19}.
The SHK statement \ref{cp:it-semiSHK} needs the following additional
argument.

\begin{proof}[Proof of \ref{cp:it-semiSHK}]
Put $A_\delta=\calP(X)\cap\calM_\delta(X)$.
Step~1 proves semiconcavity of $\tfrac12\HK^2(\cdot,w)$ along the
corresponding HK geodesics, for every observer
$w\in(A_\delta)^{\rm g}_{\SHK}$. Step~2 places the endpoints and observers
in a common SHK ball of uniform radius less than $\pi/2$. Step~3 applies
\cite[Prop.~2.27(B)]{LaM19} with unit radii to transfer the estimate to
SHK.

\emph{Step 1: the cover.} Let $(\mu_s)$ be an HK geodesic between
$\eta_0,\eta_1\in A_\delta$. By \eqref{prelim:eq-mass} and $\HK^2\le2$ on
$\calP(X)$, its mass $M(s)=1-s(1-s)\HK^2(\eta_0,\eta_1)$ lies in $[1/2,1]$,
and the SHK geodesics between $\eta_0$ and $\eta_1$ are the reparametrised
normalisations $\mu_s/M(s)$ of such geodesics
(Proposition~\ref{prelim:cone}). For a point $w=\mu_s/M(s)$, the curve
$r\mapsto\mu_r/M(s)$ is an HK geodesic through $w$ by
\eqref{prelim:eq-scaling}, and its endpoints $\eta_i/M(s)$ have densities in
$[\delta,2/\delta]$. Therefore
\begin{equation}\label{cp:eq-cover}
 (A_\delta)^{\rm g}_{\SHK}\subset
 (\calM_{\delta/2}(X))^{\rm g}_{\HK}.
\end{equation}
By \ref{cp:it-semiHK}, $\frac12\HK^2(\cdot,w)$ is
$\kappa(X,\delta/2)$-concave along every HK geodesic between points of
$A_\delta\subset\calM_{\delta/2}(X)$, for every $w$ in the set on the left
of \eqref{cp:eq-cover}.

\emph{Step 2: a common ball of radius less than $\pi/2$.}
Let $\nu_*=\Leb|_X/V$.
For any $\eta\in\calP(X)$, choose a relative ball $B$ of diameter below
$\pi/2$ around a point of $\supp\eta$. It has masses $a=\nu_*(B)>0$,
$b=\eta(B)>0$.
The LET plan $t(\nu_*|_B/a)\otimes(\eta|_B/b)$ has cost
$2+2t\log t-t\log(ab)-2t+t\ell_0<2$ for small $t>0$, where $\ell_0$ is
the finite mean transport cost. Thus $\HK^2(\nu_*,\eta)<2$, that is,
$\SHK(\nu_*,\eta)<\pi/2$.
By compactness and continuity,
$D_0=\max_{\eta\in\calP(X)}\SHK(\nu_*,\eta)<\pi/2$.

\emph{Step 3: transfer to the base.}
Choose $D_0<D<\pi/2$; all endpoints and observers lie in the SHK ball of
radius $D$ around $\nu_*$ and have mass one. Apply
\cite[Prop.~2.27(B)]{LaM19} to the cone $(\calM(X),\HK)$ over
$(\calP(X),\SHK)$ with radial coordinates all equal to $1$: the concavity of
Step~1 along the HK geodesic implies $\kappa'$-concavity of $\frac12\SHK^2$
along the corresponding SHK geodesic, with $\kappa'$ depending only on $D$
and $\kappa(X,\delta/2)$. The cited proposition requires distinct endpoints
and equal radii; for a constant geodesic there is nothing to prove. All SHK
geodesics arise by normalisation, so the conclusion
holds for every geodesic and every stated observer. LAC transfers from
convex $X\subset\R^d$ by \cite[Thm.~4.2]{LaM19}.
\end{proof}

On $\calM(X)$ the energy may tend to $-\infty$ linearly in the mass
(Example~\ref{cp:ex-sublevel}), and the distance term has to dominate; this
is the reason for the restriction on $\tau$ in the next lemma.

\begin{lemma}[Existence of minimisers]\label{cp:exist}
Let $a,b$ be as in \eqref{cp:eq-affine}. For every $y\in\calM(X)$ and every
$\tau>0$ with $2a\tau<1$ the functional
$\nu\mapsto\frac1{2\tau}\HK^2(y,\nu)+\sfE(\nu)$ has a minimiser in $\calM(X)$.
For every $y\in\calP(X)$ and every $\tau>0$ the functional
$\nu\mapsto\frac1{2\tau}\SHK^2(y,\nu)+\sfE(\nu)$ has a minimiser in $\calP(X)$.
If $\sfE(y)<\infty$, every minimiser has finite energy.
\end{lemma}

\begin{proof}
Write $M_0=y(X)$ and $M=\nu(X)$. By \eqref{prelim:eq-massbound} and
\eqref{cp:eq-affine},
\[
 \frac{\HK^2(y,\nu)}{2\tau}+\sfE(\nu)
 \ge\Big(\frac1{2\tau}-a\Big)M-\frac{\sqrt{M_0M}}\tau+\frac{M_0}{2\tau}-bV ,
\]
which is coercive in $M$ if $2a\tau<1$. A minimising sequence therefore has
bounded mass and a narrowly convergent subsequence, and the lower
semicontinuity of $\sfE$ and of $\HK$ gives a minimiser. On $\calP(X)$
compactness gives the minimiser directly. The last statement holds because
$y$ is a competitor.
\end{proof}

To reach all data and all observers, every measure of finite energy has to
be a limit of measures in $\bigcup_jA_j$ in distance \emph{and in energy}
\ref{cp:A3}. By Example~\ref{cp:ex-slope} density bounds alone do not give
finite slope, so a regularisation of the density is not enough. We average
over the cells of a fine partition, mix with a small multiple of the
uniform measure, and take one step of the scheme with a small $\tau$, which
gives finite slope and keeps the density bounds.

\begin{lemma}[Recovery by regular measures]\label{cp:recovery}
Let $\mu\in\calM(X)$ with $\sfE(\mu)<\infty$. Then there are $\delta_j\in(0,1)$
and measures $z_j\in\calM_{\delta_j}(X)$ with $\slope{\sfE}(z_j)<\infty$ such that
\begin{equation}\label{cp:eq-recovery}
 \HK(z_j,\mu)\to0\qquad\text{and}\qquad\sfE(z_j)\to\sfE(\mu).
\end{equation}
If $\mu\in\calP(X)$, the same holds with $z_j\in\calP(X)\cap\calM_{\delta_j}(X)$,
with the slope taken in $(\calP(X),\SHK)$ and with $\SHK(z_j,\mu)\to0$; if
in addition $\mu=\rho\Leb$ with $\rho\le M$ a.e., the densities of the
$z_j$ are at most $M$.
Consequently, in both spaces the measures with positive bounded density and
finite slope are dense in $\overline{\dom\sfE}$, and
$\overline{\dom\sfE}$ is the whole space $\calM(X)$, respectively
$\calP(X)$.
\end{lemma}

\begin{proof}
Partition $X$ into finitely many positive-volume cells $A_{j,i}$ with
maximal diameter tending to zero. For every $\mu\in\calM(X)$ the cell
averages
$v_j=\sum_i\mu(A_{j,i})|A_{j,i}|^{-1}\Leb|_{A_{j,i}}$ have the mass of
$\mu$, bounded densities, hence finite energy, and converge narrowly
to $\mu$; in particular $\overline{\dom\sfE}$ is the whole space. Let now
$\sfE(\mu)<\infty$, $\mu=\rho\Leb+\mu^s$. Jensen's inequality and
$E(c+s)\le E(c)+\Einf s$ give
\[
 |A_{j,i}|E(\mu(A_{j,i})/|A_{j,i}|)
 \le\int_{A_{j,i}}E(\rho)\dd x+\Einf\mu^s(A_{j,i}).
\]
If $\Einf=\infty$, finite energy forces $\mu^s=0$ and the last term is
absent. Hence $\sfE(v_j)\le\sfE(\mu)$ and the energies converge by lower
semicontinuity. Mix $v_j$ with a vanishing amount of the uniform measure:
$y_j=(1-\eps_j)v_j+\eps_j\Leb|_X$ in the HK case and
$y_j=(1-\eps_j)v_j+\eps_j\Leb|_X/V$ in the SHK case, with
$\eps_j\downarrow0$.
These measures have positive bounded densities and converge to $\mu$;
convexity of $E$ gives $\limsup_j\sfE(y_j)\le\sfE(\mu)$ and lower
semicontinuity gives $\liminf_j\sfE(y_j)\ge\sfE(\mu)$.

Let $z_j$ be a minimiser of one step of the scheme from $y_j$, with
$\tau_j\downarrow0$ small
enough for Lemma~\ref{cp:exist} and the HK upper-density bound.
Writing $r_j=\HK(y_j,z_j)$, $M_j=y_j(X)$, minimality,
\eqref{cp:eq-affine} and $\sqrt{z_j(X)}\le\sqrt{M_j}+r_j$ give
\[
 r_j^2/(2\tau_j)\le\sfE(y_j)+a(\sqrt{M_j}+r_j)^2+bV,
 \quad r_j^2\le4\tau_j(\sfE(y_j)+2aM_j+bV)\to0
\]
when $a\tau_j\le1/8$. In SHK instead
$r_j^2\le2\tau_j(\sfE(y_j)+a+bV)\to0$.
Minimality and lower semicontinuity give energy convergence, and the
slope estimate for a minimiser of one step, which is classical
\cite{AGS05} and uses minimality only, gives
$\slope{\sfE}(z_j)\le r_j/\tau_j<\infty$.
Propositions~\ref{cp:density-hk} and~\ref{cp:density-shk} give positive
lower and finite upper density bounds. For HK, fix
$e_1\in\partial E(1)$, require $1+2\tau_j\min\{e_1,0\}>0$, and apply
the upper bound at a threshold at least both $1$ and the upper bound of
$y_j$; monotonicity supplies a subgradient at least $e_1$.
If $\mu\in\calP(X)$ has density $\rho\le M$, then $MV\ge1$, so the cell
averages $v_j$ and the mixtures $y_j$ have densities at most $M$, and
Proposition~\ref{cp:density-shk} gives the same bound for $z_j$.
A diagonal argument through finite-energy points proves density in the
closure as well.
\end{proof}

With the sets $A_j$ of the next lemma, the regular measures are those with
density between two positive constants and finite slope.

\begin{lemma}[Invariant sets]\label{cp:invariance}
Assume that $\sfE$ is weakly geodesically $\lambda$-convex on
$(\calM(X),\HK)$, respectively on $(\calP(X),\SHK)$, and define for $j\in\N$,
$j\ge2$ $($the numbering in Theorem~\ref{cp:abstract} starts at $1$; shift
the index$)$,
\begin{equation}\label{cp:eq-Aj}
 A_j:=\calM_{1/j}(X)\cap\{\slope{\sfE}\le j\},\qquad\text{respectively}\qquad
 A_j:=\calP(X)\cap\calM_{1/j}(X)\cap\{\slope{\sfE}\le j\},
\end{equation}
with the slope taken in the respective space. Then the sets $A_j$ are nested and
satisfy \ref{cp:A1}, with $\kappa_j:=\kappa(X,1/j)$ from
Proposition~\ref{cp:semiconc}, as well as \ref{cp:A2} and \ref{cp:A3} of
Theorem~\ref{cp:abstract}.
\end{lemma}

\begin{proof}
The sets are nested. \ref{cp:A1} holds by Proposition~\ref{cp:semiconc},
because $A_j$ and its geodesic cover are contained in $\calM_{1/j}(X)$ and
its geodesic cover; \ref{cp:A3} is Lemma~\ref{cp:recovery}.

For \ref{cp:A2} let $\mu_0\in A_{j_0}$, $T>0$ and
$\ell=-\min\{\lambda,0\}$. In the HK case let $l=1/j_0$, $u=j_0$ and $L,U$
as in Proposition~\ref{cp:density-hk}\ref{cp:it-iter}, $a$ as in
\eqref{cp:eq-affine}, and
\[
 \tau_*:=\min\Big\{1,\ \frac1{4a+1},\ \frac1{2\ell+1},\ \frac1{4U+1}\Big\}.
\]
The four entries give $\tau\le1$, $2a\tau<1$ for Lemma~\ref{cp:exist},
$\ell\tau\le\frac12$ for the standing assumptions of
Section~\ref{cp:ss-error}, and $\tau U\le\frac14$ for
\eqref{cp:eq-iterexp}.
For $\tau\le\tau_*$ minimisers exist at every step by
Lemma~\ref{cp:exist}, all iterates have positive densities by
Proposition~\ref{cp:density-hk}\ref{cp:it-low}, and for
$n\le N=\lceil T/\tau\rceil$, hence $n\tau<T+1$,
\eqref{cp:eq-iterexp} gives
$l\,\ee^{-4L(T+1)}\le\rho_n\le u\,\ee^{8U(T+1)}$. By \eqref{cp:eq-incr},
which holds under weak convexity,
$\slope{\sfE}(\mu_n)\le\ee^{2\ell(T+1)}j_0$. Hence all iterates up to the
index $N$, for every choice of the minimisers, lie in $A_j$ as soon as
$j\ge j_0\max\{\ee^{8U(T+1)},\ee^{4L(T+1)},\ee^{2\ell(T+1)}\}$.
In the SHK case minimisers exist for all $\tau$, the densities stay in
$[1/j_0,j_0]$ by Proposition~\ref{cp:density-shk}, and the slopes are
bounded in the same way for $\tau\le\min\{1,(2\ell+1)^{-1}\}$.
\end{proof}

\subsection{The main theorem on compact domains}\label{cp:ss-main}

\begin{theorem}[Flows and minimizing movements on compact domains]\label{cp:main}
Let $X\subset\R^d$ be compact and convex with nonempty interior, let
$E:[0,\infty)\to\R$ be finite, convex and lower semicontinuous with $E(0)=0$,
and let $\sfE$ be the relaxed internal energy of Definition~\ref{prelim:energy}.
Let $\lambda\in\R$ and assume that $\sfE$ is weakly geodesically
$\lambda$-convex on $(\calM(X),\HK)$ $[$respectively, on
$(\calP(X),\SHK)$$]$, that is, any two measures of finite energy are joined
by at least one geodesic along which $\sfE$ is $\lambda$-convex.
\begin{enumerate}[label=\textup{(\roman*)}]
\item\label{cp:it-flow} There is a unique $\EVI_\lambda$ flow
$(S_t)_{t\ge0}$ of $\sfE$ on the closure $\overline{\dom\sfE}$ of the
finite-energy domain in $(\calM(X),\HK)$ $[$in $(\calP(X),\SHK)$$]$, which
is the whole space: for every
$\mu_0\in\overline{\dom\sfE}$ the curve $t\mapsto S_t\mu_0$ is the unique
complete $\EVI_\lambda$ solution starting at $\mu_0$, with values in $\calM(X)$
$[$in $\calP(X)$$]$. The flow has the properties listed in
Proposition~\ref{prelim:EVIprop}; in particular
$\sfd(S_t\mu_0,S_t\nu_0)\le\ee^{-\lambda t}\sfd(\mu_0,\nu_0)$ for
$\sfd=\HK$ $[\SHK]$.
\item\label{cp:it-mm} Let $\mu_0\in\dom\sfE$. For all sufficiently small
$\tau>0$ $($for all $\tau>0$ in the $\SHK$ case$)$ exact minimizing movement
schemes started at $\mu_0$ exist, and for every choice of the minimisers the
piecewise constant interpolants and all geodesic interpolants converge to
$t\mapsto S_t\mu_0$ locally uniformly on $[0,\infty)$ as $\tau\downarrow0$.
For every exact scheme with step $0<\tau<1$ and $12\ell\tau\le1$, where
$\ell=-\min\{\lambda,0\}$, the
piecewise constant interpolant satisfies \eqref{prelim:eq-MSerror} with
$\phi=\sfE$ and convexity constant $-\ell$.
\end{enumerate}
\end{theorem}

\begin{proof}
\ref{cp:it-flow} $(\calM(X),\HK)$ and $(\calP(X),\SHK)$ are complete geodesic
spaces, $\sfE$ is proper, lower semicontinuous and weakly geodesically
$\lambda$-convex, and by Lemma~\ref{cp:invariance} the sets \eqref{cp:eq-Aj}
satisfy \ref{cp:A1}--\ref{cp:A3}. Theorem~\ref{cp:abstract} gives the flow,
Lemma~\ref{cp:recovery} identifies $\overline{\dom\sfE}$, and
Proposition~\ref{prelim:EVIprop} gives the properties.

\ref{cp:it-mm} Minimisers exist by Lemma~\ref{cp:exist} and have finite energy,
so that the scheme can be continued. For an initial measure of infinite
slope the error estimate \eqref{cp:eq-error} is void, and convergence
follows from the existence of the flow. With $\ell=-\min\{\lambda,0\}$ as in
Section~\ref{cp:ss-error}, the flow of \ref{cp:it-flow} is an $\EVI_{-\ell}$
flow on the whole space and $\sfE$ is weakly geodesically $(-\ell)$-convex,
so Theorem~\ref{prelim:MS-error}(a) applies with the constant $-\ell$ and
$U^0=u_0=\mu_0$; this is the bound \eqref{prelim:eq-MSerror}. Its
right-hand side tends to zero locally uniformly in time by
Lemma~\ref{prelim:envelope}: the slope inequality of
Lemma~\ref{cp:euler}\ref{cp:it-slope} at a regular measure $z$, which is
the lower bound \eqref{cp:eq-lowerphi} with $z$ in place of $x_0$, together with
$\sfd(z,y)\le\sfd(z,\mu_0)+\sfd(\mu_0,y)$ and Young's inequality gives
$\sfE(y)\ge-a'-(\ell+1)\sfd^2(\mu_0,y)$ for all $y$, with a constant
$a'\ge0$, so that $\sfE_\sigma(\mu_0)\uparrow\sfE(\mu_0)$ as
$\sigma\downarrow0$. The geodesic interpolants converge as stated in
Theorem~\ref{prelim:MS-error}(a).
\end{proof}

Compared with \cite[Thm.~6.6]{LaM22}, the hypothesis is convexity along one
geodesic for each pair of measures, as in Assumption~A of \cite{LaM22}; the
next corollary shows that this is the same as convexity along every
geodesic here. Moreover the scheme converges from every initial measure of
finite energy, whereas in \cite{LaM22} the initial density is bounded above
and below.

\begin{corollary}[Weak and strong convexity coincide]\label{cp:weakstrong}
Let $X$, $E$ and $\lambda$ be as in Theorem~\ref{cp:main}. If $\sfE$ is
weakly geodesically $\lambda$-convex on $(\calM(X),\HK)$ $[$on
$(\calP(X),\SHK)$$]$, then it is strongly geodesically $\lambda$-convex
there: $\sfE$ is $\lambda$-convex along every geodesic joining two measures
of finite energy.
\end{corollary}

\begin{proof}
Theorem~\ref{cp:main} gives the $\EVI_\lambda$ flow on
$\overline{\dom\sfE}$, which is the whole space by
Lemma~\ref{cp:recovery}. Corollary~\ref{cp:strong} applies.
\end{proof}

\begin{remark}[Scope and quoted results]\label{cp:rem-scope}
(a) \emph{Convexity.} The hypothesis is weak convexity, so that the theorem
covers Assumption~A of \cite{LaM22} as written, for lower semicontinuous
$E$; the normalisation $E(0)=0$ changes $\sfE$ by a constant on a compact
domain. Convexity of the energy is
used only along geodesics that can be selected: along the geodesic from a
node to the observer in Lemma~\ref{cp:euler}\ref{cp:it-dir}; along one
geodesic between a point and a comparison point, or between consecutive
nodes, in Lemma~\ref{cp:euler}\ref{cp:it-slope}--\ref{cp:it-abc}; along the
interpolating geodesics, which are chosen admissible, in
Lemma~\ref{cp:error} and Proposition~\ref{cp:approxEVI}; and through
\cite[Cor.~3.12, Thm.~5.9]{MuS20}, which are formulated for one geodesic
for each pair, respectively for approximate convexity. The geometric
hypothesis $(\mathrm G_\kappa)$ concerns every geodesic and is therefore
available along the selected ones. If $\phi$ is strongly geodesically
$\lambda$-convex, every geodesic joining two points of $\dom\phi$ is
admissible and no selection is needed. Convexity that holds only on
sublevels, with a constant depending on the sublevel, requires a separate
argument; see Remark~\ref{cp:rem-local} and
Theorem~\ref{sub:abstract}.

(b) \emph{Energies.} The theorem covers the energies of
Theorem~\ref{prelim:LMS} and those spherical examples of
Sections~\ref{sec:transfer} and~\ref{sec:examples} that are weakly
geodesically $\lambda$-convex on the whole of $(\calP(X),\SHK)$ for some
$\lambda\in\R$ $($for $\int\rho^m$, $m>1$, on a domain of diameter larger
than $\frac\pi2$ there is no such $\lambda$ by Theorem~\ref{tr:necessary};
see Section~\ref{sec:sublevel}$)$, including energies
with infinite slope at density-bounded measures or noncompact sublevels
(Examples~\ref{cp:ex-slope} and~\ref{cp:ex-sublevel}). Unlike the original
proofs, the proofs given here do not use differentiability or strict
convexity of $E$, and scheme convergence holds for every finite-energy
initial measure.

(c) \emph{Quoted results.} The following results are used and not proved
here: the metric facts on $\HK$ and $\SHK$ of
Propositions~\ref{prelim:HKfacts} and~\ref{prelim:cone} (completeness, the
geodesic property, metrisation of narrow convergence, subadditivity,
scaling, the cone structure) \cite{LMS18,LaM19}; the transport--growth
representation of $\HK$ for absolutely continuous sources
(Lemma~\ref{prelim:tg}) \cite{LMS23}; lower semicontinuity of $\sfE$
(Lemma~\ref{prelim:lsc}); the semiconcavity of $\frac12\HK^2$ on
density-bounded sets \cite[Thm.~6.1]{LaM22}, whose proof rests on
\cite[Lemmas~4.9--4.11]{LaM19}, the transfer of semiconcavity from a cone to
its base \cite[Prop.~2.27]{LaM19} and the local angle condition
\cite[Thm.~4.2]{LaM19}; and the results of \cite{MuS20} on $\EVI$ flows:
the integral characterisation and the properties in
Proposition~\ref{prelim:EVIprop}, the extension by density
\cite[Cor.~3.8]{MuS20}, the change of the constant
\cite[Cor.~3.12]{MuS20}, the convexity implied by a flow
(Theorem~\ref{prelim:MS-convex}) and the error estimates of
Theorem~\ref{prelim:MS-error}. Of these, the existence, uniqueness and
contractivity of the flow in Theorem~\ref{cp:main}\ref{cp:it-flow} need only
the integral characterisation, \cite[Cor.~3.8, Cor.~3.12]{MuS20} and the
contraction estimate \eqref{prelim:eq-contraction}: two solutions with the
same initial point coincide by \eqref{prelim:eq-contraction} with
$s\downarrow0$. In this section Theorem~\ref{prelim:MS-convex} is used only for
Corollaries~\ref{cp:strong} and~\ref{cp:weakstrong}, and
Theorem~\ref{prelim:MS-error} only for
Theorem~\ref{cp:main}\ref{cp:it-mm}. The first
variation of the squared distance, the estimate under the local angle
condition, the derivative estimate \eqref{cp:eq-deriv}, the discrete Euler
inequalities and the approximate EVI for the interpolants are stated
without proof and belong to the approach of \cite{Sav07}; their proofs are
part of \cite{MuSII}, as explained in Section~\ref{cp:ss-points}; the
corrected error estimate is proved from them. Everything else used in the proof of
Theorem~\ref{cp:main} is proved in this section.
\end{remark}
\section{The Hellinger--Kantorovich flow on unbounded domains}\label{sec:hk-unbdd}

The theory of Section~\ref{sec:compact} works on sets of measures whose
density lies between two positive constants, and on an unbounded domain
these sets are empty. Nor can one let $R\to\infty$ in the flows on
$\Om\cap\bar B_R$ by a general stability result: such results need compact
sublevels of the energy, and a bump translated to infinity keeps its energy
and has no convergent subsequence.

We construct the flow as the limit of the flows on the balls and prove the
two facts that compactness would have provided. First, no mass escapes to
infinity, uniformly in the radius (Proposition~\ref{hk:tight}); this rests
on a bound for the growth factor in one step of the scheme
(Lemma~\ref{hk:growth}). Second, the flows on two large balls are close
(Proposition~\ref{hk:cauchy}); the difficulty is that the flow on the
larger ball cannot serve as observer for the flow on the smaller one, and
Lemma~\ref{hk:comparison} deals with it. The last subsection shows that the
minimizing movement scheme on $\calM(\Om)$ itself is well defined and
converges to this flow.

Assume throughout that $\Om$ is closed and convex with nonempty interior.
The balls $\Om_R$ and the truncation $P_R$ are those of
Section~\ref{sec:prelim} and Lemma~\ref{prelim:trunc}. Fix $R_{\min}>0$
such that $\Om_{R_{\min}}$ has nonempty interior, and use the compact flows
only for $R\ge R_{\min}$. Such a radius exists by the interior assumption.

Hypothesis \ref{hk:H1} bounds the energy that can be gained by creating
mass; \ref{hk:H2} gives convexity on every ball $\Om_R$ with the same
$\lambda$, so that the flows $S^R$ contract at a rate independent of $R$.

\begin{theorem}[HK flow on unbounded domains]\label{hk:main}
Let $\Om\subset\R^d$ be closed and convex with nonempty interior, possibly
unbounded. Let $E:[0,\infty)\to\R$ be convex and continuous at $0$, with $E(0)=0$ and assume
\begin{enumerate}[label=\textup{(H\arabic*)}]
\item\label{hk:H1} $\ezero=E'(0^+)>-\infty$;
\item\label{hk:H2} $E$ satisfies the condition of Theorem~\ref{prelim:LMS} with
 constant $\lambda\in\R$, so that $\sfE$ is geodesically $\lambda$-convex on
 $(\calM(\Om'),\HK)$ for every closed convex $\Om'\subset\R^d$ with nonempty
 interior.
\end{enumerate}
Then for every $\mu_0\in\overline{\dom\sfE}\subset\calM(\Om)$ there exists a unique
$\EVI_\lambda$ solution $u:[0,\infty)\to\calM(\Om)$ of $(\calM(\Om),\sfE,\HK)$ with
$u(0)=\mu_0$. The solutions form a $\lambda$-contractive semigroup $(S_t)_{t\ge0}$
on $\overline{\dom\sfE}$. Moreover, if $S^R_t$ denotes the $\EVI_\lambda$ flow of
$(\calM(\OmR),\sfE,\HK)$ given by Theorem~\ref{cp:main}, then for every $T>0$
\begin{equation}\label{hk:eq-exhaust}
 \sup_{t\in[0,T]}\HK\big(S^R_t(P_R\mu_0),S_t\mu_0\big)\longrightarrow0
 \qquad\text{as }R\to\infty .
\end{equation}
If $\mu_0\in\dom\sfE$ has support in $\Om_{R_0}$, then $P_R\mu_0=\mu_0$ for
$R\ge R_0$ and the convergence in \eqref{hk:eq-exhaust} has the rate
$C\,R^{-1/6}$ with $C=C(T,\lambda,\ezero,\mu_0(\Om),\sfE(\mu_0))$, for
$R\ge\max\{2R_0,1,R_{\min}\}$.
\end{theorem}

Bounded densities of compact support have finite energy and are dense in
$\calM(\Om)$ (truncate, then average over finite partitions). Thus
$\overline{\dom\sfE}=\calM(\Om)$, without moment assumptions.
Examples include the finite sums $\sum_i\alpha_ic^{p_i}+\beta c$ with
$\alpha_i\ge0$, $p_i>1$ and $\beta\in\R$, for which $\lambda=2\beta$, the
Bernstein energies, and the shifted powers in the parameter ranges of
Theorem~\ref{ex:shifted}. The sign condition matters: convexity of a signed
sum of powers does not imply~\ref{hk:H2}. The hypothesis $\ezero> -\infty$ is substantive:
for $E=-c^q$, $0<q<1$, spreading fixed mass $m$ over volume $V$ gives energy
$-m^qV^{1-q}\to-\infty$ at bounded HK cost. Adding arbitrarily small such
mass makes the energy unbounded below in every HK ball on an unbounded
convex domain, and every step functional has infimum $-\infty$.

Three facts from Section~\ref{sec:prelim} are used throughout. Truncation
to $\bar B_R$ costs at most the tail mass $\mu(\Om\setminus\bar B_R)$ in
squared distance and $\ezero^-$ times the tail mass in energy
\eqref{prelim:eq-truncHK}, and the truncated energies converge for
finite-energy measures (Lemma~\ref{prelim:trunc}). The energy is bounded
below by the mass, $\sfE(\mu)\ge\ezero\mu(\Om)\ge-\ezero^-\mu(\Om)$, by
\eqref{prelim:eq-energybounds}. Restricted HK distances agree with the
ambient distance (Lemma~\ref{prelim:normalise}). We write
$\lambda_0=\min\{\lambda,0\}$ in this section and the next; in the notation
of Section~\ref{cp:ss-error}, $\lambda_0=-\ell$.

\subsection{One step cannot grow mass by more than a fixed factor}

In the transport--growth description of Lemma~\ref{prelim:tg}, a step from
$\mu_0$ to $\mu_1$ moves the mass at $x$ to $\bfT(x)$ and multiplies it by
the square of the growth factor $q(x)$. Transport is paid for in the
distance, so the only way the scheme could send mass to infinity for free
is by growing it. The following lemma bounds $q$ by a constant that depends
on $\tau$ and $\ezero$ only, not on the domain. The idea is that a step
which grows mass too much can be improved by growing less: replacing $q$ by
$tq$, $t<1$, where $q$ is large shortens the distance to $\mu_0$ to first
order in $1-t$, and removes mass from the target at an energy price of at
most $-\ezero$ per unit. This is where $\ezero>-\infty$ enters: without it
the price of removing a thin layer of mass is not proportional to the mass
removed.

\begin{lemma}[One-step growth bound]\label{hk:growth}
Let $X\subset\R^d$ be compact and convex with nonempty interior, let
$E:[0,\infty)\to\R$ be convex and continuous at $0$, with $E(0)=0$ and $\ezero=E'(0^+)>-\infty$, and
let $\tau>0$ satisfy $1+2\tau\ezero>0$. Let $\mu_0=\rho_0\Leb\in\calM(X)$ with
$\rho_0>0$ a.e.\ in $X$, let $\mu_1$ be a minimiser of
\[
 J(\mu)=\frac1{2\tau}\HK^2(\mu_0,\mu)+\sfE(\mu)\qquad\text{on }\calM(X),
\]
and let $(\bfT,q)$ be an optimal transport--growth pair from $\mu_0$ to
$\mu_1$ $($Lemma~\textup{\ref{prelim:tg}}$)$, with transport cosine $\chi$.
Then
\begin{equation}\label{hk:eq-growth}
 q\le k\chi\le k,\qquad k:=\frac1{1+2\tau\ezero},
 \qquad\mu_0\text{-a.e.\ in }X.
\end{equation}
In particular $|\bfT(x)-x|<\pi/2$ for $\mu_0$-a.e.\ $x$ with $q(x)>0$, and
$\mu_1(X)\le k^2\mu_0(X)$.
\end{lemma}

\begin{proof}
If the conclusion fails, there is a positive-measure
set $A\subset\{q>0\}$ on which $k^{-1}-\chi/q\ge\eps>0$.
Let $\beta=\bfT\pfw(q^2\mu_0|_A)$, $n=\beta(X)>0$, and reduce $q$ to $tq$
on $A$, with $t<1$. The new target is
$\mu_1^t=\mu_1-(1-t^2)\beta$. The transport--growth upper bound
\eqref{prelim:eq-tg-upper} and \eqref{prelim:eq-tg-opt} give
\[
 \HK^2(\mu_0,\mu_1)-\HK^2(\mu_0,\mu_1^t)
 \ge(1-t)\int_Aq^2(1+t-2\chi/q)\dd\mu_0.
\]
Since $E(c)-\ezero c$ is nondecreasing and $\Einf\ge\ezero$,
$\sfE(\mu_1)-\sfE(\mu_1^t)\ge\ezero(1-t^2)n$, including singular mass.
Thus
$J(\mu_1)-J(\mu_1^t)\ge(1-t)n[2\eps-(1-t)k^{-1}]/(2\tau)>0$
for $t$ close to $1$, a contradiction. The transport and mass conclusions
follow from $q\le k\chi$.
\end{proof}

The constant $k$ does not depend on $X$. Besides the mass bound, the lemma
says that mass which survives the step travels less than $\pi/2$; this is
used in the next subsection and again in the spherical case.

\subsection{Tightness, uniformly in the size of the domain}

We now show that the flows $S^R$ keep their mass in a bounded region,
uniformly in $R$. Mass can appear far out by transport or by creation.
Transport is paid for in the HK cost, and the sum of the squared distances
over the steps is bounded by the drop of the energy; creation multiplies the
mass already present by at most $k^2$ per step. Inequality
\eqref{hk:eq-onestep} below records both effects for a Lipschitz cutoff.
The constant $K_T$ bounds the growth of mass up to time $T$, and $D_T$ the
total drop of the energy.
For $T>0$ and $\mu_0\in\dom\sfE$ with $M_0:=\mu_0(\Om)$ we set
\begin{equation}\label{hk:eq-KD}
 K_T:=\ee^{5\ezero^-T}\ge1,\qquad
 D_T:=\sfE(\mu_0)+\ezero^-K_TM_0\ge0 ,
\end{equation}
where $D_T\ge0$ by the lower bound in \eqref{prelim:eq-energybounds}.

\begin{proposition}[Uniform bounds and tightness]\label{hk:tight}
Assume \ref{hk:H1} and \ref{hk:H2}. Let $\mu_0\in\dom\sfE$ have support in
$\Om_{R_0}$, let $T>0$, and for $R\ge\max\{R_0,R_{\min}\}$ let
$u^R(t)=S^R_t\mu_0$ be the
$\EVI_\lambda$ solution of $(\calM(\OmR),\sfE,\HK)$ starting at $\mu_0$. Then
for all $R\ge\max\{R_0,R_{\min}\}$ and $t\in[0,T]$:
\begin{enumerate}[label=\textup{(\roman*)}]
\item $u^R(t)(\Om)\le K_TM_0$;
\item $-\ezero^-K_TM_0\le\sfE(u^R(t))\le\sfE(\mu_0)$;
\item for every $\rho>0$ with $\rho\ge2R_0$,
 \begin{equation}\label{hk:eq-tail}
  u^R(t)\big(\Om\setminus B_\rho\big)\le\frac{C_1}{\rho},\qquad
  C_1:=\pi\,K_T^{5/2}\sqrt{M_0\,T\,D_T}.
 \end{equation}
\end{enumerate}
The constants do not depend on $R$.
\end{proposition}

\begin{proof}
First take regular initial data on $\Om_R$, that is, with density between
two positive constants and finite slope (Section~\ref{cp:ss-abstract}).
Density bounds persist by Proposition~\ref{cp:density-hk}.
For $\tau\ezero^-\le1/10$ and $\tau\le T$, Lemma~\ref{hk:growth} gives
$k^{2n}\le\ee^{5\ezero^-n\tau}\le K_T$ and $\mu_n(\Om)\le K_TM_0$ for
$n\tau\le T$. Minimality and the lower bound in
\eqref{prelim:eq-energybounds} give the action bound
\begin{equation}\label{hk:eq-action}
 \sum_{j=1}^n\HK^2(\mu_{j-1},\mu_j)\le2\tau D_T.
\end{equation}
For any $L$-Lipschitz cutoff $0\le\varphi\le1$, put
$m_j=\int\varphi\dd\mu_j$. A useful estimate, also used below on the
sphere, is
\begin{equation}\label{hk:eq-onestep}
 m_j\le k^2m_{j-1}+\frac\pi{2\sqrt2}Lk^{3/2}
 \mu_{j-1}(\Om)^{1/2}\HK(\mu_{j-1},\mu_j).
\end{equation}
Indeed, $m_j\le k^2m_{j-1}+L\int q^2|\bfT-x|\dd\mu_{j-1}$.
Cauchy--Schwarz bounds the latter integral by
$k^{3/2}\mu_{j-1}(\Om)^{1/2}(\int q|\bfT-x|^2\dd\mu_{j-1})^{1/2}$.
Since transport has length below $\pi/2$ where $q>0$,
$\HK^2\ge(8/\pi^2)\int q|\bfT-x|^2\dd\mu_{j-1}$.
Iteration, $k^{2(n-j)},k^{3/2}\le K_T$, and \eqref{hk:eq-action} yield
\begin{equation}\label{hk:eq-cutoff}
 \int\varphi\dd\mu_n\le K_T\int\varphi\dd\mu_0
 +\tfrac\pi2 L K_T^{5/2}\sqrt{M_0TD_T}.
\end{equation}
Theorem~\ref{cp:main}\ref{cp:it-mm} on $\Om_R$ passes this bound to the flow:
for $t\le T$ use the steps $\tau=t/m$, $m\to\infty$, for which $t$ is a node
with $m\tau=t\le T$.
For general
finite-energy initial data use the energy-recovery sequence of
Lemma~\ref{cp:recovery}; contraction passes the curves to the limit, while
mass, cutoff integrals and initial energies converge. No uniform density
bound in $R$ is needed. Finally take
$\varphi(x)=\min\{1,(2|x|/\rho-1)^+\}$, $L=2/\rho$.
The energy bounds follow from monotonicity and the mass bound.
\end{proof}

\subsection{The flows on different balls form a Cauchy family}

Let $R'>R$. For two EVI solutions on the same space, closeness is the
contraction estimate, proved by using each solution as observer in the EVI
of the other. Here this works in one direction only: $u^R(t)$ is a measure
on the larger ball, but $u^{R'}(s)$ may charge $\Om_{R'}\setminus\Om_R$ and
is then not a point of the space on which $u^R$ satisfies its EVI. We use
the truncation $P_Ru^{R'}(s)$ as observer instead. By
\eqref{prelim:eq-truncHK} this costs an error of order $\sqrt\eps$ in the distance
and raises the energy by at most a multiple of $\eps$, where $\eps\sim1/R$ bounds the mass of
$u^{R'}(s)$ outside $B_R$. The error in the distance enters once for each
time interval on which the EVI is applied, however short, so it cannot be
absorbed in a differential inequality. We therefore apply the EVI on a grid
of $N$ time steps: the distance errors add up to $N^2\eps$, and the factors
$1+1/N$ from Young's inequality cost $D/N$ through the energy terms. With
$N\sim R^{1/3}$ both are of order $R^{-1/3}$ in the squared distance, which
is the origin of the rate $R^{-1/6}$ in Theorem~\ref{hk:main}. Only the
EVI and the two error bounds are used, so we state the argument in a metric
space; on the sphere $P_R$ is replaced by the normalised restriction.
\begin{lemma}[Comparison with truncated observers]\label{hk:comparison}
Let $(\calX,\sfd)$ be a metric space, $\phi:\calX\to(-\infty,\infty]$, and
let $u,v:[0,T]\to\dom\phi$ be curves with the same initial value such that
$t\mapsto\phi(u(t))$ and $t\mapsto\phi(v(t))$ are nonincreasing with values
in one interval of length $D$. Let $P$ map the points $v(s)$, $s\le T$, into
$\dom\phi$, and let $\eps,\delta,K\ge0$ be fixed numbers with
$\sfd(Pv(s),v(s))\le\delta\sqrt\eps$ and
$\phi(Pv(s))\le\phi(v(s))+K\eps$ for all $s\le T$. Suppose $u$ satisfies
the exponential EVI~\eqref{prelim:eq-expEVI} with constant $\lambda$ for
the observers $Pv(s)$, and $v$ satisfies it for the observers $u(t)$, for
all $s,t\le T$.
Then for every integer $N\ge1$,
\begin{equation}\label{hk:eq-aN}
 \sup_{t\le T}\sfd^2(u(t),v(t))
 \le C(T,\lambda,\delta,K)(N^2\eps+D/N+\eps).
\end{equation}
\end{lemma}
\begin{proof}
Fix $t_*>0$, put $h=t_*/N$, $a_k=\sfd(u(kh),v(kh))$,
$M=M_\lambda(h)$ as in Proposition~\ref{prelim:EVIprop}, and
$\theta=\ee^{-\lambda h}(1+1/N)$. Apply the EVI for $u$ with $Pv_k$, and
then for $v$ with $u_{k+1}$. At both ends of the first comparison use
$(x+\delta\sqrt\eps)^2\le(1+1/N)x^2+(1+N)\delta^2\eps$.
Substitution gives
\begin{align*}
 \tfrac12a_{k+1}^2\le{}&\tfrac12\theta^2a_k^2
 +\tfrac12\ee^{-\lambda h}(1+\theta)(1+N)\delta^2\eps\\
 &+M[\phi(v_k)-\phi(v_{k+1})]
 +M[|\theta-1|D+\theta K\eps].
\end{align*}
All remainders are nonnegative. The weights in iteration are bounded by
$\ee^{2|\lambda|T+2}$; also $M\le\ee^{|\lambda|T}T/N$ and
$|\theta-1|\le(2|\lambda|T\ee^{|\lambda|T}+1)/N$.
The energy drops telescope to at most $D$. Summing from $a_0=0$ proves the
bound, uniformly in $t_*\le T$.
\end{proof}

\begin{proposition}[Cauchy estimate]\label{hk:cauchy}
Assume \ref{hk:H1} and \ref{hk:H2}. Let $\mu_0\in\dom\sfE$ have support in
$\Om_{R_0}$ and let $T>0$. Then there is a constant
$C_2=C_2(T,\lambda,\ezero,\mu_0(\Om),\sfE(\mu_0))$ such that, with
$u^R(t)=S^R_t\mu_0$,
\begin{equation}\label{hk:eq-cauchy}
 \sup_{t\in[0,T]}\HK\big(u^R(t),u^{R'}(t)\big)\le C_2\,R^{-1/6}
 \qquad\text{for all }R'>R\ge\max\{2R_0,1,R_{\min}\}.
\end{equation}
\end{proposition}

\begin{proof}
Let $u=u^R$, $v=u^{R'}$. Use Lemma~\ref{hk:comparison} with $\phi=\sfE$,
$P=P_R$, $\delta=1$, $K=\ezero^-$, $D=D_T$ and $\eps=C_1/R$. By
Proposition~\ref{hk:tight}\,(iii) with $\rho=R\ge2R_0$ the mass of $v(s)$
outside $B_R$ is at most $\eps$ for $s\le T$, and by
Proposition~\ref{hk:tight}\,(ii) the energies of both flows lie in one
interval of length $D_T$. Then \eqref{prelim:eq-truncHK} gives
the two bounds on $P_Rv(s)$; in particular $P_Rv(s)$ has finite energy and
is supported in $\OmR$, so it can be used as observer in the EVI of the
flow on $\OmR$, and $u(t)$ as observer in the EVI of the flow on
$\Om_{R'}$.
Taking $N=\lceil R^{1/3}\rceil$, and using $R\ge1$ for
$N^2\eps\le4C_1R^{-1/3}$, gives the stated rate.
\end{proof}

\subsection{Proof of Theorem~\ref{hk:main}}

It remains to show that the limit curve given by the Cauchy estimate
satisfies the EVI against every observer, to extend the construction to all
initial data by contraction, and to prove \eqref{hk:eq-exhaust} for initial
data without compact support. The energy term of the EVI is passed to the
limit by Fatou's lemma, using the lower bound in
\eqref{prelim:eq-energybounds}.

\begin{proof}[Proof of Theorem~\ref{hk:main}]
For compactly supported finite-energy $\mu_0$, Proposition~\ref{hk:cauchy}
and completeness give a locally uniform limit $u$ of $S^R_\cdot\mu_0$.
Its mass is at most $K_TM_0$, and its energy lies in
$[-\ezero^-K_TM_0,\sfE(\mu_0)]$. With a compactly supported finite-energy
observer $w$, the integrated EVI \eqref{prelim:eq-intEVI} of the flows
$S^R$ $($Theorem~\ref{cp:main}$)$ passes to this limit: distances
converge and the energy term satisfies Fatou's lemma with that common
lower bound. Truncating any $w\in\dom\sfE$ gives
$P_Rw\to w$, $\sfE(P_Rw)\to\sfE(w)$, so the inequality holds for all
observers. Proposition~\ref{prelim:EVIprop} gives the EVI solution,
uniqueness and contraction.

For arbitrary $\mu_0\in\overline{\dom\sfE}$, approximate by compactly
supported finite-energy $\mu_0^j$. Contraction makes their flows uniformly
Cauchy on $[0,T]$. The masses are uniformly bounded by metric convergence,
so the lower bound in \eqref{prelim:eq-energybounds} again permits Fatou.
The positive-time bound
\begin{equation}\label{hk:eq-regular}
 \sfE(S_t\mu_0^j)\le\sfE(w)
 +\frac{\HK^2(\mu_0^j,w)}{2\int_0^t\ee^{\lambda r}\dd r}
\end{equation}
passes by lower semicontinuity and ensures finite energy for $t>0$ even
when the initial energy is infinite. Thus the limit is the unique EVI
solution. Uniqueness also gives the semigroup property.

For exhaustion with arbitrary initial data, fix compactly supported
$y\in\dom\sfE$. Contraction on both domains gives, with
$C_T=\ee^{-\lambda_0T}$,
\[
 \sup_{t\le T}\HK(S^R_tP_R\mu_0,S_t\mu_0)
 \le C_T[\HK(P_R\mu_0,\mu_0)+2\HK(\mu_0,y)]
 +\sup_{t\le T}\HK(S^R_ty,S_ty).
\]
Let $R\to\infty$ and then $y\to\mu_0$. The rate for compactly supported
finite-energy data follows from \eqref{hk:eq-cauchy} by letting
$R'\to\infty$.
\end{proof}

\subsection{The minimizing movement scheme on the unbounded domain}\label{hk:ss-mm}

The flow was obtained through the balls $\Om_R$; we now show that the scheme
on $\calM(\Om)$ itself is well defined and converges to it. Existence of a
minimiser is not automatic on an unbounded domain, because a minimising
sequence could lose mass at infinity. The following lemma excludes this:
far mass is too expensive. Let $\mu_0$ have all but $\eps$ of its mass in
$B_\rho$ and let $\mu$ have mass $a$ beyond distance $\pi/2$ from $B_\rho$.
Only the part of $\mu_0$ outside $B_\rho$ can be transported there, and by
Lemma~\ref{prelim:split} deleting the far mass changes the squared distance
to $\mu_0$ by at most $-a+2\sqrt{a\eps}$, while it raises the energy by at
most $-\ezero a$ (Lemma~\ref{prelim:dominated}). The distance therefore
strictly decreases when $a>4\eps$, and for $1+2\tau\ezero>0$ the gain
outweighs the energy as soon as $a$ is large compared with $\eps$; for
smaller $a$ the cut need not help, and then $a$ is already controlled by
$\eps$. \ref{hk:H1} is used as in Lemma~\ref{hk:growth}.

\begin{lemma}[Existence of exact minimisers]\label{hk:minimiser}
Let $\Om\subset\R^d$ be closed and convex with nonempty interior, let
$E:[0,\infty)\to\R$ be convex and continuous at $0$, with $E(0)=0$ and $\ezero=E'(0^+)>-\infty$, and
let $\tau>0$ satisfy $1+2\tau\ezero>0$. Then for every
$\mu_0\in\calM(\Om)$ the functional
\[
 J(\mu)=\frac1{2\tau}\HK^2(\mu_0,\mu)+\sfE(\mu)
\]
attains its infimum on $\calM(\Om)$.
\end{lemma}

\begin{proof}
Let $k=(1+2\tau\ezero)^{-1}$ as in \eqref{hk:eq-growth}. Writing
$M=\mu(\Om)$, $M_0=\mu_0(\Om)$, one has
$J(\mu)\ge[k^{-1}M-2\sqrt{MM_0}+M_0]/(2\tau)$ by
\eqref{prelim:eq-energybounds} and \eqref{prelim:eq-massbound}.
Thus $\inf J$ is finite and minimizing sequences have bounded mass.
Let $F=\Om\setminus B_{\rho+\pi/2}$,
$\eps=\mu_0(\Om\setminus B_\rho)$ and $a=\mu(F)$.
The far-cut estimate of Lemma~\ref{prelim:split} and
$\sfE(\mu)-\sfE(\mu|_{F^c})\ge\ezero a$ give
\begin{equation}\label{hk:eq-cutJ}
 J(\mu|_{F^c})-J(\mu)\le[2\sqrt{a\eps}-k^{-1}a]/(2\tau).
\end{equation}
For a minimizing sequence $\mu^j$, let
$\delta_j=J(\mu^j)-\inf J\to0$. Solving the resulting quadratic yields
$a_j\le4k^2(\eps+\tau k^{-1}\delta_j)$ for $a_j=\mu^j(F)$.
Consequently
\[
 \limsup_j\mu^j(\Om\setminus B_{\rho+\pi/2})
 \le4k^2\mu_0(\Om\setminus B_\rho)\longrightarrow0
 \qquad(\rho\to\infty).
\]
The finitely many initial members cause no tightness problem. Prokhorov's
theorem and lower semicontinuity of HK and $\sfE$ give a minimizer.
\end{proof}

Convergence of the scheme now follows from the existence of the flow. The
error estimates of Muratori and Savar\'e \cite{MuS20}, recalled in
Theorem~\ref{prelim:MS-error}, compare any exact scheme with the EVI flow
and assume only that the flow and exact steps exist and that the energy is
convex on a geodesic space. What has to be checked here is that the Moreau
envelope $\sfE_\sigma$ in the first estimate is finite at the scale
$\sigma$ which that estimate prescribes; Lemma~\ref{prelim:envelope} gives
finiteness only for small $\sigma$.

\begin{theorem}[Exact minimizing movements]\label{hk:mm}
Under the hypotheses of Theorem~\ref{hk:main}, let $\lambda_0=\min\{\lambda,0\}$ and assume
$0<\tau<1$, $-12\lambda_0\tau\le1$, $1+2\tau\ezero>0$.
Every finite-energy initial value admits exact steps on $\calM(\Om)$.
Every resulting scheme satisfies the two estimates of
Theorem~\ref{prelim:MS-error}(a), with distance HK and $\lambda=\lambda_0$,
and converges locally uniformly to the EVI flow when started at a fixed
finite-energy datum. For infinite-energy initial data, finite-energy
initial values can be chosen converging to that datum so that every scheme
from those values converges locally uniformly to the flow.
\end{theorem}

\begin{proof}
Lemma~\ref{hk:minimiser} gives exact steps, and minimality keeps their
energies finite. The space $(\calM(\Om),\HK)$ is geodesic and $\sfE$ is
geodesically $\lambda_0$-convex by \ref{hk:H2}, so
Theorem~\ref{prelim:MS-error}(a) applies to the EVI flow, which is an
$\EVI_{\lambda_0}$ flow as well. It remains to bound the Moreau envelope
from below at the scale $\sigma$ defined below. We will use
$1+2\sigma\ezero>0$. First,
convexity along the ray $t^2c\mathbf1_A\Leb$, for a Borel set $A\subset\Om$ with
$0<|A|<\infty$, first as $t\downarrow0$ and then
$c\downarrow0$, gives $\lambda\le2\ezero$. Let
$\sigma=\int_0^{3\sqrt{\tau t_\tau}}\ee^{\lambda_0r}\dd r$. Then
$1+2\sigma\ezero>0$: if $\lambda_0<0$, then $\sigma<1/|\lambda_0|$ and
$2\ezero\ge\lambda_0$, so $2\sigma\ezero>-1$; if $\lambda_0=0$, then
$\ezero\ge\lambda/2\ge0$. Using $\sfE(\nu)\ge\ezero\nu(\Om)$,
\eqref{prelim:eq-massbound} and minimisation over the mass of $\nu$,
$\sfE_\sigma(\mu)\ge\ezero\mu(\Om)/(1+2\sigma\ezero)>-\infty$.
That the envelope increases to the energy is Lemma~\ref{prelim:envelope}:
by \eqref{prelim:eq-energybounds} and \eqref{prelim:eq-massbound},
$\sfE(\nu)\ge-\ezero^-\nu(\Om)
\ge-2\ezero^-\mu(\Om)-2\ezero^-\HK^2(\mu,\nu)$ for all $\nu$, so
$\sfE_\sigma(\mu)\uparrow\sfE(\mu)$ as $\sigma\downarrow0$. Hence the
right-hand side of \eqref{prelim:eq-MSerror} tends to zero uniformly on
$[0,T]$.
For infinite-energy data choose finite-energy $\mu_j\to\mu_0$, and
thresholds $\tau_j\downarrow0$ making the scheme error from $\mu_j$ at
most $1/j$ on $[0,j]$ whenever $\tau\le\tau_j$. Use $\mu_j$ as initial
value for $\tau\in(\tau_{j+1},\tau_j]$ and apply the initial-distance term
in the same estimate.
\end{proof}

\begin{remark}\label{hk:rem-scope}
A first minimizer exists even from infinite-energy data. Convergence
started directly at such data is not asserted: the last clause of
Theorem~\ref{hk:mm} uses prepared finite-energy initial values.
\end{remark}
\section{The spherical flow on unbounded domains}\label{sec:shk-unbdd}

We now carry out the construction of Section~\ref{sec:hk-unbdd} on
$(\calP(\Om),\SHK)$: flows on the balls $\Om_R$, tightness uniformly in $R$,
and a Cauchy estimate with truncated observers. What changes is the
competitor. In Section~\ref{sec:hk-unbdd} we twice improved a measure by
removing mass from it. On the sphere mass cannot simply be removed: the
competitor has to be a probability measure, so we must divide by the mass
that is left, and this raises the density everywhere else. The energy cost
is at most the mass moved times the oscillation $\omega$ of $E'$
(Lemma~\ref{shk:energy}), which is why we assume that $E'$ is bounded.

Let $\Om$ be closed and convex with nonempty interior. As in
Section~\ref{sec:hk-unbdd}, choose $R_{\min}>0$ with
$\operatorname{int}\Om_{R_{\min}}\ne\emptyset$, and use compact flows only
for $R\ge R_{\min}$. Throughout this section, $E:[0,\infty)\to\R$ is
convex and continuous at $0$, satisfies $E(0)=0$, and obeys
\begin{equation}\label{shk:eq-bdd}
 -\infty<\ezero\le E'(c)\le\Einf<\infty,\qquad
 \omega:=\Einf-\ezero.
\end{equation}
Here either one-sided derivative may be used. The extra upper bound permits
normalization after removing mass. It costs nothing for an $\EVI_\lambda$
theory: by Theorem~\ref{tr:necessary}, on a domain of diameter larger than
$\pi/2$ an energy with $\ezero>-\infty$ that is $\lambda$-convex on the
sphere has bounded derivative in any case.

The one-step cutoff inequality \eqref{hk:eq-onestep}, the comparison lemma
(Lemma~\ref{hk:comparison}) and the passage from the Cauchy estimate to the
flow are taken over from Section~\ref{sec:hk-unbdd}, with $\omega$ in place
of $\ezero^-$ in the constants. By \eqref{prelim:eq-energybounds} the
energy of a probability measure lies between $\ezero$ and $\Einf$, so every
probability measure has finite energy and the energy varies by at most
$\omega$ along the flows. As in Section~\ref{sec:hk-unbdd},
$\lambda_0=\min\{\lambda,0\}$. We prove
the growth bound in Section~\ref{shk:ss-elem}, the flow in
Section~\ref{shk:ss-tight}, and treat the scheme in
Sections~\ref{shk:ss-min} and~\ref{shk:ss-mm}.

\begin{theorem}[Spherical flow on unbounded domains]\label{shk:main}
Let $\Om\subset\R^d$ be closed and convex with nonempty interior, possibly
unbounded. Let $E:[0,\infty)\to\R$ be convex and continuous at $0$, with $E(0)=0$ and with bounded
derivative in the sense of \eqref{shk:eq-bdd}, and put
$\omega=\Einf-\ezero$. Assume that $E$ satisfies the condition of Liero, Mielke
and Savar\'e, Theorem~\ref{prelim:LMS}, with constant $\Lambda\in\R$, and put
\[
 \lambda:=\Lambda-2\Einf .
\]
Then for every $\nu_0\in\calP(\Om)$ there is a unique $\EVI_\lambda$ solution
$u:[0,\infty)\to\calP(\Om)$ of $(\calP(\Om),\sfE,\SHK)$ with $u(0)=\nu_0$.
Moreover, if $u^R$ denotes the $\EVI_\lambda$ solution of
$(\calP(\OmR),\sfE,\SHK)$ starting at the normalised restriction
$\nu_0|_{\bar B_R}/\nu_0(\bar B_R)$ for $R\ge R_{\min}$ with
$\nu_0(\bar B_R)>0$ (both conditions hold for all sufficiently large $R$;
the restriction is $\nu_0$ itself if
$\supp\nu_0\subset\bar B_R$), then $u^R\to u$ in $\SHK$, uniformly on bounded
time intervals, as $R\to\infty$. If $\supp\nu_0\subset\bar B_{R_0}$, then
for every $T>0$ the convergence on $[0,T]$ has the rate $C\,R^{-1/6}$ with
$C=C(T,\lambda,\omega)$, for $R\ge\max\{R_1(T,\omega,R_0),R_{\min}\}$.
\end{theorem}

This includes the bounded-slope energies in the rows of
Table~\ref{ex:table} above the rule, in the parameter ranges stated there;
bounded slopes alone do not imply the condition of
Theorem~\ref{prelim:LMS} (Example~\ref{ex:ex-d2}(c)).
Positive powers require the sublevel theory below: on domains of diameter
greater than $\pi/2$, they have no global spherical convexity constant.
The constant $\lambda=\Lambda-2\Einf$ is the one given by the transfer from
the cone to the sphere (Corollary~\ref{tr:global}); it is the same on every
ball, which is what the exhaustion needs. The corollary applies on $\Om$
itself as well: bounded $E'$ gives $\Einf<\infty$ and $\ezero>-\infty$, and
the LMS condition gives strong $\Lambda$-convexity on $(\calM(\Om),\HK)$ by
Theorem~\ref{prelim:LMS}. Hence $\sfE$ is strongly geodesically
$\lambda$-convex on $(\calP(\Om),\SHK)$; the construction of the flow uses
this on the balls only, the scheme on $\calP(\Om)$ in
Section~\ref{shk:ss-mm} uses it on $\Om$.

\subsection{Energy and one-step estimates}\label{shk:ss-elem}

The following lemma is the basic estimate: if mass $m$ is
taken away in one place and put back in another, each unit put back costs
at most $\Einf$ and each unit taken away returns at least $\ezero$, so the
energy rises by at most $\omega m$. The second part applies this to a
measure that is cut and renormalised.

\begin{lemma}[Energy comparison]\label{shk:energy}
Let $E$ be convex and continuous at $0$, with $E(0)=0$ and \eqref{shk:eq-bdd}. Let
$\nu,\mu\in\calM(\Om)$ have equal mass and let
$m:=(\nu-\mu)^+(\Om)=(\nu-\mu)^-(\Om)$. Then
\begin{equation}\label{shk:eq-energy}
 \sfE(\nu)-\sfE(\mu)\le\omega\,m .
\end{equation}
In particular, if $\nu\in\calP(\Om)$, $F\subset\Om$ is a Borel set with
$a:=\nu(F)<1$ and $\nu':=\nu|_{\Om\setminus F}/(1-a)\in\calP(\Om)$, then
\begin{equation}\label{shk:eq-cutE}
 \sfE(\nu')\le\sfE(\nu)+\omega a .
\end{equation}
\end{lemma}

\begin{proof}
Integrate $E(b)-E(a)\le\Einf(b-a)^+-\ezero(b-a)^-$ and use the same
inequality for singular masses; the two variations have equal mass. For
\eqref{shk:eq-cutE} the positive and negative variations of $\nu'-\nu$
both have mass $a$.
\end{proof}

For $F=\Om\setminus\bar B_R$ the measure $\nu'$ is the normalised
restriction $\tilde P_R\nu$ of Lemma~\ref{prelim:trunc}, and
\eqref{shk:eq-cutE} is \eqref{prelim:eq-truncEs}.

Renormalisation is estimated on the cosine of $\SHK$, while the step
functional contains $\tfrac12\SHK^2$. To pass between them we repeatedly use
\begin{equation}\label{shk:eq-G}
 G(z)=\tfrac12\arccos^2z,\qquad
 G(y)-G(x)\le-(y-x)\quad(0\le x\le y\le1),
\end{equation}
since $G'(z)=-\arccos(z)/\sqrt{1-z^2}\le-1$ for $0\le z<1$ and $G$ is
continuous at $1$. In particular $G$ is
decreasing. Note that $\tfrac12\SHK^2=G(\cos\SHK)$.
\label{shk:ss-growth}

The growth bound has the same form as Lemma~\ref{hk:growth}. We again
reduce the growth factor from $q$ to $tq$ where it is too large, but now
the target must be divided by its remaining mass $M_t$. To first order in
$1-t$ this increases the cosine of the distance to $\mu_0$, and it raises
the energy by at most $\omega$ times the mass moved; the threshold in
\eqref{shk:eq-growth} is the value of $q$ above which the gain wins. The
condition $3\tau\omega<1$ plays the role of $1+2\tau\ezero>0$. The argument
uses the energy only through the two bounds \ref{shk:G1} and \ref{shk:G2}
below, and we state the lemma with these as hypotheses, because for the
power energies of Section~\ref{sec:sublevel} they hold with a different
constant $\omega$ (Lemma~\ref{sub:growth}).

\begin{lemma}[Growth bound for one $\SHK$ step]\label{shk:growth}
Let $X\subset\R^d$ be compact and convex with nonempty interior and let
$E:[0,\infty)\to\R$ be convex and continuous at $0$ with $E(0)=0$. Let
$\tau>0$, let $\mu_0=\rho_0\Leb\in\calP(X)$ with $\rho_0>0$ a.e.\ in $X$ and
$\sfE(\mu_0)<\infty$, and let $\mu_1$ be a minimiser of
\[
 J(\nu):=\frac1{2\tau}\SHK^2(\mu_0,\nu)+\sfE(\nu),\qquad\nu\in\calP(X).
\]
Let $(\bfT,q)$ be an optimal transport--growth pair from $\mu_0$ to $\mu_1$
$($Lemma~\textup{\ref{prelim:tg}}$)$, with transport cosine $\chi$. For a
Borel set $A$ with $n_A:=\int_Aq^2\dd\mu_0>0$ and $t\in(0,1)$ put
\begin{equation}\label{shk:eq-comp}
 \beta_A:=\bfT\pfw(q^2\mu_0|_A),\qquad
 \nu_A^t:=\frac{\mu_1-(1-t^2)\beta_A}{1-(1-t^2)n_A}\in\calP(X).
\end{equation}
Suppose that a number $\omega\ge0$ (in this lemma any such number; the
standing assumption \eqref{shk:eq-bdd} is not used) satisfies
$3\tau\omega<1$ and
\begin{enumerate}[label=\textup{(G\arabic*)}]
\item\label{shk:G1} $\sfE(\mu_0)-\sfE(\mu_1)\le\omega$;
\item\label{shk:G2} $\displaystyle\limsup_{t\uparrow1}
 \frac{\sfE(\nu_A^t)-\sfE(\mu_1)}{1-t}\le2\omega\,n_A$ for every such $A$.
\end{enumerate}
Then, with $\theta_1:=\SHK(\mu_0,\mu_1)$,
\begin{equation}\label{shk:eq-growth}
 q\le\frac{\chi}{\cos\theta_1-2\tau\omega}\le\frac1{1-3\tau\omega}
 \qquad\mu_0\text{-a.e.}
\end{equation}
In particular $|\bfT(x)-x|<\frac\pi2$ for $\mu_0$-a.e.\ $x$ with $q(x)>0$.
If $E$ satisfies \eqref{shk:eq-bdd}, then \ref{shk:G1} and \ref{shk:G2}
hold with $\omega=\Einf-\ezero$.
\end{lemma}

\begin{proof}
One step is short by \ref{shk:G1}; we then argue by contradiction with the
renormalised competitor \eqref{shk:eq-comp}.

\emph{Step 1: the step is short.} Put $I=\int_X q\chi\dd\mu_0$. Since
$\int_X q^2\dd\mu_0=\mu_1(X)=1$, \eqref{prelim:eq-tg-opt} gives
$\HK^2(\mu_0,\mu_1)=2-2I$, and the cone identity \eqref{prelim:eq-HKSHK}
gives $I=\cos\theta_1$. Comparing $\mu_1$ with the competitor $\mu_0$ and
using \ref{shk:G1}, $\theta_1^2\le2\tau(\sfE(\mu_0)-\sfE(\mu_1))
\le2\tau\omega$, hence $I\ge1-\theta_1^2/2\ge1-\tau\omega$ and
$I-2\tau\omega\ge1-3\tau\omega>0$. This gives the second inequality in
\eqref{shk:eq-growth}, because $\chi\le1$, and the last assertion follows
from the first inequality, because $q>0$ then forces $\chi>0$.

\emph{Step 2: the competitor.} Suppose that $q\le\chi/(I-2\tau\omega)$
fails. Then there are $k'>(I-2\tau\omega)^{-1}$ and a set $A$ with
$\mu_0(A)>0$ on which $q>0$ and $q\ge k'\chi$. Write $\beta=\beta_A$,
$n=n_A>0$ and
$J_A=\int_Aq\chi\dd\mu_0\le n/k'<n(I-2\tau\omega)$. For $t\in(0,1)$ let
\[
 \mu^t=\mu_1-(1-t^2)\beta,\quad M_t=1-(1-t^2)n,
 \quad\nu^t=\nu_A^t=\mu^t/M_t.
\]
Reducing $q$ to $tq$ on $A$ gives a pair with target $\mu^t$, so the
transport--growth bound \eqref{prelim:eq-tg-upper} gives
$\HK^2(\mu_0,\mu^t)\le1+M_t-2(I-(1-t)J_A)$, and by
\eqref{prelim:eq-normalise}
\begin{equation}\label{shk:eq-Gamma}
 \cos\SHK(\mu_0,\nu^t)\ge
 \Gamma(t):=\frac{I-(1-t)J_A}{\sqrt{M_t}}.
\end{equation}

\emph{Step 3: first-order comparison.} We have $\Gamma(1)=I$ and
$\Gamma'(1)=J_A-In<-2\tau\omega n\le0$, so $\Gamma(t)>I$ for $t<1$ close to
$1$. Since $\Gamma(t)\le1$, this is already impossible if $I=1$. Otherwise
\eqref{shk:eq-G} gives
$\tfrac12\SHK^2(\mu_0,\nu^t)\le G(\Gamma(t))\le\tfrac12\theta_1^2-(\Gamma(t)-I)$,
and minimality of $\mu_1$ gives
\[
 0\le J(\nu^t)-J(\mu_1)\le-\frac{\Gamma(t)-I}{\tau}
 +\sfE(\nu^t)-\sfE(\mu_1).
\]
Divide by $1-t$ and let $t\uparrow1$. The first term on the right tends to
$\Gamma'(1)/\tau=(J_A-In)/\tau$, and by \ref{shk:G2} the upper limit of the
second is at most $2\omega n$. Hence $In-J_A\le2\tau\omega n$, which
contradicts $J_A<n(I-2\tau\omega)$.

\emph{Step 4: bounded $E'$.} Let $E$ satisfy \eqref{shk:eq-bdd} and
$\omega=\Einf-\ezero$. Then \ref{shk:G1} holds because the energies of
probability measures lie in $[\ezero,\Einf]$. For \ref{shk:G2}, the
measure $\nu^t-\mu_1=(M_t^{-1}-1)\mu_1-(1-t^2)\beta/M_t$ has positive
variation at most $(M_t^{-1}-1)\mu_1$, of mass $(1-t^2)n/M_t$, so
Lemma~\ref{shk:energy} gives
$\sfE(\nu^t)-\sfE(\mu_1)\le\omega(1-t^2)n/M_t$, and $(1+t)/M_t\to2$.
\end{proof}

\subsection{Exhaustion}\label{shk:ss-tight}
In this subsection the hypotheses of Theorem~\ref{shk:main} are in force,
in particular the condition of Theorem~\ref{prelim:LMS} with constant
$\Lambda$. On each $\Om_R$, Corollary~\ref{tr:global} and Theorem~\ref{cp:main} give the
flow $S^R$ on $\calP(\OmR)$, with the common constant
$\lambda=\Lambda-2\Einf$.
Tightness and the Cauchy estimate follow as in Section~\ref{sec:hk-unbdd}.
Since all masses equal $1$ and the energy drops by at most $\omega$ along a
scheme, the constants depend on the initial measure only through $R_0$.

\begin{proposition}[Tightness, uniformly in $R$]\label{shk:tight}
Let $T>0$ and $R_0>0$. There is a constant $C_1=C_1(T,\omega)$ such that for
every $\nu_0\in\calP(\Om)$ with $\supp\nu_0\subset\bar B_{R_0}$, every
$\rho\ge2R_0$ and every $R\ge\max\{R_0,R_{\min}\}$,
\begin{equation}\label{shk:eq-tight}
 \sup_{0\le t\le T}\big(S^R_t\nu_0\big)(\Om\setminus B_\rho)\le
 \frac{C_1}{\rho}.
\end{equation}
\end{proposition}

\begin{proof}
For regular initial data, Proposition~\ref{cp:density-shk} preserves positive
upper and lower density bounds. With $21\tau\omega\le1$ and $\tau\le T$,
Lemma~\ref{shk:growth} gives $q\le k:=(1-3\tau\omega)^{-1}$ and
$k^{2n}\le K_{T,\omega}:=\ee^{7\omega T}$ for $n\tau\le T$.
All masses equal $1$ and
$\sum_{j\le n}\HK^2(\mu_{j-1},\mu_j)
\le\sum_{j\le n}\SHK^2(\mu_{j-1},\mu_j)\le2\tau\omega$.
Apply the one-step cutoff inequality \eqref{hk:eq-onestep}, which requires
only $q\le k$ and transport below $\pi/2$, and iterate to obtain
\begin{equation}\label{shk:eq-mt}
 \int\varphi\dd\mu_n\le K_{T,\omega}\int\varphi\dd\mu_0
 +\tfrac\pi2\operatorname{Lip}(\varphi)K_{T,\omega}^2\sqrt{T\omega}.
\end{equation}
Theorem~\ref{cp:main}\ref{cp:it-mm} passes the estimate to the flow
$($for $0<t\le T$ use the steps $\tau=t/m$, $m\to\infty$, for which $t$ is a
node$)$, and
Lemma~\ref{cp:recovery} with contraction removes the regularity assumption.
A cutoff zero on $B_{\rho/2}$, one outside $B_\rho$, with Lipschitz constant
$2/\rho$, proves the result with
$C_1=\pi K_{T,\omega}^2\sqrt{T\omega}$.
\end{proof}

For the comparison of the flows on two balls the truncated observer is the
normalised restriction $\tilde P_R$, which Lemma~\ref{hk:comparison}
covers. The radius $R_1$ below ensures that the tail mass is at most $1/2$,
so that $\tilde P_R$ is defined.

\begin{proposition}[Cauchy property]\label{shk:cauchy}
Let $T>0$, $R_0>0$ and $\nu_0\in\calP(\Om)$ with
$\supp\nu_0\subset\bar B_{R_0}$. There are constants
$C_2=C_2(T,\lambda,\omega)$ and $R_1=R_1(T,\omega,R_0)\ge2R_0$ such that
\begin{equation}\label{shk:eq-cauchy}
 \sup_{0\le t\le T}\SHK\big(S^R_t\nu_0,S^{R'}_t\nu_0\big)\le C_2\,R^{-1/6}
 \qquad\text{for all } R'>R\ge\max\{R_1,R_{\min}\} .
\end{equation}
\end{proposition}

\begin{proof}
Put $R_1:=\max\{2R_0,2C_1,1\}$ with $C_1$ from
Proposition~\ref{shk:tight}, and let $R'>R\ge\max\{R_1,R_{\min}\}$. By
Proposition~\ref{shk:tight} with $\rho=R\ge2R_0$ the mass of
$S^{R'}_s\nu_0$ outside $B_R$ is at most $\eps:=C_1/R\le1/2$ for $s\le T$.
Use the normalised restriction $P=\tilde P_R$.
Lemma~\ref{prelim:trunc} supplies the two bounds on $\tilde P_R$ required
in Lemma~\ref{hk:comparison}: \eqref{prelim:eq-truncSHK} with
$\delta=\pi/2$ and \eqref{prelim:eq-truncEs} with $K=\omega$. The energies
lie in an interval of length $D=\omega$ by \eqref{prelim:eq-energybounds}.
Taking $N=\lceil R^{1/3}\rceil$, and using $R\ge1$ for
$N^2\eps\le4C_1R^{-1/3}$, proves the claim.
\end{proof}

\begin{proof}[Proof of Theorem~\ref{shk:main}]\label{shk:ss-limit}
The probability space is complete by Proposition~\ref{prelim:cone}.
For compactly supported initial data, Proposition~\ref{shk:cauchy} gives
a locally uniform limit. The integrated EVI \eqref{prelim:eq-intEVI} of
the flows $S^R$ $($Theorem~\ref{cp:main}$)$ passes to the limit
for compactly supported observers by Fatou, since $\sfE\ge\ezero$.
For an arbitrary observer $w$, Lemma~\ref{prelim:trunc} and lower
semicontinuity give
$\tilde P_Rw\to w$ and $\sfE(\tilde P_Rw)\to\sfE(w)$, so the limiting
inequality holds for $w$ as well. This identifies the EVI solution.

For arbitrary initial data use normalized truncations and contraction;
their global flows converge uniformly on $[0,T]$, and the same energy
lower bound passes the EVI to this limit. Finally, fixing a compactly
supported approximation $y$ to $\nu_0$,
\[
 \sup_{t\le T}\SHK(S^R_t\tilde P_R\nu_0,S_t\nu_0)
 \le\ee^{|\lambda|T}[\SHK(\tilde P_R\nu_0,\nu_0)+2\SHK(\nu_0,y)]
 +\sup_{t\le T}\SHK(S^R_ty,S_ty).
\]
First $R\to\infty$, then $y\to\nu_0$, proves exhaustion. Uniqueness and
the semigroup property follow from EVI contraction. For initial data with
support in $\bar B_{R_0}$ the rate follows from \eqref{shk:eq-cauchy} by
letting $R'\to\infty$.
\end{proof}

\subsection{Exact minimisers of the step functional}\label{shk:ss-min}

As on the cone, a minimising sequence of the step functional could lose
mass at infinity, and we exclude this by cutting off the far mass $a$ of a
near-minimiser; the cut measure must now be divided by $1-a$. Cutting and
then normalising raises the energy by at most $\omega a$ relative to the
original measure, and under the condition $aI\ge2\sqrt{a\eps}$ of
Lemma~\ref{shk:farcut} the distance to $\nu_0$ does not increase; it
strictly decreases when the cut is nonzero and $\cos\SHK(\nu_0,\nu)>0$. Either $a$ is
already bounded by a multiple of the tail $\eps$ of $\nu_0$, or the cut is
profitable by an amount proportional to $a$, which must then be small
because $\nu^j$ is almost a minimiser.

The gain in distance does not involve the energy. We state it separately,
because it is used again for the power energies in
Theorem~\ref{sub:resolvent}.

\begin{lemma}[Normalised far cut]\label{shk:farcut}
Let $\nu_0,\nu\in\calP(\Om)$ and $\rho>0$, and put
$F=\Om\setminus B_{\rho+\pi/2}$, $\eps=\nu_0(\Om\setminus B_\rho)$,
$a=\nu(F)$ and $I=\cos\SHK(\nu_0,\nu)$. Assume $a<1$ and let
$\nu'=\nu|_{F^c}/(1-a)$. Then
\begin{equation}\label{shk:eq-farcut}
 \cos\SHK(\nu_0,\nu')\ge\frac{I-\sqrt{a\eps}}{\sqrt{1-a}} .
\end{equation}
If moreover $aI\ge2\sqrt{a\eps}$, then
\begin{equation}\label{shk:eq-fargain}
 \tfrac12\SHK^2(\nu_0,\nu')\le\tfrac12\SHK^2(\nu_0,\nu)
 -\Big(\frac{aI}2-\sqrt{a\eps}\Big),
\end{equation}
and $\SHK(\nu_0,\nu')<\SHK(\nu_0,\nu)$ if $a>0$ and $I>0$.
\end{lemma}

\begin{proof}
By the cone identity \eqref{prelim:eq-HKSHK}, $\HK^2(\nu_0,\nu)=2-2I$. The
far-cut estimate \eqref{prelim:eq-farcut} gives
$\HK^2(\nu_0,\nu|_{F^c})\le2-2I-a+2\sqrt{a\eps}$, and
\eqref{prelim:eq-normalise} with $M=1-a$ gives
\[
 \cos\SHK(\nu_0,\nu')
 =\frac{2-a-\HK^2(\nu_0,\nu|_{F^c})}{2\sqrt{1-a}}
 \ge\Gamma:=\frac{I-\sqrt{a\eps}}{\sqrt{1-a}},
\]
which is \eqref{shk:eq-farcut}. Let now $aI\ge2\sqrt{a\eps}$. If $a=0$,
then $\nu'=\nu$ and there is nothing to prove. If $a>0$, then $I\ge0$.
Put $s=\sqrt{1-a}\in(0,1)$, so that $1-s=a/(1+s)>a/2$. Then
\[
 s(\Gamma-I)=(1-s)I-\sqrt{a\eps}\ge\frac{aI}2-\sqrt{a\eps}\ge0,
\]
and dividing the nonnegative left-hand side by $s\le1$ gives
$\Gamma-I\ge\frac{aI}2-\sqrt{a\eps}$; the first inequality is strict if
$I>0$. Thus $0\le I\le\Gamma\le\cos\SHK(\nu_0,\nu')\le1$, and
\eqref{shk:eq-G} gives
$\tfrac12\SHK^2(\nu_0,\nu')\le G(\Gamma)\le G(I)-(\Gamma-I)$, which is
\eqref{shk:eq-fargain}; for $I>0$ we get $\Gamma>I$ and a strict decrease.
\end{proof}

\begin{remark}[The constant $2$ is optimal]\label{shk:farcut-sharp}
On $\Om=\R$ take $\rho=1$ and, for $0<u<1$,
\[
 \nu_0=z^2\delta_0+w^2\delta_4,\qquad \nu=x^2\delta_0+y^2\delta_4,
\]
with $(x,y)=\big(\tfrac{1-u^2}{1+u^2},\tfrac{2u}{1+u^2}\big)$ and
$(z,w)=\big(\tfrac{4-u^2}{4+u^2},\tfrac{4u}{4+u^2}\big)$. Then
$a=y^2$, $\eps=w^2\le a<1$ and $\nu'=\delta_0$. The two sites are more than
$\pi/2$ apart, so by Corollary~\ref{prelim:far}, applied with $\delta_0$ as
the first measure and $F=\{4\}$, and by \eqref{prelim:eq-scaling},
$\HK^2(\nu_0,\delta_0)=(1-z)^2+w^2=2-2z$, that is $\cos\SHK(\nu_0,\nu')=z$.
Subadditivity \eqref{prelim:eq-subadd} and \eqref{prelim:eq-scaling} give
$\HK^2(\nu_0,\nu)\le(z-x)^2+(w-y)^2$, that is $I\ge xz+yw$, and
\eqref{shk:eq-farcut} gives the reverse inequality. Hence
\[
 I-\cos\SHK(\nu_0,\nu')=\frac{2u^4}{(4+u^2)(1+u^2)}>0,
 \qquad
 \frac{aI}{\sqrt{a\eps}}=\frac{4+3u^2+u^4}{2(1+u^2)^2}\to2
 \quad(u\downarrow0).
\]
So for every $C<2$ there is a normalised far cut with
$aI>C\sqrt{a\eps}$ that strictly increases the distance to $\nu_0$.
\end{remark}

\begin{lemma}[Existence of minimisers]\label{shk:minimiser}
Let $\Om\subset\R^d$ be closed and convex with nonempty interior, let $E$ be
convex and continuous at $0$, with $E(0)=0$ and \eqref{shk:eq-bdd}, and let $\tau>0$ with
$3\tau\omega<1$. Then for every $\nu_0\in\calP(\Om)$ the functional
\[
 J(\nu):=\frac1{2\tau}\SHK^2(\nu_0,\nu)+\sfE(\nu)
\]
attains its minimum on $\calP(\Om)$.
\end{lemma}

\begin{proof}
For a minimizing sequence $\nu^j$, let
$\delta_j=J(\nu^j)-\inf J\to0$ and $I_j=\cos\SHK(\nu_0,\nu^j)$.
Energy bounds and comparison with $\nu_0$ give
$\SHK^2(\nu_0,\nu^j)\le2\tau(\omega+\delta_j)$, hence
$I_j\ge1-\tau(\omega+\delta_j)$.
Let $F=\Om\setminus B_{\rho+\pi/2}$,
$\eps=\nu_0(\Om\setminus B_\rho)<1/4$, and $a=\nu^j(F)$; we may assume
$a\ge\eps$.
For large $j$, $a<1$: otherwise the splitting of Lemma~\ref{prelim:split}
gives $\HK^2(\nu_0,\nu^j)\ge1-\eps>3/4$, whereas it is at most
$2\tau(\omega+\delta_j)<3/4$. Let
$\nu'=\nu^j|_{F^c}/(1-a)$ and apply Lemma~\ref{shk:farcut} with
$\nu=\nu^j$.
If $aI_j<2\sqrt{a\eps}$, then $a<4\eps/I_j^2$.
Otherwise \eqref{shk:eq-fargain} bounds the gain in distance, the energy
increases by at most $\omega a$ by \eqref{shk:eq-cutE}, and
$J(\nu')\ge J(\nu^j)-\delta_j$, so that
\[
 \kappa_ja-\sqrt{a\eps}\le\tau\delta_j,
 \qquad \kappa_j=\tfrac12(1-3\tau\omega-\tau\delta_j)>0
\]
for large $j$, where we used $I_j/2-\tau\omega\ge\kappa_j$.
Young's inequality yields
$a\le\eps/\kappa_j^2+2\tau\delta_j/\kappa_j$.
Including $a<\eps$ proves $\limsup_j\nu^j(F)\le C\eps$ with finite $C$.
Hence the sequence is tight; Prokhorov and lower semicontinuity yield a
probability minimizer.
\end{proof}

\subsection{Exact minimizing movements}\label{shk:ss-mm}

As in Section~\ref{hk:ss-mm}, convergence of the scheme follows from the
error estimates of Theorem~\ref{prelim:MS-error}. The envelope term is
easier than on the cone: the energy varies by at most $\omega$ on
$\calP(\Om)$, so the difference between the energy and its Moreau envelope
is bounded by $\omega$ and no preparation of the initial data is needed.

\begin{theorem}[Exact minimizing movements]\label{shk:mm}
Under the hypotheses of Theorem~\ref{shk:main}, let $\lambda_0=\min\{\lambda,0\}$.
If $0<\tau<1$, $3\tau\omega<1$ and $-12\lambda_0\tau\le1$, exact schemes
exist from every probability measure and satisfy
Theorem~\ref{prelim:MS-error}(a) with distance SHK and constant $\lambda_0$.
Here $0\le\sfE(\nu_0)-\sfE_\sigma(\nu_0)\le\omega$, and the difference
tends to zero as $\sigma\downarrow0$. Thus every scheme converges locally
uniformly to the flow, with error $O(\tau^{1/4})$, or $O(\tau^{1/2})$ for
finite initial slope.
\end{theorem}
\begin{proof}
Existence is Lemma~\ref{shk:minimiser}. By Corollary~\ref{tr:global} on
$\Om$, $\sfE$ is strongly geodesically $\lambda$-convex on the geodesic
space $(\calP(\Om),\SHK)$, hence weakly geodesically $\lambda_0$-convex,
and the flow of Theorem~\ref{shk:main} is an $\EVI_{\lambda_0}$ flow on
$\calP(\Om)$. This verifies the hypotheses of
Theorem~\ref{prelim:MS-error}(a). Since $\ezero\le\sfE\le\Einf$ on
$\calP(\Om)$ by \eqref{prelim:eq-energybounds}, the envelope satisfies
$\ezero\le\sfE_\sigma(\nu_0)\le\sfE(\nu_0)\le\Einf$ for every $\sigma>0$,
and Lemma~\ref{prelim:envelope} with $b=0$ gives
$\sfE_\sigma(\nu_0)\uparrow\sfE(\nu_0)$.
\end{proof}
\section{Power energies on the sphere: EVI on energy sublevels}\label{sec:sublevel}

The power energies $\int\rho^m\dd x$, $m>1$, are geodesically convex on the
cone, and their flow on $(\calM(\Om),\HK)$ is given by
Theorem~\ref{hk:main}. On the sphere, by Theorem~\ref{tr:necessary}, a
power energy on a domain of diameter larger than $\pi/2$ is not
$\lambda$-convex for any $\lambda\in\R$, so no $\EVI_\lambda$ theory can
apply to it. This section shows that the flow exists nevertheless. By
Theorem~\ref{tr:powers} a power energy is $\lambda_c$-convex between two
probability measures of energy at most $c$, with $\lambda_c\to-\infty$ as
$c\to\infty$. Since the energy decreases along the minimising movement
scheme, the iterates of a scheme that starts at energy $c$ stay at energy
at most $c$, and the existence proof of Section~\ref{sec:compact} can be
run with a constant tied to the energy level. The limit satisfies the
evolution variational inequality with constant $\lambda_c$ against all
observers of energy at most $c$; we call this a sublevel solution.

Sections~\ref{sub:ss-def}--\ref{sub:ss-mm} are written for a metric space
and contain the definitions, the generation theorem
(Theorem~\ref{sub:abstract}) and the convergence of the scheme from every
initial datum of finite energy (Theorem~\ref{sub:mm}).
Section~\ref{sub:ss-compact} states the result for sums of powers and
proves it on compact domains; Section~\ref{sub:ss-unbdd} contains the proof
on unbounded domains.

Semiconvexity on sublevels and the corresponding EVI are not new as
notions. Sublevel semiconvexity is studied in a finite-dimensional setting
in \cite{LMTZ25}, where a theory of sublevel EVI in the spherical Hellinger
geometry is announced \cite{MST25}; that work was not available to us. The
notion of solution used here, an EVI whose constant depends on the common
energy level of the solution and the observer, is the one of Mielke and
Sumners \cite[Sec.~5.3]{MiSu26}. They treat one-dimensional viscoelasticity
with a dissipation distance comparable to the Bhattacharya (spherical
Hellinger) distance on probability densities, where there is no transport;
solutions are constructed by spatial discretisation, and the sublevel EVI is
proved for the discrete problems and passed to the limit under a uniform
bound on the convexity constant of each sublevel
\cite[Prop.~5.6, Thm.~5.7]{MiSu26}. What is done here is different in three
respects: the solution is generated by the minimising movement scheme in a
metric space whose squared distance is semiconcave only on a family of
sets, under hypothesis (SL); the space is the Hellinger--Kantorovich sphere,
where transport is present; and the domain may be unbounded.

\subsection{Sublevel semiconvexity and sublevel solutions}\label{sub:ss-def}

Convexity along the geodesics between two points of a sublevel, with a
constant that depends on the level, is not quite enough for the existence
proof: a geodesic between two points of energy at most $c$ may pass through
points of higher energy, and in the proof such points are in turn endpoints
of geodesics along which convexity is used. Hypothesis (SL) below therefore
also bounds the energy along the geodesic, by a function $g(c)$ of the
level of its endpoints. Theorem~\ref{tr:powers} verifies (SL) for sums of
powers on $(\calP(\Om),\SHK)$.

Throughout Sections~\ref{sub:ss-def}--\ref{sub:ss-mm}, $(\calX,\sfd)$ is a
metric space and $\phi:\calX\to(-\infty,+\infty]$ is proper and lower
semicontinuous. For $c\in\R$ we write
\[
 \calX_c:=\{x\in\calX:\ \phi(x)\le c\}.
\]

\begin{definition}[Sublevel semiconvexity]\label{sub:SL}
The functional $\phi$ satisfies \emph{(SL)} if there are a nonincreasing
function $c\mapsto\lambda_c\in(-\infty,0]$ and a nondecreasing function
$g:\R\to\R$ with $g(c)\ge c$ such that for every $c\in\R$, all
$x,y\in\calX_c$ and \emph{every} geodesic $\gamma$ from $x$ to $y$ the
following holds: $\phi(\gamma_s)\le g(c)$ for all $s\in[0,1]$, and for all
$0\le s_0<s_1\le1$ and $\vartheta\in[0,1]$,
\begin{equation}\label{sub:eq-SL}
 \phi\big(\gamma_{(1-\vartheta)s_0+\vartheta s_1}\big)
 \le(1-\vartheta)\phi(\gamma_{s_0})+\vartheta\phi(\gamma_{s_1})
 -\frac{\lambda_c}2\,\vartheta(1-\vartheta)\,(s_1-s_0)^2\,\sfd^2(x,y).
\end{equation}
\end{definition}

Inequality \eqref{sub:eq-SL} for all $s_0<s_1$ and $\vartheta$ says that
$\phi$ is $\lambda_c$-convex along $\gamma$ in the sense of
Definition~\ref{prelim:convexity}(a); this is what makes every such
geodesic admissible in the proof of Theorem~\ref{sub:abstract}.
No continuity of $c\mapsto\lambda_c$ or of $g$ is assumed. The geodesic
$\gamma$ need not stay in $\calX_c$; it stays in $\calX_{g(c)}$. With
$s_0=0$, $s_1=1$, \eqref{sub:eq-SL} is the endpoint inequality
\begin{equation}\label{sub:eq-SL-end}
 \phi(\gamma_s)\le(1-s)\phi(x)+s\phi(y)-\frac{\lambda_c}2\,s(1-s)\,\sfd^2(x,y).
\end{equation}

In an $\EVI_\lambda$ the constant $\lambda$ is the convexity constant along
the geodesic from the solution to the observer. Here it depends on the
energy of both, so the inequality is asked for each level separately, as in
\cite[Thm.~5.7]{MiSu26}.

\begin{definition}[Sublevel solution]\label{sub:def-sol}
Let $\phi$ satisfy (SL). A continuous curve $u:[0,\infty)\to\dom\phi$ is a
\emph{sublevel solution} of the evolution variational inequality for
$(\calX,\phi,\sfd)$ if $\phi(u(t))\le\phi(u(0))$ for all $t\ge0$ and, for
every $c\ge\phi(u(0))$, every observer $v\in\calX_c$ and all $0\le s<t$,
\begin{equation}\label{sub:eq-EVI}
 \tfrac12\sfd^2(u(t),v)-\tfrac12\sfd^2(u(s),v)
 +\int_s^t\Big(\phi(u(r))+\frac{\lambda_c}2\,\sfd^2(u(r),v)\Big)\dd r
 \le(t-s)\,\phi(v).
\end{equation}
\end{definition}

Two sublevel solutions that start in $\calX_c$ stay in $\calX_c$, so each
can be used as observer for the other and they contract with the constant
$\lambda_c$. The next lemma makes this precise: on the metric space
$(\calX_C,\sfd)$ the inequality \eqref{sub:eq-EVI} for one level is the
integrated EVI of Proposition~\ref{prelim:EVIprop}. We use the theory of
\cite{MuS20} on $\calX_C$ only through statements that do not require a
geodesic or a length space, because a sublevel need not be one.

\begin{lemma}[Restriction to a sublevel]\label{sub:restrict}
Let $C,\kappa\in\R$ and let $u:[0,\infty)\to\calX_C$ be a continuous curve
such that for all $v\in\calX_C$ and $0\le s<t$
\begin{equation}\label{sub:eq-base}
 \tfrac12\sfd^2(u(t),v)-\tfrac12\sfd^2(u(s),v)
 +\int_s^t\Big(\phi(u(r))+\frac{\kappa}2\,\sfd^2(u(r),v)\Big)\dd r
 \le(t-s)\,\phi(v).
\end{equation}
Then $u$ is an $\EVI_\kappa$ solution starting at $u(0)$, in the sense of
Definition~\ref{prelim:EVI}, for the functional $\phi|_{\calX_C}$ on the
metric space $(\calX_C,\sfd)$. Consequently the conclusions of
Proposition~\ref{prelim:EVIprop} hold in that space with $\lambda=\kappa$:
the differential inequality \eqref{prelim:eq-EVI} and the exponential
estimate \eqref{prelim:eq-expEVI} for all observers $v\in\calX_C$, the
latter for $0\le s<t$; monotonicity of $t\mapsto\phi(u(t))$ on
$[0,\infty)$; and the contraction \eqref{prelim:eq-contraction}, for
$0\le s<t$, between $u$ and any second continuous curve
$\tilde u:[0,\infty)\to\calX_C$ satisfying \eqref{sub:eq-base}. Moreover
$u$ is locally Lipschitz on $(0,\infty)$ and belongs to
$AC^2([0,T];\calX)$ for every $T>0$, with
$\int_0^T|\dot u|^2(r)\dd r\le\phi(u(0))-\phi(u(T))$.
\end{lemma}

\begin{proof}
On the metric space $\calX_C$, the restricted energy is proper and lower
semicontinuous, and its domain is $\calX_C$. Along a continuous curve its
values are bounded below on compact time intervals (the image is compact)
and above by $C$. Thus \eqref{sub:eq-base} is the integrated inequality
\eqref{prelim:eq-intEVI} in that space, and
Proposition~\ref{prelim:EVIprop} applies; the solution starts at $u(0)$ by
continuity, so the statements for $s=0$ hold, in particular
$\phi(u(s))\le\phi(u(0))$, which is not a consequence of lower
semicontinuity (see the remark after Proposition~\ref{prelim:EVIprop}).
What is not contained in that proposition is the regularity. Local
Lipschitz continuity at positive times and the energy identity are
\cite[Thm.~3.5]{MuS20}. Integrating the latter from $s>0$ to $T$ gives
$\int_s^T|\dot u|^2\le\phi(u(s))-\phi(u(T))
\le\phi(u(0))-\phi(u(T))$; let $s\downarrow0$.
No length-space property of $\calX_C$ is used.
\end{proof}

\subsection{The abstract generation theorem}\label{sub:ss-abstract}

We now replace the $\lambda$-convexity of $\phi$ in
Theorem~\ref{cp:abstract} by (SL). We retain its solvability, invariance
and recovery hypotheses and use only the derivative estimate supplied by
its geometric hypothesis. Let
$(\calX,\sfd)$ be complete and geodesic, and let
$A_1\subset A_2\subset\cdots\subset\dom\slope\phi$ be sets with the
following properties.
\begin{enumerate}[label=\textup{(A\arabic*)}]
\item\label{sub:A1} For every $j$ the geometric derivative estimate
 \eqref{cp:eq-deriv} of Lemma~\ref{cp:deriv} holds along geodesics with endpoints in $A_j$,
 with a finite constant $\kappa_j\ge0$ and with observers in the geodesic
 cover $A_j^{\rm g}$ of $A_j$, for almost every $t\in(0,1)$
 $($Section~\textup{\ref{cp:ss-geom}}$)$.
\item\label{sub:A2} For every $x\in\bigcup_jA_j$ and every $T>0$ there are
 $j$ and $\tau_0>0$ such that for all $\tau\in(0,\tau_0)$ exact minimisers
 of the scheme of Definition~\ref{prelim:MM} starting at $x$ exist, and all
 iterates $x_n$ with $n\le\lceil T/\tau\rceil$ lie in $A_j$, for every
 choice of the minimisers.
\item\label{sub:A3} For every $v\in\dom\phi$ there are
 $v_m\in\bigcup_jA_j$ with $\sfd(v_m,v)\to0$ and $\phi(v_m)\to\phi(v)$.
\end{enumerate}

The full condition $(\mathrm G_{\kappa_j})$ implies \ref{sub:A1} by
Lemma~\ref{cp:deriv}, but is not needed here. Squared distance to a fixed
observer along a geodesic is Lipschitz on $[0,1]$, hence absolutely
continuous. Integrating \ref{sub:A1} therefore gives the corresponding
secant estimate in Lemma~\ref{cp:deriv}, which is the geometric input to
the approximate EVI and the construction below.

\begin{theorem}[Generation under sublevel semiconvexity]\label{sub:abstract}
Let $(\calX,\sfd)$ be a complete geodesic metric space, let
$\phi:\calX\to(-\infty,+\infty]$ be proper, lower semicontinuous and satisfy
\textup{(SL)}, and let $(A_j)_j$ be nested subsets of $\dom\slope\phi$ with
\ref{sub:A1}--\ref{sub:A3}. Then there is a unique family of maps
$S_t:\dom\phi\to\dom\phi$, $t\ge0$, with the following properties.
\begin{enumerate}[label=\textup{(\roman*)}]
\item For every $x\in\dom\phi$ the curve $u(t)=S_tx$ is a sublevel solution
 in the sense of Definition~\ref{sub:def-sol} with $u(0)=x$; that is, $u$
 is continuous, $\phi(S_tx)\le\phi(x)$, and \eqref{sub:eq-EVI} holds with
 the constant $\lambda_c$ for every $c\ge\phi(x)$ and every observer
 $v\in\calX_c$.
\item For every $c\in\R$ and all $x,y\in\calX_c$,
 \begin{equation}\label{sub:eq-contr}
  \sfd(S_tx,S_ty)\le\ee^{-\lambda_ct}\,\sfd(x,y)\qquad(t\ge0).
 \end{equation}
\item $S_{t+s}=S_tS_s$ for $s,t\ge0$; the map $t\mapsto\phi(S_tx)$ is
 nonincreasing; every trajectory is locally Lipschitz on $(0,\infty)$ and
 belongs to $AC^2([0,T];\calX)$ for every $T>0$.
\item If $c\in\R$ and $\tilde u:[0,\infty)\to\calX_c$ is a continuous curve
 which satisfies \eqref{sub:eq-EVI} for this one value of $c$ and all
 $v\in\calX_c$, then $\tilde u(t)=S_t\tilde u(0)$ for all $t\ge0$.
\item For every $x\in\bigcup_jA_j$, all geodesic interpolants of the
 minimising movement scheme starting at $x$ converge to $t\mapsto S_tx$ as
 $\tau\to0$, uniformly on bounded time intervals, for every choice of the
 minimisers and of the interpolating geodesics.
\end{enumerate}
\end{theorem}

As in Section~\ref{cp:ss-abstract}, the points of $\bigcup_jA_j$ are called
regular. The proof has three steps. For regular initial data of energy at
most $c$ we use the construction of Section~\ref{sec:compact} in the
localised form of Remark~\ref{cp:rem-local}. The
interpolants
of the scheme reach the level $g(c)$, and the regular observers that
approximate an observer of energy at most $c$ have energy up to $c+1$, so
we work at the level $g(c+1)$ and obtain the EVI on $\calX_c$ only with the
worse constant
\begin{equation}\label{sub:eq-Lambda}
 \Lambda(c)=2\lambda_{g(c+1)}-2.
\end{equation}
Contraction on sublevels then extends the flow to all data of finite
energy. In the last step we test the differential EVI with an observer
$\gamma_s$ close to $u(t)$ on the geodesic from $u(t)$ to $v$; the term
with $\Lambda$ is of order $s^2$ and disappears as $s\downarrow0$. This
restores $\lambda_c$ without a continuity assumption on $\lambda$ or $g$.

\begin{proof}
\emph{Regular data.} Fix $x_0\in\bigcup_jA_j$, $T>0$, $c\ge\phi(x_0)$, and
put $\bar c=g(c+1)$, $\ell=-\lambda_{\bar c}$ and
$\hat\lambda=-2(\ell+1)=\Lambda(c)$, in the notation of
Section~\ref{cp:ss-error}. The iterates lie in a common $A_j$ by
\ref{sub:A2} and have energy at most $c$; their interpolants have energy
at most $g(c)\le\bar c$. Lemmas~\ref{cp:euler} and~\ref{cp:error},
Proposition~\ref{cp:approxEVI} and Steps~1 and~2 of the proof of
Theorem~\ref{cp:abstract} are stated for a weakly geodesically convex
functional; we use them in the localised form of
Remark~\ref{cp:rem-local}, with $K=\calX_{\bar c}$ and this $\ell$. Its
hypotheses follow from (SL) and $c\le c+1\le g(c+1)=\bar c$:
$\{\phi\le\phi(x_0)\}\subset\calX_c\subset\calX_{\bar c}$, which is
\ref{cp:K1}; any two points of $\calX_{\bar c}$ are joined by a geodesic,
and by \eqref{sub:eq-SL} at the level $\bar c$ every such geodesic is
admissible in the sense of Section~\ref{cp:ss-error} with this $\ell$,
which is \ref{cp:K2}; every geodesic between two points of $\calX_c$ has
energy at most $g(c)\le\bar c$, which is \ref{cp:K3}; and observers are
taken in $\calX_{\bar c}$, as \ref{cp:K4} asks. The level $\bar c$ is
$g(c+1)$ and not $g(c)$ because of the last step below: the regular points
that approximate an observer of energy at most $c$ have energy up to
$c+1$, and they have to lie in $K$.
With this $\ell$ and $\hat\lambda$ we therefore have: the lower bound
\eqref{cp:eq-lowerphi} for $y\in\calX_{\bar c}$, in particular at the
points of the interpolants and of the limit curve; the bounds of
Lemma~\ref{cp:euler}\ref{cp:it-abc} and \eqref{cp:eq-incr} for the slopes
and the increments, which also hold with $\ell=-\lambda_{\phi(x_0)}$, since
they use only consecutive nodes; the energy bound in \eqref{cp:eq-phiu}
along the interpolants; the approximate EVI \eqref{cp:eq-approxEVI} for
observers $v\in A_j^{\rm g}\cap\calX_{\bar c}$, with the error function
\eqref{cp:eq-Er}; and the error estimate \eqref{cp:eq-error} with
$\kappa=\kappa_j$.
Two interpolants can be used as observers for each other, and the
comparison argument of Step~1 of the proof of Theorem~\ref{cp:abstract}
gives \eqref{cp:eq-cauchy}.
Completeness gives a locally uniform, Lipschitz limit $u$ with
$\phi(u(t))\le\phi(x_0)$. For a regular observer in $\calX_{\bar c}$
enlarge $A_j$ if needed
and integrate \eqref{cp:eq-approxEVI}; \eqref{cp:eq-lowerphi} supplies the
common lower bound for Fatou. For an arbitrary $v\in\calX_c$, the recovery
sequence in \ref{sub:A3} eventually has energy at most $c+1\le\bar c$,
hence lies in $K$, and passes this inequality to $v$. Thus $u$ satisfies \eqref{sub:eq-base}
with $\kappa=\Lambda(c)$ on $\calX_c$.

\emph{Finite-energy data.} For $x\in\dom\phi$, take a recovery sequence
$x_m$ and fix $C>\phi(x)$. Eventually $\phi(x_m)\le C$.
Lemma~\ref{sub:restrict} with $\Lambda(C)$ makes their flows uniformly
Cauchy on $[0,T]$. Their limit starts at $x$ and has energy at most
$\phi(x)$. The union of the images of a uniformly convergent sequence of
continuous curves on $[0,T]$ and its limit is compact. Hence lower
semicontinuity supplies a common energy lower bound, and Fatou passes
\eqref{sub:eq-base} to the limit on every $\calX_C$, $C>\phi(x)$.
Contraction makes the limit independent of the recovery sequence.

\emph{Restoring the constant.} Fix $t>0$, an observer $v$ and a level $c$
containing both $u(t)$ and $v$. Take $C$ larger than $\phi(x)$ and $g(c)$,
and put
$\kappa=\Lambda(C)$. Let $\gamma$ join $u(t)$ to $v$ and
$D=\sfd(u(t),v)$. Every $\gamma_s$ belongs to $\calX_C$.
Writing $a_s=\frac12\rmd^+_r\sfd^2(u(r),\gamma_s)|_{r=t}$,
\cite[Lemma~3.13]{MuS20} gives $a_1\le a_s/s$ (this is also the argument of
\cite[Prop.~5.6]{MiSu26}). The EVI with observer
$\gamma_s$, followed by (SL), gives
\[
 a_1\le\frac{a_s}s
 \le\phi(v)-\phi(u(t))-\tfrac{\lambda_c}2(1-s)D^2
                         -\tfrac\kappa2sD^2.
\]
Let $s\downarrow0$. This proves the differential EVI with $\lambda_c$
at every positive time whenever both energies are at most $c$ (the case
$D=0$ follows directly). For $c\ge\phi(x)$ the entire trajectory lies in
$\calX_c$, because $\phi(u(r))\le\phi(x)$ for all $r\ge0$.
On the metric space $(\calX_c,\sfd)$ with the
restricted energy this is the definition of an $\EVI_{\lambda_c}$ solution
\cite[Def.~3.1]{MuS20}, so \cite[Thm.~3.3]{MuS20} gives the integrated
inequality \eqref{sub:eq-EVI} for $0<s<t$, and continuity at $0$ extends it
to $s=0$. Lemma~\ref{sub:restrict} then gives contraction, monotonicity,
regularity and uniqueness. Applying
uniqueness after any time shift gives the semigroup property. Convergence
of the scheme from regular data was proved in the first step.
\end{proof}

\begin{remark}\label{sub:rem-scope}
The flow is defined on $\dom\phi$, with contraction on each sublevel;
continuity along sequences with unbounded energies is not asserted.
A sublevel need not be a length space, so the general finite-energy
estimate of \cite[Thm.~5.9]{MuS20} cannot simply be applied there.
The next theorem uses the estimate for data of finite slope after a short
discrete initial layer. It is conditional on exact steps and a sublevel
lower bound.
\end{remark}

\subsection{Minimising movements from finite-energy data}\label{sub:ss-mm}

Theorem~\ref{sub:abstract} gives convergence of the scheme only from
regular initial data, that is, from points of $\bigcup_jA_j$, and by Remark~\ref{sub:rem-scope} the estimate of
\cite{MuS20} for general data is not available on a sublevel. We use the
estimate for data of finite slope and start it a few steps late. The
scheme from $x$ dissipates at most $B=\phi(x)-\inf_{\calX_c}\phi$, so among
its first $N$ iterates there is one whose slope is at most
$\sqrt{2B/(N\tau)}$, and it is within $\sqrt{2N\tau B}$ of $x$. The choice
$N\approx\tau^{-1/2}$ balances the two errors and gives the rate
$\tau^{1/4}$. The lower bound for the energy on the sublevel and the
existence of exact minimisers are assumed.

The piecewise constant interpolant $\bar U_\tau$ of a sequence $(U^n)$ is
that of Definition~\ref{prelim:MM}.

\begin{theorem}[Convergence of the scheme from finite-energy data]\label{sub:mm}
Let $(\calX,\sfd)$ be a geodesic metric space and let
$\phi:\calX\to(-\infty,+\infty]$ be proper, lower semicontinuous and satisfy
\textup{(SL)}. Let $x\in\dom\phi$, $c:=\phi(x)$, $\ell:=-\lambda_c$, and
assume:
\begin{enumerate}[label=\textup{(\alph*)}]
\item there are continuous maps $S_t:\calX_c\to\calX_c$, $t\ge0$, with
 $S_{t+s}=S_tS_s$, $S_0=\mathrm{id}$, such that every trajectory
 $t\mapsto S_ty$, $y\in\calX_c$, is a continuous curve satisfying
 \eqref{sub:eq-base} with $C=c$ and $\kappa=\lambda_c$ for all
 $v\in\calX_c$;
\item $m:=\inf_{\calX_c}\phi>-\infty$; put $B:=c-m\ge0$;
\item $0<\tau\le1$, $12\,\ell\,\tau\le1$, and $(U^n)_{n\ge0}$ is an exact
 minimising movement scheme in $\calX$ with step $\tau$ starting at
 $U^0=x$:
 \begin{equation}\label{sub:eq-scheme}
  U^n\in\operatorname*{argmin}_{y\in\calX}
  \Big\{\frac1{2\tau}\sfd^2(U^{n-1},y)+\phi(y)\Big\}\qquad(n\ge1).
 \end{equation}
\end{enumerate}
Then for every $T>0$,
\begin{equation}\label{sub:eq-rate}
 \sup_{0\le t\le T}\sfd\big(S_tx,\bar U_\tau(t)\big)\le
 C_{T,\ell}\,\sqrt B\,\tau^{1/4},\qquad
 C_{T,\ell}:=2\big(\ee^{\ell T}+1\big)
 +\sqrt2\,\big(\sqrt{T+1}+1\big)\,\ee^{3\ell(T+1)},
\end{equation}
and every geodesic interpolant $\tilde U_\tau$ of $(U^n)$ satisfies
$\sup_{t\le T}\sfd(S_tx,\tilde U_\tau(t))\le
C_{T,\ell}\sqrt B\,\tau^{1/4}+\sqrt{2\tau B}$.
\end{theorem}

\begin{proof}
Put $\phi_c=\phi|_{\calX_c}$. By Lemma~\ref{sub:restrict} and assumption
(a), $(S_t)$ is an $\EVI_{\lambda_c}$ flow in $\calX_c$, in the sense of
Definition~\ref{prelim:EVI}, for $\phi_c$ on the metric space
$(\calX_c,\sfd)$. We apply Theorem~\ref{prelim:MS-regular} in this space,
with $\lambda=\lambda_c=-\ell$ and $D=\calX_c$: if a sequence $(V^n)$ in
$\calX_c$ of points of finite slope satisfies \eqref{prelim:eq-MSreg} for
$\phi_c$, then \eqref{prelim:eq-MSslope} holds for the flow started at any
$z\in\calX_c$ and the scheme $(V^n)$. That theorem needs neither a
geodesic sublevel nor approximate convexity. It does need the flow as a
semigroup, which is why (a) asks for the semigroup property: in the proof
of \cite[Thm.~5.5]{MuS20} the point $S_{n\tau}V^0$ is compared with
$S_\tau V^{n-1}$. The condition $12\ell\tau\le1$ in (c) is its hypothesis
$-12\lambda\tau\le1$.

For the actual scheme let $d_n=\sfd(U^{n-1},U^n)$ and
$\varsigma_n=\slope{\phi_c}(U^n)$. Minimality and the triangle inequality give
\begin{equation}\label{sub:eq-sum}
 \sum_{n=1}^N d_n^2/\tau\le2B,\qquad
 \varsigma_n\le d_n/\tau\quad(n\ge1).
\end{equation}
Using points on the geodesic from $U^n$ to $U^{n-1}$ as competitors in the
\emph{ambient} minimization and letting the interpolation parameter tend
to zero gives the second inequality of \eqref{prelim:eq-MSreg}; that geodesic
may leave $\calX_c$, but (SL) applies to its endpoints. The first
inequality of \eqref{prelim:eq-MSreg} is the slope bound in
\eqref{sub:eq-sum}.

Choose $N=\lceil\tau^{-1/2}\rceil$ and
$j\in\{1,\ldots,N\}$ with $\varsigma_j^2\le2B/(N\tau)$.
Cauchy--Schwarz and the action bound give, for every $n\ge0$,
$\sfd(U^{n+j},U^n)\le\sqrt{2N\tau B}$.
Apply \eqref{prelim:eq-MSslope} to $V^n=U^{n+j}$ and $z=x$, then compare back
to $U^{\lceil t/\tau\rceil}$. For $t\le T$ this yields
\[
 \sfd(S_tx,\bar U_\tau(t))\le
 (\ee^{\ell T}+1)\sqrt{2N\tau B}
 +(\sqrt{\tau(T+\tau)}+\tau)\ee^{3\ell(T+\tau)}\sqrt{2B/(N\tau)}.
\]
Since $\sqrt\tau\le N\tau\le2\sqrt\tau$, this is the stated constant.
The interpolation error is at most $\max_n d_n\le\sqrt{2\tau B}$.
\end{proof}

\subsection{Sums of powers}\label{sub:ss-compact}

In the rest of this section,
\begin{equation}\label{sub:eq-E}
 E(c)=\sum_{i=1}^k\alpha_i\,c^{p_i},\qquad k\ge1,\quad\alpha_i>0,\quad p_i>1,\qquad
 p:=\max_ip_i ,
\end{equation}
and $\sfE$ is the functional of Definition~\ref{prelim:energy}. Since
$\Einf=\infty$, $\sfE(\mu)=\int E(\rho)\dd x$ if $\mu=\rho\Leb$ and
$\sfE(\mu)=+\infty$ otherwise; $\sfE\ge0$; and for a probability density
$\rho$ on a closed convex set, $\sfE(\rho\Leb)<\infty$ if and only if
$\rho\in L^p$, because $c^{q}\le c+c^p$ for $1<q\le p$. We put
\begin{equation}\label{sub:eq-lam}
 \lambda_c:=-p\,2^{p+1}\,c,\qquad g(c):=2^p\,c\qquad(c\ge0).
\end{equation}
We shall use three elementary properties of $E$. First, $E(0)=0$, $E$ is
convex and $E'\ge0$ is nondecreasing, so $E(a)\le aE'(a)$ and
\begin{equation}\label{sub:eq-Mbound}
 \sfE(\rho\Leb)\le\int\rho\,E'(\rho)\dd x\le E'(M)
 \qquad\text{for a probability density }\rho\le M .
\end{equation}
Secondly, for $a\in[0,\frac12]$ the convexity of $a\mapsto(1-a)^{-p}$ gives
$(1-a)^{-p}\le1+2(2^p-1)a\le1+2^{p+1}a$. Hence, if $\nu\in\calP(\Om)$ has
finite energy, $F\subset\Om$ is a Borel set with $a:=\nu(F)\le\frac12$, and
$\nu':=\nu|_{\Om\setminus F}/(1-a)$, then
\begin{equation}\label{sub:eq-cut}
 \sfE(\nu')\le(1-a)^{-p}\sfE(\nu)\le\big(1+2(2^p-1)a\big)\,\sfE(\nu)
 \le2^p\,\sfE(\nu).
\end{equation}
In particular the energy increase is at most $2^{p+1}a\sfE(\nu)$.
We also use the scalar distance bound \eqref{shk:eq-G}.

The theorem is stated for bounded and unbounded domains together. The
conclusions (i)--(iv) are the same in both cases, except that on an
unbounded domain the step of the scheme is restricted in terms of the
energy of the initial measure; (v) concerns unbounded domains only. The
proof for bounded $\Om$ follows the statement. The proof for unbounded
$\Om$ uses the bounded case on the balls $\OmR$ and is given in
Section~\ref{sub:ss-unbdd}.

\begin{theorem}[Sums of powers]\label{sub:powers-compact}\label{sub:powers-unbdd}
Let $\Om\subset\R^d$ be closed and convex with nonempty interior, bounded or
unbounded, let $E$ be as in \eqref{sub:eq-E} and
$\lambda_c=-p\,2^{p+1}c$. Then there is a unique family of maps $S_t$,
$t\ge0$, on
$\dom\sfE=\{\rho\Leb\in\calP(\Om):\rho\in L^p(\Om)\}$ with the following
properties.
\begin{enumerate}[label=\textup{(\roman*)}]
\item For every $\nu_0\in\dom\sfE$ the curve $u(t)=S_t\nu_0$ is continuous
 with $u(0)=\nu_0$, $t\mapsto\sfE(u(t))$ is nonincreasing, and for every
 $c\ge\sfE(\nu_0)$, every $w\in\calP(\Om)$ with $\sfE(w)\le c$ and all
 $0\le s<t$,
 \begin{equation}\label{sub:eq-EVIpow}
  \tfrac12\SHK^2(u(t),w)-\tfrac12\SHK^2(u(s),w)
  +\int_s^t\Big(\sfE(u(r))+\frac{\lambda_c}2\,\SHK^2(u(r),w)\Big)\dd r
  \le(t-s)\,\sfE(w).
 \end{equation}
\item $\SHK(S_t\nu_0,S_t\nu_1)\le\ee^{-\lambda_ct}\,\SHK(\nu_0,\nu_1)$
 whenever $\sfE(\nu_0),\sfE(\nu_1)\le c$.
\item $S_{t+s}=S_tS_s$; trajectories are locally Lipschitz on $(0,\infty)$
 and of class $AC^2$ on every $[0,T]$ with respect to $\SHK$; a continuous curve in
 $\{\sfE\le c\}$ that satisfies \eqref{sub:eq-EVIpow} for this $c$ and all
 $w$ with $\sfE(w)\le c$ is a trajectory of $(S_t)$.
\item Let $\nu_0\in\dom\sfE$, $c:=\sfE(\nu_0)$, $\ell:=p\,2^{p+1}c$, and let
 $\tau\in(0,1]$ satisfy $12\ell\tau\le1$ and, if $\Om$ is unbounded, also
 $2^{p+6}\tau c\le1$. Then exact minimising movement schemes $(U^n)$ on
 $\calP(\Om)$ with step $\tau$ starting at $\nu_0$ exist, and each of them
 satisfies, for every $T>0$,
 \[
  \sup_{0\le t\le T}\SHK\big(S_t\nu_0,\bar U_\tau(t)\big)\le
  C_{T,\ell}\,\sqrt c\;\tau^{1/4},
 \]
 with $C_{T,\ell}$ from \eqref{sub:eq-rate}; the geodesic interpolants
 converge likewise.
\item Let $\Om$ be unbounded. If $\nu_0=\rho_0\Leb$ with
 $\rho_0\in L^\infty$ and $\supp\nu_0\subset\bar B_{R_0}$, and $S^R$ denotes
 the flow on $\calP(\OmR)$ given by this theorem for the compact domain
 $\OmR$, then for every $T>0$ there are constants $C_2$ and $R_1$,
 depending on $T$, $p$, $R_0$ and $E'(\|\rho_0\|_\infty)$, such that
 \[
  \sup_{0\le t\le T}\SHK\big(S^R_t\nu_0,S_t\nu_0\big)\le C_2\,R^{-1/6}
  \qquad(R\ge R_1).
 \]
\end{enumerate}
\end{theorem}

\begin{proof}[Proof for bounded $\Om$]
Write $X=\Om$. Properness and lower semicontinuity hold, and compactness of
$\calP(X)$ gives exact steps for every $\tau$ (Lemma~\ref{cp:exist}).
Theorem~\ref{tr:powers} is (SL) with
\eqref{sub:eq-lam}; put $\lambda_c=0$, $g(c)=c$ on empty levels $c<0$.
Take $A_j=\{\rho\Leb\in\calP(X):1/j\le\rho\le j,
\slope{\sfE}(\rho\Leb)\le j\}$.
Proposition~\ref{cp:semiconc} and Lemma~\ref{cp:deriv} prove \ref{sub:A1}.
The exact maximum principle (Proposition~\ref{cp:density-shk}) preserves
$[1/j,j]$. For a scheme started at $\mu_0\in A_j$, \eqref{cp:eq-incr}
with $\ell=-\lambda_{\sfE(\mu_0)}$, which holds because it uses only
geodesics between consecutive nodes (see the proof of
Theorem~\ref{sub:abstract}), bounds all slopes
up to $T$ by $\ee^{2\ell(T+1)}j$ for small steps. This proves \ref{sub:A2}.
The energy recovery of Lemma~\ref{cp:recovery} uses only ordinary convexity
of $E$, density bounds and one step of the scheme; it proves \ref{sub:A3}
without global SHK convexity. Theorem~\ref{sub:abstract} now gives the flow
and (i)--(iii).
Finally $\sfE\ge0$, so Theorem~\ref{sub:mm} applies with $B\le c$ and gives
(iv).
\end{proof}

\begin{remark}[The formal equation]\label{sub:rem-pde}
Formally, the flow of Theorem~\ref{sub:powers-compact} is the equation
\[
 \partial_t\rho=\operatorname{div}\big(\rho\,\nabla E'(\rho)\big)
 -4\rho\,\Big(E'(\rho)-\int_\Om\rho\,E'(\rho)\dd x\Big)
\]
with no-flux boundary conditions, a porous-medium equation with a nonlocal
reaction term that preserves the total mass, in any dimension $d$.
Theorem~\ref{sub:powers-compact} is a statement about the metric flow. We
do not prove that the trajectories are weak solutions of this equation.
\end{remark}

\subsection{Sums of powers on unbounded domains}\label{sub:ss-unbdd}

On an unbounded domain we follow Section~\ref{sec:shk-unbdd}, where
tightness and the comparison of the flows on different balls rest on the
bound for $E'$, which a power energy does not have. Two observations
replace it. The exact maximum principle keeps the density of every iterate
below the bound $M$ of the initial density, so $E'(M)$ plays the role of
the bound (Lemma~\ref{sub:growth}); this gives the flow for bounded,
compactly supported data. Contraction on sublevels, which does not see the
density bound, then extends the flow to all data of finite energy.

Let now $\Om$ be closed and convex with nonempty interior, and unbounded.
The ambient probability space is complete and geodesic by
Proposition~\ref{prelim:cone}. For the normalised restriction $\tilde P_R$
of Lemma~\ref{prelim:trunc} and a tail mass
$\eps=\nu(\Om\setminus\bar B_R)\le1/2$ we have the distance bound
\eqref{prelim:eq-truncSHK}, and \eqref{sub:eq-cut} with
$F=\Om\setminus\bar B_R$ bounds the energy of $\tilde P_R\nu$.

A minimising sequence for one step could lose mass at infinity. The next
theorem excludes this for small steps: unless the far mass is already
bounded by a multiple of the tail of $\nu_0$, cutting it off and
renormalising gains more in distance than it costs in energy. The
conditional gain in distance is Lemma~\ref{shk:farcut}. Since \eqref{sub:eq-cut} bounds the
cost in terms of the energy itself, the step restriction
\eqref{sub:eq-tau} involves the energy $c$ of $\nu_0$.

\begin{theorem}[Exact minimisers on unbounded domains]\label{sub:resolvent}
Let $\Om\subset\R^d$ be closed and convex with nonempty interior, let $E$ be
as in \eqref{sub:eq-E}, and let $\nu_0\in\calP(\Om)$ with
$c:=\sfE(\nu_0)<\infty$. Then $c>0$, and for every $\tau>0$ with
\begin{equation}\label{sub:eq-tau}
 2^{p+6}\,\tau\,c\le1
\end{equation}
the functional
\[
 J(\mu):=\frac1{2\tau}\SHK^2(\nu_0,\mu)+\sfE(\mu)
\]
attains its minimum on $\calP(\Om)$. Every minimiser $\mu_1$ satisfies
$\sfE(\mu_1)\le c$; consequently, for every $\tau$ with \eqref{sub:eq-tau}
the minimising movement scheme with step $\tau$ starting at $\nu_0$ exists
for all steps, for every choice of the minimisers.
\end{theorem}

\begin{proof}
Since a probability density cannot vanish a.e., $c>0$.
Take a minimizing sequence $\mu_j$ with
$\delta_j=J(\mu_j)-\inf J\to0$, $\delta_j\le c$.
Then $\sfE(\mu_j)\le2c$, $\theta_j^2=\SHK^2(\nu_0,\mu_j)\le4\tau c$,
and $I_j=\cos\theta_j\ge1-2\tau c\ge7/8$.
For $F=\Om\setminus B_{\rho+\pi/2}$ and
$\eps=\nu_0(\Om\setminus B_\rho)\le1/256$, put $a=\mu_j(F)$.
If $a=1$, Lemma~\ref{prelim:split} would give
$\HK^2(\nu_0,\mu_j)\ge1-\eps$, contradicting $\HK^2\le4\tau c\le1/4$.
Let $\eps\le a<1$ and $\nu'=\mu_j|_{F^c}/(1-a)$. By \eqref{shk:eq-farcut},
\[
 \frac{I_j-\sqrt{a\eps}}{\sqrt{1-a}}\le\cos\SHK(\nu_0,\nu')\le1.
\]
The numerator is at least $13/16$, so
$a\le87/256<1/2$. Only now apply \eqref{sub:eq-cut}: the energy increase
from $\mu_j$ to $\nu'$ is at most $\beta a$, where $\beta=2^{p+2}c$ and
$\tau\beta\le1/16$.
By Lemma~\ref{shk:farcut}, either $aI_j<2\sqrt{a\eps}$, and then
$a<4\eps/I_j^2<16\eps$, or \eqref{shk:eq-fargain} holds. In the latter
case $J(\nu')\ge J(\mu_j)-\delta_j$ gives
\[
 \kappa a-\sqrt{a\eps}\le\tau\delta_j,
 \quad\kappa=I_j/2-\tau\beta\ge3/8,
 \quad a\le16\eps+\tfrac{16}3\tau\delta_j.
\]
This last bound covers the other cases too. Hence the sequence is tight,
and Prokhorov and lower semicontinuity give a minimizer. Its energy is at
most $c$, so the same step restriction permits iteration.
\end{proof}

The step restriction is sufficient; it is not claimed to be sharp.
For tightness of flows we also need a growth bound for bounded densities.
It is Lemma~\ref{shk:growth} with $E'(M)$ in the role of the oscillation of
$E'$.

\begin{lemma}[Growth bound for bounded densities]\label{sub:growth}
Let $X\subset\R^d$ be compact and convex with nonempty interior, let $E$ be
as in \eqref{sub:eq-E}, let $0<m\le M<\infty$ and
\[
 \omega:=E'(M)
\]
(this number takes the place of the oscillation of $E'$ in
Lemma~\ref{shk:growth}; for powers that oscillation is infinite),
and let $\tau>0$ with $3\tau\omega<1$. Let $\mu_0=\rho_0\Leb\in\calP(X)$ with
$m\le\rho_0\le M$ a.e., let $\mu_1$ be a minimiser of
$J(\nu)=\frac1{2\tau}\SHK^2(\mu_0,\nu)+\sfE(\nu)$ on $\calP(X)$, and let
$(\bfT,q)$ be an optimal transport--growth pair from $\mu_0$ to $\mu_1$
$($Lemma~\textup{\ref{prelim:tg}}$)$, with transport cosine $\chi$. Then
$\mu_1=\rho_1\Leb$ with $m\le\rho_1\le M$, and
\begin{equation}\label{sub:eq-growth}
 q\le\frac{\chi}{\cos\SHK(\mu_0,\mu_1)-2\tau\omega}\le\frac1{1-3\tau\omega}
 \qquad\mu_0\text{-a.e.}
\end{equation}
In particular $|\bfT(x)-x|<\frac\pi2$ for $\mu_0$-a.e.\ $x$ with $q(x)>0$.
\end{lemma}

\begin{proof}
The density bounds are Proposition~\ref{cp:density-shk}. The rest is
Lemma~\ref{shk:growth}, once \ref{shk:G1} and \ref{shk:G2} are verified
with $\omega=E'(M)$. For \ref{shk:G1}, $\sfE(\mu_0)\le E'(M)$ by
\eqref{sub:eq-Mbound} and $\sfE(\mu_1)\ge0$. For \ref{shk:G2}, let $A$,
$n=n_A$ and $\nu^t=\nu_A^t$ be as in \eqref{shk:eq-comp}, and write
$M_t=1-(1-t^2)n$ and $\nu^t=\rho^t\Leb$. Since $\beta_A\le\mu_1$,
\[
 \rho^t\le\rho_1/M_t\le M/M_t,\qquad
 \int_X(\rho^t-\rho_1)^+\dd x\le(M_t^{-1}-1)\int_X\rho_1\dd x
 =(1-t^2)n/M_t .
\]
Convexity gives $E(\rho^t)-E(\rho_1)\le E'(\rho^t)(\rho^t-\rho_1)$, and
since $E'\ge0$ is nondecreasing,
\[
 \sfE(\nu^t)-\sfE(\mu_1)
 \le E'(M/M_t)\int_X(\rho^t-\rho_1)^+\dd x
 \le E'(M/M_t)(1-t^2)n/M_t.
\]
As $t\uparrow1$, $M_t\to1$ and $E'(M/M_t)\to\omega$ by continuity of $E'$,
so the upper limit in \ref{shk:G2} is at most $2\omega n$. For $t<1$ the
factor $E'(M/M_t)$ is larger than $\omega$, so this is not the bound with
the fixed constant $\omega$ that Step~4 of the proof of
Lemma~\ref{shk:growth} gives for bounded $E'$; the constant is reached
only in the limit, and \ref{shk:G2} is formulated as an upper limit for
this reason.
\end{proof}

\begin{proof}[Proof of Theorem~\ref{sub:powers-unbdd} for unbounded $\Om$]
For bounded, compactly supported data we argue as in
Section~\ref{sec:shk-unbdd}, but the exhaustion argument of
Sections~\ref{sec:hk-unbdd} and~\ref{sec:shk-unbdd} cannot be quoted as it
stands. There the flows $S^R$ satisfy the $\EVI$ with one constant, and
truncation raises the energy by at most a multiple of the tail mass. Here
the constant in the inequality for $S^R$ depends on the energy level of
the solution and of the observer, and by \eqref{sub:eq-cut} the normalised
truncation multiplies the energy by a factor up to $(1-\eps)^{-p}$, so a
truncated observer lies on a higher level. The levels therefore have to be
carried along, and the comparison is made at a
common level $c_*$. General data are then reached by contraction. The
flows $S^R$ on $\calP(\OmR)$ are those of the bounded case.

\emph{Bounded, compactly supported data.} Suppose $\rho_0\le M$ and
$\supp\nu_0\subset\bar B_{R_0}$; put $\omega=E'(M)$ and, as in
Proposition~\ref{shk:tight}, $K_{T,\omega}=\ee^{7\omega T}$. Consider
first a scheme on $\calP(\OmR)$ from a regular measure $\mu_0$ with density
at most $M$. Lemma~\ref{sub:growth} gives the growth bound and keeps the
density bounds, and by \eqref{sub:eq-Mbound}
\[
 \sum_j\SHK^2(\mu_{j-1},\mu_j)\le2\tau\sfE(\mu_0)\le2\tau\omega.
\]
These are precisely the inputs to the cutoff iteration in
Proposition~\ref{shk:tight}, so \eqref{shk:eq-mt} holds for such schemes,
and part (iv) on $\OmR$ passes it to the trajectory $S^R_\cdot\mu_0$ $($for
$0<t\le T$ use the steps $\tau=t/m$, $m\to\infty$, for which $t$ is a
node$)$. For
$\nu_0$ itself, Lemma~\ref{cp:recovery} on $\OmR$ gives regular measures
converging to $\nu_0$ with densities at most $M$. Their energies are at
most $\omega$, so contraction at level $\omega$ passes the cutoff bound to
$S^R_\cdot\nu_0$. Thus
\begin{equation}\label{sub:eq-tight}
 \sup_{t\le T}(S^R_t\nu_0)(\Om\setminus B_\rho)
 \le C_1/\rho,\qquad C_1=\pi K_{T,\omega}^2\sqrt{T\omega},\quad\rho\ge2R_0,
\end{equation}
for every $R\ge R_0+1$. For such $R$ the set $\OmR$ has nonempty interior:
$\nu_0$ is absolutely continuous with support in $\bar B_{R_0}$, so $\Om$
has interior points in $\bar B_{R_0+1}$.

Put $R_1:=\max\{R_0+1,2R_0,2C_1\}$ and let $R'\ge R\ge R_1$. Then
$\nu_0\in\calP(\OmR)\subset\calP(\Om_{R'})$, $R\ge1$, and by
\eqref{sub:eq-tight} with $\rho=R$ the mass of $S^{R'}_s\nu_0$ outside
$B_R$ is at most $\eps:=C_1/R\le1/2$ for $s\le T$.
Let $u=S^R_\cdot\nu_0$ and $v=S^{R'}_\cdot\nu_0$. Their energies lie in $[0,\omega]$, while
\eqref{sub:eq-cut} gives
\[
 \sfE(\tilde P_Rv)\le\sfE(v)+2^{p+1}\omega\eps,
 \qquad \sfE(\tilde P_Rv)\le2^p\omega=:c_*.
\]
Both curves therefore satisfy the exponential EVI with the common
constant $\lambda_{c_*}$ for all observers used in
Lemma~\ref{hk:comparison}. Apply that lemma with
$\delta=\pi/2$, $K=2^{p+1}\omega$, $D=\omega$, and
$N=\lceil R^{1/3}\rceil$. It yields
\begin{equation}\label{sub:eq-cauchy}
 \sup_{t\le T}\SHK(S^R_t\nu_0,S^{R'}_t\nu_0)\le C_2R^{-1/6}.
\end{equation}
Completeness gives the locally uniform limit $u$. For compactly supported
observers, pass the integrated EVI by Fatou, using $\sfE\ge0$.
For $w$ with $\sfE(w)\le c$, use $w_R=\tilde P_Rw$ at levels
$c_R=(1-\eps_R)^{-p}c\downarrow c$. Then $w_R\to w$,
$\sfE(w_R)\to\sfE(w)$ and $\lambda_{c_R}\to\lambda_c$, proving the
EVI for all observers.

\emph{General data.} By normalized truncation and cell averaging choose
bounded, compactly supported $\nu_0^j\to\nu_0$ with
$e_j=\sfE(\nu_0^j)\to\sfE(\nu_0)$. Contraction at a fixed larger energy
level makes their flows uniformly Cauchy. To pass the EVI with observer
energy at most $c$, use levels $\max\{c,e_j\}\to c$; Fatou and continuity
of $\lambda_c$ give the asserted constant in the limit.
Lemma~\ref{sub:restrict} gives energy monotonicity, regularity and
contraction. The constant-recovery argument in Theorem~\ref{sub:abstract},
using (SL) on $(\calP(\Om),\SHK)$, which is Theorem~\ref{tr:powers}, gives
the differential EVI at any level
containing the current point and observer. Consequently time shifts satisfy
the same sublevel definition, and uniqueness gives the semigroup property.
This proves (i)--(iii), and (v) follows from \eqref{sub:eq-cauchy} by
letting $R'\to\infty$.

Finally, Theorem~\ref{sub:resolvent} supplies every exact step, and
Theorem~\ref{sub:mm} applies with (SL), $m\ge0$ and $B\le c$. This is (iv).
\end{proof}

\begin{remark}\label{sub:rem-unbdd}
The tightness constant uses a density bound only for the preliminary
construction. The final flow and scheme convergence require finite energy,
without a moment bound. For arbitrary initial probabilities see
Corollary~\ref{sub:powers-measure-data}. Identification with weak PDE
solutions is not treated.
\end{remark}

\subsection{Extension to measure data}\label{sub:ss-extension}

Contraction on a fixed finite sublevel gives no uniform estimate for
approximations whose initial energies tend to infinity. Test the sublevel
EVI with a finite-energy observer $w$. While the solution stays close to
$w$, the linear bound on $-\lambda_c$ allows half of the energy above the
lower bound $m$ of $\phi$ to be absorbed. The resulting integral inequality
bounds the accumulated energy uniformly on a short interval for nearby
initial data. Monotonicity also gives the positive-time bound
$\phi(S_tx)-m\le2M+\sfd^2(x,w)/t$. The integral bound controls the product
of the sublevel contraction factors, giving the local Lipschitz estimate
needed to extend the flow to $\overline{\dom\phi}$. The sums of powers
satisfy the required linear bound by \eqref{sub:eq-lam}.

\begin{theorem}[Extension under a linear bound on the sublevel constants]
\label{sub:linear-extension}
Let $(\calX,\sfd)$ be a complete metric space and let
$\phi:\calX\to(-\infty,+\infty]$ be proper, lower semicontinuous and
bounded below by $m\in\R$. Let $c\mapsto\lambda_c\in(-\infty,0]$ be
nonincreasing, with $-\lambda_c\le A(c-m)+B$ for $c\ge m$, where $A>0$
and $B\ge0$. Let $(S_t)_{t\ge0}$ be a semigroup on $\dom\phi$, with
$S_0=\mathrm{id}$. For each $x\in\dom\phi$ assume that $u(t)=S_tx$ is
continuous, $\phi(u(t))\le\phi(x)$, and \eqref{sub:eq-EVI} holds for
every $c\ge\phi(x)$, $v\in\calX_c$ and $0\le s<t$.

Then $(S_t)$ extends uniquely by continuity to a semigroup $(\bar S_t)$
on $D:=\overline{\dom\phi}$. For each $x\in D$ its trajectory is the unique
continuous curve $u$ with $u(0)=x$ such that $u(s)\in\dom\phi$ and
$u(s+h)=S_hu(s)$ for $s>0$, $h\ge0$. Each positive-time shift satisfies
the sublevel EVI. For $w\in\dom\phi$ put
$M=\tfrac32(\phi(w)-m)+\tfrac{B}{2A}$. Then
\begin{equation}\label{sub:eq-extension-energy}
 \phi(\bar S_t x)-m\le 2M+\frac{\sfd^2(x,w)}{t}
 \qquad\bigl(x\in D,\ \sfd^2(x,w)<1/A,\ t>0\bigr).
\end{equation}
For $\sfd^2(x,w)\le1/(2A)$ the right-hand side is at most
$(2At)^{-1}+2M$, hence at most $6M+4$ once $t\ge1/(8A(M+1))$.
On $B_w:=\{x\in D:\sfd^2(x,w)\le1/(2A)\}$ a Lipschitz constant for
$\bar S_t$ is $K_w(t)=\exp(7/4-\lambda_{m+6M+4}(t-\tau)_+)$,
where $\tau=1/(8A(M+1))$.
\end{theorem}

\begin{proof}
First let $x\in\dom\phi$. Write $e_x(r)=\phi(S_rx)-m$,
$e_w=\phi(w)-m$ and $f(r)=\sfd^2(S_rx,w)/2$. By the semigroup law
$e_x$ is nonincreasing. On any interval $[a,b]$ where $f\le1/(2A)$,
partition the interval and apply \eqref{sub:eq-EVI} on each piece with
observer $w$, restarting at its left endpoint $r_i$ at level
$c(r_i)=m+\max\{e_x(r_i),e_w\}$. The left sums for $-\lambda_{c(r)}$
converge to the integral by monotonicity. Summing and using
$-\lambda_{c(r)}\le A(e_x(r)+e_w)+B$ gives
\[
 f(b)-f(a)+\tfrac12\int_a^b e_x(r)\dd r\le M(b-a).
\]
Suppose $f(0)<1/(2A)$ and fix $t>0$. If $e_x(t)>2M$, then
$e_x(r)\ge e_x(t)>2M$ on $[0,t]$. At a first hitting time
$b\le t$ of $1/(2A)$ the integral estimate would give
$f(b)<f(0)$. Hence $2f(t)+t(e_x(t)-2M)\le2f(0)$, which implies
\eqref{sub:eq-extension-energy}. The low-energy case is immediate.

For $x\in B_w\cap\dom\phi$, a first hit at $b\le\tau$ would give
$f(b)\le f(0)+M\tau<3/(8A)$. Thus $\int_0^\tau e_x(r)\dd r\le3/(4A)$
and $e_x(\tau)\le6M+4$. For $x,y\in B_w\cap\dom\phi$, apply
Lemma~\ref{sub:restrict} on each partition interval at its initial
common level $m+\max\{e_x,e_y\}$. Up to $\tau$ the limiting logarithm
of the product of the contraction factors is bounded by
$A\int_0^\tau(e_x+e_y)\dd r+B\tau\le7/4$. Thereafter both energies
are at most $m+6M+4$. This gives the stated $K_w$.

The open balls $\{x\in D:\sfd^2(x,w)<1/(2A)\}$ cover $D$. On each
of them the estimate gives uniformly Cauchy trajectories on $[0,T]$
for every convergent finite-energy approximation. The common limit
defines $\bar S_t$. Spatial continuity gives its uniqueness and the
semigroup law. Lower semicontinuity passes \eqref{sub:eq-extension-energy}
to the limit for $\sfd^2(x,w)<1/A$, so $\bar S_sx$ has finite energy
for $s>0$ and the restart identity follows. A second such curve $v$
from $x$ has $\sfd(v(t),\bar S_tx)\le K_w(t)\sfd(v(s),\bar S_sx)\to0$
as $s\downarrow0$, choosing an open ball containing $x$.

It remains to include the boundary of $B_w$. The stronger distance
inequality above also passes to the limit whenever
$e(t):=\phi(\bar S_tx)-m>2M$, since the approximating terminal
energies eventually exceed this threshold. If $\phi(x)=+\infty$,
lower semicontinuity gives $e(s)\to+\infty$ as $s\downarrow0$.
Thus $\bar S_sx$ lies strictly closer to $w$ for all sufficiently
small $s>0$. Approximate each infinite-energy endpoint by these
finite-energy points, keep each finite-energy endpoint fixed, and
pass the bound with $K_w(t)$ to the limit by the semigroup law.
\end{proof}

\begin{corollary}[Sums of powers from arbitrary probabilities]
\label{sub:powers-measure-data}
Let $\Om\subset\R^d$ be closed and convex with nonempty interior, bounded
or unbounded. Let $E$ be given by \eqref{sub:eq-E}, put $p=\max_i p_i$
and $A=p\,2^{p+1}$. The flow of Theorem~\ref{sub:powers-compact}
extends uniquely by continuity to $(\calP(\Om),\SHK)$. For each initial
probability $\nu_0$, the measure $u(t)=\bar S_t\nu_0$ has an $L^p(\Om)$
density at every $t>0$. For $w\in\dom\sfE$ with $\SHK^2(\nu_0,w)<1/A$,
\[
 \sfE(u(t))\le3\sfE(w)+\frac{\SHK^2(\nu_0,w)}{t}\qquad(t>0).
\]
For $s>0$ the shift $h\mapsto u(s+h)$ satisfies \eqref{sub:eq-EVIpow}
for every $c\ge\sfE(u(s))$ and every observer of energy at most $c$.
\end{corollary}
\begin{proof}
Theorem~\ref{sub:powers-compact} supplies the semigroup and
\eqref{sub:eq-lam} gives the linear bound with $m=B=0$. Hence
Theorem~\ref{sub:linear-extension} applies with $M=\tfrac32\sfE(w)$, and
it remains to show that $\dom\sfE$ is dense in $\calP(\Om)$. Finite sums
of point masses are narrowly dense in $\calP(\Om)$. Replace each point
mass at $z$ by normalised Lebesgue measure on $B_\eps(z)\cap\Om$,
which has positive volume by convexity and nonempty interior. These
bounded, compactly supported densities have finite energy and converge
narrowly to the atomic probability as $\eps\downarrow0$.
Proposition~\ref{prelim:cone} gives completeness and identifies narrow
convergence with $\SHK$ convergence.
\end{proof}

\section{Open problems}\label{sec:outlook}

We collect the questions that the present work leaves open, and say for
each where the methods of this paper stop. After the results of
Sections~\ref{sec:transfer}--\ref{sec:sublevel} were complete, answers or
partial answers to items~1 to~5 were obtained in the same way as
the results of this paper. They are not part of it and will be presented
elsewhere; we say at each of these items what is concerned.
Appendix~\ref{sec:verification} does not cover them.

\begin{enumerate}[leftmargin=2em]
\item \emph{Characterisation on the sphere.} Theorem~\ref{tr:T1} gives a
sufficient condition for $\lambda$-convexity on $(\calP(\Om),\SHK)$ and
Theorem~\ref{tr:necessary} a necessary one on large domains; for energies
with bounded derivative they differ only through the constant $\Lambda$ on
the cone. Is every internal energy that is $\lambda$-convex on the sphere
also semiconvex on the cone? The necessary condition is obtained from
geodesics without transport, between measures with supports more than
$\pi/2$ apart. Along these the energy depends on $E$ only through its
values at two densities, which yields the bound on $E'$ but says nothing
about the condition of \cite{LMS23}; we have no family of test geodesics
with transport that would detect that condition for a general $E$. For
domains of diameter at most $\pi/2$ the supports of two measures are never
that far apart, so neither the necessary condition nor the rigidity theorem
applies, and the methods of this paper give nothing. \emph{Added.} A characterisation of
semiconvexity on the sphere in terms of semiconvexity on the cone, without
a restriction on the diameter, is among the results to be presented
elsewhere.

\item \emph{The scalar condition in dimension $d\ge3$.} Condition
\eqref{ex:eq-Cd} is exact for convex $E$ with $G''>0$, but we know its
solutions only through the criterion of Theorem~\ref{ex:Kd}, the Bernstein
class of Theorem~\ref{ex:bernstein} and the shifted powers of
Theorem~\ref{ex:shifted}. A description of all convex $G$ satisfying it,
comparable to Theorem~\ref{ex:d2} for $d=2$, is missing: for $d\ne2$ the
condition contains a quadratic term in addition to convexity of $G$ and a
bound on $G'$, and we have no way of solving it in general. A concrete
instance is the optimal constant for the shifted powers $1-(1+c)^q$.
Theorem~\ref{ex:shifted} gives a constant for every exponent above the
critical one, but Proposition~\ref{ex:counter} shows that for $d\ge3$ the
natural constant $-2q$ fails near the critical exponent, and the best
constant there is not known. \emph{Added.} Results on the exponents and
dimensions for which $-2q$ is the sharp constant are among those to be
presented elsewhere; an explicit description of all solutions of
\eqref{ex:eq-Cd}, comparable to Theorem~\ref{ex:d2}, remains open.

\item \emph{Energies with $E'(0^+)=-\infty$ on unbounded domains.}
Theorem~\ref{hk:main} needs $E'(0^+)$ finite. The hypothesis enters in the
growth bound for one step (Lemma~\ref{hk:growth}), where a step that
creates too much mass is compared with one that creates less, and removing
mass must not cost too much energy. For $-\int\rho^q$ the difficulty is
more basic: the energy itself is unbounded below on finite measures of an
unbounded domain, so a theory with moment or weight conditions would be
needed. \emph{Added.} Flows of $-\int\rho^q$ with a confining potential
are among the results to be presented elsewhere.

\item \emph{Sublevel flows.} For the sums of powers in \eqref{sub:eq-E},
Corollary~\ref{sub:powers-measure-data} gives a flow from every initial
probability, with an $L^p$ density at positive times and a local
$O(t^{-1})$ energy bound. Do these trajectories solve the no-flux
reaction--diffusion equation of Remark~\ref{sub:rem-pde} in the weak
sense? On an unbounded convex domain without a global Sobolev
inequality, does the short-time energy bound $O(t^{-\alpha})$, with
$\alpha=d(p-1)/(d(p-1)+2)$, hold uniformly over initial probabilities?
For $E(c)=c^p$, the transport part of the equation of
Remark~\ref{sub:rem-pde} is $\partial_t\rho=(p-1)\Delta\rho^p$. Its
mass-preserving scaling on $\mathbb R^d$ gives the exponent $\alpha$ for
$\int\rho^p$.
One would also like quantitative bounds with the same uniformity for the
density in $L^\infty$ and the pressure $\rho E'(\rho)-E(\rho)$ in
$L^2_{\rm loc}((0,\infty);H^1(\operatorname{int}\Om))$. \emph{Added.} The
identification of these trajectories as weak solutions of the no-flux
equation is among the results to be presented elsewhere.

\item \emph{Reference measures other than Lebesgue.} All density bounds are
with respect to $\Leb$. On an unbounded domain no finite measure has a
density bounded below, which is why we had to argue through bounded
subdomains. Relative bounds with respect to a log-concave reference measure
would allow confining potentials and would be the natural setting on
unbounded domains. The one-step density bounds of
Section~\ref{sec:compact} and the semiconcavity estimates of \cite{LaM19}
are formulated for Lebesgue densities, and we have not examined whether
they survive a change of reference measure. \emph{Added.} A pointwise
criterion for semiconvexity with respect to a weighted reference measure
is among the results to be presented elsewhere; a general existence theory
with such a reference measure is not available.

\item \emph{Other geometries.} The transfer principle of
Section~\ref{sec:transfer} uses only the cone structure. It applies verbatim
to any metric cone, along geodesics of the base of length less than $\pi$
(Remark~\ref{tr:rem-scope}), for instance to cones over Wasserstein-type spaces with
other reaction terms. The rest of the theory does not transfer in this
way: the existence proofs depend on the density bounds for one step and on
the semiconcavity of the squared distance, which are known only for $\HK$
and $\SHK$. Whether analogous statements hold on other cones is open.
\end{enumerate}

\appendix
\section{Computer-checked components}\label{sec:verification}

The results of Sections~\ref{sec:prelim}--\ref{sec:sublevel} and of
Appendices~\ref{nc:sec-engine}--\ref{nc:sec-rigidity} have been formalised in the
Lean~4 proof assistant with its library Mathlib. The development has two
parts. The first, described in the next three subsections, works at the
level of abstract metric spaces and of functions of real variables: the
facts about $\HK$, optimal plans and integral functionals that the proofs
use enter as \emph{interfaces}, that is, named hypotheses, and a theorem
that uses one is proved for every object satisfying it. Such a theorem is
not a formal proof of the corresponding theorem of the paper. The second
part proves facts about actual measures and
constructs instances of the interfaces; it is described in the
subsection ``Formalised for measures'' below, and when a result is
formalised both relative to interfaces and for measures, that subsection
gives the statement for measures. The development
is organised in layers; a layer imports only Mathlib, earlier layers and,
without cycles, files of its own layer, which is checked mechanically. It has no unproved statements and no axioms
beyond the three standard ones of Lean with classical logic.

\emph{Where each result is described.} The following list names, for each
numbered result of the paper that is discussed in this appendix, the
subsections and paragraphs below that treat it.
\begin{itemize}[leftmargin=2em,itemsep=0pt]
\item Definition~\ref{prelim:convexity}, Remark~\ref{prelim:rem-convexity},
  Lemma~\ref{prelim:envelope}, Proposition~\ref{prelim:EVIprop},
  Theorems~\ref{prelim:MS-error} and~\ref{prelim:MS-convex}: ``Proved
  without interface'', Section~\ref{sec:prelim}.
\item Definition~\ref{prelim:LET}: the opening of ``Formalised for
  measures''.
\item Proposition~\ref{prelim:HKfacts}: ``The metric''.
\item Proposition~\ref{prelim:cone}: ``Geodesics in the spherical
  space''.
\item Remark~\ref{prelim:rem-dependence}: ``Interfaces''.
\item Lemma~\ref{prelim:tg}: ``Transport--growth''.
\item Lemmas~\ref{prelim:dominated} and~\ref{shk:energy}: ``The
  energy''.
\item Lemma~\ref{prelim:trunc}: ``Interfaces''; ``The energy''.
\item Lemma~\ref{prelim:lsc}: ``Lower semicontinuity''.
\item Theorem~\ref{prelim:LMS}: ``Interfaces''; ``Theorem~\ref{prelim:LMS},
  sufficiency''; ``Theorem~\ref{prelim:LMS}, sufficiency: the parts of the
  proof in \cite{LMS23}''; ``Theorem~\ref{prelim:LMS}, converse, and
  Theorem~\ref{nc:LMS-pointwise}''.
\item Theorem~\ref{nc:LMS-pointwise}: ``Theorem~\ref{prelim:LMS}, converse,
  and Theorem~\ref{nc:LMS-pointwise}''.
\item Lemma~\ref{prelim:normalise}: ``The metric''; ``The energy''.
\item Lemma~\ref{prelim:split}, Corollary~\ref{prelim:far},
  Remark~\ref{shk:farcut-sharp}: ``Far cuts''.
\item Lemma~\ref{shk:farcut}: ``Interfaces''; ``Far cuts''.
\item Lemma~\ref{tr:planar}: ``Proved without interface'',
  Sections~\ref{sec:transfer} and~\ref{sec:examples}.
\item Theorem~\ref{tr:T1}: ``Proved without interface'',
  Sections~\ref{sec:transfer} and~\ref{sec:examples}; the opening of
  ``Formalised for measures''.
\item Lemma~\ref{tr:lift}: ``Proved without interface'',
  Sections~\ref{sec:transfer} and~\ref{sec:examples}; ``Geodesics in the
  spherical space''.
\item Corollary~\ref{tr:global}: ``Proved without interface'',
  Sections~\ref{sec:transfer} and~\ref{sec:examples}; ``Lower
  semicontinuity''.
\item Theorem~\ref{tr:powers}: ``Proved without interface'',
  Sections~\ref{sec:transfer} and~\ref{sec:examples}; ``Sums of powers:
  convexity on sublevels''; ``Consequences for measures''.
\item Proposition~\ref{tr:negpowers}: ``Proved without interface'',
  Sections~\ref{sec:transfer} and~\ref{sec:examples}; ``Consequences for
  measures''.
\item Lemma~\ref{tr:geo-unique}, Corollary~\ref{tr:geo-unique-shk},
  Remark~\ref{tr:rem-strict}: ``Geodesics between separated measures''.
\item Theorem~\ref{tr:necessary}, Remark~\ref{tr:rem-consequences}: ``The
  necessary condition''.
\item Corollary~\ref{tr:consequences}: ``The necessary condition'';
  ``Consequences for measures''.
\item Theorem~\ref{tr:rigidity}: ``Proved without interface'',
  Sections~\ref{sec:transfer} and~\ref{sec:examples}; ``Rigidity''.
\item Proposition~\ref{nc:engine}, Lemmas~\ref{nc:shk-asymp}
  and~\ref{nc:dilution}: ``Rigidity''.
\item Propositions~\ref{ex:reduced} and~\ref{ex:Cd}: ``Proved without
  interface'', Sections~\ref{sec:transfer} and~\ref{sec:examples};
  ``Computer algebra''.
\item Proposition~\ref{ex:counter}, Lemma~\ref{ex:stieltjes},
  Remark~\ref{ex:rem-counter}: ``Proved without interface'',
  Sections~\ref{sec:transfer} and~\ref{sec:examples}.
\item Theorems~\ref{ex:d2}, \ref{ex:Kd} and~\ref{ex:shifted},
  Proposition~\ref{ex:d1}: ``Proved without interface'',
  Sections~\ref{sec:transfer} and~\ref{sec:examples}; ``Consequences for
  measures''.
\item Theorem~\ref{ex:bernstein}, Proposition~\ref{ex:mixture}:
  ``Consequences for measures''.
\item Proposition~\ref{ex:kernel}: ``Consequences for measures'';
  ``Computer algebra''.
\item Example~\ref{ex:ex-d2}: ``Proved without interface'',
  Sections~\ref{sec:transfer} and~\ref{sec:examples}; ``Consequences for
  measures''.
\item Example~\ref{ex:ex-bern}: ``Consequences for measures''.
\item Lemmas~\ref{cp:deriv}, \ref{cp:euler} and~\ref{cp:error},
  Proposition~\ref{cp:approxEVI}, Remark~\ref{cp:rem-dini},
  Theorem~\ref{cp:abstract}, Remarks~\ref{cp:rem-abstract}
  and~\ref{cp:rem-local}, Corollary~\ref{cp:strong},
  Example~\ref{cp:ex-error}: ``Proved without interface'',
  Section~\ref{sec:compact}.
\item Lemma~\ref{cp:nouphill}: ``No transport towards higher density''.
\item Proposition~\ref{cp:density-hk}: ``Density bounds for $\HK$''.
\item Proposition~\ref{cp:density-shk}: ``Density bounds for $\SHK$''.
\item Remark~\ref{cp:rem-sharp}: the end of ``Density bounds for $\SHK$''.
\item Proposition~\ref{cp:semiconc}: ``Interfaces''; ``Semiconcavity and the local angle
  condition''.
\item Lemma~\ref{cp:exist}: ``Minimisers on compact domains''.
\item Lemma~\ref{cp:recovery}: ``Interfaces''; ``Recovery''.
\item Lemma~\ref{cp:invariance}: ``Interfaces''; ``Invariant sets''.
\item Theorem~\ref{cp:main}: ``Interfaces''; ``Minimisers on compact
  domains''; ``Flows on compact domains''.
\item Corollary~\ref{cp:weakstrong}: ``Interfaces''; ``Flows on compact
  domains''.
\item Theorem~\ref{hk:main}: ``Proved without interface'',
  Sections~\ref{sec:hk-unbdd}--\ref{sec:sublevel}; ``Interfaces''; ``Flows
  on unbounded domains''.
\item Theorems~\ref{hk:mm}, \ref{shk:main} and~\ref{shk:mm}:
  ``Interfaces''; ``Flows on unbounded domains''.
\item Lemmas~\ref{hk:growth} and~\ref{shk:growth}: ``Proved without
  interface'', Sections~\ref{sec:hk-unbdd}--\ref{sec:sublevel}; ``Growth in
  one step''.
\item Lemma~\ref{sub:growth}: ``Growth in one step''.
\item Propositions~\ref{hk:tight} and~\ref{shk:tight}: ``Interfaces'';
  ``Uniform bounds on bounded domains''.
\item Lemma~\ref{hk:comparison}: ``Proved without interface'',
  Sections~\ref{sec:hk-unbdd}--\ref{sec:sublevel}.
\item Propositions~\ref{hk:cauchy} and~\ref{shk:cauchy}: ``Interfaces''.
\item Lemma~\ref{hk:minimiser}: ``Interfaces''; ``Minimisers on unbounded
  domains''.
\item Lemma~\ref{shk:minimiser}, Theorem~\ref{sub:resolvent}: ``Minimisers
  on unbounded domains''.
\item Definition~\ref{sub:def-sol}: ``Flows on subdomains''.
\item Lemma~\ref{sub:restrict}, Theorem~\ref{sub:abstract}: ``Proved
  without interface'', Sections~\ref{sec:hk-unbdd}--\ref{sec:sublevel};
  ``Regularity on sublevels''.
\item Theorem~\ref{sub:mm}: ``Proved without interface'',
  Sections~\ref{sec:hk-unbdd}--\ref{sec:sublevel}.
\item The estimate \eqref{sub:eq-cut}: ``Interfaces''; ``Truncation for power energies''.
\item Theorem~\ref{sub:powers-compact}: ``Interfaces''; ``Sums of powers:
  the flow''.
\item Theorem~\ref{sub:linear-extension},
  Corollary~\ref{sub:powers-measure-data}: ``Extension to measure data''.
\end{itemize}

\subsection*{Proved without interface}

\emph{Section~\ref{sec:prelim}.} The notions of convexity in
Definition~\ref{prelim:convexity} and their relations in
Remark~\ref{prelim:rem-convexity}(i); Lemma~\ref{prelim:envelope} on the
Moreau envelope; and the
results of \cite{MuS20} quoted in Section~\ref{prelim:ss-EVI}, for the
integral form \eqref{prelim:eq-intEVI} and lower semicontinuous
functionals: contraction, uniqueness, monotonicity of the energy and
\eqref{prelim:eq-expEVI}; the extension of a flow to the closure and the
improvement of the constant (Corollaries~3.8 and~3.12 of \cite{MuS20});
Theorem~\ref{prelim:MS-convex}; and both estimates of
Theorem~\ref{prelim:MS-error} with their constants, together with its
statement on convergence for $U^0=u_0\in\dom\phi$: uniformly on bounded
time intervals, for the piecewise constant and for every geodesic
interpolant, with the $\EVI$ in integral form. The proof of
\eqref{prelim:eq-MSerror} avoids \cite[Prop.~5.8]{MuS20} and so differs
from the published one. For a proper lower semicontinuous energy,
Proposition~\ref{prelim:EVIprop} gives the equivalence with the differential
$\EVI$. If the energy along the curve and the squared distances to every
finite-energy observer are locally integrable on $(0,\infty)$,
\eqref{prelim:eq-intEVI} implies continuity and finite energy at every
positive time.

\emph{Sections~\ref{sec:transfer} and~\ref{sec:examples}.}
Lemma~\ref{tr:planar}. On an abstract metric cone, that is, for two metric
spaces related by the cone identity in \eqref{prelim:eq-HKSHK}:
Theorem~\ref{tr:T1}, the uniqueness assertion of Lemma~\ref{tr:lift}, and
Corollary~\ref{tr:global}, Theorem~\ref{tr:powers} and
Proposition~\ref{tr:negpowers} for functionals with the properties of the
energy used in their proofs. At the level of the scalar condition
\eqref{prelim:eq-LMScond}: Proposition~\ref{ex:reduced} in both
directions, which is the equivalence with the matrix condition;
the kernel and Bernstein families, shifted powers, and the arctangent and
square-root densities of Proposition~\ref{ex:d1}; the differentiable
criteria of Proposition~\ref{ex:Cd} and Theorems~\ref{ex:d2} and~\ref{ex:Kd}; and
\eqref{tr:eq-qrange} with Step~2 of the proof of
Theorem~\ref{tr:rigidity} in Appendix~\ref{nc:sec-rigidity}. The best scalar constants
stated in Section~\ref{sec:examples} are formalised as well: that of
Theorem~\ref{ex:d2} (see below), also in its nonsmooth form, that of
Example~\ref{ex:ex-d2}(c), and that of Remark~\ref{ex:rem-counter}(b):
whenever $-2q$ is admissible for the shifted power, it is the best constant.
Lemma~\ref{ex:stieltjes} is formalised, with both total masses (part (c) of
Theorem~\ref{ex:shifted} is proved directly, without it). The optimal scalar constant of Theorem~\ref{ex:d2}, its recession-slope
identity and its affine equality case are formalised. The nonsmooth convex Lipschitz case of Theorem~\ref{ex:d2}
is formalised in a derivative-free form: the scalar condition is equivalent
to convexity of $G$ and a lower bound on its secant slopes, and the
constants on $\HK$ and $\SHK$ are expressed through the infimum and the
supremum of these slopes. For
Proposition~\ref{ex:counter}, the limit $F_{-2q}(u)\to-\infty$ at
$q=(d-1)/d$ and the interval of failing exponents are formalised for every
$d\ge3$, together with the explicit formula and the negative value in
$d=3$. For Proposition~\ref{tr:negpowers}, the abstract
argument gives the endpoint-energy maximum bound and the corresponding
convexity constant along a geodesic; its application to the actual
negative-power energies and their sublevels is given under ``Consequences
for measures'' below.

\emph{Section~\ref{sec:compact}.} The first variation of the squared
distance and the estimate of Cauchy--Schwarz type under the local angle
condition that enter the derivative estimate; Lemmas~\ref{cp:deriv},
\ref{cp:euler} and~\ref{cp:error},
Proposition~\ref{cp:approxEVI} in the form of Remark~\ref{cp:rem-dini},
Example~\ref{cp:ex-error} as a statement about the scheme on $[0,1]$,
Theorem~\ref{cp:abstract} with the convergence of all geodesic
interpolants, and Corollary~\ref{cp:strong}. Steps~1 and~2 of the proof of
Theorem~\ref{cp:abstract} are formalised in the variant of
Remark~\ref{cp:rem-abstract}(ii), under the hypothesis of
Remark~\ref{cp:rem-local}.

\emph{Sections~\ref{sec:hk-unbdd}--\ref{sec:sublevel}.}
Lemma~\ref{hk:comparison}; the proof of Theorem~\ref{hk:main}, as a
theorem about flows on an increasing family of closed subsets of a
complete metric space, given the estimates of Lemma~\ref{prelim:trunc} and
Proposition~\ref{hk:tight}; the first-order arguments of
Lemmas~\ref{hk:growth} and~\ref{shk:growth}, with the inequalities for the
competitor as hypotheses; the real-variable parts of the tightness and
minimiser statements of these sections, with their constants; and
Lemma~\ref{sub:restrict}, Theorem~\ref{sub:abstract} and
Theorem~\ref{sub:mm}, including the regularity assertions (see
``Regularity on sublevels'' below).
For Theorem~\ref{sub:abstract}, hypothesis~\ref{sub:A1} is formalised in
its stated derivative-only form. Absolute continuity of squared distance
along a geodesic and the passage from the almost-everywhere derivative
bound to the secant estimate are also proved in Lean; the stronger
condition $(\mathrm G_\kappa)$ is not assumed in this theorem.

\subsection*{Interfaces}

In the first part four groups of facts enter as interfaces.
(a) \emph{The metric $\HK$}: a metric space with a mass, a zero element and
a scaling $\mu\mapsto r^2\mu$ that satisfy $\HK^2(\mu,0)=\mu(\Om)$, the
upper bound in \eqref{prelim:eq-massbound} and \eqref{prelim:eq-scaling};
$\SHK$ is defined by \eqref{prelim:eq-SHK}, and its properties listed in
Remark~\ref{prelim:rem-dependence} are derived.
(b) \emph{The criterion of \cite{LMS23}}: the implication from
\eqref{prelim:eq-LMScond} to $\lambda$-convexity of the functional on
measures (Theorem~\ref{prelim:LMS}), with its lower semicontinuity. This
implication is itself formalised for the functional
\eqref{prelim:eq-energy}, along every given geodesic (see
``Theorem~\ref{prelim:LMS}, sufficiency'' below), and so is the lower semicontinuity. The statements of
Sections~\ref{sec:transfer} and~\ref{sec:examples} that need only this
criterion have been combined with its proof (see ``Consequences for
measures'' below); so have the existence theorems of
Sections~\ref{sec:compact}--\ref{sec:sublevel} (see ``Flows on unbounded
domains'' and ``Sums of powers: the flow'' below).
(c) \emph{Radial growth of internal energies}: the bound
$\sfE(r^2\nu)\le\sfE(s^2\nu)+\Einf(r^2-s^2)$, $0<s\le r$, from the proof
of Corollary~\ref{tr:global}, and its analogues for sums of powers in the
proof of Theorem~\ref{tr:powers}.
(d) \emph{The measure-level estimates of
Sections~\ref{sec:compact}--\ref{sec:sublevel}}: existence of minimisers
of the step functionals; the density bounds, as invariance of
density-bounded sets under the scheme; semiconcavity on these sets
(Proposition~\ref{cp:semiconc}); recovery (Lemma~\ref{cp:recovery});
truncation (Lemma~\ref{prelim:trunc} and \eqref{sub:eq-cut}); the far-cut
inequality \eqref{prelim:eq-farcut}; tightness
(Propositions~\ref{hk:tight} and~\ref{shk:tight}, \eqref{sub:eq-tight});
and, in Sections~\ref{sec:hk-unbdd}--\ref{sec:sublevel}, the flows on the
bounded domains $\OmR$.

For the densities of Section~\ref{sec:examples} for which the scalar LMS
condition was proved above under ``Proved without interface'', (a)--(c)
give the stated convexity constants on $\HK$ and $\SHK$ under the
respective transfer hypotheses. Relative to (a)--(d), the
remaining steps are formalised for Lemma~\ref{cp:invariance},
Theorem~\ref{cp:main}, Corollary~\ref{cp:weakstrong},
Theorems~\ref{hk:main} and~\ref{shk:main} with
Propositions~\ref{hk:cauchy} and~\ref{shk:cauchy},
Lemmas~\ref{hk:minimiser} and~\ref{shk:farcut}, Theorems~\ref{hk:mm}
and~\ref{shk:mm} as estimates of order $\tau^{1/4}$, and
Theorem~\ref{sub:powers-compact}, whose regularity assertion follows from
its assertion (i) as described under ``Regularity on sublevels''.
Proposition~\ref{cp:semiconc}, Lemma~\ref{cp:invariance},
Theorem~\ref{cp:main}, Corollary~\ref{cp:weakstrong},
Theorems~\ref{hk:main}, \ref{hk:mm}, \ref{shk:main}, \ref{shk:mm}
and~\ref{sub:powers-compact} are also formalised for measures with no
interface; the results for actual measures and their hypotheses are
described under ``Formalised for measures'' below.

\subsection*{Scope of the first part}

The sharpness statements of Sections~\ref{sec:transfer}
and~\ref{sec:examples}, including those of the examples and of
Remarks~\ref{tr:rem-strict}, \ref{cp:rem-sharp}
and~\ref{shk:farcut-sharp}, are formalised. The abstract
treatment of Sections~\ref{sec:hk-unbdd} and~\ref{sec:shk-unbdd} takes the
inequalities for competitors as hypotheses and does not use
transport--growth pairs; Lemma~\ref{prelim:tg} itself is formalised in
full (see ``Transport--growth'' below).

\subsection*{Formalised for measures}

The following is proved in Lean for finite Borel measures on a complete
second-countable metric space with its Borel sets, without interfaces. The
distance is defined as the square root of the infimum of the cost over
finite plans on the cone, as in \eqref{prelim:eq-coneOT} but with the
cutoff $\pi$ and an infimum in place of a minimum. That this infimum
equals the one with the cutoff $\pi/2$ of \eqref{prelim:eq-coneOT} is
proved in Lean, and so is its equality with the infimum of the
entropy--transport objective \eqref{prelim:eq-LET} over the plans of
Definition~\ref{prelim:LET}; both infima are attained: a finite
plan attains the minimum in \eqref{prelim:eq-LET}, and a probability
measure on the product of the two cones attains the minimum in
\eqref{prelim:eq-coneOT}.

Theorem~\ref{tr:T1} is formally proved for measures by applying the
abstract cone theorem to the cone of finite measures over the
probabilities. Theorems~\ref{tr:necessary}, \ref{tr:powers},
\ref{ex:Kd} and~\ref{ex:shifted} and Corollary~\ref{tr:global} are also
formalised for measures under their printed hypotheses.
Theorem~\ref{ex:d2} is formalised under its printed hypotheses.
So is Theorem~\ref{ex:bernstein}. The measure-level existence and convergence
results include Theorem~\ref{cp:main}, Corollary~\ref{cp:weakstrong} and
Theorems~\ref{hk:main}, \ref{hk:mm}, \ref{shk:main}, \ref{shk:mm} and
\ref{sub:powers-compact}, in the form described below.

\emph{The metric.} This distance is a metric on finite measures, and
\eqref{prelim:eq-SHK} defines a metric on probability measures. The
finite measures satisfy interface (a), and the probabilities are its
sphere of mass one, so that interface is instantiated (the cone structure
over the probabilities for a nonempty base). Proved for
them: \eqref{prelim:eq-subadd}, \eqref{prelim:eq-massbound},
\eqref{prelim:eq-scaling} for nonnegative factors, the mass identity
\eqref{prelim:eq-mass} along a given geodesic, the cone identity in
\eqref{prelim:eq-HKSHK} and \eqref{prelim:eq-normalise}; that both
distances metrise narrow convergence, that both spaces are complete and
separable, and that the probabilities over a compact base form a compact
space; that any two finite measures on a closed convex subset of a
complete separable normed space are joined by an $\HK$ geodesic, so that
all assertions of Proposition~\ref{prelim:HKfacts} are formalised for
such a subset; the isometry and
closedness assertions of Lemma~\ref{prelim:normalise} for closed subsets;
and, for the objective \eqref{prelim:eq-LET}, that a minimising plan
multiplied by $r_0r_1$ is minimising for the scaled measures, and that
the minimiser is unique and explicit when the two measures are
proportional or are carried by two sets at distance at least $\pi/2$.

\emph{The energy.} The relaxed energy \eqref{prelim:eq-energy} is defined
in Lean, with \eqref{prelim:eq-energybounds},
Lemma~\ref{prelim:dominated}, all of Lemma~\ref{prelim:trunc} including
the convergence of the truncations, Lemma~\ref{shk:energy}, the energy
and $\EVI$ assertions of Lemma~\ref{prelim:normalise} (the latter for the
integral form), and the radial growth bound of interface (c) for finite
recession slope on a nonempty base. The reference measure is
$\sigma$-finite, and the statements about the energy carry the hypothesis
under which its negative part is finite, that is, $\ezero>-\infty$ or a
finite reference measure.

\emph{Lower semicontinuity.} Lemma~\ref{prelim:lsc} is formalised in full
for the relaxed energy \eqref{prelim:eq-energy}, in both of its cases:
for a finite reference measure $m$ and every convex $E$ with $E(0)=0$ that is
continuous at $0$, including $\ezero=-\infty$; and for a reference measure
that is finite on closed balls and finite $\ezero$. The energy is lower
semicontinuous for narrow convergence of finite measures on a metric space,
and for $\HK$ and $\SHK$ on a complete separable one. The formal proof does
not use \cite[Thm.~2.34]{AFP00}: after subtracting an affine minorant, the
energy is written as a supremum of functionals
$\mu\mapsto\int\varphi\dd\mu+\int\psi\dd m$ with $\varphi,\psi$ bounded
and continuous, and an infinite reference measure is reached by continuous
cut-off functions. Only the forward implication, from convergence in $\HK$
or $\SHK$ to narrow convergence, enters. With it,
Corollary~\ref{tr:global} is formalised for the relaxed energy on measures
in both cases.

\emph{Transport--growth.} The two upper bounds of Lemma~\ref{prelim:tg}
are formalised: \eqref{prelim:eq-tg-upper} for a finite source measure on
any measurable space, measurable maps and square-integrable factors, with
target a separable metric space, and the split-source bound
\eqref{prelim:eq-tg-split} for finitely many pieces. For the optimal pair
\eqref{prelim:eq-tg-opt}, the step from optimal potentials to the maps is
formalised, on a separable real Hilbert space and for every $\tau>0$, under
a regularity hypothesis: let $(\varphi_0,\varphi_1)$ be an optimal pair of
Borel potentials for $\mu_0,\mu_1$, and let $D$ be a Borel set such that
$\varphi_0$ is finite near every point of $D$ and differentiable there, and
such that $D$ carries the part of $\mu_0$ on the points of its support at
distance less than $\pi/2$ from the support of $\mu_1$. Then the minimiser
of \eqref{prelim:eq-LET} is unique and concentrated on the graph of
$\bfT(x)=x+\boldsymbol{\arctan}(\tau\nabla\varphi_0(x))$, where
$\boldsymbol{\arctan}(0)=0$ and
$\boldsymbol{\arctan}(v)=\arctan(|v|)\,v/|v|$ for $v\ne0$. If moreover
$\mu_1$ gives no mass to the points of its support at distance at least
$\pi/2$ from the support of $\mu_0$, then $\bfT$ and
$q=\ee^{2\tau\varphi_0}\sqrt{1+\tau^2|\nabla\varphi_0|^2}$ on $D$, extended by
$\bfT(x)=x$ and $q=0$, satisfy \eqref{prelim:eq-tg-opt}, and every pair
satisfying \eqref{prelim:eq-tg-opt} has the same $q$ $\mu_0$-almost
everywhere and the same $\bfT$ $\mu_0$-almost everywhere on $\{q>0\}$. The
existence statement \eqref{prelim:eq-tg-opt} of the lemma is formalised as
stated, with no interface: for a compact convex set $X$ with nonempty
interior in a finite-dimensional real inner product space, $\mu_0$ with a
density positive almost everywhere on $X$ and any finite $\mu_1$ on $X$,
there are Borel maps $\bfT:X\to X$ and $q:X\to[0,\infty)$ satisfying
\eqref{prelim:eq-tg-opt}, where $\HK$ is the distance of $\calM(X)$; the
case $\mu_1=0$ is included. The regularity hypothesis is obtained there
from a Lipschitz optimal potential and Rademacher's theorem. The
second-order regularity of optimal potentials of \cite{LMS23} is not
formalised and is not used.

\emph{Far cuts.} The far-cut inequality \eqref{prelim:eq-farcut} and the
identity \eqref{prelim:eq-far}. The objective \eqref{prelim:eq-LET} is
defined for finite plans, and for every plan of finite objective the
measure $\mu_0''$ of \eqref{prelim:eq-splitdef} is concentrated within
$\pi/2$ of $F$ and satisfies $\mu_0''(\Om)\le\eps$ in the situation of
Lemma~\ref{prelim:split}; for a given optimal plan on the cone, the
squared distances of the two pieces obtained by restricting the plan over
$F$ add up. For a minimising plan of \eqref{prelim:eq-LET} the objective splits
exactly over the pieces \eqref{prelim:eq-splitdef}, both restricted plans
are minimising, and the identity \eqref{prelim:eq-split} holds. Thus
Lemma~\ref{prelim:split} is formalised for the minimising plan in its
statement, and Corollary~\ref{prelim:far} under its stated separation
hypothesis. Further, Lemma~\ref{shk:farcut} with its
strictness statement, in the form of an equivalence; and the optimality
of the constant in Remark~\ref{shk:farcut-sharp}.

\emph{Theorem~\ref{prelim:LMS}, sufficiency.} The sufficiency
part of Theorem~\ref{prelim:LMS} itself is formalised, in the following
form. Let $\Om\subset\R^d$, $d\ge1$, be closed, let $E:[0,\infty)\to\R$ be
convex and continuous at $0$ with $E(0)=0$, and suppose that $\Om$ has
finite Lebesgue measure or that $\ezero>-\infty$. If
\eqref{prelim:eq-LMScond} holds with a constant $\lambda\in\R$, then for
every constant-speed geodesic $(\mu_t)_{t\in[0,1]}$ of $(\calM(\Om),\HK)$
with $\sfE(\mu_0)<\infty$ and $\sfE(\mu_1)<\infty$,
\[
  \sfE(\mu_t)\le(1-t)\sfE(\mu_0)+t\,\sfE(\mu_1)
    -\frac\lambda2\,t(1-t)\,\HK^2(\mu_0,\mu_1),\qquad t\in[0,1],
\]
where $\sfE$ is \eqref{prelim:eq-energy} with its singular term and $\HK$
is built from the Euclidean distance restricted to $\Om$. Convexity of
$\Om$ and a nonempty interior are not used; measures with a singular part,
a zero endpoint and $\Einf=+\infty$ are included. The statement concerns
given geodesics; that any two finite measures are joined by a geodesic is
formalised when $\Om$ is closed and convex, as stated above. The formal proof follows \cite[Theorem~7.2]{LMS23} up to the
density formula; it then bounds the energy per unit mass, $E(c)/c$ on the
absolutely continuous part and $\Einf$ on the singular part, at
$\mu_t$-almost every point $x$ by the convex combination of its values at
the two endpoints of the characteristic through $x$, weighted by the
squared growth factors $q_{t\to0}(x)^2$ and $q_{t\to1}(x)^2$, minus
$\frac\lambda2t(1-t)$ times the squared speed
$4\xi_t(x)^2+|\nabla\xi_t(x)|^2$, and integrates against $\mu_t$; it does
not treat the absolutely continuous and the
singular part of the geodesic as two separate geodesics, as that proof
does.

\emph{Theorem~\ref{prelim:LMS}, sufficiency: the parts of the proof in
\cite{LMS23}.} Of the proof of the sufficiency part in \cite{LMS23}, the
scalar step is formalised, for every real-valued density satisfying the condition of
Theorem~\ref{prelim:LMS} and with no smoothness assumed: along twice
differentiable positive curves $(\rho,\gamma)$ obeying the curvature
conditions \cite[(7.3)]{LMS23} (an equality in dimension one), the function
$t\mapsto N_{\lambda,E}(\rho(t),\gamma(t))$ is convex. The pointwise geometric
part is formalised as well, on a finite-dimensional real inner-product space
of positive dimension: for
a pair of potentials $\xi_0\ge-\frac12$ and $\bar\xi_1$ not exceeding the
forward evolution of $\xi_0$ at time one, the forward and backward
generalised Hopf--Lax evolutions, for each $s\in(0,1)$ their contact sets
$\Xi_s$ with the common Lipschitz
gradient, the base maps $\bfT_{s\to t}$, the radius factors $q_{s\to t}$ and
the spatial Jacobians $\delta_s(t,x)=\det\mathrm{D}\bfT_{s\to t}(x)$, the derivative
being taken within $\Xi_s$, as in \cite[Sect.~4]{LMS23} with the corrections
listed in Appendix~\ref{sec:errata}; for each $s$ there is one Borel set
$D\subseteq\Xi_s$ with $\Leb(\Xi_s\setminus D)=0$ such that, for every $x\in D$,
the curvature conditions \cite[(5.12)]{LMS23} hold at every $t\in(0,1)$ and,
for every $m>0$, the function
$t\mapsto\delta_s(t,x)\,E\big(mq_{s\to t}(x)^2/\delta_s(t,x)\big)-\frac\lambda2mq_{s\to t}(x)^2$ is
convex on $(0,1)$. The formal proof uses neither the convex blow-up of
\cite[Lemma~4.6]{LMS23} nor the invariance of density points under
bi-Lipschitz maps. Of the statements about measures that the proof uses,
the first is formalised as well, for finite Borel measures on a complete
second-countable metric space and every $\tau>0$: the dual problem for
$\HK^2/(2\tau)$ of \cite[Theorem~2.12]{LMS23}, posed over pairs of Borel
potentials with values in $[-\infty,\infty]$, has a maximiser; all
minimisers of \eqref{prelim:eq-LET} have the same two marginals; and for
every optimal pair the two exponential factors are the densities of these
marginals, almost everywhere for the respective given measures. For every
minimising plan the finite contact identity \cite[(2.45a)]{LMS23} holds
almost everywhere for that plan. Among the optimal pairs there is one
in which each potential is the transform of the other,
$\varphi_0(x)=\sup_y\big(\varphi_1(y)+\tau^{-1}\log\cos|x-y|\big)$ and
$\varphi_1(y)=\inf_x\big(\varphi_0(x)-\tau^{-1}\log\cos|x-y|\big)$ with the
supremum and the infimum over $|x-y|<\pi/2$, so that $\varphi_0$ is lower
and $\varphi_1$ upper semicontinuous, as in \cite[Theorem~3.3(1)]{LMS23};
here $|x-y|$ stands for the distance of the space. The inclusions of
the support of $\mu_0$ in $\{\varphi_0<\infty\}$ and of the support of
$\mu_1$ in $\{\varphi_1>-\infty\}$ stated there are formalised with the
closures of these two sets; without the closures they are false.
\cite[Theorem~5.1]{LMS23} is formalised on a separable real Banach space: let
$(\mu_t)_{t\in[0,1]}$ be a given constant-speed $\HK$ geodesic and
$(\varphi_0,\varphi_1)$ such a pair for $\mu_0,\mu_1$ with $\tau=1$, and
let $\xi_t$ and $\bar\xi_t$ be the forward and backward evolutions of
$\xi_0=\frac12(\ee^{2\varphi_0}-1)$ and
$\bar\xi_1=\frac12(1-\ee^{-2\varphi_1})$. Then
$\int\xi_t\dd\mu_t-\int\xi_s\dd\mu_s=\HK^2(\mu_s,\mu_t)/(2(t-s))$ for
$0\le s<t\le1$, the same holds for $\bar\xi$, and the support of $\mu_t$
lies in $\{\xi_t=\bar\xi_t\}$ for every $t\in[0,1]$. The
representation of the geodesic by the maps above,
\cite[Theorem~5.2]{LMS23}, is formalised on a finite-dimensional real
inner-product space, with no assumption of absolute continuity: for $0<s<1$ and $t\in[0,1]$,
\[
 \mu_t=(\bfT_{s\to t})\pfw\big(q_{s\to t}^2\mu_s\big),\qquad
 \HK^2(\mu_s,\mu_t)=\int_{\Xi_s}\big[1+q_{s\to t}^2
 -2q_{s\to t}\cospi|\bfT_{s\to t}(x)-x|\big]\dd\mu_s ;
\]
for every pair $(\mathbf{S},p)$ satisfying these two identities, with $p\ge0$
and both maps measurable for the completion of $\mu_s$, the points
$[\mathbf{S},p]$ and $[\bfT_{s\to t},q_{s\to t}]$ of the cone agree
$\mu_s$-almost everywhere; and the minimiser $\eta$ of
\eqref{prelim:eq-LET} for $\mu_s,\mu_t$ is unique and satisfies
$\eta=(\mathrm{id},\bfT_{s\to t})\pfw\eta_0$ for its first marginal
$\eta_0$. The source times $s=0$ and $s=1$ are not covered. The
geodesic assertion of \cite[Theorem~5.4]{LMS23} is formalised on the same
spaces: for $0<s<1$ and every finite Borel measure $\nu$ concentrated on
$\Xi_s$, the curve $\nu_t=(\bfT_{s\to t})\pfw(q_{s\to t}^2\nu)$,
$t\in[0,1]$, is a constant-speed $\HK$ geodesic with
$\HK^2(\nu_0,\nu_1)=\int(4\xi_s^2+|\nabla\xi_s|^2)\dd\nu$, the gradient
being the common gradient on $\Xi_s$; this holds for any lower
semicontinuous $\xi_0\ge-\frac12$ and upper semicontinuous $\bar\xi_1$ not
exceeding the forward evolution of $\xi_0$ at time one, and in particular
for the restriction of $\mu_s$ to any Borel set when $(\mu_t)$ is a given
geodesic, with the tight optimal pair above.
\cite[Corollary~5.5]{LMS23} is formalised on the same spaces in the
following form, for a given geodesic $(\mu_t)$ with the tight optimal pair
above. If $\mu_s\ll\Leb$ for one $s\in(0,1)$, or if $\mu_0\ll\Leb$ and the
restriction of $\mu_1$ to the points of its support at distance at least
$\pi/2$ from the support of $\mu_0$ is absolutely continuous, then
$\mu_t\ll\Leb$ for every $t\in(0,1)$. For $s,t\in(0,1)$ and Borel densities
$c_s$, $c_t$ of $\mu_s$, $\mu_t$,
\[
  c_t\big(\bfT_{s\to t}(x)\big)=\frac{c_s(x)\,q_{s\to t}(x)^2}{\delta_s(t,x)}
  \quad\text{and}\quad\delta_s(t,x)>0
  \qquad\text{for $\Leb$-almost every }x\in\Xi_s ;
\]
the exceptional set depends on $s$, $t$ and the two densities, whereas
\cite[(5.6a)]{LMS23} is printed for every $x$ in the differentiability
set of \cite[(4.57)]{LMS23} and every $t\in(0,1)$. For
$s,t\in(0,1)$ the restrictions to $\{\bar\xi_1=\frac12\}$ satisfy
$\mu_t=(t/s)^2\mu_s$ and those to $\{\xi_0=-\frac12\}$ satisfy
$\mu_t=((1-t)/(1-s))^2\mu_s$, also when the measures are singular; this
gives \cite[(5.7)]{LMS23} for the densities, almost everywhere for each
pair of times. At an endpoint $i\in\{0,1\}$ the formal statement is a sum
over preimages: let $0<s<1$, $\mu_s=c_s\Leb$ and $\mu_i=c_i\Leb$ with Borel
densities, and let $S$ be the set of the points $x$ of a certain Borel
subset of $\Xi_s\cap\{c_s>0\}$ of full $\mu_s$-measure at which
$q_{s\to i}(x)>0$ and $\delta_s(i,x)>0$. Then the points of that subset
with $q_{s\to i}>0$ and $\delta_s(i,\cdot)=0$ form a Lebesgue null set,
every $y$ has at most countably many preimages in $S$ under
$\bfT_{s\to i}$,
\[
  c_i(y)=\sum_{x\in S,\ \bfT_{s\to i}(x)=y}
    \frac{c_s(x)\,q_{s\to i}(x)^2}{\delta_s(i,x)}
  \quad\text{for $\Leb$-almost every }y,
\]
and consequently
$c_s\,q_{s\to i}^2/\delta_s(i,\cdot)\le c_i\circ\bfT_{s\to i}$ almost
everywhere on $S$. The identity of \cite[Corollary~5.5(3)]{LMS23} with a
single preimage is not formalised, nor is part~(4) of that corollary.
\cite[Corollary~5.6]{LMS23} is formalised on the same spaces: let
$(\mu_t)$ be a given geodesic with the tight optimal pair above,
$0<s<1$, and $\mu_s=\mu_s^{\mathrm a}+\mu_s^\perp$ the Lebesgue
decomposition of $\mu_s$ with respect to $\Leb$. Then
$a_t=(\bfT_{s\to t})\pfw(q_{s\to t}^2\mu_s^{\mathrm a})$ and
$b_t=(\bfT_{s\to t})\pfw(q_{s\to t}^2\mu_s^\perp)$, $t\in[0,1]$, are
constant-speed $\HK$ geodesics; $\mu_t=a_t+b_t$ and $b_t\perp\Leb$ for
every $t\in[0,1]$; and for $0<t<1$, $a_t\ll\Leb$, so that $a_t$ and $b_t$
are the absolutely continuous and the singular part of $\mu_t$. Absolute
continuity of $a_0$ and $a_1$ is not asserted.

\emph{Theorem~\ref{prelim:LMS}, converse, and
Theorem~\ref{nc:LMS-pointwise}.}
Lean also proves the converse of Theorem~\ref{prelim:LMS} under its stated
hypotheses in every dimension. For $d\ge1$, the isotropic test of
Appendix~\ref{nc:sec-lms} gives $B_\lambda(c)\ge0$ from weak geodesic
$\lambda$-convexity. Lean verifies the density integral on the compactly
supported test measures even when $\Om$ has infinite volume. It also
proves the shear clause of Theorem~\ref{nc:LMS-pointwise} for $d\ge2$.
For $d=0$, the unique geodesic interpolates the square roots of the masses
linearly, including at zero mass. The global criterion for densities
taking $+\infty$ in \cite{LMS23} remains open in Lean. So does the
regularity theorem of \cite{LMS23} for general tight optimal pairs of
potentials, which includes their second-order differentiability almost
everywhere; neither is used in this paper. The first-order
differentiability needed for Lemma~\ref{prelim:tg} is formalised (see
``Transport--growth'' above).

\emph{Geodesics between separated measures.} For measures carried by two
sets at distance at least $\pi/2$: the values
$\HK^2(\mu_0,\mu_1)=\mu_0(X)+\mu_1(X)$ and $\SHK(\nu_0,\nu_1)=\pi/2$ of
Lemma~\ref{tr:geo-unique} and Corollary~\ref{tr:geo-unique-shk}, and that
the curves \eqref{tr:eq-reaction} and \eqref{tr:eq-reaction-shk} are
geodesics; and, for two Dirac masses at the endpoints of a geodesic of the
base of length exactly $\pi/2$, the two different geodesics of
Remark~\ref{tr:rem-strict}. If every point of the support of one measure
is at distance larger than $\pi/2$ from every point of the support of the
other, the curves \eqref{tr:eq-reaction} and \eqref{tr:eq-reaction-shk}
are the only geodesics; this gives Lemma~\ref{tr:geo-unique} and
Corollary~\ref{tr:geo-unique-shk} in full. Lemma~\ref{tr:geo-unique}
also holds when one of the two finite measures is zero.

\emph{Geodesics in the spherical space.} Geodesicity of $(\calM(X),\HK)$
implies that of $(\calP(X),\SHK)$, in particular for closed convex
subsets of complete separable normed spaces.
For distinct probabilities $\nu_0,\nu_1$, the
geodesics of the two spaces correspond to each other by the
reparametrisation and rescaling of Proposition~\ref{prelim:cone}.
This bijection extends to $\HK$ endpoints $r_i^2\nu_i$ for $r_i>0$.
Normalisation at the reparametrised time recovers the $\SHK$ geodesic \cite[Thm.~2.7]{LaM19}.
Lean proves Lemma~\ref{tr:lift} in full, including both identities in
\eqref{tr:eq-lift} at every affine time.

\emph{The necessary condition.} Theorem~\ref{tr:necessary} is formalised
in full, for the relaxed energy \eqref{prelim:eq-energy} and the metric
$\SHK$ on the probability measures, with a Haar measure as reference, over
a compact convex set with nonempty interior and diameter larger than
$\pi/2$ in a finite-dimensional normed space, and over a closed set with
nonempty interior in a finite-dimensional real inner-product space of
positive dimension: for $E$ convex with $E(0)=0$ and
finite $\ezero=E'_+(0)$, weak geodesic $\lambda$-convexity
implies $\Einf<\infty$, the bound \eqref{tr:eq-necessary}, and $E(c)=\ezero
c$ if $\lambda=0$. With it, the assertion of
Corollary~\ref{tr:consequences}(iii) that no constant larger than
$2(\ezero-\Einf)$ is possible is formalised. The three examples of
Remark~\ref{tr:rem-consequences} are formalised on the compact convex
domains just described: no
$\lambda\in\R$ is possible for $E(c)=c^m$ with $m>1$, for $E(c)=c\log c$,
and for $E(c)=-c^q$ with $\frac12<q<1$, together with the bound
$\lambda\le-2mh^{m-1}$ and the two expansions printed there; for
$0<q\le\frac12$ the profile \eqref{tr:eq-profile} is convex, so these test
geodesics give no obstruction.

\emph{Sums of powers: convexity on sublevels.} For $E(c)=\sum_i\alpha_ic^{p_i}$ with $\alpha_i\ge0$,
$p_i\ge1$ and at least one term with $\alpha_i>0$, $p_i>1$, the relaxed
energy is identified with the functional of Theorem~\ref{tr:powers} (the
integral on absolutely continuous measures, $+\infty$ otherwise), and the
two bounds used in its proof are formalised for measures: $\sfE(a\mu)\le
a^p\sfE(\mu)$ for $a\ge1$, and the radial derivative at a probability
measure of finite energy exists as a right limit and equals
$2\sum_ip_i\alpha_i\int\rho^{p_i}\le2p\,\sfE$. The scalar condition of
Theorem~\ref{prelim:LMS} with $\lambda=0$ is proved for these densities in
every dimension. Both conclusions of Theorem~\ref{tr:powers} then follow
for measures, by the lower semicontinuity above and the sufficiency part of
Theorem~\ref{prelim:LMS}, which is formalised separately (see
``Theorem~\ref{prelim:LMS}, sufficiency'' above); the combined statement is given under ``Consequences for
measures''.

\emph{Consequences for measures.} With the sufficiency part of
Theorem~\ref{prelim:LMS} the following are formal statements about measures
with no interface, on a closed $\Om\subset\R^d$, $d\ge1$, with Lebesgue
measure, and along every given geodesic. Convexity of $\Om$ and a nonempty
interior are used only for the sharpness statements.
Theorem~\ref{tr:powers}: for $\alpha_i\ge0$, $1\le p_i\le p$ with at least
one $\alpha_i>0$, $p_i>1$, every $\SHK$ geodesic $\nu$ in $\calP(\Om)$ and
every real $c\ge\max\{\sfE(\nu_0),\sfE(\nu_1)\}$, both conclusions hold
with the constants $2^pc$ and $-p\,2^{p+1}c$; $\Om$ may be unbounded.
Proposition~\ref{tr:negpowers}, for $\Om$ of finite Lebesgue measure and
$(d,q)$ as in \eqref{tr:eq-qrange}: the bound of the energy along the
geodesic by the larger endpoint value, the constant
$-2q\max\{\sfE_q(\nu_0),\sfE_q(\nu_1)\}\ge0$, strong $0$-convexity on
$\HK$ and on $\SHK$, and the statement on the sublevels
$\{\sfE_q\le-c\}$, $c>0$. Corollary~\ref{tr:consequences}: (i) and (ii) for $\Om$
of finite Lebesgue measure or $\ezero>-\infty$, part (ii) with any real
number $L\ge\Einf$ in place of $\Einf$; and (iii) on the domains of the
formal necessary condition above, for
$\ezero>-\infty$, $\Einf<\infty$ and $\lambda=2\ezero$, where it says that
no constant larger than $2(\ezero-\Einf)$ is possible. The constants on $\HK$ and on $\SHK$
stated in Theorems~\ref{ex:Kd} and~\ref{ex:shifted} and in
Propositions~\ref{ex:mixture} and~\ref{ex:d1} are formalised in the same
way, with the sharpness on $\SHK$ of the constants $2(\ezero-\Einf)$ in
Theorem~\ref{ex:Kd} for $\Einf<\infty$, $-2q$ in
Theorem~\ref{ex:shifted} and $-2$ in Proposition~\ref{ex:d1} on compact
convex sets with nonempty interior and diameter greater than $\pi/2$ (on
closed convex sets that are not compact the sharpness follows from the
formal statement of Theorem~\ref{tr:necessary}); so is the corresponding statement for the kernel
$K_d$ of Proposition~\ref{ex:kernel}, whose functional has the constant
$-2$ on $\HK$ and the sharp constant $-2$ on $\SHK$. For Theorem~\ref{ex:bernstein} the two constants and
the sharpness on $\SHK$ are formalised, for every real $\beta$, and so are
the assertions $E\in C^1([0,\infty))$, $E'(0^+)=-\beta-|m|$ and
$-\beta-|m|\le E'\le-\beta$.
For Theorem~\ref{ex:d2} the constants are formalised in the form
$2G'(0^+)$ on $\HK$ and $2(G'(0^+)-G'(\infty))$ on $\SHK$ for convex
Lipschitz $G$, under the differentiability hypotheses of that theorem and
for $\Om$ of finite Lebesgue measure or $\ezero>-\infty$, as in
Corollary~\ref{tr:consequences}(i). The derivative limits equal $\inf G'$
at zero and $\sup G'=\Einf$ at infinity. The stated $\SHK$ constant is
zero exactly for $E(c)=-\alpha\sqrt{c}+\beta c$, with $\alpha\ge0$.
Example~\ref{ex:ex-d2}(b) is formalised on the
actual measure spaces, with the sharpness of its constant on $\SHK$ on
compact convex domains with nonempty interior and diameter greater than
$\pi/2$. Example~\ref{ex:ex-bern}(b), (c) and (d) are formalised for
the printed densities on the actual measure spaces, in every dimension
$d\ge1$: the constants on $\HK$ and on $\SHK$, the effect of a drift, and
the sharpness of the constant on $\SHK$ on compact convex domains with
nonempty interior and diameter greater than $\pi/2$.

\emph{Rigidity.} For Theorem~\ref{tr:rigidity}, Lean proves all three implications on a
compact convex domain with nonempty interior and diameter greater than
$\pi/2$, for every $0<q<1$. Outside the range \eqref{tr:eq-qrange},
the negative midpoint defect at the centre of the test ball remains strictly
negative after integration. Proposition~\ref{nc:engine} gives the exact
squared $\HK$ distance and the common midpoint of all $\HK$ geodesics
after a separated Dirac probability is added. A geodesic witnessing weak
$\SHK$ convexity would then satisfy the chord inequality at the
normalisation of that midpoint, contrary to
Lemmas~\ref{nc:shk-asymp} and~\ref{nc:dilution}.

\emph{Minimisers on compact domains.} Lemma~\ref{cp:exist} is formalised
for measures with no interface, on a compact metric space with a finite
reference measure, for $E$ convex and continuous at $0$ with $E(0)=0$ and
constants $a,b\ge0$ with $E(c)\ge-ac-b$ for $c\ge0$: for every finite
measure $y$ and every $\tau>0$ with $2a\tau<1$ the $\HK$ step functional has
a minimiser; for every probability $y$ and every $\tau>0$ the $\SHK$ step
functional has a minimiser; and in both cases every minimiser has finite
energy if $\sfE(y)<\infty$. Consequently, for every $\mu_0$ of finite
energy an exact minimizing movement scheme started at $\mu_0$ exists, for
every $\tau>0$ on $\SHK$ and for every $\tau>0$ with $2a\tau<1$ on $\HK$;
this is the existence assertion of Theorem~\ref{cp:main}\ref{cp:it-mm};
its convergence assertions are described under ``Flows on compact
domains'' below. This
concerns the compact domains of
Section~\ref{sec:compact} only.

\emph{No transport towards higher density.} Lemma~\ref{cp:nouphill} is
formalised for measures with no interface, on a complete separable metric
space with a finite reference measure $m$, in the following form. Let
$\tau>0$, let $\mu_0\ll m$ be a nonzero finite measure
$[$a probability$]$, and let $\mu_1$ be a minimiser of the $\HK$
$[\SHK]$ step functional. Then $\mu_1\ll m$; and for every pair
$(\bfT,q)$ satisfying \eqref{prelim:eq-tg-opt} and every representative
$\rho_1$ of the density of $\mu_1$, \eqref{cp:eq-nouphill} holds. The
hypothesis $\rho_0>0$ almost everywhere of the lemma is used only to make
$\mu_0$ nonzero. The formal proof does not use transport--growth pairs:
it shows the inequality $\rho_1(x_1)\le\rho_1(x_0)$ for almost every pair
$([x_0,r_0],[x_1,r_1])$ of points of the cone with $r_0>0$ and $r_1>0$,
with respect to every finite measure on the product of two cones whose
homogeneous marginals, that is, the images of $r_i^2$ times the measure
under $[x_i,r_i]\mapsto x_i$, are $\mu_0$ and $\mu_1$ and whose cost with
$\cospi$ equals $\HK^2(\mu_0,\mu_1)$; every optimal plan on the cone is of
this kind. It obtains
\eqref{cp:eq-nouphill} from it; absolute continuity is proved with an
optimal plan of \eqref{prelim:eq-LET}. The existence of a pair
$(\bfT,q)$ is not needed and not asserted.

\emph{Density bounds for $\HK$.} Proposition~\ref{cp:density-hk} is
formalised for measures with no interface, in the same setting (a complete
separable metric space, a finite reference measure $m$ in place of $\Leb$,
densities with respect to $m$) and for every minimiser of the step
functional that is supplied; existence of minimisers is not part of the
statements. In parts (i) and (ii) the hypothesis $\rho_0>0$ almost
everywhere is replaced by $\mu_0\ne0$, and the assertion $\mu_1\ne0$ of
part (i) requires $m(X)>0$. Part (iii) is proved for $m(X)>0$ and every
finite run of successive minimisers, with $n=0$ included; that every
$\mu_n$ is absolutely continuous is proved, not assumed, and the elements
$e_l$, $e_u$ may differ when $l=u$. The bound \eqref{cp:eq-iterexp} includes
the endpoint $\tau U=\frac14$. As for Lemma~\ref{cp:nouphill}, the formal
proof works with a plan on the cone that attains $\HK^2(\mu_0,\mu_1)$
instead of a pair $(\bfT,q)$, and it covers the case $\mu_1=0$ in the
lower bound directly. Also formalised, and not stated in this paper: if
$m(X)>0$, $0<l\le\rho_0\le u$ and the right derivative of $E$ satisfies
$E'_+(c)\ge-K/\sqrt c$ for all $c\ge c_0$, with $c_0>0$ and $K\ge0$, then
for every $\tau>0$ and every finite run of successive minimisers
$\rho_n\le\big(\sqrt{\max\{u,c_0\}}+2Kn\tau\big)^2$ almost everywhere.

\emph{Density bounds for $\SHK$.} Proposition~\ref{cp:density-shk} is
formalised for measures with no interface, in the same setting and again
for every supplied minimiser: if $\tau>0$, $\mu_0\ll m$ is a probability
whose density with respect to $m$ lies in $[l,u]$ almost everywhere,
$0<l\le u<\infty$, and $\mu_1$ minimises the $\SHK$ step functional over
the probabilities, then $\mu_1\ll m$ and its density lies in $[l,u]$
almost everywhere. The statement on the iterates is proved for every
finite run of successive minimisers, with step sizes that may change from
step to step. The formal proof differs from the printed one. Suppose that
a bound fails. Then there are two levels $a<b$ such that $\mu_1$ has less
mass than $\mu_0$ on $\{\rho_1\le a\}$ and more mass than $\mu_0$ on
$\{\rho_1\ge b\}$. On a plan on the cone that attains $\HK^2(\mu_0,\mu_1)$,
these two strict inequalities show that adding a part of $\mu_0$ on the
first set, and removing a part of $\mu_1$ on the second, each strictly
decrease the $\HK$ distance from $\mu_0$; a convex combination of the two
targets has mass one and is still strictly closer, because
$\HK^2(\mu_0,\cdot)$ is convex, and after a further convex combination
with $\mu_1$ it has no larger energy. The theorem needs neither a pair
$(\bfT,q)$ nor a prescribed element of $\partial E$; for the comparison of
the energies the proof chooses a supporting slope of $E$ at a density
between $a$ and $b$, and the two linear contributions cancel because the
added and the removed mass are equal.
Remark~\ref{cp:rem-sharp} is formalised: the identity for the functional, the unique minimiser, and equality in \eqref{cp:eq-low} for constant initial density.

\emph{Recovery.} Lemma~\ref{cp:recovery} is formalised for measures with
no interface, on a compact metric space with a finite reference measure
$m$ of full support, that is, with $m(U)>0$ for every nonempty open set
$U$, in place of $\Leb$ on a compact convex set. For every
$\mu\in\calM(X)$ there are $z_j$, each with a density with respect to $m$
between $\delta_j$ and $\delta_j^{-1}$ for some $\delta_j\in(0,1)$ and with
$\slope{\sfE}(z_j)<\infty$, such that \eqref{cp:eq-recovery} holds; the
same in $(\calP(X),\SHK)$; and if $\mu\in\calP(X)$ has a density at most
$M$, the densities of the $z_j$ are at most $M$. The measure $\mu$ need
not have finite energy: if $\sfE(\mu)=\infty$ then $\sfE(z_j)\to\infty$.
Consequently the measures with positive bounded density and finite slope
are dense in both spaces, and $\overline{\dom\sfE}$ is the whole space;
both statements are formalised. The proof follows the printed one, with
two changes. The cells consist of the points nearest to the centres of a
finite cover by small balls, with ties given to the least index; each cell
contains a neighbourhood of its centre and so has positive $m$-measure,
which is where full support is used. Without full support the lemma
fails. And the smallness of the
distance from $y_j$ to $z_j$ is obtained by comparing two minimisers $z$
and $w$ of steps from $y_j$ with step sizes $t$ and $s\ge2t$, which gives
$\HK^2(y_j,z)\le4t\,(\sfE(y_j)-\sfE(w))$ $[$the same with $\SHK]$ without
a lower bound for the energy; the three displacement estimates printed in
the proof, including the case $a\tau_j=\frac18$, are formalised
separately. The inequality $\sfE(v_j)\le\sfE(\mu)$ for the cell averages is
proved for every finite measure $\mu$, with singular part and with
$\Einf=\infty$.

\emph{Semiconcavity and the local angle condition.}
Proposition~\ref{cp:semiconc} is formalised for measures with no
interface. For the concavity assertion of part~\ref{cp:it-semiHK}: for a compact subset $X$ of a finite-dimensional Euclidean
space, a finite reference measure $m$ of full support and $0<\delta\le1$
there is $\kappa\ge0$ such that for every $\HK$ geodesic $(\mu_t)$ whose
endpoints have densities between $\delta$ and $\delta^{-1}$ with respect to
$m$, and every $w$ in the $\HK$ geodesic cover of this class, the function
$t\mapsto\frac12\HK^2(\mu_t,w)-\frac\kappa2\HK^2(\mu_0,\mu_1)t^2$ is concave
on $[0,1]$. Convexity of $X$ is not used, since the statement concerns
given geodesics. The formal proof does not follow the proof quoted in the
text: it uses no doubling property of the reference measure. The densities
of the marginals of every optimal entropy-transport plan between two
measures of the geodesic cover, and the distance over which it transports,
are bounded by means of the optimality conditions, a positive lower bound
for the mass of balls of a fixed radius and the bound for the total mass; a geodesic is represented from an interior time
by transport and growth along characteristics; and the pointwise inequality
for the cone distance is integrated. The concavity assertion of
part~\ref{cp:it-semiSHK} is formalised in the same generality, with $\SHK$
geodesics between probability measures of the class and observers in its
$\SHK$ geodesic cover. As in the printed proof it is derived from the $\HK$
statement for $\delta/2$ by passing from the cone to its base; the formal
proof does not use the common ball of Step~2 of the printed proof, because
the transfer is carried out with the bound $\SHK\le\pi/2$, which holds for
all pairs, and it gives the constant
$\frac\pi2\bigl((1+2\kappa)/\sqrt2+2\kappa\bigr)+1$, where $\kappa$ is the
$\HK$ constant for $\delta/2$. The local angle condition is formalised at
every finite measure with $\HK$ and at every probability measure with
$\SHK$, on every closed subset of a finite-dimensional Euclidean space, in
the following form: for every finite family of nonconstant geodesics
$\gamma_1,\dots,\gamma_n$ from the point and all $b_1,\dots,b_n>0$ one has
$\sum_{i,j}b_ib_j\cos\sphericalangle(\gamma_i,\gamma_j)\ge0$ for the upper
angles; for three geodesics this gives that the three upper angles add up
to at most $2\pi$, which is the condition used in Section~\ref{sec:compact}.
For $\HK$ the formal proof lifts the geodesics to paths in the cone that
start over a common law on the initial measure, obtained by gluing along
that measure, with separate labels for the mass that each geodesic creates
from zero; it bounds the distance between two geodesics at two small times
by the cost of the resulting coupling and passes to the limit under the
integral. For $\SHK$ the inequality is transferred from the cone to its
base.

\emph{Invariant sets.} Of Lemma~\ref{cp:invariance}, the nestedness of the
sets $A_j$ and the conditions \ref{cp:A2} and~\ref{cp:A3} are formalised
for measures with no interface, in $(\calM(X),\HK)$ and in
$(\calP(X),\SHK)$, on a compact metric space with a finite reference
measure of full support, under the hypothesis of the lemma that $\sfE$ is
weakly geodesically $\lambda$-convex in the respective space. For
\ref{cp:A2} in the HK case the statement is the printed one: with $\tau_*$
as in the proof, for every $\tau\le\tau_*$, every $n\le\lceil T/\tau\rceil$
and every choice of minimisers of the first $n$ steps, the $n$-th iterate
lies in a set $A_j$ whose index depends only on $j_0$, $T$ and the
constants, and a minimiser of the next step exists; in the SHK case the
density bounds are preserved for every step size and the cutoff is
$\min\{1,(2\ell+1)^{-1}\}$. Condition \ref{cp:A3} is the recovery lemma
above. Condition \ref{cp:A1} follows from Proposition~\ref{cp:semiconc},
which is formalised for measures (previous paragraph), with constants
$\kappa_j$ that serve for the density classes and for the sets $A_j$; with
it the lemma is formalised completely, for a compact subset of a
finite-dimensional Euclidean space.

\emph{Flows on compact domains.} For a compact subset $X$ of a
finite-dimensional Euclidean space, a finite reference measure of full
support, and $E$ convex on $[0,\infty)$ with $E(0)=0$ and continuous at
$0$, the following is formalised for measures with no interface, in
$(\calM(X),\HK)$ and in $(\calP(X),\SHK)$, under the hypothesis that
$\sfE$ is weakly geodesically $\lambda$-convex in the respective space.
There is an $\EVI_\lambda$ flow of $\sfE$ on the whole space; every solution
of the $\EVI_\lambda$ in integral form starting at any point coincides with
it; and it contracts distances by the factor $\ee^{-\lambda t}$: this is
Theorem~\ref{cp:main}\ref{cp:it-flow}. By the equivalence of the two
forms of the $\EVI$, the trajectories of this flow are also solutions in
the differential form, and a solution in the differential form with the
same initial value coincides with them. Consequently $\sfE$ is
$\lambda$-convex along every geodesic, which is
Corollary~\ref{cp:weakstrong}. Theorem~\ref{cp:main}\ref{cp:it-mm} is
formalised as follows. From every $\mu_0\in\dom\sfE$ schemes of exact
minimisers exist for sufficiently small positive step sizes in the HK case
and for every positive step size in the SHK case. Weak convexity provides
geodesic interpolants between consecutive iterates. For every such scheme
with $12\ell\tau\le1$, where $\ell=-\min\{\lambda,0\}$, its piecewise
constant interpolant satisfies \eqref{prelim:eq-MSerror} for the same flow
with convexity constant $-\ell$. For every $T\ge0$ there is a function
$\eta$ with $\tau^{-1/4}\eta(\tau)\to0$ as $\tau\downarrow0$ such that, for
all small $\tau$, every scheme of exact minimisers from $\mu_0$, its
piecewise constant interpolant and each of its geodesic interpolants stay
within $\eta(\tau)$ of the flow on $[0,T]$; this is more than the locally
uniform convergence stated in the theorem. For initial data in
$\bigcup_jA_j$ the distance is at most $C\sqrt\tau$. The formal proof of
the convergence does not use the lower bound of the printed proof: the
values of the Moreau--Yosida regularisation in \eqref{prelim:eq-MSerror}
are bounded below by the energy along the flow. The hypothesis of
weak convexity is assumed, as in the theorem; that it holds for concrete
energies is the content of Theorem~\ref{prelim:LMS}, which is formalised
as described above.

\emph{Growth in one step.} Lemma~\ref{hk:growth} is formalised for measures
with no interface, on a complete separable metric space with a
$\sigma$-finite reference measure, for an arbitrary finite measure $\mu_0$
and any real $\ezero$ not larger than the slope of $E$ at $0$ with
$1+2\tau\ezero>0$. Here $\tau>0$ and $\mu_1$ minimises
$\mu\mapsto\sfE(\mu)+\HK^2(\mu_0,\mu)/(2\tau)$ over $\calM(X)$.
The formal statement is about optimal plans on the cone
instead of transport--growth pairs: for every plan that realises
$\HK^2(\mu_0,\mu_1)$, the radius at the target is at most $k$ times the
radius at the source times the transport cosine, with
$k=(1+2\tau\ezero)^{-1}$; in particular almost every pair of the plan
with positive radius at the target has positive radius at the source and
base points at distance strictly below $\pi/2$, and
$\mu_1(X)\le k^2\mu_0(X)$. This does not say that a $\mu_0$-null set stays
$\mu_1$-null: mass can be transported into it. For every optimal transport--growth pair this
gives \eqref{hk:eq-growth} as printed. The hypotheses $\mu_0=\rho_0\Leb$
with $\rho_0>0$ and the existence of an optimal pair
(Lemma~\ref{prelim:tg}) are not used. For Lemma~\ref{shk:growth}, take $\mu_0$ of finite
energy, $\tau>0$, and a minimiser $\mu_1$ of $\nu\mapsto\sfE(\nu)+\SHK^2(\mu_0,\nu)/(2\tau)$
over $\calP(X)$. Suppose $\omega\ge0$, $3\tau\omega<1$, and \ref{shk:G1} holds.
If an optimal cone plan satisfies \ref{shk:G2} on the parts of the target obtained by
restricting that same plan, it obeys the radius bound with the factor
$(\cos\theta_1-2\tau\omega)^{-1}$ and the transport cosine, and with the
factor $(1-3\tau\omega)^{-1}$, and almost every pair with positive radius
at the target has positive radius at the source and base points at
distance below $\pi/2$; for a given optimal transport--growth pair the
printed \ref{shk:G2}, on subsets of the source, suffices and the conclusion
is \eqref{shk:eq-growth}; and if the slopes of $E$ at $0$ and at infinity lie
between two real numbers $\ezero\le L$, then \ref{shk:G1} and \ref{shk:G2}
hold with $\omega=L-\ezero$ for every initial probability measure and
every optimal plan. For
Lemma~\ref{sub:growth}: with a finite reference measure, an initial
density between $m$ and $M$ and $3\tau E'(M)<1$, the minimiser has a
density between $m$ and $M$ and \eqref{sub:eq-growth} holds, again for
every optimal plan. The one-step
inequality \eqref{hk:eq-onestep} is formalised for every optimal plan with
such a radius bound, with the constant $\pi/(2\sqrt2)$, in $\HK$ and in
$\SHK$.

\emph{Uniform bounds on bounded domains.} Proposition~\ref{hk:tight} is
formalised for measures: for $\Om$ closed and convex, $R$ so large that
$\OmR$ has nonempty interior, and under \ref{hk:H1} and \ref{hk:H2}, the
flow on $\calM(\OmR)$ satisfies (i) and (ii) for every initial measure of finite energy.
If $\supp\mu_0\subseteq\Om_{R_0}$, it also satisfies (iii) for every $\rho>0$ with
$\rho\ge2R_0$. The constants $K_T,D_T$ in \eqref{hk:eq-KD} and $C_1$ in
\eqref{hk:eq-tail} are independent of $R$. The formal proof iterates
the one-step inequality along exact schemes started at $\mu_0$ itself and
passes to the flow by the convergence in
Theorem~\ref{cp:main}\ref{cp:it-mm}. The construction of that compact flow
uses the density bounds and recovery by regular measures
(Lemma~\ref{cp:recovery}). Proposition~\ref{shk:tight} is
formalised in the same way, for every initial probability measure on
$\OmR$ with support in $\bar B_{R_0}$ and every $\rho>0$ with
$\rho\ge2R_0$, with $C_1=\pi K_{T,\omega}^2\sqrt{T\omega}$.
For sums of powers the formalisation proves \eqref{sub:eq-tight} on every
admissible compact truncation, with $C_1=\pi K_{T,\omega}^2\sqrt{T\omega}$
and $\omega=E'(M)$ for initial densities bounded by $M$ and supported in
$\bar B_{R_0}$. The same conditions on $\rho$ apply.

\emph{Minimisers on unbounded domains.} Lemmas~\ref{hk:minimiser}
and~\ref{shk:minimiser} and Theorem~\ref{sub:resolvent} are formalised for
measures, on a proper separable metric space with a $\sigma$-finite
reference measure that is finite on closed balls; this includes every
closed subset of $\R^d$ with Lebesgue measure, so that convexity and a
nonempty interior are not needed. For Lemma~\ref{hk:minimiser} the
hypotheses are those printed, and in addition every minimiser is shown to
have finite energy, whatever the energy of $\mu_0$. For
Lemma~\ref{shk:minimiser} they are those printed. For
Theorem~\ref{sub:resolvent} all four assertions are formalised: $c>0$;
existence of a minimiser under \eqref{sub:eq-tau}; $\sfE(\mu_1)\le c$ for
every minimiser; and the existence of the scheme for every choice of the
minimisers, in the form that every finite sequence of exact steps extends
to an infinite scheme.

\emph{Flows on subdomains.} For an isometry $i$ into a pseudometric space,
a curve solves the $\EVI$ for $\phi\circ i$ if and only if its image
solves it for the functional equal to $\phi$ on the range of $i$ and
$+\infty$ elsewhere; when the domain of $i$ is a metric space, an $\EVI$
flow is transported in the same way. The same equivalence holds for the
sublevel solutions of Definition~\ref{sub:def-sol}, with the same
constants $\lambda_c$.
For a closed subset $A$ of a complete separable metric space, a
$\sigma$-finite reference measure and an energy whose negative part is
finite on the measures concerned, this gives:
an $\EVI_\lambda$ flow of the relaxed energy on $\calM(A)$ or $\calP(A)$
is, after extension by zero, an $\EVI_\lambda$ flow for the ambient energy
restricted to the measures carried by $A$; and for sums of powers, a
family of sublevel solutions on $\calP(A)$ starting from every
probability of finite energy gives, after extension by zero, such a family
for the ambient energy restricted to the probabilities carried by $A$.
This is how the compact flows on $\OmR$ enter
Sections~\ref{sec:hk-unbdd} to~\ref{sec:sublevel}.

\emph{Flows on unbounded domains.} Theorems~\ref{hk:main}
and~\ref{shk:main} are formalised for measures, with no interface: for
$\Om$ closed and convex with nonempty interior in a Euclidean space of any
finite dimension and under the printed hypotheses on $E$, the formal
statements give a flow on all of $\calM(\Om)$, respectively $\calP(\Om)$,
which is an $\EVI_\lambda$ flow; uniqueness of $\EVI_\lambda$ solutions;
the contraction with factor $\ee^{-\lambda t}$; the convergence
\eqref{hk:eq-exhaust} of the flows on $\OmR$, respectively of the flows
from the normalised restrictions, uniformly on $[0,T]$; and the rate
$C\,R^{-1/6}$ for data supported in $\bar B_{R_0}$ and
$R\ge\max\{2R_0,1,R_{\min}\}$. In the formal statement $R_{\min}$ is a fixed function of $\Om$.
In $\HK$ the rate requires $\sfE(\mu_0)<\infty$, with $C=C(T,\lambda,\ezero,\mu_0(\Om),\sfE(\mu_0))$.
In $\SHK$ it holds for every such probability with $C=C(T,\lambda,\omega)$. In Theorem~\ref{shk:main} the threshold
$\max\{2R_0,1,R_{\min}\}$ is thus explicit and independent of $T$ and
$\omega$. Under the step-size hypotheses of
Theorems~\ref{hk:mm} and~\ref{shk:mm}, the formalisation gives the two
estimates of Theorem~\ref{prelim:MS-error}(a) for exact schemes and their
locally uniform convergence. In $\HK$ it includes preparation by
finite-energy initial values for an initial measure of infinite energy.
In $\SHK$ it applies directly to every initial probability.

\emph{Truncation for power energies.} The estimate \eqref{sub:eq-cut} is
formalised for measures: for a $\sigma$-finite reference measure, a sum of
powers with nonnegative coefficients and exponents in $[1,p]$, at least one
term with positive coefficient and exponent above $1$, a probability
measure $\nu$ of finite energy and a measurable set $F$ with
$a=\nu(F)\le\frac12$, the normalised restriction $\nu'$ has finite energy,
the three inequalities of \eqref{sub:eq-cut} hold, and
$\sfE(\nu')-\sfE(\nu)\le2^{p+1}a\,\sfE(\nu)$. The first inequality is also
proved for every $a<1$ and without finiteness of the energy. The
characterisation of finite energy stated after \eqref{sub:eq-E} is
formalised as well: if $p>1$ is the largest exponent and its coefficient
is positive, a finite measure has finite energy if and only if it is
absolutely continuous with a density in $L^p$.

\emph{Regularity on sublevels.} The last sentence of
Lemma~\ref{sub:restrict} is formalised in an arbitrary metric space, for
$\phi$ lower semicontinuous with values in $(-\infty,\infty]$: a
continuous curve $u:[0,\infty)\to\calX_C$ satisfying \eqref{sub:eq-base}
for all $v\in\calX_C$ and $0\le s<t$, as in the lemma, is locally Lipschitz
on $(0,\infty)$; the
limit $|\dot u|(t)=\lim_{h\to0}\sfd(u(t+h),u(t))/|h|$ exists for almost
every $t>0$; $\sfd(u(a),u(b))\le\int_a^b|\dot u|$ for $0\le a\le b$, so that
$u$ is absolutely continuous on every $[0,T]$; every locally integrable
function with this property dominates $|\dot u|$ almost everywhere; and
$\int_s^T|\dot u|^2\le\phi(u(s))-\phi(u(T))$ for $0\le s\le T$. Consequently
the regularity assertions of Theorem~\ref{sub:abstract}(iii) and of the
theorems on sums of powers follow from their assertion (i).

\emph{Sums of powers: the flow.} The bounded and unbounded cases of
Theorem~\ref{sub:powers-compact} are formalised for measures under their
printed hypotheses. In the exhaustion proof every initial probability
admits a sequence of compactly supported probabilities of bounded density
converging to it in $\SHK$ and in extended energy. For
$c=\sfE(\nu_0)<\infty$ and $\ell=p\,2^{p+1}c$, the scheme estimates
in~(iv) hold for $0<\tau\le1$ and $12\ell\tau\le1$, with the additional
condition $2^{p+6}\tau c\le1$ on unbounded domains. In~(v) the
$R^{-1/6}$ estimate holds uniformly for initial densities with a common
bound $M$ and support in $\bar B_{R_0}$, with the constant and threshold
depending only on $T$, $p$, $R_0$ and $E'(M)$.

\emph{Extension to measure data.} Theorem~\ref{sub:linear-extension} is
formalised for a complete metric space under its printed hypotheses: the
extension to the closure $D$ of the domain with its semigroup law, the
continuity of each $\bar S_t$ on $D$ and its uniqueness as a continuous
extension of $S_t$, the uniqueness of each trajectory, the sublevel EVI
for positive-time shifts, the estimate~\eqref{sub:eq-extension-energy},
and the Lipschitz constant $K_w$ on the closed ball.
Corollary~\ref{sub:powers-measure-data} is formalised for measures on
bounded and unbounded closed convex domains with nonempty interior,
including the density of the finite-energy probabilities.

\subsection*{What the formalisation returned to the paper}

Several hypotheses turned out not to be needed, or statements to hold in a
sharper form; where this concerns the paper it is recorded at the
statement, for example in Remarks~\ref{prelim:rem-dependence},
\ref{cp:rem-dini} and~\ref{cp:rem-abstract}. The proof of
Theorem~\ref{sub:abstract} uses Section~\ref{sec:compact} under a
hypothesis weaker than the one stated there; the precise hypothesis is
recorded in Remark~\ref{cp:rem-local}. In three places the first formal
statement of a result quoted from \cite{MuS20} omitted standing
assumptions of that paper (lower semicontinuity, or the exclusion of the
value $-\infty$) and was false as formulated; counterexamples were checked
formally and the statements were corrected. This concerned the formal
statements, not the published theorems.

The formal proof recovers Theorem~6.6 of \cite{LaM22} from
Theorem~\ref{cp:main} for the original energy, including its recession term.
On a compact convex $X\subset\R^d$ with nonempty interior and $d\ge1$,
weak geodesic $\lambda$-convexity in $\HK$ or in $\SHK$ forces $E$ to be
continuous at zero, as one sees by testing the energy on midpoints of
explicit pairs of measures. The $\HK$ flow takes values in $\calM(X)$ and
the $\SHK$ flow in $\calP(X)$. The restriction $d\ge1$ cannot be dropped for
$\HK$: on $\R^0$, the density $E(0)=1$, $E(c)=0$ for $c>0$ gives a strongly
geodesically $0$-convex energy for which the zero measure has finite energy
but admits neither a first minimising step nor an EVI solution starting
there.

\subsection*{Computer algebra}

The identities for $B$, $Q$ and $F$ in Propositions~\ref{ex:reduced},
\ref{ex:Cd} and~\ref{ex:kernel} were also verified symbolically, for general
dimension $d$, with a computer algebra system (exact arithmetic). Numerical
experiments were used only to find conjectures; no statement of this paper
relies on them, and one conjecture suggested by numerics is refuted in
Proposition~\ref{ex:counter}.

\section{Explicit test measures on the cone}\label{nc:sec-engine}

This appendix and the next two prove the converse part of
Theorem~\ref{prelim:LMS} and the implication $(1)\Rightarrow(3)$ of
Theorem~\ref{tr:rigidity}. Both proofs compare the energy at three times of
explicit test measures supported in a small ball, built from a quadratic
potential. A strict contact inequality shows that an explicit plan is
optimal and that every $\HK$ geodesic between the two endpoints passes
through the third measure at its midpoint (Proposition~\ref{nc:engine}). No
solutions of a Hamilton--Jacobi equation and no non-branching property of
$\HK$ geodesics are used.

Throughout, $d\ge1$, $|\cdot|$ and $\langle\cdot,\cdot\rangle$ denote the Euclidean norm and inner product on $\R^d$, $\|M\|_{\mathrm{op}}=\sup_{|v|=1}|Mv|$ is the operator norm of a $d\times d$ matrix $M$, and $B(x_*,r)=\{x\in\R^d:|x-x_*|<r\}$ and $\bar B(x_*,r)=\{x\in\R^d:|x-x_*|\le r\}$ denote open and closed balls centred at $x_*$. On the metric cone $\mathfrak C(\Om)$ over a closed set $\Om\subset\R^d$, with points written $z=[x,r]$ ($x\in\Om$, $r\ge0$) and all points $[x,0]$ identified with the apex $\mathfrak o$, we use the $\pi$-cutoff cone distance
\[
 d_C^2([x_0,r_0],[x_1,r_1])=r_0^2+r_1^2-2r_0r_1\cos\big(\min\{|x_0-x_1|,\pi\}\big),
\]
which is the restriction of the geodesic cone metric on $\mathfrak C(\R^d)$ and agrees with the cone formulation \eqref{prelim:eq-coneOT} of $\HK^2$ \cite[Thms.~2.2, 2.5]{LMS23} (the endpoint support constraints for positive radii keep the minimum over $\mathfrak C(\Om)$ equal to the ambient minimum over $\mathfrak C(\R^d)$; no convexity of $\Om$ is required for this metric).

At time $s=0$ the test map displaces $v$ by $-\arctan(|\nabla\phi(v)|)\nabla\phi(v)/|\nabla\phi(v)|$, interpreted as zero at a critical point. With the displaced point $x$ fixed, the function $v'\mapsto\phi(v')+\log\cos|x-v'|$ has a strict maximum at $v$ on the test ball. The factor $\tilde\Upsilon$ expresses the displacement as a $C^1$ function of the gradient also at zero.

\begin{lemma}[Smooth radial arctangent factor]\label{nc:upsilon}
Define $\tilde\Upsilon,\tilde\Upsilon':\R\to\R$ by
\[
 \tilde\Upsilon(w):=
 \begin{cases}
  \dfrac{\arctan\sqrt{w}}{\sqrt{w}}, & w>0,\\[2mm]
  1-\dfrac{w}{3}, & w\le0,
 \end{cases}
 \qquad
 \tilde\Upsilon'(w):=
 \begin{cases}
  \dfrac{1}{2w}\!\left(\dfrac{1}{1+w}-\dfrac{\arctan\sqrt{w}}{\sqrt{w}}\right), & w>0,\\[2mm]
  -\dfrac{1}{3}, & w\le0.
 \end{cases}
\]
Then $\tilde\Upsilon\in C^1(\R)$ with derivative $\tilde\Upsilon'$, $\tilde\Upsilon(0)=1$, $\tilde\Upsilon'(0)=-\frac13$, and for every $v\in\R^d$,
\begin{equation}\label{nc:eq-upsilon-norm}
 \big|\tilde\Upsilon(|v|^2)\,v\big|=\arctan|v|.
\end{equation}
\end{lemma}

\begin{proof}
For $z\ge0$, the function $k(z):=\arctan z-z+\frac13z^3$ satisfies $k(0)=0$ and $k'(z)=\frac{z^4}{1+z^2}\in[0,z^4]$, whence $0\le\arctan z-z+\frac13z^3\le\frac15z^5$. Substituting $z=\sqrt{w}$ ($w>0$) and dividing by $w\sqrt{w}$ yields
\[
 0\le\frac{\tilde\Upsilon(w)-(1-\frac13w)}{w}\le\frac{w}{5}\qquad(w>0).
\]
Thus the right difference quotient of $\tilde\Upsilon$ at $0$ tends to $-\frac13$, matching the left derivative on $(-\infty,0]$. On $(0,\infty)$ the chain and quotient rules give the displayed formula for $\tilde\Upsilon'$. Moreover, for $w>0$ we have $1-w\le(1+w)^{-1}\le1-w+w^2$ and $0\le\tilde\Upsilon(w)-(1-\frac13w)\le\frac15w^2$, so
\[
 -\frac13-\frac{w}{10}\le\tilde\Upsilon'(w)\le-\frac13+\frac{w}{2}\qquad(w>0).
\]
Hence $\tilde\Upsilon'$ is continuous at $0$, and therefore on all of $\R$. Identity \eqref{nc:eq-upsilon-norm} is immediate for $v=0$ and follows from $\tilde\Upsilon(|v|^2)=\arctan(|v|)/|v|>0$ for $v\ne0$.
\end{proof}

The next two lemmas give the contact inequality behind the optimality of the test plans, and Lemma~\ref{nc:cone-mid} identifies the midpoint of a single cone geodesic; both are used in Proposition~\ref{nc:engine}.

\begin{lemma}[$1$-strong concavity of $\log\cos|\cdot|$]\label{nc:logcos}
For all $w_0,w_1\in B(0,\pi/2)\subset\R^d$ and all $t\in[0,1]$,
\begin{align*}
 (1-t)\Big(\log\cos|w_0|+\frac{|w_0|^2}{2}\Big)
 &+t\Big(\log\cos|w_1|+\frac{|w_1|^2}{2}\Big)\\
 &\le\log\cos\big|(1-t)w_0+tw_1\big|+\frac{|(1-t)w_0+tw_1|^2}{2}.
\end{align*}
\end{lemma}

\begin{proof}
Consider $f(s):=-\log\cos\sqrt{2s}-s$ on $S:=[0,\pi^2/8)$. Then $f$ is continuous on $S$ and differentiable on $(0,\pi^2/8)$ with $f'(s)=\frac{\tan\sqrt{2s}}{\sqrt{2s}}-1$. For $r\in(0,\pi/2)$, we have $(\tan r-r)'=\sec^2r-1>0$, so $\frac{\tan r}{r}\ge1$, and
\[
 \frac{\rmd}{\rmd r}\Big(\frac{\tan r}{r}\Big)=\frac{2r-\sin(2r)}{2r^2\cos^2r}\ge0
\]
since $\sin u\le u$ for $u\ge0$. Hence $f'\ge0$ is nondecreasing on $(0,\pi^2/8)$, so $f$ is both nondecreasing and convex on $S$. Put $s_0=|w_0|^2/2\in S$, $s_1=|w_1|^2/2\in S$, and $s_t=|(1-t)w_0+tw_1|^2/2$. The inner-product identity
\[
 |(1-t)w_0+tw_1|^2=(1-t)|w_0|^2+t|w_1|^2-t(1-t)|w_1-w_0|^2
\]
implies $0\le s_t\le(1-t)s_0+ts_1\in S$, and therefore $f(s_t)\le f((1-t)s_0+ts_1)\le(1-t)f(s_0)+tf(s_1)$. Negating gives the lemma.
\end{proof}

\begin{lemma}[Strict quadratic contact]\label{nc:contact}
Let $M\in\R^{d\times d}$ be symmetric with $\|M\|_{\mathrm{op}}\le\frac14$, let $r_*>0$ satisfy $9r_*<2\pi$, and let $\alpha\in\R$. Put $\phi(v):=\alpha+\frac12\langle Mv,v\rangle$, $c(v):=\ee^{-2\phi(v)}$, and for $v_0\in\bar B(0,r_*)$ set $x:=v_0-\tilde\Upsilon(|Mv_0|^2)\,Mv_0$. Then for every $v_1\in\bar B(0,r_*)$,
\begin{equation}\label{nc:eq-contact-log}
 |x-v_0|<\frac\pi2,\qquad |x-v_1|<\frac\pi2,\qquad
 \phi(v_1)+\log\cos|x-v_1|\le\phi(v_0)+\log\cos|x-v_0|-\frac38|v_1-v_0|^2.
\end{equation}
In particular, with $a(x):=\cos^2|x-v_0|/c(v_0)$ and $b(v_1):=c(v_1)$,
\begin{equation}\label{nc:eq-contact-cos}
 \cos^2|x-v_1|\le a(x)\,b(v_1)\,\ee^{-\frac34|v_1-v_0|^2}\le a(x)\,b(v_1),
\end{equation}
with equality if and only if $v_1=v_0$.
\end{lemma}

\begin{proof}
By \eqref{nc:eq-upsilon-norm}, $|x-v_0|=\arctan|Mv_0|\le|Mv_0|\le\frac14r_*<\frac\pi2$, and hence
\[
 |x-v_1|\le|x-v_0|+|v_0|+|v_1|\le\frac94r_*<\frac\pi2.
\]
Set $h:=v_1-v_0$, $v_t:=(1-t)v_0+tv_1$, and $H(v):=\phi(v)+\log\cos|x-v|$. Applying Lemma~\ref{nc:logcos} to $w_0=x-v_0$ and $w_1=x-v_1$ yields
\begin{align*}
 \log\cos|x-v_t|&-(1-t)\log\cos|x-v_0|-t\log\cos|x-v_1|\\
 &\ge\frac12\big((1-t)|w_0|^2+t|w_1|^2-|w_t|^2\big)
 =\frac12t(1-t)|h|^2.
\end{align*}
By symmetry of $M$ and $\|M\|_{\mathrm{op}}\le\frac14$,
\[
 \phi(v_t)-(1-t)\phi(v_0)-t\phi(v_1)=-\frac12t(1-t)\langle Mh,h\rangle\ge-\frac18t(1-t)|h|^2.
\]
Adding the two inequalities and dividing by $t\in(0,1)$ gives
\[
 H(v_1)-H(v_0)+\frac38(1-t)|h|^2\le\frac{H(v_t)-H(v_0)}{t}\qquad(t\in(0,1)).
\]
We verify that $\left.\frac{\rmd}{\rmd t}\right|_{t=0}H(v_t)=0$. First, $\left.\frac{\rmd}{\rmd t}\right|_{t=0}\phi(v_t)=\langle Mv_0,h\rangle$. For the second term, if $Mv_0=0$ then $x=v_0$ and $|x-v_t|=t|h|$, so $\left.\frac{\rmd}{\rmd t}\right|_{t=0}\log\cos(t|h|)=0=-\langle Mv_0,h\rangle$. If $Mv_0\ne0$, then $x-v_0=-\tilde\Upsilon(|Mv_0|^2)Mv_0\ne0$ has norm $\arctan|Mv_0|\in(0,\pi/2)$, so
\begin{align*}
 \left.\frac{\rmd}{\rmd t}\right|_{t=0}\log\cos|x-v_t|
 &=\frac{\tan|x-v_0|}{|x-v_0|}\langle x-v_0,h\rangle\\
 &=\frac{|Mv_0|}{\arctan|Mv_0|}\Big(-\frac{\arctan|Mv_0|}{|Mv_0|}\langle Mv_0,h\rangle\Big)
 =-\langle Mv_0,h\rangle.
\end{align*}
Thus $\left.\frac{\rmd}{\rmd t}\right|_{t=0}H(v_t)=0$. Letting $t\downarrow0$ gives \eqref{nc:eq-contact-log}; multiplying by $2$ and exponentiating gives \eqref{nc:eq-contact-cos}.
\end{proof}

\begin{lemma}[Midpoint of a cone segment]\label{nc:cone-mid}
Let $z_0=[x_0,r_0]$ and $z_1=[x_1,r_1]$ in $\mathfrak C(\R^d)$ satisfy $r_0>0$, $r_1>0$, and $\theta:=|x_0-x_1|<\pi/2$, and put
\[
 \theta_1:=\arctan\!\left(\frac{r_0\sin\theta}{r_0\cos\theta+r_1}\right)\in[0,\theta].
\]
If $z_m\in\mathfrak C(\R^d)$ satisfies $2\big(d_C^2(z_0,z_m)+d_C^2(z_m,z_1)\big)=d_C^2(z_0,z_1)$, then $z_m=[x_m,r_m]$ with $r_m>0$,
\[
 r_m^2=\frac{r_0^2+r_1^2+2r_0r_1\cos\theta}{4},\qquad
 |x_m-x_1|=\theta_1,\qquad |x_0-x_m|+|x_m-x_1|=\theta,
\]
and $x_m=x_1$ if $\theta=0$, while $x_m=(1-\theta_1/\theta)x_1+(\theta_1/\theta)x_0$ if $\theta>0$. The same holds in $\mathfrak C(\Om)$ for any subset $\Om\subset\R^d$ containing $x_0,x_1$.
\end{lemma}

\begin{proof}
Since $d_C(z_0,z_1)\le d_C(z_0,z_m)+d_C(z_m,z_1)$, the identity
\[
 \big(d_C(z_0,z_m)-d_C(z_m,z_1)\big)^2+\big((d_C(z_0,z_m)+d_C(z_m,z_1))^2-d_C^2(z_0,z_1)\big)=0
\]
forces $d_C(z_0,z_m)=d_C(z_m,z_1)=\frac12d_C(z_0,z_1)$. If $z_m=\mathfrak o$, then $2(r_0^2+r_1^2)=r_0^2+r_1^2-2r_0r_1\cos\theta$, contradicting $r_0^2+r_1^2+2r_0r_1\cos\theta>0$. Hence $z_m=[x_m,r_m]$ with $r_m>0$.

Set $\alpha_0:=\min\{|x_0-x_m|,\pi\}$ and $\alpha_1:=\min\{|x_m-x_1|,\pi\}$. Since $\theta<\pi/2<\pi$, the triangle inequality on $\R^d$ gives $|\theta-\alpha_1|\le\alpha_0\le\pi$. Embed the points into the complex plane via $v_0:=r_0\ee^{\mathrm i\theta}$, $v_1:=r_1$, $w_m:=r_m\ee^{\mathrm i\alpha_1}$, and $v_m:=\frac12(v_0+v_1)$. Then $|v_0-v_1|=d_C(z_0,z_1)$ and $|w_m-v_1|=d_C(z_m,z_1)=\frac12|v_0-v_1|$, while monotonicity of $\cos$ on $[0,\pi]$ gives
\[
 |w_m-v_0|^2=r_0^2+r_m^2-2r_0r_m\cos|\theta-\alpha_1|
 \le r_0^2+r_m^2-2r_0r_m\cos\alpha_0=d_C^2(z_0,z_m)=\frac14|v_0-v_1|^2.
\]
By the parallelogram identity in $\mathbb C$,
\[
 |w_m-v_m|^2=\frac12|w_m-v_0|^2+\frac12|w_m-v_1|^2-\frac14|v_0-v_1|^2\le0.
\]
Thus $w_m=v_m$ and $|w_m-v_0|=d_C(z_0,z_m)$. Comparing real and imaginary parts of $w_m=v_m$ gives
\[
 r_m\cos\alpha_1=\frac{r_0\cos\theta+r_1}{2}>0,\qquad r_m\sin\alpha_1=\frac{r_0\sin\theta}{2}\ge0.
\]
Summing their squares yields the formula for $r_m^2$. Since $\cos\alpha_1>0$ and $\alpha_1\in[0,\pi]$, we have $\alpha_1\in[0,\pi/2)$, so $|x_m-x_1|=\alpha_1$ and $\tan\alpha_1=\frac{r_0\sin\theta}{r_0\cos\theta+r_1}$, whence $\alpha_1=\theta_1$. Moreover, $r_m\sin(\theta-\alpha_1)=\frac12r_1\sin\theta\ge0$ and $\theta-\alpha_1\in(-\pi/2,\pi/2)$ imply $\alpha_1\le\theta$. Equality $|w_m-v_0|^2=d_C^2(z_0,z_m)$ forces $\cos|\theta-\alpha_1|=\cos\alpha_0$, and strict decrease of $\cos$ on $[0,\pi]$ gives $\alpha_0=\theta-\alpha_1<\pi/2$, hence $|x_0-x_m|=\theta-\alpha_1$. Therefore $|x_0-x_m|+|x_m-x_1|=\theta$. If $\theta=0$, then $|x_m-x_1|=0$, so $x_m=x_1$. If $\theta>0$, put $\lambda:=\theta_1/\theta$, $a:=x_0-x_m$, and $b:=x_m-x_1$; then $|a|=(1-\lambda)\theta$, $|b|=\lambda\theta$, and $|a+b|=\theta$, whence $\langle a,b\rangle=\lambda(1-\lambda)\theta^2$ and $|(1-\lambda)b-\lambda a|^2=0$, giving $x_m=(1-\lambda)x_1+\lambda x_0$. The final sentence holds because $d_C$ on $\mathfrak C(\Om)$ is the restriction of $d_C$ on $\mathfrak C(\R^d)$.
\end{proof}

We now define the test measures. The measure $\bar\Gamma_1$ is the constant density $c_*$ on a small ball $K_1$ around $x_*$. For $s\in\{0,\frac12\}$, the mass at $y=x_*+v$ is moved to $\bfT_s(y)$ and multiplied by $|w_s(v)|^2$. In the planar picture of the cone, $w_s(v)=\big((1-s)c(v)+s,(1-s)c(v)|Mv|\big)$ runs along a segment from $w_0(v)$ to $(1,0)$; its polar angle $\arctan(p_s(v)|Mv|)$ is the distance moved (Remark~\ref{nc:rem-strong-shortcut}). The potential sets $c(v)=\ee^{-2\phi(v)}$. Its gradient $Mv$ determines the displacement, directed along $-Mv$, and $|w_s(v)|^2$ gives the mass multiplier. The conditions at $v=0$ say that the transport maps are locally invertible ($J_s(0)>0$) and that the three densities at the centre lie in the interval $(c_1,c_2)$ where $E$ is controlled; by continuity this persists on a small ball. Part (3) computes $\HK^2(\Gamma_0^\eps,\Gamma_1^\eps)$ exactly and classifies all optimal cone plans, and part (4) shows that every $\HK$ geodesic from $\Gamma_0^\eps$ to $\Gamma_1^\eps$ passes through $\Gamma^\eps_{1/2}$ at its midpoint. This is what allows the convexity inequality at $t=\frac12$ on the cone to be evaluated exactly. The same proposition serves without a reservoir ($\eta=0$, Appendix~\ref{nc:sec-lms}) and with one (Appendix~\ref{nc:sec-rigidity}).

\begin{proposition}[Test measures, optimal plans and midpoints]\label{nc:engine}
Let $d\ge1$, let $\Om\subset\R^d$ be closed with nonempty interior, and fix $x_*\in\operatorname{int}\Om$ and $r_*>0$ with $B(x_*,2r_*)\subseteq\operatorname{int}\Om$ and $9r_*<2\pi$. Fix $c_*>0$, an open interval $(c_1,c_2)\subset(0,\infty)$ with $0\le c_1<c_*<c_2\le\infty$, $\alpha\in\R$, and a symmetric matrix $M\in\R^{d\times d}$ with $\|M\|_{\mathrm{op}}\le\frac14$. For $v\in\R^d$ and $s\in[0,1]$, define
\begin{gather*}
 \phi(v):=\alpha+\frac12\langle Mv,v\rangle,\qquad
 c(v):=\ee^{-2\phi(v)}>0,\qquad
 p_s(v):=\frac{(1-s)c(v)}{(1-s)c(v)+s}\in[0,1],\\
 \bfT_s(x_*+v):=x_*+v-p_s(v)\,\tilde\Upsilon\big(p_s(v)^2|Mv|^2\big)\,Mv,\\
 |w_s(v)|^2:=\big((1-s)c(v)+s\big)^2+(1-s)^2c(v)^2|Mv|^2>0.
\end{gather*}
Assume that at $v=0$, for each $s\in\{0,\frac12,1\}$,
\[
 J_s(0):=\det\big(I_d-p_s(0)M\big)>0,\qquad
 \rho_s^{\mathrm{ctr}}:=c_*\frac{\big((1-s)\ee^{-2\alpha}+s\big)^2}{J_s(0)}\in(c_1,c_2).
\]
Then there exists $r_0\in(0,r_*)$ such that for every $r_1\in(0,r_0]$, putting $K_1:=\bar B(x_*,r_1)$ and $K_s:=\bfT_s(K_1)$ for $s\in\{0,\frac12,1\}$:
\begin{enumerate}[label=\textup{(\arabic*)}]
 \item \emph{Local $C^1$ diffeomorphism and common interior carrier.} For each $s\in\{0,\frac12,1\}$, $\bfT_s$ is $C^1$ on $\R^d$ with $D\bfT_s(x_*)=I_d-p_s(0)M$ and $\|D(\bfT_s-\mathrm{id})(y)\|_{\mathrm{op}}\le\frac12$ on an open ball $B(x_*,\tilde r_0)$ with $r_0<\tilde r_0\le r_*$. For each $s\in\{0,\frac12,1\}$, $\bfT_s$ is a $C^1$ diffeomorphism from $B(x_*,\tilde r_0)$ onto an open subset of $B(x_*,2r_*)$, with $J_s(y-x_*):=\det(D\bfT_s(y))>0$ on $B(x_*,\tilde r_0)$, and $K_s=\bfT_s(K_1)$ is a compact subset satisfying
 \[
  K_s\subseteq\bar B\big(x_*,\tfrac32r_1\big)\subset B(x_*,2r_*)\subset\operatorname{int}\Om.
 \]
 \textup{(}Note that $K_0$ may extend outside $K_1$, for example when $M=-\frac14I_d$.\textup{)}

 \item \emph{Pushforward densities and finite energy formula.} For $s\in\{0,\frac12,1\}$, define the measure $\bar\Gamma_s:=(\bfT_s)\pfw\big(c_*|w_s(\cdot-x_*)|^2\Leb|_{K_1}\big)\in\calM(\Om)$ and set
 \[
  \tilde\rho_s(x):=c_*\frac{|w_s(\bfT_s^{-1}(x)-x_*)|^2}{J_s(\bfT_s^{-1}(x)-x_*)}\quad(x\in K_s),\qquad
  \tilde\rho_s(x):=0\quad(x\in\Om\setminus K_s).
 \]
 Then $\tilde\rho_s$ is continuous on the compact set $K_s$ with compact range $K_E^{(s)}:=\tilde\rho_s(K_s)\subset(c_1,c_2)$, and $\bar\Gamma_s=\tilde\rho_s\Leb$. Assuming in addition that $E:[0,\infty)\to(-\infty,\infty]$ satisfies $E(0)=0$, that its restriction $E|_{(c_1,c_2)}$ is real-valued and continuous \textup{(}or real-valued, locally bounded, and Borel\textup{)}, and that $\sfE:\calM(\Om)\to(-\infty,\infty]$ agrees with $\int_\Om E(\rho)\dd x$ on every $\rho\Leb|_\Om\in\calM(\Om)$ whose density $\rho$ is bounded, compactly supported in $\operatorname{int}\Om$, and takes values in $\{0\}\cup K_E$ for a compact $K_E\subset(c_1,c_2)$, we have
 \begin{equation}\label{nc:eq-energy-pullback}
  \sfE(\bar\Gamma_s)=\int_{K_s}E(\tilde\rho_s(x))\dd x
  =\int_{K_1}E\!\left(c_*\frac{|w_s(y-x_*)|^2}{J_s(y-x_*)}\right)J_s(y-x_*)\dd y\in\R.
 \end{equation}

 \item \emph{Exact $\HK^2$ and classification of optimal cone plans \textup{(}with or without a separated reservoir\textup{)}.} Let $\eta\in\calM(\Om)$ be either $\eta=0$ \textup{(}setting $V:=\emptyset$ and $\eps:=1$\textup{)} or a nonzero measure supported on a compact set $V\subset\Om$ satisfying $|x-y|>\pi/2$ for all $x\in B(x_*,2r_*)$ and $y\in V$, and let $\eps>0$. Put $\Gamma_s^\eps:=\eta+\eps\bar\Gamma_s$ for $s\in\{0,\frac12,1\}$. Then
 \begin{equation}\label{nc:eq-hk-exact}
  \HK^2(\Gamma_0^\eps,\Gamma_1^\eps)=\eps\,D^2,\qquad
  D^2:=c_*\int_{K_1}\Big((c(y-x_*)-1)^2+c(y-x_*)^2|M(y-x_*)|^2\Big)\dd y.
 \end{equation}
 Moreover, every optimal finite Borel cone plan $\Pi$ for $(\Gamma_0^\eps,\Gamma_1^\eps)$ in \eqref{prelim:eq-coneOT} \textup{(}with the $\pi$-cutoff distance $d_C$\textup{)} is concentrated on the union of three disjoint sets:
 \begin{enumerate}[label=\textup{(\roman*)}]
  \item $\{(\mathfrak o,\mathfrak o)\}$;
  \item $\{([x,r],[x,r]):x\in V,\ r>0\}$;
  \item $\{([\bfT_0(y),|w_0(y-x_*)|r],[y,r]):y\in K_1,\ r>0\}$.
 \end{enumerate}

 \item \emph{Single-midpoint $\HK$ uniqueness.} Every $\HK$ geodesic $\tilde\Gamma:[0,1]\to\calM(\Om)$ from $\Gamma_0^\eps$ to $\Gamma_1^\eps$ satisfies $\tilde\Gamma_{1/2}=\Gamma_{1/2}^\eps$.
\end{enumerate}
\end{proposition}

\begin{proof}
\emph{(1)} By Lemma~\ref{nc:upsilon}, $q_s(v):=-p_s(v)\tilde\Upsilon(p_s(v)^2|Mv|^2)$ is $C^1$ on $\R^d$ with $q_s(0)=-p_s(0)$. Hence $\Phi_s(v):=\bfT_s(x_*+v)-(x_*+v)=q_s(v)Mv$ is $C^1$ on $\R^d$ with $D\Phi_s(0)=-p_s(0)M$ (since $M(0)=0$), so $D\bfT_s(x_*)=I_d-p_s(0)M$ and $\|D\Phi_s(0)\|_{\mathrm{op}}\le\|M\|_{\mathrm{op}}\le\frac14$. By continuity of $D\Phi_s$, $\det$, $|w_s|^2$, and $J_s$ at $v=0$, there exists $\tilde r_0\in(0,r_*]$ such that on the open ball $B(0,\tilde r_0)$, for all $s\in\{0,\frac12,1\}$,
\[
 \|D\Phi_s(v)\|_{\mathrm{op}}\le\frac12,\qquad J_s(v)>0,\qquad
 c_*\frac{|w_s(v)|^2}{J_s(v)}\in(c_1,c_2).
\]
Pick any $r_0\in(0,\tilde r_0)$. On $B(0,\tilde r_0)$, the mean-value inequality gives $|\Phi_s(v)-\Phi_s(v')|\le\frac12|v-v'|$, whence $|\bfT_s(x_*+v)-\bfT_s(x_*+v')|\ge\frac12|v-v'|$ (injectivity) and, since $\Phi_s(0)=0$, $|\bfT_s(x_*+v)-x_*|\le\frac32|v|$. By the inverse function theorem, $\bfT_s$ is a $C^1$ diffeomorphism from $B(x_*,\tilde r_0)$ onto an open subset of $B(x_*,2r_*)$. For every $r_1\in(0,r_0]$, compactness of $K_1=\bar B(x_*,r_1)$ implies that $K_s=\bfT_s(K_1)\subseteq\bar B(x_*,\frac32r_1)\subset B(x_*,2r_*)\subset\operatorname{int}\Om$ is compact.

\emph{(2)} By (1), $\tilde\rho_s$ is continuous on the compact set $K_s$ with compact image $K_E^{(s)}=\tilde\rho_s(K_s)\subset(c_1,c_2)$. Since $|\det D\bfT_s(y)|=J_s(y-x_*)>0$ on $B(x_*,\tilde r_0)\supset K_1$, the $C^1$ change-of-variables formula gives, for every Borel set $S\subseteq\Om$,
\begin{align*}
 \bar\Gamma_s(S)&=\int_{K_1}\mathbf 1_S(\bfT_s(y))\,c_*|w_s(y-x_*)|^2\dd y\\
 &=\int_{K_1}\mathbf 1_S(\bfT_s(y))\,\tilde\rho_s(\bfT_s(y))\,J_s(y-x_*)\dd y
 =\int_{K_s}\mathbf 1_S(x)\,\tilde\rho_s(x)\dd x,
\end{align*}
so $\bar\Gamma_s=\tilde\rho_s\Leb$. Under the stated assumptions on $E$ and $\sfE$, the restriction $E|_{K_E^{(s)}}$ is real-valued, bounded, and Borel, so $E\circ\tilde\rho_s$ is real-valued, bounded, and Borel on the compact set $K_s$ (and vanishes off $K_s$ since $E(0)=0$), hence integrable; applying the same change of variables to $E(\tilde\rho_s)$ gives \eqref{nc:eq-energy-pullback}.

\emph{(3)} To prove optimality we build a dual bound on the endpoint pairs allowed by the weighted marginals. Integrating it against any admissible cone plan gives the cost of the explicit plan as a lower bound. The strict contact inequality identifies the test pairs where equality holds. On the reservoir, equality pairs a point with itself at the same radius. The remaining equality case is the pair of apices. Note that $\bfT_1=\mathrm{id}$, $|w_1|^2=1$, and $J_1=1$, so $\bar\Gamma_1=c_*\Leb|_{K_1}$. Set $S_0:=K_0\cup V$ and $S_1:=K_1\cup V$. Since $K_0,K_1\subset B(x_*,2r_*)$ and $|x-y|>\pi/2$ on $B(x_*,2r_*)\times V$, both unions are disjoint unions of compact sets, and $\Gamma_0^\eps(\Om\setminus S_0)=\Gamma_1^\eps(\Om\setminus S_1)=0$. Define bounded nonnegative Borel functions $a^\eps,b^\eps:\Om\to[0,\infty)$ by
\[
 a^\eps(x):=
 \begin{cases}
  \dfrac{\cos^2|x-\bfT_0^{-1}(x)|}{c(\bfT_0^{-1}(x)-x_*)}, & x\in K_0,\\[2mm]
  1, & x\in V,\\
  0, & x\notin S_0,
 \end{cases}
 \qquad
 b^\eps(y):=
 \begin{cases}
  c(y-x_*), & y\in K_1,\\
  1, & y\in V,\\
  0, & y\notin S_1.
 \end{cases}
\]
For $y=x_*+v\in K_1$, \eqref{nc:eq-upsilon-norm} gives $|\bfT_0(y)-y|=\arctan|Mv|<\pi/2$, so $\cos^2|\bfT_0(y)-y|=(1+|Mv|^2)^{-1}$ and $|w_0(v)|^2=c(v)^2(1+|Mv|^2)$, whence
\[
 a^\eps(\bfT_0(y))\,c_*|w_0(v)|^2=c_*\,c(v)\qquad\text{and}\qquad
 |w_0(v)|\cos|\bfT_0(y)-y|=c(v).
\]
Integrating yields
\[
 \int_\Om a^\eps\dd\Gamma_0^\eps=\int_\Om b^\eps\dd\Gamma_1^\eps
 =\eta(V)+\eps c_*\int_{K_1}c(y-x_*)\dd y=:I_\eps.
\]
For the upper bound on $\HK^2(\Gamma_0^\eps,\Gamma_1^\eps)$, push forward $\eta$ by $y\mapsto([y,1],[y,1])$ and $\eps c_*\Leb|_{K_1}$ by $y\mapsto([\bfT_0(y),|w_0(y-x_*)|],[y,1])$. Their sum is an admissible cone plan for $(\Gamma_0^\eps,\Gamma_1^\eps)$ whose cost is
\[
 \eps c_*\int_{K_1}\big(|w_0(v)|^2+1-2|w_0(v)|\cos|\bfT_0(y)-y|\big)\dd y
 =\Gamma_0^\eps(\Om)+\Gamma_1^\eps(\Om)-2I_\eps=\eps D^2,
\]
where $v=y-x_*$ and we used $|w_0(v)|^2+1-2c(v)=(c(v)-1)^2+c(v)^2|Mv|^2$. Thus $\HK^2(\Gamma_0^\eps,\Gamma_1^\eps)\le\eps D^2$.

For the lower bound and classification, let $\Pi$ be any admissible finite Borel cone plan for $(\Gamma_0^\eps,\Gamma_1^\eps)$ in \eqref{prelim:eq-coneOT}. The marginal conditions force $r_0=0$ or $x_0\in S_0$, and $r_1=0$ or $x_1\in S_1$, for $\Pi$-a.e.\ $([x_0,r_0],[x_1,r_1])$. We claim that on this full-measure set,
\begin{equation}\label{nc:eq-amgm-contact}
 a^\eps(x_0)r_0^2+b^\eps(x_1)r_1^2-2r_0r_1\cos\big(\min\{|x_0-x_1|,\pi\}\big)\ge0,
\end{equation}
with equality if and only if $(z_0,z_1)$ lies in one of (i)--(iii). Indeed:
\begin{itemize}
 \item If $r_0=0$ or $r_1=0$, the left-hand side of \eqref{nc:eq-amgm-contact} is $a^\eps(x_0)r_0^2+b^\eps(x_1)r_1^2\ge0$, and since $a^\eps>0$ on $S_0$ and $b^\eps>0$ on $S_1$, it vanishes if and only if $r_0=r_1=0$, i.e.\ (i).
 \item If $r_0>0$ and $r_1>0$, then $x_0\in S_0$ and $x_1\in S_1$, so $a^\eps(x_0)>0$ and $b^\eps(x_1)>0$. We verify that $\cos(\min\{|x_0-x_1|,\pi\})\le\sqrt{a^\eps(x_0)b^\eps(x_1)}$ on $S_0\times S_1$:
 \begin{itemize}
  \item on $(K_0\times V)\cup(V\times K_1)$, $|x_0-x_1|>\pi/2$, so $\cos(\min\{|x_0-x_1|,\pi\})<0<\sqrt{a^\eps(x_0)b^\eps(x_1)}$ (strict inequality);
  \item on $V\times V$, $a^\eps(x_0)b^\eps(x_1)=1$, and $\cos(\min\{|x_0-x_1|,\pi\})\le1$ with equality iff $x_0=x_1$;
  \item on $K_0\times K_1$, writing $v_0=\bfT_0^{-1}(x_0)-x_*\in\bar B(0,r_*)$ and $v_1=x_1-x_*\in\bar B(0,r_*)$, Lemma~\ref{nc:contact} gives $|x_0-x_1|<\pi/2$ and $\cos^2|x_0-x_1|\le a^\eps(x_0)b^\eps(x_1)$ with equality iff $v_1=v_0$, i.e.\ $x_0=\bfT_0(x_1)$.
 \end{itemize}
 Combining this with $2\sqrt{a^\eps(x_0)}\,r_0\cdot\sqrt{b^\eps(x_1)}\,r_1\le a^\eps(x_0)r_0^2+b^\eps(x_1)r_1^2$ proves \eqref{nc:eq-amgm-contact}, with equality iff either $x_0=x_1\in V$ and $r_0=r_1$ (case (ii)), or $x_1\in K_1$, $x_0=\bfT_0(x_1)$, and $r_0=\sqrt{b^\eps(x_1)/a^\eps(x_0)}\,r_1=|w_0(x_1-x_*)|r_1$ (case (iii)).
\end{itemize}
Integrating \eqref{nc:eq-amgm-contact} against $\Pi$ and using $\int r_0^2\dd\Pi=\Gamma_0^\eps(\Om)$, $\int r_1^2\dd\Pi=\Gamma_1^\eps(\Om)$, $\int a^\eps(x_0)r_0^2\dd\Pi=I_\eps$, and $\int b^\eps(x_1)r_1^2\dd\Pi=I_\eps$ yields
\begin{align*}
 \int d_C^2(z_0,z_1)\dd\Pi&-(\Gamma_0^\eps(\Om)+\Gamma_1^\eps(\Om)-2I_\eps)\\
 &=\int\Big(a^\eps(x_0)r_0^2+b^\eps(x_1)r_1^2-2r_0r_1\cos\big(\min\{|x_0-x_1|,\pi\}\big)\Big)\dd\Pi\ge0.
\end{align*}
Taking the infimum over $\Pi$ yields $\HK^2(\Gamma_0^\eps,\Gamma_1^\eps)\ge\eps D^2$, proving \eqref{nc:eq-hk-exact}. Furthermore, for every optimal $\Pi$ the nonnegative integrand vanishes $\Pi$-a.e., which yields the stated concentration on (i)--(iii).

\emph{(4)} For any HK midpoint, normalise two optimal cone plans to have the same middle marginal and glue them. The geodesic midpoint identity implies that the endpoint marginal of the glued plan is optimal and that almost every middle cone point is a metric midpoint of its endpoints. Part (3) classifies the endpoint pairs. Lemma~\ref{nc:cone-mid} then determines the positive-radius midpoints, and the apex pair stays at the apex. Let $\tilde\Gamma$ be any $\HK$ geodesic from $\Gamma_0^\eps$ to $\Gamma_1^\eps$, and let $\Pi_{0m}$ and $\Pi_{m1}$ be optimal finite Borel cone plans for $(\Gamma_0^\eps,\tilde\Gamma_{1/2})$ and $(\tilde\Gamma_{1/2},\Gamma_1^\eps)$ in \eqref{prelim:eq-coneOT}. By 2-homogeneity of the cone constraints and of $d_C^2$, replacing $\Pi_{0m}$ on $\{r_m>0\}$ by the pushforward of $r_m^2\Pi_{0m}|_{\{r_m>0\}}$ under $([x_0,r_0],[x_m,r_m])\mapsto([x_0,r_0/r_m],[x_m,1])$ (and similarly for $\Pi_{m1}$) preserves admissibility, finiteness, and optimality, and makes both plans share the same $z_m$-marginal on $\{r_m>0\}$, namely $(x\mapsto[x,1])\pfw\tilde\Gamma_{1/2}$. On $\{z_m=\mathfrak o\}$, adding a nonnegative multiple of $\delta_{(\mathfrak o,\mathfrak o)}$ to whichever plan has smaller mass at $z_m=\mathfrak o$ adds the same apex mass to both ordinary marginals of that plan while preserving its $r_i^2$-weighted marginals and cost and yields two optimal finite Borel plans $P,Q$ on $\mathfrak C(\Om)^2$ with identical $z_m$-marginals on the entire Polish cone $\mathfrak C(\Om)$. Since $\Gamma_0^\eps(\Om)>0$, their common total mass $L:=P(\mathfrak C(\Om)^2)=Q(\mathfrak C(\Om)^2)$ is finite and strictly positive. Applying the probability gluing lemma on Polish spaces \cite[Lem.~5.3.2]{AGS05} to $P/L$ and $Q/L$ yields a probability measure $\gamma$ on $\mathfrak C(\Om)^3$, and setting $\Pi^{(3)}:=L\gamma$ gives a finite Borel measure on $\mathfrak C(\Om)^3$ whose $(0,m)$- and $(m,1)$-marginals are $P$ and $Q$.

Let $\Pi_{01}$ be the $(0,1)$-marginal of $\Pi^{(3)}$, which is admissible for $(\Gamma_0^\eps,\Gamma_1^\eps)$. By $d_C^2(z_0,z_1)\le2d_C^2(z_0,z_m)+2d_C^2(z_m,z_1)$ and the geodesic midpoint property $\HK(\Gamma_0^\eps,\tilde\Gamma_{1/2})=\HK(\tilde\Gamma_{1/2},\Gamma_1^\eps)=\frac12\HK(\Gamma_0^\eps,\Gamma_1^\eps)$,
\begin{align*}
 \HK^2(\Gamma_0^\eps,\Gamma_1^\eps)
 &\le\int d_C^2(z_0,z_1)\dd\Pi^{(3)}
 \le\int 2\big(d_C^2(z_0,z_m)+d_C^2(z_m,z_1)\big)\dd\Pi^{(3)}\\
 &=2\HK^2(\Gamma_0^\eps,\tilde\Gamma_{1/2})+2\HK^2(\tilde\Gamma_{1/2},\Gamma_1^\eps)
 =\HK^2(\Gamma_0^\eps,\Gamma_1^\eps).
\end{align*}
Hence $\Pi_{01}$ is optimal for $(\Gamma_0^\eps,\Gamma_1^\eps)$, and $2\big(d_C^2(z_0,z_m)+d_C^2(z_m,z_1)\big)=d_C^2(z_0,z_1)$ for $\Pi^{(3)}$-a.e.\ $(z_0,z_m,z_1)$. Combining this with part (3) and Lemma~\ref{nc:cone-mid} identifies $z_m$ $\Pi^{(3)}$-a.e.:
\begin{itemize}
 \item in case (i), $z_0=z_1=\mathfrak o$, so $d_C(\mathfrak o,z_m)=0$, hence $z_m=\mathfrak o$;
 \item in case (ii), $z_0=z_1=[x_1,r_1]$ ($x_1\in V$, $r_1>0$), so $\theta=0$ gives $z_m=[x_1,r_1]$;
 \item in case (iii), $z_1=[y,r_1]$ ($y=x_*+v\in K_1$, $r_1>0$) and $z_0=[\bfT_0(y),|w_0(v)|r_1]$ with $\theta=|\bfT_0(y)-y|=\arctan|Mv|<\pi/2$: if $Mv=0$, then $\theta=0$, $x_m=y=\bfT_{1/2}(y)$, and $r_m=\frac12(c(v)+1)r_1=|w_{1/2}(v)|r_1$; if $Mv\ne0$, then $|w_0(v)|\cos\theta=c(v)$ and $|w_0(v)|\sin\theta=c(v)|Mv|$, so
 \[
  r_m^2=\frac{c(v)^2(1+|Mv|^2)+1+2c(v)}{4}r_1^2=|w_{1/2}(v)|^2r_1^2,
  \qquad \theta_1=\arctan\big(p_{1/2}(v)|Mv|\big),
 \]
 and $x_m=y+\frac{\theta_1}{\theta}(\bfT_0(y)-y)=y-p_{1/2}(v)\tilde\Upsilon\big(p_{1/2}(v)^2|Mv|^2\big)Mv=\bfT_{1/2}(y)$.
\end{itemize}
Therefore, for every bounded Borel function $\varphi:\Om\to\R$, since $V$ and $K_1$ are disjoint compact subsets of $\Om$ and $\bfT_{1/2}(K_1)=K_{1/2}\subset\Om$, we obtain a well-defined bounded Borel function $\psi:\Om\to\R$ by setting
\[
 \psi(y):=
 \begin{cases}
  \varphi(y), & y\in V,\\
  |w_{1/2}(y-x_*)|^2\varphi(\bfT_{1/2}(y)), & y\in K_1,\\
  0, & y\in\Om\setminus(V\cup K_1),
 \end{cases}
\]
and compute
\begin{align*}
 \int_\Om\varphi\dd\tilde\Gamma_{1/2}
 &=\int r_m^2\varphi(x_m)\dd\Pi^{(3)}
 =\int r_1^2\psi(x_1)\dd\Pi^{(3)}
 =\int_\Om\psi\dd\Gamma_1^\eps\\
 &=\int_V\varphi\dd\eta+\eps c_*\int_{K_1}|w_{1/2}(y-x_*)|^2\varphi(\bfT_{1/2}(y))\dd y
 =\int_\Om\varphi\dd\Gamma_{1/2}^\eps.
\end{align*}
Hence $\tilde\Gamma_{1/2}=\Gamma_{1/2}^\eps$.
\end{proof}

\begin{remark}[Strong-convexity shortcut]\label{nc:rem-strong-shortcut}
If only \emph{strong} geodesic $\lambda$-convexity is assumed (as in \cite[Thm.~7.5]{LMS23}), part (4) of Proposition~\ref{nc:engine} (gluing and Lemma~\ref{nc:cone-mid}) is not needed at all. Extending $\Gamma_s^\eps:=\eta+\eps(\bfT_s)\pfw(c_*|w_s(\cdot-x_*)|^2\Leb|_{K_1})$ to all $s\in[0,1]$, each $\bfT_s(y)$ is a convex combination of $y\in K_1$ and $\bfT_0(y)\in K_0$ inside the convex open ball $B(x_*,2r_*)\subset\operatorname{int}\Om$, and in the planar cone representation $w_s(v)=((1-s)c(v)+s,(1-s)c(v)|Mv|)$ interpolates linearly from $w_0(v)$ to $w_1(v)=(1,0)$ with polar angle $\phi_s=\arctan(p_s(v)|Mv|)\in[0,\theta]\subset[0,\pi/2)$. Consequently $d_C([\bfT_s(y),|w_s(v)|],[\bfT_t(y),|w_t(v)|])=|t-s|\,d_C([\bfT_0(y),|w_0(v)|],[y,1])$ for all $s,t\in[0,1]$, so pushing forward $\eta$ and $\eps c_*\Leb|_{K_1}$ yields $\HK(\Gamma_s^\eps,\Gamma_t^\eps)\le|t-s|\sqrt\eps\,D$. Since $\HK(\Gamma_0^\eps,\Gamma_1^\eps)=\sqrt\eps\,D$ by part (3), the triangle inequality forces $\HK(\Gamma_s^\eps,\Gamma_t^\eps)=|t-s|\sqrt\eps\,D$, so $s\mapsto\Gamma_s^\eps$ is itself an $\HK$ geodesic in $\calM(\Om)$. Keeping part (4) upgrades the necessity theorem to \emph{weak} geodesic $\lambda$-convexity.
\end{remark}

\section{Converse of Theorem~\ref{prelim:LMS} under weak HK \texorpdfstring{$\lambda$}{lambda}-convexity}\label{nc:sec-lms}

In \cite[Thm.~7.5]{LMS23} the converse part of Theorem~\ref{prelim:LMS} is
proved under strong geodesic $\lambda$-convexity. We prove it under weak
convexity, as stated in Theorem~\ref{prelim:LMS}, in the following
pointwise form, which requires neither convexity nor boundedness of $\Om$,
nor convexity of $E$ for the monotonicity clause.

\begin{theorem}[Pointwise necessity of the LMS condition under weak HK $\lambda$-convexity]\label{nc:LMS-pointwise}
Let $d\ge1$, let $\Om\subset\R^d$ be closed with nonempty interior, and let $E:[0,\infty)\to(-\infty,\infty]$ satisfy $E(0)=0$. Let $\sfE:\calM(\Om)\to(-\infty,\infty]$ be any functional that agrees with $\int_\Om E(\rho)\dd x$ on every absolutely continuous measure $\mu=\rho\Leb|_\Om\in\calM(\Om)$ whose density $\rho$ is bounded, compactly supported in $\operatorname{int}\Om$, and takes values in $\{0\}\cup K_E$ for a compact subset $K_E\subset\operatorname{int}(\dom E)\subset(0,\infty)$. Assume that $\sfE$ is \emph{weakly} geodesically $\lambda$-convex on $(\calM(\Om),\HK)$ for some $\lambda\in\R$.

Then for every $c_*\in(0,\infty)$ and every open interval $(c_1,c_2)\subset(0,\infty)$ containing $c_*$ on which $E$ is finite-valued:
\begin{enumerate}[label=\textup{(\arabic*)}]
 \item \emph{Hessian clause \textup{(}all $d\ge1$\textup{)}.} If $E$ is differentiable on $(c_1,c_2)$ and $E'$ is differentiable at $c_*$, then the matrix $B_\lambda(c_*)$ of Definition~\ref{prelim:matrix} is positive semidefinite:
 \[
  \forall x_0,y_0\in\R,\qquad (x_0,y_0)\,B_\lambda(c_*)\,(x_0,y_0)^\top\ge0.
 \]
 \item \emph{Dimension-monotonicity clause \textup{(}$d\ge2$, no convexity of $E$ assumed\textup{)}.} If $d\ge2$, $E$ is continuous on $(c_1,c_2)$, and $E$ is differentiable at $c_*$, then
 \[
  R(c_*):=c_*E'(c_*)-E(c_*)\ge0.
 \]
\end{enumerate}
\end{theorem}

We prove Theorem~\ref{nc:LMS-pointwise} with Proposition~\ref{nc:engine} for $\eta=0$ and $\eps=1$. The number $\Delta_{\mathrm{ctr}}$ below is the midpoint convexity defect, per unit volume, at the centre of the test measures. Lemma~\ref{nc:center-contradiction} shows that it cannot be negative. Lemmas~\ref{nc:clause1} and~\ref{nc:clause2} choose the potential so that $\Delta_{\mathrm{ctr}}$ becomes negative for a small parameter whenever $B_\lambda(c_*)$ has a negative direction (an isotropic $M$ together with $\alpha$), respectively whenever $R(c_*)<0$ (a traceless shear).

\begin{lemma}[Contradiction from negative center excess]\label{nc:center-contradiction}
Let $\Om$, $E$ and $\sfE$ be as in Theorem~\ref{nc:LMS-pointwise}, without the convexity assumption, let $\lambda\in\R$, and let the data of Proposition~\ref{nc:engine} be given \textup{(}with $\eta=0$ and $\eps=1$\textup{)}, with $E$ finite-valued on $(c_1,c_2)$. If $E$ is continuous on $(c_1,c_2)$ and
\[
 \Delta_{\mathrm{ctr}}:=\frac12E(\rho_0^{\mathrm{ctr}})J_0(0)+\frac12E(c_*)-E(\rho_{1/2}^{\mathrm{ctr}})J_{1/2}(0)-\frac\lambda8\,c_*(\ee^{-2\alpha}-1)^2<0,
\]
then $\sfE$ cannot be weakly geodesically $\lambda$-convex on $(\calM(\Om),\HK)$.
\end{lemma}

\begin{proof}
Define on $\bar B(0,r_0)$ the continuous function
\begin{align*}
 \tilde F_\lambda(v):={}&\frac12E\!\left(c_*\frac{|w_0(v)|^2}{J_0(v)}\right)J_0(v)+\frac12E(c_*)-E\!\left(c_*\frac{|w_{1/2}(v)|^2}{J_{1/2}(v)}\right)J_{1/2}(v)\\
 &-\frac\lambda8\,c_*\Big((c(v)-1)^2+c(v)^2|Mv|^2\Big).
\end{align*}
Since $c(0)=\ee^{-2\alpha}$ and $M(0)=0$, we have $\tilde F_\lambda(0)=\Delta_{\mathrm{ctr}}<0$. By continuity of $\tilde F_\lambda$, choose $r_1\in(0,r_0]$ such that $\tilde F_\lambda(v)\le\frac12\Delta_{\mathrm{ctr}}<0$ for all $v\in\bar B(0,r_1)$. Suppose for contradiction that $\sfE$ is weakly geodesically $\lambda$-convex on $(\calM(\Om),\HK)$. Since $\bar\Gamma_0,\bar\Gamma_1\in\dom\sfE$ by Proposition~\ref{nc:engine}(2), there exists an $\HK$ geodesic $\tilde\Gamma$ from $\bar\Gamma_0$ to $\bar\Gamma_1$ along which \eqref{prelim:eq-convex} holds; at $t=\frac12$ this reads
\[
 0\le\frac12\sfE(\bar\Gamma_0)+\frac12\sfE(\bar\Gamma_1)-\sfE(\tilde\Gamma_{1/2})-\frac\lambda8\,\HK^2(\bar\Gamma_0,\bar\Gamma_1).
\]
By Proposition~\ref{nc:engine}(2)--(4), $\tilde\Gamma_{1/2}=\bar\Gamma_{1/2}$, and the right-hand side equals
\[
 \int_{K_1}\tilde F_\lambda(y-x_*)\dd y\le\frac12\Delta_{\mathrm{ctr}}\,\Leb(K_1)<0,
\]
a contradiction.
\end{proof}

\begin{lemma}[Clause (1): isotropic perturbation and exact center second derivative]\label{nc:clause1}
Let $d\ge1$, $c_*>0$, $(c_1,c_2)\subset(0,\infty)$ with $0\le c_1<c_*<c_2\le\infty$, and $\lambda\in\R$. Suppose $E:(c_1,c_2)\to\R$ is differentiable on $(c_1,c_2)$ and $E'$ is differentiable at $c_*$. For $x_0,y_0\in\R$, put $\beta:=-x_0/d$, $\alpha_1:=y_0/4$, and $\tau_{\max}:=\big(4(1+|\beta|)\big)^{-1}>0$, and define
\begin{equation}\label{nc:eq-F0-general}
 F_{0,\lambda,E}(\tau):=\frac12E\!\left(c_*\frac{\ee^{-4\alpha_1\tau}}{(1-\tau\beta)^d}\right)(1-\tau\beta)^d
 +\frac12E(c_*)-E\!\left(c_*\frac{\mathcal R(\tau)^2}{\mathcal H(\tau)^d}\right)\mathcal H(\tau)^d
 -\frac\lambda8\,c_*(\ee^{-2\alpha_1\tau}-1)^2,
\end{equation}
where $\mathcal R(\tau):=\frac12(\ee^{-2\alpha_1\tau}+1)$ and $\mathcal H(\tau):=1-\tau\beta\,\frac{\ee^{-2\alpha_1\tau}}{\ee^{-2\alpha_1\tau}+1}$. Then $F_{0,\lambda,E}(0)=F'_{0,\lambda,E}(0)=0$, $F'_{0,\lambda,E}$ is differentiable at $0$, and
\begin{equation}\label{nc:eq-F0-second}
 F''_{0,\lambda,E}(0)=\frac14\,(x_0,y_0)\,B_\lambda(c_*)\,(x_0,y_0)^\top.
\end{equation}
Consequently, if $(x_0,y_0)B_\lambda(c_*)(x_0,y_0)^\top<0$, then for every $\tau_{\mathrm{hi}}\in(0,\tau_{\max}]$ there exists $\tau_0\in(0,\tau_{\mathrm{hi}}]$ such that $F_{0,\lambda,E}(\tau_0)<0$ and both arguments $c_*\ee^{-4\alpha_1\tau_0}/(1-\tau_0\beta)^d$ and $c_*\mathcal R(\tau_0)^2/\mathcal H(\tau_0)^d$ lie in $(c_1,c_2)$.
\end{lemma}

\begin{proof}
Write $a(\tau):=\ee^{-2\alpha_1\tau}=\ee^{-\frac12y_0\tau}$, $P(\tau):=\frac{a(\tau)}{a(\tau)+1}$, $J_0(\tau):=(1-\tau\beta)^d=(1+\frac{x_0}{d}\tau)^d$, $u_0(\tau):=c_*a(\tau)^2J_0(\tau)^{-1}$, $J_{1/2}(\tau):=\mathcal H(\tau)^d$, and $u_{1/2}(\tau):=c_*\mathcal R(\tau)^2J_{1/2}(\tau)^{-1}$. For $\tau$ in an open interval around $0$ small enough that $1-\tau\beta>0$, $\mathcal H(\tau)>0$, and $u_0(\tau),u_{1/2}(\tau)\in(c_1,c_2)$, $E$ is differentiable at $u_0(\tau)$ and $u_{1/2}(\tau)$, so the chain rule gives
\[
 \frac{\rmd}{\rmd\tau}\big(J_s(\tau)E(u_s(\tau))\big)=J'_s(\tau)E(u_s(\tau))+J_s(\tau)E'(u_s(\tau))u'_s(\tau)\qquad(s\in\{0,\tfrac12\}).
\]
Since $E'$ is differentiable at $c_*=u_s(0)$ and $J_s(0)=1$, this derivative is differentiable at $\tau=0$ with
\begin{equation}\label{nc:eq-JE-second}
 \big(J_sE(u_s)\big)''\Big|_0
 =J''_s(0)E(c_*)+\big(2J'_s(0)u'_s(0)+u''_s(0)\big)E'(c_*)+E''(c_*)(u'_s(0))^2.
\end{equation}
We evaluate the derivatives at $\tau=0$:
\begin{itemize}
 \item For $s=0$: $J'_0(0)=x_0$, $J''_0(0)=\frac{d-1}{d}x_0^2$, and $\ell_0(\tau):=\log(u_0(\tau)/c_*)=-y_0\tau-d\log(1+\frac{x_0}{d}\tau)$ has $\ell'_0(0)=-(x_0+y_0)$, $\ell''_0(0)=\frac{x_0^2}{d}$. Hence $u'_0(0)=-c_*(x_0+y_0)$ and $u''_0(0)=c_*(\frac{x_0^2}{d}+(x_0+y_0)^2)$, giving $2J'_0(0)u'_0(0)+u''_0(0)=c_*(-\frac{d-1}{d}x_0^2+y_0^2)$. Substituting into \eqref{nc:eq-JE-second} with $A=c_*^2E''(c_*)$ and $R=c_*E'(c_*)-E(c_*)$ yields
 \[
  \big(J_0E(u_0)\big)''\Big|_0=-\frac{d-1}{d}R\,x_0^2+c_*E'(c_*)\,y_0^2+A\,(x_0+y_0)^2.
 \]
 \item For $s=\frac12$: $a(0)=1$, $a'(0)=-\frac12y_0$, $a''(0)=\frac14y_0^2$; $\mathcal R(0)=1$, $\mathcal R'(0)=-\frac14y_0$, $\mathcal R''(0)=\frac18y_0^2$, so $(\log\mathcal R)'(0)=-\frac14y_0$ and $(\log\mathcal R)''(0)=\frac1{16}y_0^2$; $P(0)=\frac12$ and $P'(0)=-\frac18y_0$; $\mathcal H(0)=1$, $\mathcal H'(0)=\frac{x_0}{2d}$, $\mathcal H''(0)=-\frac{x_0y_0}{4d}$, so $(\log\mathcal H)'(0)=\frac{x_0}{2d}$ and $(\log\mathcal H)''(0)=-\frac{x_0y_0}{4d}-\frac{x_0^2}{4d^2}$. Therefore $J'_{1/2}(0)=\frac12x_0$, $J''_{1/2}(0)=\frac{d-1}{4d}x_0^2-\frac14x_0y_0$, and $\ell_{1/2}(\tau):=\log(u_{1/2}(\tau)/c_*)=2\log\mathcal R(\tau)-d\log\mathcal H(\tau)$ satisfies $\ell'_{1/2}(0)=-\frac12(x_0+y_0)$ and $\ell''_{1/2}(0)=\frac18y_0^2+\frac14x_0y_0+\frac{x_0^2}{4d}$. Thus $u'_{1/2}(0)=-\frac{c_*}{2}(x_0+y_0)$, $u''_{1/2}(0)=c_*(\ell''_{1/2}(0)+\frac14(x_0+y_0)^2)$, and $2J'_{1/2}(0)u'_{1/2}(0)+u''_{1/2}(0)=c_*(-\frac{d-1}{4d}x_0^2+\frac14x_0y_0+\frac38y_0^2)$. Substituting into \eqref{nc:eq-JE-second} yields
 \[
  \big(J_{1/2}E(u_{1/2})\big)''\Big|_0
  =\Big(-\frac{d-1}{4d}x_0^2+\frac14x_0y_0\Big)R+\frac38c_*E'(c_*)\,y_0^2+\frac14A\,(x_0+y_0)^2.
 \]
\end{itemize}
At $\tau=0$, $(J_{1/2}E(u_{1/2}))'(0)=\frac12(J_0E(u_0))'(0)$ and $((a-1)^2)'(0)=0$, so $F_{0,\lambda,E}(0)=F'_{0,\lambda,E}(0)=0$. Moreover, $-\frac\lambda8c_*((a-1)^2)''|_0=-\frac\lambda4c_*(a'(0))^2=-\frac{\lambda c_*}{16}y_0^2$. Combining $\frac12(J_0E(u_0))''|_0-(J_{1/2}E(u_{1/2}))''|_0-\frac{\lambda c_*}{16}y_0^2$ gives
\[
 F''_{0,\lambda,E}(0)
 =\frac14\!\left[\Big(A-\frac{d-1}{d}R\Big)x_0^2+2\Big(A-\frac12R\Big)x_0y_0+\Big(A+\frac{c_*}{4}(2E'(c_*)-\lambda)\Big)y_0^2\right],
\]
which is \eqref{nc:eq-F0-second} by \eqref{ex:eq-B}.

If $(x_0,y_0)B_\lambda(c_*)(x_0,y_0)^\top<0$, then $F'_{0,\lambda,E}(\tau)/\tau\to F''_{0,\lambda,E}(0)<0$ as $\tau\downarrow0$, so there exists $\delta\in(0,\tau_{\mathrm{hi}}]$ such that $u_0(\tau),u_{1/2}(\tau)\in(c_1,c_2)$ and $F'_{0,\lambda,E}(\tau)<0$ for all $\tau\in(0,\delta]$. By the mean value theorem, $F_{0,\lambda,E}$ is strictly decreasing on $[0,\delta]$, hence $F_{0,\lambda,E}(\tau_0)<F_{0,\lambda,E}(0)=0$ for $\tau_0=\delta$.
\end{proof}

\begin{lemma}[Clause (2): traceless shear perturbation for $d\ge2$]\label{nc:clause2}
Let $c_*>0$ and $(c_1,c_2)\subset(0,\infty)$ with $0\le c_1<c_*<c_2\le\infty$. Suppose $E:(c_1,c_2)\to\R$ is continuous on $(c_1,c_2)$ and differentiable at $c_*$ with $R(c_*)=c_*E'(c_*)-E(c_*)<0$. Define
\[
 F_0^{\mathrm{shear}}(\tau):=\frac12E\!\left(\frac{c_*}{1-\tau^2}\right)(1-\tau^2)
 +\frac12E(c_*)-E\!\left(\frac{c_*}{1-\tau^2/4}\right)\!\left(1-\frac{\tau^2}{4}\right).
\]
Then there exists $\tau_0\in(0,\frac14]$ such that $F_0^{\mathrm{shear}}(\tau_0)<0$ and both $c_*/(1-\tau_0^2)$ and $c_*/(1-\tau_0^2/4)$ lie in $(c_1,c_2)$.
\end{lemma}

\begin{proof}
Set $\psi(t):=t\,E(c_*/t)$ for $t>0$ such that $c_*/t\in(c_1,c_2)$. This domain is an open neighbourhood of $1$. Differentiability of $E$ at $c_*$ implies that $\psi$ is differentiable at $1$ with $\psi'(1)=E(c_*)-c_*E'(c_*)=-R(c_*)>0$. Consider the extension
\[
 g(h):=\frac12\psi(1-h)+\frac12\psi(1)-\psi\big(1-\tfrac14h\big),
\]
defined on a neighbourhood of $0$ and differentiable \emph{at $0$} (not necessarily on a neighbourhood of $0$), with $g(0)=0$, $g(\tau^2)=F_0^{\mathrm{shear}}(\tau)$ for small $\tau\ge0$, and
\[
 g'(0)=-\frac12\psi'(1)+\frac14\psi'(1)=-\frac14\psi'(1)=\frac14R(c_*)<0.
\]
(Equivalently, $h\mapsto F_0^{\mathrm{shear}}(\sqrt{h})$ has right derivative $\frac14R(c_*)<0$ at $0^+$.) Since $g(h)/h\to g'(0)<0$ as $h\downarrow0$, there exists $h_0\in(0,\frac{1}{16}]$ such that $g(h_0)<0$ and $c_*/(1-h_0),c_*/(1-h_0/4)\in(c_1,c_2)$. Putting $\tau_0:=\sqrt{h_0}\in(0,\frac14]$ yields the lemma.
\end{proof}

\begin{proof}[Proof of Theorem~\ref{nc:LMS-pointwise} and of the converse part of Theorem~\ref{prelim:LMS}]
Fix $x_*\in\operatorname{int}\Om$ and $r_*>0$ with $B(x_*,2r_*)\subseteq\operatorname{int}\Om$ and $9r_*<2\pi$.

\emph{Proof of Theorem~\ref{nc:LMS-pointwise}(1).} Suppose for contradiction that $(x_0,y_0)B_\lambda(c_*)(x_0,y_0)^\top<0$ for some $x_0,y_0\in\R$. Put $\beta:=-x_0/d$, $\alpha_1:=y_0/4$, and $\tau_{\max}:=(4(1+|\beta|))^{-1}$. Lemma~\ref{nc:clause1} yields $\tau_0\in(0,\tau_{\max}]$ with $F_{0,\lambda,E}(\tau_0)<0$ and $c_*\ee^{-4\alpha_1\tau_0}/(1-\tau_0\beta)^d,c_*\mathcal R(\tau_0)^2/\mathcal H(\tau_0)^d\in(c_1,c_2)$. Apply Proposition~\ref{nc:engine} (with $\eta=0$, $\eps=1$) and Lemma~\ref{nc:center-contradiction} to $\alpha:=\tau_0\alpha_1$ and the symmetric matrix $M:=(\tau_0\beta)I_d$:
\begin{itemize}
 \item $\|M\|_{\mathrm{op}}=\tau_0|\beta|\le\frac14$;
 \item $D\bfT_s(x_*)=(1-\tau_0\beta p_s(0))I_d$, so $J_s(0)=(1-\tau_0\beta p_s(0))^d\ge(\frac34)^d>0$ for $s\in\{0,\frac12,1\}$;
 \item $\rho_0^{\mathrm{ctr}}=c_*\ee^{-4\alpha_1\tau_0}/(1-\tau_0\beta)^d\in(c_1,c_2)$, $\rho_{1/2}^{\mathrm{ctr}}=c_*\mathcal R(\tau_0)^2/\mathcal H(\tau_0)^d\in(c_1,c_2)$, $\rho_1^{\mathrm{ctr}}=c_*\in(c_1,c_2)$;
 \item $\Delta_{\mathrm{ctr}}=F_{0,\lambda,E}(\tau_0)<0$.
\end{itemize}
Since $E$ is differentiable, hence continuous, on $(c_1,c_2)$, Lemma~\ref{nc:center-contradiction} contradicts weak geodesic $\lambda$-convexity of $\sfE$.

\emph{Proof of Theorem~\ref{nc:LMS-pointwise}(2).} Assume $d\ge2$ and suppose for contradiction that $R(c_*)<0$. Lemma~\ref{nc:clause2} provides $\tau_0\in(0,\frac14]$ with $F_0^{\mathrm{shear}}(\tau_0)<0$ and $c_*/(1-\tau_0^2),c_*/(1-\tau_0^2/4)\in(c_1,c_2)$. Pick an orthonormal pair $e_1,e_2\in\R^d$ and apply Proposition~\ref{nc:engine} (with $\eta=0$, $\eps=1$) and Lemma~\ref{nc:center-contradiction} to $\alpha:=0$ and the symmetric traceless shear $Mv:=\tau_0\big(\langle e_1,v\rangle e_1-\langle e_2,v\rangle e_2\big)$:
\begin{itemize}
 \item $\|M\|_{\mathrm{op}}=\tau_0\le\frac14$ (by Bessel's inequality, $|Mv|^2=\tau_0^2(\langle e_1,v\rangle^2+\langle e_2,v\rangle^2)\le\tau_0^2|v|^2$);
 \item $p_0(0)=1$, $p_{1/2}(0)=\frac12$, $p_1(0)=0$, and in an orthonormal basis extending $\{e_1,e_2\}$, $I_d-t(e_1e_1^\top-e_2e_2^\top)=\operatorname{diag}(1-t,1+t,1,\dots,1)$ has determinant $1-t^2$, whence
 \[
  J_0(0)=1-\tau_0^2\ge\frac{15}{16}>0,\qquad J_{1/2}(0)=1-\frac{\tau_0^2}{4}\ge\frac{63}{64}>0,\qquad J_1(0)=1;
 \]
 \item since $\alpha=0$, $\rho_0^{\mathrm{ctr}}=c_*/(1-\tau_0^2)\in(c_1,c_2)$, $\rho_{1/2}^{\mathrm{ctr}}=c_*/(1-\tau_0^2/4)\in(c_1,c_2)$, $\rho_1^{\mathrm{ctr}}=c_*\in(c_1,c_2)$, and the $\lambda$-term $-\frac\lambda8c_*(\ee^{-2\alpha}-1)^2$ vanishes identically, so $\Delta_{\mathrm{ctr}}=F_0^{\mathrm{shear}}(\tau_0)<0$.
\end{itemize}
Lemma~\ref{nc:center-contradiction} again yields a contradiction.

\emph{The converse part of Theorem~\ref{prelim:LMS}.} Under the hypotheses of Theorem~\ref{prelim:LMS} with $E\in C^2((0,\infty))$ and $\sfE$ weakly geodesically $\lambda$-convex on $(\calM(\Om),\HK)$, Definition~\ref{prelim:energy} guarantees that $\sfE(\rho\Leb|_\Om)=\int_\Om E(\rho)\dd x$ for every bounded compactly supported density $\rho$ (zero singular part). Applying Theorem~\ref{nc:LMS-pointwise} at every $c_*>0$ gives $B_\lambda(c_*)\ge0$ for all $c_*>0$ and $(d-1)R(c_*)\ge0$ for all $c_*>0$ (the latter is also automatic from convexity of $E$ and $E(0)=0$). Proposition~\ref{ex:reduced} then gives \eqref{prelim:eq-LMScond} for $d\ge1$.

For $d=0$ the set $\Om$ is a point $\{x\}$ with unit measure, and the finite
measures are $r^2\delta_x$, $r\ge0$. The cone identity gives
$\HK(r^2\delta_x,s^2\delta_x)=|r-s|$, so the only geodesic from $r_0^2\delta_x$
to $r_1^2\delta_x$ is $r_t^2\delta_x$ with $r_t=(1-t)r_0+tr_1$. Since
$r_t^2-(r_1-r_0)^2t^2$ is affine in $t$, weak geodesic $\lambda$-convexity
of $\sfE(r^2\delta_x)=E(r^2)$ is convexity of $r\mapsto E(r^2)-\frac\lambda2r^2$
on $[0,\infty)$. For $d=0$ this function of $\gamma=r$ is $\Nfun(\rho,\gamma)$,
which does not depend on $\rho$; so $\Nfun$ is convex on $(0,\infty)^2$, and
the monotonicity condition on $(d-1)\Nfun$ holds trivially.
\end{proof}

\section{Proof of the rigidity theorem}\label{nc:sec-rigidity}

Here $X\subset\R^d$ is compact and convex with nonempty interior, $\diam X>\pi/2$ (which forces $d\ge1$), and $\sfE_q(\mu)=-\int_X\rho^q\dd x$ for $\mu=\rho\Leb+\mu^s$ (with $\Einf=0$).

The proof of $(1)\Rightarrow(3)$ uses Proposition~\ref{nc:engine} with a separated reservoir $\eta$; it compares energies at three times and needs no second derivatives along curves. We add to a fixed probability $\eta$ far away a small multiple $\eps$ of the test measures of Appendix~\ref{nc:sec-engine}, and normalise. Nothing is transported between the reservoir and the test part. So the negative midpoint defect of the test part enters the energy on the sphere with the factor $\eps^q$, while the normalisation and the term with $\lambda$ contribute terms of order $\eps$ (Lemma~\ref{nc:dilution}). Since $q<1$, the defect wins for small $\eps$. Lemma~\ref{nc:shk-asymp} shows that the normalised midpoint is reached at a time $t_\eps\to\frac12$. We isolate two elementary scalar limits.

\begin{lemma}[SHK midpoint asymptotics]\label{nc:shk-asymp}
Let $m_0,m_1,D\in\R$, and let $\delta_\eps\in(0,\pi/2]$ be defined for sufficiently small $\eps>0$ so that
\[
 (1+\eps m_0)+(1+\eps m_1)-2\sqrt{1+\eps m_0}\sqrt{1+\eps m_1}\cos\delta_\eps=\eps D^2.
\]
Put
\[
 t_\eps:=\frac{1}{\delta_\eps}\arctan\!\left(\frac{\sqrt{1+\eps m_1}\sin\delta_\eps}{\sqrt{1+\eps m_0}+\sqrt{1+\eps m_1}\cos\delta_\eps}\right).
\]
Then, as $\eps\downarrow0$,
\[
 \frac{\delta_\eps^2}{\eps}\to D^2\qquad\text{and}\qquad t_\eps\to\frac12.
\]
\end{lemma}

\begin{proof}
Write $a_\eps:=\sqrt{1+\eps m_0}\to1$ and $b_\eps:=\sqrt{1+\eps m_1}\to1$. Using $1-\cos\delta_\eps=2\sin^2(\delta_\eps/2)$, the hypothesis rewrites as
\[
 (a_\eps-b_\eps)^2+4a_\eps b_\eps\sin^2(\delta_\eps/2)=\eps D^2.
\]
Since $(a_\eps-b_\eps)^2/\eps=\eps(m_0-m_1)^2/(a_\eps+b_\eps)^2\to0$, dividing by $\eps a_\eps b_\eps$ yields $4\sin^2(\delta_\eps/2)/\eps\to D^2$. In particular $\sin(\delta_\eps/2)\to0$. Because $\delta_\eps/2\in(0,\pi/4]$, Jordan's inequality $\sin(\delta_\eps/2)\ge\frac{2}{\pi}\frac{\delta_\eps}{2}=\frac{\delta_\eps}{\pi}$ forces $\delta_\eps\to0$. Since $\frac{\sin u}{u}\to1$ as $u\to0$,
\[
 \frac{\delta_\eps^2}{\eps}=\frac{4\sin^2(\delta_\eps/2)/\eps}{\big(\sin(\delta_\eps/2)/(\delta_\eps/2)\big)^2}\to D^2.
\]
For the second limit, set $y_\eps:=\frac{b_\eps\sin\delta_\eps}{a_\eps+b_\eps\cos\delta_\eps}>0$. Then
\[
 \frac{y_\eps}{\delta_\eps}=\frac{b_\eps}{a_\eps+b_\eps\cos\delta_\eps}\cdot\frac{\sin\delta_\eps}{\delta_\eps}\to\frac{1}{1+1}\cdot1=\frac12.
\]
In particular $y_\eps\to0$, and since $\frac{\arctan y}{y}\to1$ as $y\to0$, we conclude $t_\eps=\frac{\arctan y_\eps}{y_\eps}\cdot\frac{y_\eps}{\delta_\eps}\to\frac12$.
\end{proof}

\begin{lemma}[Dilution contradiction]\label{nc:dilution}
Let $q<1$, $\lambda\in\R$, $A_0,A_1,A_{1/2},b_0,m_0,m_1,m_{1/2},D\in\R$, and suppose $\Delta_A:=\frac12A_0+\frac12A_1-A_{1/2}<0$. If functions $t_\eps,\delta_\eps$ satisfy $t_\eps\to\frac12$ and $\delta_\eps^2/\eps\to D^2$ as $\eps\downarrow0$, then the inequality
\begin{equation}\label{nc:eq-dil-ineq}
 0\le(1-t_\eps)G_\eps(0)+t_\eps G_\eps(1)-G_\eps(\tfrac12)-\frac\lambda2\,t_\eps(1-t_\eps)\,\delta_\eps^2,
\end{equation}
where $G_\eps(s):=(1+\eps m_s)^{-q}(\eps^qA_s+b_0)$, cannot hold for all sufficiently small $\eps>0$.
\end{lemma}

\begin{proof}
Suppose \eqref{nc:eq-dil-ineq} holds for all sufficiently small $\eps>0$. Divide by $\eps^q>0$:
\begin{align*}
 0&\le(1-t_\eps)(1+\eps m_0)^{-q}A_0+t_\eps(1+\eps m_1)^{-q}A_1-(1+\eps m_{1/2})^{-q}A_{1/2}\\
 &\quad+b_0\left[(1-t_\eps)\frac{(1+\eps m_0)^{-q}-1}{\eps^q}+t_\eps\frac{(1+\eps m_1)^{-q}-1}{\eps^q}-\frac{(1+\eps m_{1/2})^{-q}-1}{\eps^q}\right]\\
 &\quad-\frac\lambda2\,t_\eps(1-t_\eps)\,\frac{\delta_\eps^2}{\eps}\,\eps^{1-q}.
\end{align*}
Since $1-q>0$, for each $m\in\R$ we have $\frac{(1+\eps m)^{-q}-1}{\eps^q}=\eps^{1-q}\cdot\frac{(1+\eps m)^{-q}-1}{\eps}\to0\cdot(-qm)=0$ as $\eps\downarrow0$, and similarly $\frac{\delta_\eps^2}{\eps}\,\eps^{1-q}\to D^2\cdot0=0$. Passing to the limit $\eps\downarrow0$ gives $0\le\frac12A_0+\frac12A_1-A_{1/2}=\Delta_A<0$, a contradiction.
\end{proof}

\begin{proof}[Proof of Theorem~\ref{tr:rigidity}]
$(3)\Rightarrow(2)$ is contained in Proposition~\ref{tr:negpowers}, and $(2)\Rightarrow(1)$ is immediate because $(\calP(X),\SHK)$ is a geodesic space. It remains to prove $(1)\Rightarrow(3)$. Assume (1) and suppose for contradiction that (3) fails.

\emph{Step 1 (separated interior balls).} Since $X$ is convex with nonempty interior, $\overline{\operatorname{int}X}=X$ and hence $\diam(\operatorname{int}X)=\diam X>\pi/2$. Pick $x_*,y_*\in\operatorname{int}X$ with $|x_*-y_*|>\pi/2$ (which also forces $d\ge1$). Choose $r_*>0$ with $r_*\le\pi/5$ (so $9r_*<2\pi$), $B(x_*,2r_*)\subseteq\operatorname{int}X$, $B(y_*,r_*)\subseteq\operatorname{int}X$, and $4r_*\le|x_*-y_*|-\pi/2$. By the triangle inequality, every $p\in B(x_*,2r_*)$ and $z\in B(y_*,r_*)$ satisfy $|p-z|>\pi/2$. Put $V:=\bar B(y_*,r_*/2)\subset B(y_*,r_*)\subset\operatorname{int}X$ and $\eta:=\Leb(V)^{-1}\mathbf 1_V\Leb\in\calP(X)$, and write $b_0:=\sfE_q(\eta)$.

\emph{Step 2 (negative direction and fixed time $\tau_0$).} For $E(c)=-c^q$ at $c_*=1$, we have $E(1)=-1$, $E'(1)=-q$, $E''(1)=q(1-q)$, so $A(1)=q(1-q)$ and $R(1)=1-q$, and \eqref{ex:eq-B} at $\lambda=0$ gives
\[
 B_0(1)=\begin{pmatrix}
 (1-q)\big(q-\frac{d-1}d\big)&(1-q)\big(q-\frac12\big)\\[0.3em]
 (1-q)\big(q-\frac12\big)&q\big(\frac12-q\big)
 \end{pmatrix},\qquad
 \det B_0(1)=\frac{(1-q)\big(\frac12-q\big)\big((d+2)q-d\big)}{2d}.
\]
If $q>\frac12$, the $(2,2)$ entry is negative; if $q\le\frac12$ and $d\ge3$, or $q<\frac12$ and $d=2$, the $(1,1)$ entry is negative; if $d=1$ and $q<\frac13$, the determinant is negative. Hence, since (3) fails, there exist $x_0,y_0\in\R$ with $(x_0,y_0)B_0(1)(x_0,y_0)^\top<0$ (in particular $(x_0,y_0)\ne(0,0)$). Setting $\beta:=-x_0/d$, $\alpha_1:=y_0/4$, and $\tau_{\max}:=(4(1+|\beta|))^{-1}$, Lemma~\ref{nc:clause1} (with $E(c)=-c^q$, $c_*=1$, $\lambda=0$, and $(c_1,c_2)=(0,\infty)$) yields $\tau_0\in(0,\tau_{\max}]$ with $F_{0,0,E}(\tau_0)<0$.

\emph{Step 3 (localized measures, negative defect, and midpoint uniqueness).} Apply Proposition~\ref{nc:engine} on $\Om=X$ with $c_*=1$, $(c_1,c_2)=(0,\infty)$, $\alpha:=\tau_0\alpha_1$, $M:=(\tau_0\beta)I_d$, and the separated reservoir $\eta$ on $V$. Here $\|M\|_{\mathrm{op}}=\tau_0|\beta|\le\frac14$ and $J_s(0)=(1-\tau_0\beta p_s(0))^d\ge(\frac34)^d>0$. The continuous integrand
\[
 \tilde F_0(v):=\frac12E\!\left(\frac{|w_0(v)|^2}{J_0(v)}\right)J_0(v)+\frac12E(1)-E\!\left(\frac{|w_{1/2}(v)|^2}{J_{1/2}(v)}\right)J_{1/2}(v)
\]
satisfies $\tilde F_0(0)=F_{0,0,E}(\tau_0)<0$. Choose $r_1\in(0,r_0]$ so that $\tilde F_0(v)\le\frac12F_{0,0,E}(\tau_0)<0$ on $\bar B(0,r_1)$. Writing $A_s:=\sfE_q(\bar\Gamma_s)$ and $m_s:=\bar\Gamma_s(X)$ for $s\in\{0,\frac12,1\}$, Proposition~\ref{nc:engine}(2) yields
\[
 \Delta_A:=\frac12A_0+\frac12A_1-A_{1/2}=\int_{K_1}\tilde F_0(y-x_*)\dd y\le\frac12F_{0,0,E}(\tau_0)\,\Leb(K_1)<0.
\]
Moreover, in \eqref{nc:eq-hk-exact} the integrand $(c(v)-1)^2+c(v)^2(\tau_0\beta)^2|v|^2$ is continuous and nonnegative on $K_1=\bar B(x_*,r_1)$ and cannot vanish identically because $\tau_0>0$ and $(\alpha_1,\beta)\ne(0,0)$ (if $\alpha_1\ne0$ it is positive at $v=0$; if $\alpha_1=0$ then $\beta\ne0$ and it is positive for all $v\ne0$). Hence $D^2>0$.

For each $\eps>0$, put $\Gamma_s^\eps:=\eta+\eps\bar\Gamma_s$, $M_\eps(s):=\Gamma_s^\eps(X)=1+\eps m_s$, and $\nu_s^\eps:=\Gamma_s^\eps/M_\eps(s)\in\calP(X)$. By Proposition~\ref{nc:engine}(3)--(4), $\HK^2(\Gamma_0^\eps,\Gamma_1^\eps)=\eps D^2>0$ and every $\HK$ geodesic from $\Gamma_0^\eps$ to $\Gamma_1^\eps$ has midpoint $\Gamma_{1/2}^\eps$. Because $K_0,K_1$ are disjoint from $V$, $\Gamma_0^\eps$ and $\Gamma_1^\eps$ both coincide with $\eta$ on $V$ ($\eta(V)=1$) while $\Gamma_0^\eps\ne\Gamma_1^\eps$, so they are not proportional; thus $\nu_0^\eps\ne\nu_1^\eps$ and $\delta_\eps:=\SHK(\nu_0^\eps,\nu_1^\eps)\in(0,\pi/2]$. Furthermore, by homogeneity of $\sfE_q$, disjointness of $V$ and $K_s$, and $\Einf=0$,
\[
 G_\eps(s):=\sfE_q(\nu_s^\eps)=(1+\eps m_s)^{-q}\big(\eps^qA_s+b_0\big)\qquad(s\in\{0,\tfrac12,1\}).
\]

\emph{Step 4 (lifting an arbitrary $\SHK$ geodesic and dilution).} By weak geodesic $\lambda$-convexity (1), for each $\eps>0$ there exists an $\SHK$ geodesic $\sigma^\eps:[0,1]\to\calP(X)$ from $\nu_0^\eps$ to $\nu_1^\eps$ along which $\sfE_q$ is $\lambda$-convex. As in the proof of Lemma~\ref{tr:lift}, the cone identity \eqref{prelim:eq-HKSHK} implies that $\iota(r\ee^{\mathrm i\theta}):=r^2\sigma^\eps(\theta/\delta_\eps)$ is an isometry of the planar sector $S_{\delta_\eps}=\{r\ee^{\mathrm i\theta}:r>0,\ 0\le\theta\le\delta_\eps\}$ into $(\calM(X),\HK)$. Putting
\[
 z_\eps(s):=(1-s)\sqrt{M_\eps(0)}+s\sqrt{M_\eps(1)}\,\ee^{\mathrm i\delta_\eps}=|z_\eps(s)|\,\ee^{\mathrm i\theta_\eps(s)}\qquad(s\in[0,1]),
\]
the curve $s\mapsto\tilde\Gamma^\eps_s:=\iota(z_\eps(s))$ is an $\HK$ geodesic from $\Gamma_0^\eps$ to $\Gamma_1^\eps$. By single-midpoint $\HK$ uniqueness (Proposition~\ref{nc:engine}(4)), $\tilde\Gamma^\eps_{1/2}=\Gamma_{1/2}^\eps$. Normalising both sides gives $\sigma^\eps(t_\eps)=\nu_{1/2}^\eps$, where
\[
 t_\eps:=\frac{\theta_\eps(1/2)}{\delta_\eps}
 =\frac{1}{\delta_\eps}\arctan\!\left(\frac{\sqrt{M_\eps(1)}\sin\delta_\eps}{\sqrt{M_\eps(0)}+\sqrt{M_\eps(1)}\cos\delta_\eps}\right)\in(0,1).
\]
Here $t_\eps\in(0,1)$ because $\delta_\eps\in(0,\pi/2]$ ensures $\operatorname{Re}z_\eps(1/2)>0$, $\operatorname{Im}z_\eps(1/2)>0$, and $\sin(\delta_\eps-\theta_\eps(1/2))=\frac12\sqrt{M_\eps(0)}\sin\delta_\eps/|z_\eps(1/2)|>0$. Evaluating \eqref{prelim:eq-convex} along $\sigma^\eps$ at $t=t_\eps$ yields \eqref{nc:eq-dil-ineq}. On the other hand, the cone identity \eqref{prelim:eq-HKSHK} gives
\[
 (1+\eps m_0)+(1+\eps m_1)-2\sqrt{1+\eps m_0}\sqrt{1+\eps m_1}\cos\delta_\eps
 =\HK^2(\Gamma_0^\eps,\Gamma_1^\eps)=\eps D^2.
\]
Lemma~\ref{nc:shk-asymp} gives $\delta_\eps^2/\eps\to D^2$ and $t_\eps\to\frac12$ as $\eps\downarrow0$, and Lemma~\ref{nc:dilution} yields the desired contradiction.
\end{proof}

\section{Corrections to the cited papers}\label{sec:errata}

Formalising a proof means reading every cited statement literally. In
doing so we found the points listed below. They are misprints, missing
hypotheses, and a few auxiliary statements that are false as printed but
true in the form in which they are used. As far as we could determine,
none of them makes a main theorem of its paper false. We list only points
that were confirmed independently by at least two of the AI assistants or by a
counterexample checked in Lean; the mark ``(Lean)'' is given only where the
check is a file of the repository given in the note after the abstract. Unless stated
otherwise, the versions read are the arXiv versions named in each heading,
and a point may already be corrected in the published article; where we
compared with the published article we say so. Points that concern only
the presentation of a proof are omitted.

\subsection*{Liero, Mielke, Savar\'e \cite{LMS18}, arXiv:1508.07941v2}

\begin{enumerate}[label=\textup{(\arabic*)},leftmargin=2.2em]
\item \emph{Lemma~7.9(iii).} The cost identity is claimed for all
completions of a relaxed plan by a plan between the residual marginals. It
holds for optimal completions. For $\mu_1=\delta_0$, $\mu_2=\delta_\pi$ on
$\R$ the zero plan is optimal for the relaxed problem and $\HK^2=2$, but
the completion by the single pair of unit radii has cost~$4$. Also in the
published article.
\item \emph{Theorem~4.11, Step~1 of the proof.} The lower semicontinuity
needed for the minimax theorem fails on the class used there; the term in
question is upper semicontinuous. This is a gap in one step and does not
affect the duality theorem.
\item \emph{Theorem~6.3 d).} The inequality
$\sigma_1\sigma_2\ge\cos^2$ of part b) is stated there on a product
$A_1\times A_2$ of sets of full measure. The potentials of part d) are
set to $+\infty$ on $\supp\mu_i\setminus A_i$, so the inequality
between the potentials and the cost asserted in d) is obtained on
$A_1\times A_2$ and where a potential equals $-\infty$, not on the
supports.
\item \emph{Theorem~7.6.} The statement fails literally for
$\mu_1=\mu_2=0$, which has to be treated separately. Also in the published
article.
\item \emph{Theorem~7.20(i).} For an arbitrary finite plan $\alpha$ on the
cone the converse needs the hypothesis that both homogeneous marginals of
$\alpha$ are dominated by the given measures (Lean). Also in the published
article.
\item \emph{Remark~7.5.} A finite plan on the cone has a lift that is a
probability measure exactly when it has mass one. The minimisation
formula~(7.23), which is (7.24) in the published article, correctly uses
finite plans.
\item \emph{The product plan (5.28)}, which is (5.29) in the published
article. When both measures are nonzero it satisfies the homogeneous
marginal constraints exactly when $a_i^p=\mu_i(X_i)$, not for arbitrary
radii $a_1,a_2$.
\item \emph{Lemma~7.3, proof.} The middle inequality of the geometric
series estimate fails at $m=1$. The final bound $4^{-m}(1+M)$ holds for
all $m\ge1$, so the conclusion is unaffected.
\item \emph{Example~5.2~E.1, (5.4).} The infimum over positive transported
mass and its lower semicontinuous envelope are equal on the open quadrant
only; they differ on the axes (Lean).
\end{enumerate}

\subsection*{Laschos, Mielke \cite{LaM19}, arXiv:1712.01888v2}

Only the arXiv version was read.

\begin{enumerate}[label=\textup{(\arabic*)},leftmargin=2.2em]
\item \emph{After (3.11).} The class $[\nu,r]$ is identified with $r\nu$;
it is $r^2\nu$, as in the introduction and in (3.10).
\item \emph{Theorem~3.3 (scaling), proof.} The statement allows
$r_0,r_1\ge0$ and the proof divides by both radii. If one factor vanishes
the only feasible plan is the zero plan and the identity holds.
\item \emph{Before the identity for the squared radius along a geodesic.}
The second derivative of $r^2$ equals $2\,\sfd_{\mathfrak C}^2(z_0,z_1)$;
the factor $2$ is missing. The mass identity that follows is correct.
\item \emph{Lemma~2.1, (2.2).} The leading factor $r_0r_1$ should be
$\tilde r_0\tilde r_1$.
\item \emph{The formula for the angle along a cone geodesic.} The
expressions with the principal arcsine and arctangent are not valid for
all angles in $(0,\pi)$; the form with the arccosine is correct.
\item \emph{Lemma~2.23(A).} The printed constant is false; there is a
counterexample in $(\R^2,\|\cdot\|_4)$. For a distance profile with values
in $(0,R]$, $R<\pi$, the constant $1+(1+2K)R/\sin R$ works, and the
hypothesis is needed along the whole curve.
\item \emph{Proposition~2.27.} In (A) the hypotheses $r_0,r_1>0$ and
$x_0\ne x_1$ are implicit. In (B) a positive radius of the observer is
necessary.
\item \emph{Definition~2.5 and Proposition~2.8.} A geodesic is required to
have positive speed, so there is none from a point to itself, and the
rescaled curve of Proposition~2.8 can be constant (Lean). Constant curves
should be allowed.
\item \emph{Proposition~4.1 (angle condition).} The lifted cone geodesics
do not all have the same length, so the normalisation has to be weighted;
an integrable minorant is needed for Fatou's lemma; ``strictly positive''
for the vertex term should be ``nonnegative''.
\item \emph{Lemma~4.3, (4.2).} The set $A_b$ is defined with ``for all
$x\in A$'' where ``there exists'' is meant; as printed the estimate is
false (Lean). The open $b$-neighbourhood of $A$ is the right set.
\item \emph{Lemma~4.9.} The proof divides by $\sigma_1$ before showing
that it is positive, and the covering argument needs local doubling up to
a larger radius than stated. The constants depend on $\ell$ and on the
doubling data.
\item \emph{Lemma~4.11.} The radii of the lift are inverted; they are
$r_i=\sigma_i^{-1/2}$.
\item \emph{Theorem~4.12.} The inclusion of the density-bounded
probabilities in the corresponding set of measures, used in the proof,
depends on the total mass $V$ of the reference measure. Let
$0<V<\infty$. On a geodesic
space with more than one point the inclusion holds for a given $\delta$ only if $V=1$, and on a
one-point space exactly when $\delta^2\le V\le\delta^{-2}$ (Lean).
Normalising the reference measure to mass one repairs this step. (For
$V=0$, which the definition allows, both sets of probabilities are empty.)
\item \emph{Definition of $\mathsf{LET}$, (3.3).} The functional can take
the value $+\infty$.
\end{enumerate}

\subsection*{Liero, Mielke, Savar\'e \cite{LMS23}, arXiv:2208.14299}

The numbering agrees with the published article.

\begin{enumerate}[label=\textup{(\arabic*)},leftmargin=2.2em]
\item \emph{Theorem~2.14(1) and Theorem~3.3(1).} The inequality
$\sigma_0\sigma_1\ge\cos^2$ and the bounds $0<\sigma_i<\infty$ are claimed
on the supports. They hold on sets of full measure. For
$\mu_0=\sum_n\ee^{-n^3}\delta_{-2^{-n}}$ and
$\mu_1=\sum_k2^{-k}\delta_{\pi/2-2^{1-k}}$ on $\R$ they fail at the point
$0$ of the support of $\mu_0$. The full-measure form suffices for
everything that follows.
\item \emph{Corollary~4.2, (4.31).} The bounds for the backward solution
are misprinted. They should read
$-\frac1{2(1-t)}\le\bar\xi_t\le\frac1{2t}$. Also in the published article.
\item \emph{Corollary~3.5.} Uniqueness of the pair $(\bfT,q)$ fails on
$\{q=0\}$, where neither the push-forward nor the cost depends on $\bfT$.
The cone-valued map $[\bfT,q]$ is unique.
\item \emph{Theorem~5.2, proof.} Only $0<s<t\le1$ is treated; the case
$0\le t<s<1$ needs the same argument for a pair that is made tight first.
\item \emph{The first-order identity after Corollary~3.5.} The sign is
reversed: $\nabla\varphi_0(x)=\nabla\mathsf L_1(\bfT(x)-x)$.
\item \emph{Proposition~4.1(5).} The relation $r_0=Zr_2$ should be
$r_0^2=Zr_2^2$, and a positive radius is needed.
\item \emph{Proposition~4.1(6), (4.21).} The equivalence fails in a
degenerate case at $s=0$, when $\xi_0=+\infty$ on a ball of radius
$\pi/2$.
\item \emph{Before (4.6).} Attainment of the infimum in (4.4) is claimed
without a restriction on $t$; it fails at $t=\tau$ and holds for $t<\tau$.
\item \emph{Theorem~4.3(1) and (7).} The contact sets at the two endpoint
times need not be closed, and the transport maps to the endpoint times
need not be Lipschitz; both statements hold for interior times. Example:
$\mu_0=\delta_0$, $\mu_1=\delta_{\pi/2}$ (Lean). Also in the published
article.
\item \emph{Theorem~4.3(2), (3), estimates (4.40).} The two constants are
exchanged, and the hypothesis on $\bar a$ should be $\bar a\ge-\frac12$
(Lean for the counterexamples). Also in the published article.
\item \emph{Theorem~4.3(4), proof.} The sign of the elapsed time is
reversed in the case $s>t$. The statement (4.43) is correct.
\item \emph{Theorem~4.8, (4.62c,d).} The right-hand sides are evaluated at
a point $x$ that is independent of the point $y$ at which the Jacobian is
taken. The two points have to coincide.
\item \emph{Remark~2.11.} The formula loses the case of no transport: for
unit radii at distance $\pi/2$ the printed infimum is $+\infty$ and the
cone cost is $2$ (Lean). The minimum should be taken over transported mass
$\theta\ge0$.
\item \emph{Proposition~6.2, last line of the proof.} The exponent is
inverted: $s=c^{d/(2(d+2))}$.
\item \emph{Normalisation of $\lambda$.} The introduction and the text
before (6.4) use the shift $E-\lambda c$, Theorem~7.2 uses
$E-\frac\lambda2c$. Only the latter is consistent with (1.21).
\item \emph{Theorem~7.2, proof.} Steps~2 and~3 apply Proposition~6.2 to a
density that need not be smooth, while that proposition is proved for
$C^2$ densities. Applying it to smoothed densities and passing to the
limit closes the gap.
\item \emph{Theorem~7.5, proof.} It is asserted without proof that the
characteristic curves of a smooth solution of the Hamilton--Jacobi
equation are geodesics. Misprint: $\|\mathrm D^2\xi_s\|^2=4d\beta^2$. The
present paper does not use this theorem; see
Appendices~\ref{nc:sec-engine}--\ref{nc:sec-rigidity}.
\end{enumerate}

\subsection*{Laschos, Mielke \cite{LaM22}, arXiv:2207.09815v3}

The article has appeared in print \cite{LaM22}. The points below and all
our references to \cite{LaM22} follow arXiv version~3; we did not compare
with the published article.

The corrections that matter for the results are discussed in
Section~\ref{cp:ss-points}. For completeness:

\begin{enumerate}[label=\textup{(\arabic*)},leftmargin=2.2em]
\item \emph{Lemma~5.17.} The bound is false as printed; a term of order
$\kappa\tau$ is missing (Lean). See Lemma~\ref{cp:error}.
\item \emph{Proposition~4.3.} The final inequality is reversed when
$E'(c_{\rm low})>0$. The conclusion is true; see
Proposition~\ref{cp:density-hk}.
\item \emph{Proposition~4.6.} The argument gives $16\tau^2e_*^2$ where the
statement has $4\tau^2e_*^2$, and a denominator has the wrong sign.
\item \emph{Theorem~5.9.} Its hypotheses are not met in the application:
it assumes convexity along every geodesic, compact sublevels, finite slope
on the density-bounded sets, and solvability of the scheme for every step.
See Theorem~\ref{cp:abstract}.
\item \emph{Theorem~6.6 and the Main Result.} Both flows are said to take
values in $\calP(X)$; the flow on the cone takes values in $\calM(X)$.
\item \emph{Lemma~5.15, proof.} The function $q$ should be $\sfd^2/2$, and
the second difference quotient has its numerator reversed.
\item \emph{Proposition~A.1.} The radial identity should read
$r(t)^2=1-t(1-t)D^2$, and the lower bound for the radius is
$\cos(\delta/2)$.
\item \emph{Section~2, the example with $E(c)=-\sqrt c+\frac13c^{3/2}$.}
The displayed derivative of the energy and the displayed squared slope
omit the factor $\Leb(X)$, and the text gives $E'(\rho)$ an extra factor
$\frac12$. The conclusion is unaffected.
\item \emph{Lemma~6.2.} The hypothesis $0<d_1<\frac\pi2$ on the radius,
which is part of \cite[Lemma~4.9]{LaM19}, is missing; the statement fails
for $d_1=\frac\pi2$. See the proof of Proposition~\ref{cp:semiconc}.
\end{enumerate}
\section{The dependency tree}\label{sec:tree}

\begin{figure}[p]
\centering
\includegraphics[angle=90,width=\textwidth,height=0.84\textheight,keepaspectratio]{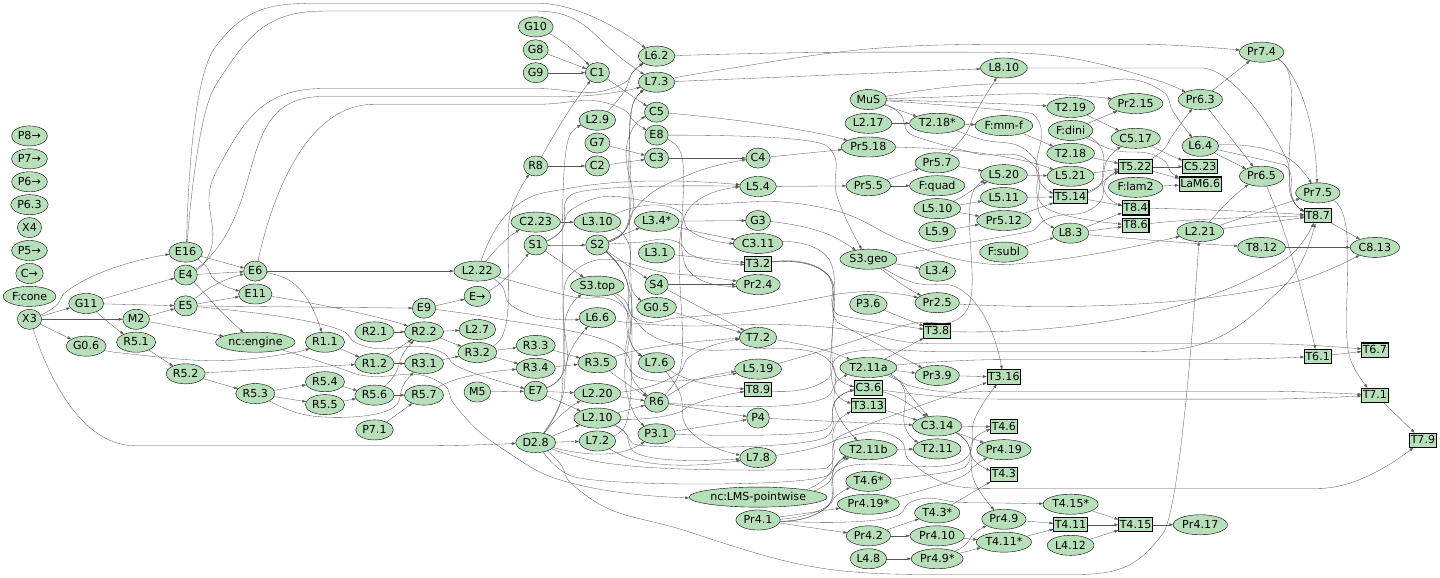}
\caption{The dependency tree, to be read from the bottom of the page
upwards. An arrow points from a statement to a statement whose proof uses
it; arrows implied by others are omitted. Every node is formalised in
Lean. Boxes are the main theorems. Labels such as T5.22, L2.7, Pr4.17, C3.14 are the numbers of
theorems, lemmas, propositions and corollaries of this paper; the other
labels name auxiliary and cited results.}\label{tree:fig}
\end{figure}

Figure~\ref{tree:fig} shows the statements of this paper and the results
of the literature they rest on, with an arrow from each statement to
every statement whose proof uses it. The tree was built together with the formalisation and is
checked mechanically against the sources: every node carries the place
where its statement is printed and, if it is formalised, the names of the
corresponding Lean theorems. At the time of writing it has
\treeNodes{} nodes, and all of them are formalised in Lean. What is
formalised for each statement, and at which scope, is described in
Appendix~\ref{sec:verification}. A version with the full statement of
every node is part of the repository given in the note after the abstract.

\end{document}